\documentclass[a4paper,10pt]{article}
\usepackage{geometry}
\usepackage{stmaryrd}
\usepackage{bbm}
\usepackage{amsfonts}
\usepackage{amsmath}
\usepackage{amssymb}
\usepackage{mathrsfs}
\usepackage{accents}
\usepackage{amscd}

\newtheorem{atheorem}{\bf \temp}[section]
\newtheorem{thm}[atheorem]{Theorem}
\newtheorem{cor}[atheorem]{Corollary}
\newtheorem{lem}[atheorem]{Lemma}
\newtheorem{prop}[atheorem]{Proposition}
\newtheorem{de}[atheorem]{Definition}
\newtheorem{rem}[atheorem]{Remark}

\numberwithin{equation}{section}

\title{\textbf{Cyclotron  Damping along a Uniform Magnetic Field }}
\author{Xixia Ma \footnote{Yau Mathematical Sciences Center, Tsinghua University. E-mail addresses:kfmaxixia@tsinghua.edu.cn} \ \ \ \ \ \ \
 \\}
\date{}

\begin{document}

\maketitle
\textbf{Abstract.}  We prove cyclotron damping  for the collisionless Vlasov-Maxwell equations  on $\mathbb{T}_{x}^{3}\times\mathbb{R}_{v}^{3}$ under the assumptions
that the electric induction is zero and  $(\mathcal{\mathbf{PSC}})$ holds. It is a crucial step  to solve the stability problem of the Vlasov-Maxwell equations. Our proof is based on
a new dynamical system of the plasma particles, originating  from Faraday Law of Electromagnetic induction and Lenz's Law. On the basis of it,
we use the improved Newton iteration scheme to show the damping mechanism.

%\textbf{Keywords.}
\begin{center}
%\item\section{Linear}
{\bf\large 1. \quad Introduction }
\end{center}

In this paper, it is assumed that the plasma system is collisionless, nonrelativistic, and hot. Cyclotron damping describes the phenomenon  that a plasma with a prescribed
zero-order distribution function, imbedded in a uniform magnetic field, which is assumed to be perturbed by an electromagnetic wave propagating parallel to the field.
A number of treatments of the problem of cyclotron damping have appeared in the literature[18,20,24,25], but there are few rigorous mathematical results except
 the recent  results of  Bedrossian and Wang [7]  . The usual method that deals with this phenomenon is via the Vlasov
equations. In this paper we study cyclotron damping  at the level of kinetic description based on the Vlasov equations
from the mathematical view point. First, we analyze the  Vlasov-Maxwell equations from perspective of both equilibrium and stability theories.

Now we give a detailed description of the Vlasov-Maxwell equations. We  denote the particles distribution function by $f=f(t,x,v),$
 and the electric and magnetic fields by $E(t,x)$ and $B(t,x),$ respectively. Then the Vlasov equation says
 \begin{align}
 \partial_{t}f+v\cdot\nabla_{x}f+(E+v\times B)\cdot\nabla_{v}f=0.
 \end{align}
 The electric and magnetic fields $E(t,x)$ and $B(t,x)$ in Eq.$(0.1)$ are determined from Maxwell's equations:
 \begin{align}
 &\nabla\cdot E(t,x)=\int_{\mathbb{R}^{3}}f(t,x,v)dv,\quad \nabla\times B(t,x)=\int_{\mathbb{R}^{3}}vf(t,x,v)dv+\frac{\partial E(t,x)}{\partial t},\notag\\
 &\frac{\partial B(t,x)}{\partial t}=-\nabla\times E(t,x),\quad \nabla\cdot B(t,x)=0.
 \end{align}
 Note that Eq.$(0.1)$ is nonlinear since $E(t,x)$ and $B(t,x)$ are determined in terms of $f(t,x,v)$ from Maxwell's equations $(0.2).$

 An equilibrium analysis of Eq.$(0.1)$ and Eqs.$(0.2)$ proceeds by setting $\frac{\partial}{\partial t}=0$ and looking for stationary solutions,
 $f^{0}(x,v),E^{0}(x),B^{0}(x),$ that satisfy the equations
 \begin{align}
 & v\cdot\nabla_{x}f^{0}(x,v)+(E^{0}+v\times B^{0})\cdot\nabla_{v}f^{0}(x,v)=0,\quad \nabla\cdot B^{0}(x)=0,\notag\\
 &\nabla\times B^{0}(x)=\int_{\mathbb{R}^{3}}vf^{0}(x,v)dv,\quad \nabla\cdot E^{0}(x)=\int_{\mathbb{R}^{3}}f^{0}(x,v)dv.\notag\\
 \end{align}
 An analysis of Eq.(0.3) reduces to a determination of the particle constants of the motion in the equilibrium fields $E^{0}(x)$ and  $B^{0}(x).$ In this paper, we
 assume that $E^{0}(x)=0,$ namely, $\int_{\mathbb{R}^{3}}f^{0}(x,v)dv=0.$ This implies that there are no deviations from charge neutrality in equilibrium,
 $B^{0}(x)$ is produced by external current sources as well as any equilibrium plasma currents.

 A stability analysis based on Eq.(0.1) and Eqs.(0.2) proceeds in the following manner. The quantities $f(t,x,v),E(t,x),$ and $B(t,x)$ are expressed as the sum
  of their equilibrium values plus a time-dependent perturbation:
  \begin{align}
  &f(t,x,v)=f^{0}(x,v)+\delta f(t,x,v),\quad E(t,x)=E^{0}(x)+\delta E(t,x),\quad B(t,x)=B^{0}(x)+\delta B(t,x).\notag\\
  \end{align}
  The quantities $f^{0}(x,v),E^{0}(x)$ and $B^{0}(x)$ satisfy (0.3). The time development of the perturbations $\delta f(t,x,v),\delta E(t,x),$ and $\delta B(t,x)$
  is studied by using Eq.(0.1) and Eqs.(0.2). For small-amplitude perturbations, the Vlasov-Maxwell equations are linearized about the equilibrium
  $f^{0}(x,v),E^{0}(x)$ and $B^{0}(x).$ This gives
  \begin{align}
 &\frac{\partial\delta f}{\partial t} + v\cdot\nabla_{x}\delta f(t,x,v)+(E^{0}+v\times B^{0})\cdot\nabla_{v}\delta f(t,x,v)
 =-(\delta E+v\times \delta B)\cdot\nabla_{v}f^{0}(x,v),\quad\notag\\
 &\nabla\cdot \delta B(t,x)=0,\nabla\times \delta B(t,x)=\int_{\mathbb{R}^{3}}v\delta f(t,x,v)dv+\frac{\partial\delta E}{\partial t},\quad\nabla\cdot
  \delta E(t,x)=\int_{\mathbb{R}^{3}}\delta f(t,x,v)dv.\notag\\
 \end{align}
If the perturbations $\delta f(t,x,v),\delta E(t,x),$ and $\delta B(t,x)$ grow, then the equilibrium distribution $f^{0}(x,v)$ is unstable.
Otherwise, the perturbations damp, so
 the system returns to equilibrium and is stable.  We assume that the equilibrium $f^{0}$ is independent of space, namely, $f^{0}(x,v)=f^{0}(v).$

From the above analysis and the form of Vlasov-Maxwell equations, it is obvious that when $B\equiv0,$ cyclotron damping is reduced to Landau damping. Hence,
the method used is similar to that employed by
Mouhot and Villani [23]. However,  compared with the electric field, a static magnetic field introduces a fascinating complication
into the motion of charged particles. And the particles trajectories become helices, spiraling around the magnetic lines of force.
This severe alteration of the orbits tends to inhibit transport across the magnetic field.
%To make the statement clear, we have to use the mathematical language different from the above section, but here we don't give the rigorous mathematical proof.
 The mechanism of Landau damping  does depend on the transfer of electric field energy to particles moving in phase with the wave. However, for cyclotron damping
 in electromagnetic  plasmas, the electric field of the wave is perpendicular to the direction of the magnetic field  and the particle drifts
  and accelerates the particle perpendicular to the drift direction.

In the following  we recall Landau damping through  gathering lots of physical literature and results of mathematical articles.
  The existence of a damping mechanism by which plasma particles
absorb wave energy  was found by L.D.Landau at the linear level,
under the condition that the plasma is not cold and the velocity distribution is of finite extent. Next in linear case, many works from mathematical aspects found in [9,15,25,28]
 gave rigorous proofs under different assumptions. Later, a ground-breaking work for Landau damping was made by Mouhot and Villani in the  nonlinear case.
 They gave the first and rigours proof of nonlinear Landau damping under the assumption of the electric field. In this paper we will extend their results and prove that
 cyclotron damping in electromagnetic fields  still occurs.

Now we will make a brief statement about  the connection and difference  between  the results of [23] and ours.
 First, in   electric field case, Mouhot and Villani proved the existence of Landau damping
 under assumption of the $(\mathbf{L})$ condition, that is expressed as follows:

 ($\mathbf{L}$) There are constants  $C_{0},\lambda,\kappa>0$ such that $|\hat{f}^{0}(\eta)|\leq C_{0}e^{-2\pi\lambda|\eta|}$ for any $\eta\in\mathbb{R}^{d};$
  and for any $\xi\in\mathbb{C}$ with $0\leq \textmd{Re}\xi<\lambda,$
 \begin{align}
  \inf_{k\in\mathbb{Z}^{d}}|\mathcal{L}(\xi,k)-1|\geq\kappa,
  \end{align}
  where we define a function
  $\mathcal{L}(\xi,k)=-4\pi^{2}\int^{\infty}_{0}e^{2\pi\xi^{\ast}|k|t}\widehat{W}(k)\hat{f}^{0}(kt)|k|^{2}tdt,$
  and $\xi^{\ast}$ is the complex conjugate to $\xi.$

  To some extent, (0.6) of the  $(\mathbf{L})$ condition is similar to the ``Small Denominators" condition in KAM theory in [1], but  is stronger. Here we will consider
  the condition of cyclotron damping from a totally different  perspective, in detail, we will give a  physical condition that we call
  the ``$\mathbf{Physical}$ $\mathbf{Stability}$ $\mathbf{Condition}$",in short, ``$\mathcal{\mathcal{\mathbf{PSC}}}$", which is stated in the following form (here we assume the background magnetic field $B_{0}$ along the  $\hat{z}$ direction):

  $(\mathbf{PSC})$ :
   for any component velocity in the $\hat{z}$ direction $v_{3}\in\mathbb{R},$
             there exists some positive constant $v_{Te}$ such that if $v_{3}=\frac{\omega}{k_{3}},\omega,k$
             are frequencies of time and space $t,x,$ respectively,  then $|v_{3}|\gg v_{Te}.$

              $(\mathbf{PSC})$  tells us that  the number of particles that the wave velocity greatly
         exceeds  their velocity is much larger than the number of particles' velocity  slower than the wave velocity. And
          we will show that cyclotron damping occurs under the above conditions, and that isn't only consistent with the physical observation, but also is the same
          with  the ``Small Denominators" condition in KAM theory in [1] in some sense.

          Second, compared with the  electric field case, it is easy to find a new term $v\times B$ in the electromagnetic field setting.
         And this brings many difficulties because of the unboundedness of $ v.$  Based on the physical facts of  Faraday Law of
          Electromagnetic
           induction and Lenz's Law, we know that $\delta B$ generates the force that inhibits the change of the electric field.
           This helps us estimate such term. And the above fact leads  us  to
           study the following
           dynamics of the particles trajectory:

 \begin{align}
   \left\{\begin{array}{l}
   \frac{d}{dt}X_{t,\tau}(x,v)=V_{t,\tau}(x,v),\quad\frac{d}{dt}V_{t,\tau}(x,v) =V_{t,\tau}(x,v)\times B_{0}
   +E[f](t,X_{t,\tau}(x,v)),\\
 X_{\tau,\tau}(x,v)=x,\quad  V_{\tau,\tau}(x,v)=v,
  \end{array}\right.
           \end{align}
where $B=B_{0}+\delta B.$ In other words, we reduce inhomogeneous  dynamical system  to homogeneous dynamical system. Hence, based on the above dynamical system (0.7),
           we call the adopted  Newton iteration as the improved Newton iteration.

 Indeed,  there are many papers that contribute to Landau damping. Here
we  only list some results.  Bedrossian, Masmoudi, and Mouhot [4] provided a new, simple  and short proof of nonlinear Landau damping on $\mathbb{T}^{d}$ in only electric field case
 that nearly obtains the `` critical "  Gevrey-$\frac{1}{s}$ regularity predicted in [23]. Although their proofs have lots of
  the same ingredients as the proof in [23]
 from a physical point of view,
 at a mathematical level, the two proofs  are quite different, they ``mod out" by the characteristics of free transport and work in the coordinates $z=x-vt$ with
 $(t,x,v)\rightarrow(t,z,v).$ The evolution equation (0.1) $ (B=0)$ becomes $$\partial_{t}f+E(t,z+vt)\cdot (\nabla_{v}-t\nabla_{z})f+E(t,z+vt)
 \cdot\nabla_{v}f^{0}=0.$$  From this formula,  it is easy to see the phase mixing mechanism.  And this coordinate shift
 is related to the notion of ``gliding regularity" used in [23].
 One of the main ingredients of their proof is to split nonlinear terms into the transport structure term and ``reaction" term in [23] by using  paradifferential calculus .
  Bedrossian, Masmoudi [3] also used this method to prove the inviscid damping
and asymptotic stability of 2-D Euler equations  and later also proved the stability threshold for the 3D Couette flow in Sobolev regularity in [5]
 and so on. They [6] also proved Landau damping for the collisionless Vlasov equation with a class of $L^{1}$ interaction potentials on
 $\mathbb{R}_{x}^{3}\times\mathbb{R}_{v}^{3}$ for localized disturbances of infinite, inhomogeneous background.
  Also,  there are counterexamples that can be found in [2,10,13], showing that there is in
  general no exponential decay without analyticity and confining. Bedrossian stated one of these counterexamples by  proving that the theorem of
   Mouhot and Villani
  on Landau damping near equilibrium for the Vlasov-Poisson equations on $\mathbb{T}_{x}\times\mathbb{R}_{v}$ cannot, in general, be extended to Sobolev spaces
 by  constructing a sequence of homogeneous
  background distributions and arbitrarily small perturbations in $H^{s}$ which deviates arbitrary far from free transport for a long time.
  Lin and Zeng [10] also showed that
  there exist nontrivial traveling wave solutions to the Vlasov equation in Sobolev space $W^{s.p}_{x,v}$ $(p>1,s<1+\frac{1}{p})$ with arbitrary traveling speed.
  This implies that nonlinear Landau dampping is not true in $W^{s.p}_{x,v}$ $(p>1,s<1+\frac{1}{p})$ space for any homogeneous equilibria and in any period box.
   In addition, Deng and Masmoudi [13] showed the instability of  the Couette flow in low Gevrey spaces.
In recent years, there are also lots of results on the stability in other setting such as those in [12,14,15,28].

 This paper is organized as follows. Section 1 mainly introduces  hybrid analytic norms. In section 2, we will prove cyclotron damping at the linear level.
 We will state the new observation  and  sketch the proof of main theorem in section 3.   Section 4 is dedicated  to  the deflection estimates of the particles
 trajectory. Section 5 is the key section, which
 states the phenomena of plasma echoes. We will control the error terms in section  6, and give the iteration in section 7.

 Before stating our main theorem,   we assume that the electric induction is zero, then  the Maxwell's equations reduce to the following forms
  \begin{align}
 &\nabla\cdot E(t,x)=\int_{\mathbb{R}^{3}}f(t,x,v)dv,\quad \nabla\times B(t,x)=0,\quad \partial_{t} B(t,x)=-\nabla\times E(t,x),\quad \nabla\cdot B(t,x)=0.
 \end{align}
Now based on  the assumption that the electric induction is zero,  we first give two results of the Vlasov equation with
the electric field $ E(t,x)$ and the magnetic field $B(t,x) $ on  both the  linear  and the nonlinear levels   satisfying
 the conditions $ E= W(x)\ast\rho(t,x),\partial_{t}B=\nabla_{x}\times E,$ where $W(x)$ is a vector function and satisfies $|\widehat{W(k)}|\leq\frac{1}{1+|k|^{\gamma}}$.
%Now we rewrite  the  linear Vlasov equation under the above assumptions:
%\begin{equation}
%\left\{\begin{array}{l}
%\partial_{t}f+v\cdot\nabla_{x}f+\frac{q}{m}(v\times B_{0})\cdot\nabla_{v}f=-\frac{q}{m}(E+v\times B)\cdot\nabla_{v}f^{0},\\
%\partial_{t}B=\nabla_{x}\times E,\quad \nabla\cdot B=0,\\
%E=W(x)\ast\rho(t,x),\rho(t,x)=\int_{\mathbb{R}^{3}}f(t,x,v)dv,\\
%f(0,x,v)=f_{0}(x,v_{1},v_{2},v_{z}),f^{0}(v)=f^{0}(v_{1},v_{2},v_{z}),
%\end{array}\right.
 %\end{equation}
 %then the corresponding result is as follows.

 Now we state our main result as follows.
 \begin{thm} Let $f^{0}:\mathbb{R}^{3}\rightarrow\mathbb{R}_{+}$ be an analytic velocity profile,
 and assume $W(x)=(W_{1}(x),W_{2}(x),0):\mathbb{T}^{3}\rightarrow\mathbb{R}^{3}$ and $ W(x)$ is an odd function on $x_{3}$  satisfying $
  |\widehat{W}(k)|\leq\frac{1}{1+|k|^{\gamma}},\gamma>1.$ Further we assume that, for some constant $\lambda_{0}$ such that  $\lambda_{0}-B_{0}>0,$
 \begin{align}
 \sup_{\eta\in\mathbb{R}^{3}}e^{2\pi(\lambda_{0}-B_{0})|\eta|}|\tilde{f}^{0}(\eta)|\leq C^{0},\quad \sum_{n\in\mathbb{N}_{0}^{3}}\frac{(\lambda_{0}-B_{0})^{n}}{n!}\|\nabla^{n}_{v}f^{0}\|_{L_{dv}^{1}}\leq C_{0}<\infty.
 \end{align}

And we consider the following system,

\begin{equation}
\left\{\begin{array}{l}
\partial_{t}f+v\cdot\nabla_{x}f+\frac{q}{m}(v\times B_{0})\cdot\nabla_{v}f=-\frac{q}{m}(E+v\times B)\cdot\nabla_{v}f,\\
\partial_{t}B=\nabla_{x}\times E,\quad \nabla\cdot B=0,\\
E=W(x)\ast\rho(t,x),\rho(t,x)=\int_{\mathbb{R}^{3}}f(t,x,v)dv,\\
f(0,x,v)=f_{0}(x,v)=f_{0}(x,v_{1},v_{2},v_{z}),f^{0}(v)=f^{0}(v_{1},v_{2},v_{z}),
\end{array}\right.
 \end{equation}
 there is $\varepsilon=\varepsilon(\lambda_{0},\mu_{0},\beta,\gamma,\lambda'_{0},\mu'_{0})$ verifying the following property:
 $f_{0}=f_{0}(x,v)$ is an initial data such that
\begin{align}
\sup_{k\in\mathbb{Z}^{3},\eta\in\mathbb{R}^{3}}e^{2\pi(\lambda_{0}-B_{0})|\eta|}e^{2\pi\mu_{0}|k|}|f^{0}-f_{0}|+\int_{\mathbb{T}^{3}}\int_{\mathbb{R}^{3}}
|f^{0}-f_{0}|e^{\beta|v|}dvdx\leq\varepsilon,
\end{align}
where any $\beta>0,\lambda_{0}>\lambda'_{0}>B_{0},\mu_{0}>\mu'_{0}>0.$

In addition, we also assume that the $(\mathbf{PSC})$ holds.
%$\quad\mathbf{ how} \quad\mathbf{ to}\quad \mathbf{make} \quad \mathbf{assumptions }$
%$\quad \mathbf{on}\quad f^{0}$ is suitable
       %\item[(i)]$f^{0}(v)$ is  Maxwellian, that is, $f^{0}(v)=f^{0}(v_{\perp},v_{z})=C_{M}exp(-\alpha(v^{2}_{\perp}+v^{2}_{z}))$ for some $\alpha>0$  ;
         %\item[(ii)] $|\hat{f}^{0}(\eta)|\leq C^{0}e^{-\lambda_{1}|\eta_{3}|}e^{-\lambda_{1}|\eta_{1}|}e^{-\lambda_{1}|\eta_{2}|},|\partial_{\eta_{3}}\hat{f}^{0}(\eta)|\leq C^{1}(\lambda_{1})e^{-\lambda_{1}|\eta_{1}|}e^{-\lambda_{1}|\eta_{2}|}e^{-\lambda_{1}|\eta_{3}|},$
            % $|f^{0}(\eta_{1},\eta_{2},v_{3})|\leq Ce^{-\lambda_{1}|\eta_{1}|}e^{-\lambda_{1}|\eta_{2}|}e^{-\alpha_{1}|v_{3}|},$
             % for some constant $\lambda_{1},\alpha_{1}>0;$

        % \item[(iii)]$|\hat{f}_{0}(k,\eta)|\leq C_{0} e^{-\lambda_{0}|\eta_{1}|}e^{-\lambda_{0}|\eta_{2}|}e^{-\lambda_{0}|\eta_{3}|};$

                  %\item[(v)] $f^{0}(\eta_{1},\eta_{2},v_{3})\leq Ce^{-\lambda_{1}|\eta_{1}|}e^{-\lambda_{1}|\eta_{2}|}e^{-\alpha_{1}|v_{3}|},$ otherwise $f^{0}(\eta_{1},\eta_{2},v_{3})<\varepsilon(v_{3}),$ for some $0<\varepsilon(v_{3})<1$ sufficient small.%$\inf_{k\in\mathbb{Z}^{3}}|\tilde{\mathcal{L}}(\omega,k)-1|>\kappa,$ for some $0<\kappa<1.$
        Then %if $f_{0}(x,v),f^{0}(v)$ are axisymetric, then the solution $f$ of (1.1) is also axisymmetric.
         there exists a unique classical solution $(f(t,x,v),E(t,x),B(t,x))$  to the non-linear Vlasov system (0.12).

          Moreover,
for any fixed $\eta_{3},k_{3},\forall r\in\mathbb{N},$ as $|t|\rightarrow\infty,$ we have
         \begin{align}
         &|\hat{f}(t,k,\eta_{1}+k_{1}\frac{\sin\Omega t}{\Omega},\eta_{2}+k_{2}\frac{\sin\Omega t}{\Omega},\eta_{3})-\hat{f}_{0}(k,\eta)|\leq  e^{-(\lambda'_{0}-B_{0})|\eta_{3}+k_{3}t|},\quad\|\rho(t,\cdot)-\rho_{0}\|_{C^{r}(\mathbb{T}^{3})}
         =O(e^{-2\pi(\lambda'_{0}-B_{0})|t|}),\notag\\
          &\|E(t,\cdot)\|_{C^{r}(\mathbb{T}^{3})}=O(e^{-2\pi(\lambda'_{0}-B_{0})|t|}),\quad \|B(t,\cdot)\|
          _{C^{r}(\mathbb{T}^{3})}=O(e^{-2\pi(\lambda'_{0}-B_{0})|t|}),\notag\\
         \end{align}
         where $\rho_{0}=\int_{\mathbb{T}^{3}}\int_{\mathbb{R}^{3}}f_{0}(x,v)dvdx.$
         %(这里可能不需要假设初值为轴对称在速度方向 %maybe not assume $f_{0}$ axisymmetric in $v$)
         \end{thm}

         %Here we sketch the difficulties  and methods that the proof of Theorem 0.3 may meet. From Faraday Law of
         % Electromagnetic induction,  the force $F$ generated by the magnetic  field is that $F=BLV,$ where $L$ is a constant and $V$ is the velocity relative to the
         % magnetic field. Since $B(t,x)$ is independent of $V$ and $V$ is a unbounded variable, we almost have no hope to estimate $F,$ but this is the key point to
          %prove cyclotron damping. The idea to solve this problem  is from a new observation on the basis of Lenz's Law, then we reduce the inhomogeneous dynamical
         % equations of particles trajectory to the homogeneous equations  that is similar to that of Landau damping[]. Another difficulty is that the equation that
         % the density
          %$\rho[f]$ function satisfies  does not form the closed equation because of the additional term $v\times B\cdot\nabla_{v} f.$ For this, we have to estimate the
          %equations of the distribution equation and the density equation that bring  different kinds of resonances in different norms that is much more complicated
         % and difficult
          %than Landau damping in[].

Now we simply analyze the relation among the Vlasov-Poisson equations, the Vlasov-Maxwell equations,and our model. If we assume that both the electric induction and the
magnetic field are  zero, then  Vlasov-Maxwell equations reduce to the  Vlasov-Poisson equations; if we only  assume that the electric induction is zero,
then the  Vlasov-Maxwell equations reduce to our case. In other words, our case is the generalized case of the Vlasov-Poisson equations, and provides a new observation
 from
the physical viewpoint
 to solve the corresponding problem of the Vlasov-Maxwell equations. However, we cannot still solve
the Vlasov-Maxwell equations  completely only through this new observation. Therefore, the stability's or unstability's problem of the Vlasov-Maxwell equations is still open.

           In the following  we sketch the difficulties  and methods in our paper's setting. The crucial estimates of
            cyclotron damping  include the following two inequalities:

           $\bullet$ a control of $\rho=\int_{\mathbb{R}^{3}}fdv$ in $\mathcal{F}^{\lambda\tau+\mu}$ norm, that is,
            $\sup_{\tau\geq0}\|\rho_{\tau}\|_{\mathcal{F}^{\lambda\tau+\mu}}<\infty,$

            $\bullet$ a control of $f_{\tau}\circ\Omega_{t,\tau}$ in $\mathcal{Z}^{\lambda'(1+b),\mu';1}_{\tau-\frac{bt}{1+b}}$ norm, where
             $\lambda'<\lambda,\mu'<\mu.$

           However, during the iteration scheme, for cyclotron damping, from the view of the original Newton iteration, the
            characteristics are not only determined by the density $\rho^{n},$  but also related with  the velocity at  stage $n.$ However,   $\rho^{n}$
            is independent of the velocity and the key difficulty  is that
            the velocity is unbounded. This makes that we obtain the estimates of  the associated deflection $\Omega^{n}$ more difficult.  But this difficulty doesn't
            exist for Landau damping in [23].  To overcome this difficulty,
            on the basis of a new observation from Lenz's Law,
            we reduce the classical dynamical system to the improved dynamical system (0.7). And the corresponding  equation of the density $\rho[h^{n+1}]$ at stage
            $n+1$  becomes
            $$ \rho[h^{n+1}](t,x)=\int^{t}_{0}\int_{\mathbb{R}^{3}}-\bigg[(E[h^{n+1}]\circ\Omega^{n}_{s,t}(x,v)\cdot G_{s,t}^{n})-
            (B[h^{n+1}]\circ\Omega^{n}_{s,t}(x,v)
            \cdot G_{s,t}^{n,v})$$
           $$-(B[f^{n}]\circ\Omega^{n}_{s,t}(x,v))\cdot(\nabla'_{v}h^{n+1}\times V^{0}_{s,t}(x,v))\circ\Omega^{n}_{s,t}(x,v)\bigg](s, X^{0}_{s,t}(x,v),V^{0}_{s,t}(x,v))dvds
            +(\textmd{terms}\quad \textmd{from}\quad \textmd{stage}\quad n),$$
            From the above equation, we see that, comparing with  that  in [23], there is a new term $(B[f^{n}]\circ\Omega^{n}_{s,t}(x,v))\cdot$
             $
            (\nabla'_{v}h^{n+1}\times V^{0}_{s,t}(x,v)\circ\Omega^{n}_{s,t}(x,v))$  that have the information of the stage $n+1.$ Of course,
            this is due to the reason that we regard the perturbation of  the magnetic field as a negligible term.
            To get a self-consistent estimate,
             we have to deal with this term and  have little choice but to come back the
            equation of $h^{n+1}$. This leads to  different kinds of resonances (in term of different norms), for example,
            $v\times B[f^{n}]\cdot\nabla_{v}h^{n+1}$ may generate resonance in $\mathcal{Z}^{\lambda,\mu;1}_{\tau}$ norm  on $h^{n+1},$
            except in $\mathcal{F}^{\lambda\tau+\mu}$ norm  on $\rho[h^{n+1}],$
            because  both $B[f^{n}]$ and $h^{n+1}$ contain the space variable.  But there are no these problems in [23], and Landau damping in [23]
            only has the resonances in  $\mathcal{F}^{\lambda\tau+\mu}$ norm  on $\rho[h^{n+1}].$

          \begin{rem} $\gamma>1$ of Theorems 0.1 and 0.3 can be extended to $\gamma\geq1,$ the difference between $\gamma>1$ and $\gamma=1$ is the proof of the growth
          integral in section 7. The proof of $\gamma=1$ is similar to section 7 in [23], here we omit this case.
          \end{rem}
          \begin{rem} Indeed, we don't need to assume that $ W(x)=(W_{1}(x),W_{2}(x),0)$ and
            $B_{0}=(0,0,B_{0}),$ only
            need to assume that the electric field $E(t,x)$ is perpendicular to the background magnetic field $B_{0}$.
          \end{rem}
         \begin{rem}  From the physics viewpoint, the condition of the damping  is that the number of particles that the wave velocity greatly
         exceeds  their velocity
         is much larger than the number of particles whose velocity is slower than the wave velocity.
          The $(\mathbf{PSC})$ condition  in above theorem is in consistent with  the  statement of the physics viewpoint. In fact,
           the  $(\mathbf{PSC})$ condition and (0.11)
           imply that the number of resonant particles  is exponentially small and their effect corresponding is weak.
            From the conclusion of  Theorem (0.3), we know that,  under this condition, most particles absorb energy from the wave, and then the damping occurs.
         \end{rem}
         \begin{rem} From the definition of the hybrid analytic norms and the proof of the following sections, there is
         still the phenomena of cyclotron damping for the nonlinear Vlasov-Maxwell equations,  which is the same
         to that in the nonlinear Vlasov-Poisson equations  but only in $\hat{z}$ direction, the position of the corresponding resonances translates 0 into $B_{0},$ and in the horizon direction,
         the action of the particles moves along circle that is the same to the linear case.
         \end{rem}

         \begin{center}
\item\section{Linear Cyclotron Damping}
%{\bf\large 1. \quad Introduction }
\end{center}

In this section, let us consider the following linear Vlasov equations in a uniform magnetic field, and recall the equations:
\begin{equation}
\left\{\begin{array}{l}
\partial_{t}f+v\cdot\nabla_{x}f+\frac{q}{m}(v\times B_{0})\cdot\nabla_{v}f=-\frac{q}{m}(E+v\times B)\cdot\nabla_{v}f^{0},\\
\partial_{t}B=\nabla_{x}\times E,\quad \nabla\cdot B=0,\\
E=W(x)\ast\rho(t,x),\rho(t,x)=\int_{\mathbb{R}^{3}}f(t,x,v)dv,\\
f(0,x,v)=f_{0}(x,v_{1},v_{2},v_{z}),f^{0}(v)=f^{0}(v_{1},v_{2},v_{z}),
\end{array}\right.
 \end{equation}
where the distribution function  $f=f(t,x,v):\mathbb{R}^{+}\times\mathbb{T}^{3}\times\mathbb{R}^{3}\rightarrow\mathbb{R},W(x)=(W_{1}(x),W_{2}(x),0):
\mathbb{T}^{3}\rightarrow\mathbb{T}^{3},$ $ B_{0}$ is a constant magnetic field along the $\hat{z}$ direction, $E=E(t,x)$ is the electric field,
 $ B_{0}+B$ is the magnetic field. %And in this paper we consider a simple case in which the thermal energy of particles perpendicular to the magnetic field is zero and  the velocity of particles parallel to the magnetic field $B_{0}$ is $ V, $ that is, $V=V\hat{z},v\times B_{1}=V\times B_{1}.$

%\begin{align}
%&\mathcal{L}(t,k)=\frac{q}{m}\int^{t}_{0}k_{3}(\hat{W}_{2}(i\eta_{k_{1}})\partial_{\eta_{3}}
%\hat{f}^{0}(\eta_{k_{1}},\eta_{k_{2}},k_{3}t)+\hat{W}_{1}(i\eta_{k_{2}})\partial_{\eta_{3}}\hat{f}^{0}(\eta_{k_{1}},\eta_{k_{2}},k_{3}t))dt\notag\\
%&-\frac{q}{m}
%(\hat{W}_{1}(i\eta_{k_{1}})\hat{f}^{0}(\eta_{k_{1}},\eta_{k_{2}},k_{3}t)+\hat{W}_{2}(i\eta_{k_{2}})\hat{f}^{0}(\eta_{k_{1}},\eta_{k_{2}},k_{3}t)),\notag\\
%\end{align}
%where $\eta_{k_{1}}=\frac{1}{\Omega}(-k_{2}\cos\Omega t+k_{2}-k_{1}\sin\Omega t),\eta_{k_{2}}=\frac{1}{\Omega}(-k_{2}\sin\Omega t-k_{1}+k_{1}\cos\Omega t).$

\begin{thm}
 For any $\eta,v\in\mathbb{R}^{3},k\in\mathbb{N}_{0}^{3},$  we assume that the following conditions hold in equations (0.9).
%$\quad\mathbf{ how} \quad\mathbf{ to}\quad \mathbf{make} \quad \mathbf{assumptions }$
%$\quad \mathbf{on}\quad f^{0}$ is suitable
\begin{itemize}
        \item[(i)] $ W(x)$ is an odd function on $x_{3},$  $|\widehat{W}(k)|\leq\frac{1}{1+|k|^{\gamma}},\gamma>1,$ where $ W(x)=(W_{1}(x),W_{2}(x),0);$
%\item[(i)]$f^{0}(v)$ is  Maxwellian, that is, $f^{0}(v)=f^{0}(v_{\perp},v_{z})=C_{M}exp(-\alpha(v^{2}_{\perp}+v^{2}_{z}))$ for some $\alpha>0$  ;
         \item[(ii)] $|\hat{f}^{0}(\eta)|\leq C^{0}e^{-2\pi\lambda_{0}|\eta|},|\partial_{\eta_{3}}\hat{f}^{0}(\eta)|\leq C^{0}e^{-2\pi\lambda_{0}|\eta|},$
             $|f^{0}(\cdot,v_{3})|\leq C^{0}e^{-2\pi\alpha_{0}|v_{3}|},$
              for some constants $\lambda_{0},\alpha_{0},C^{0}>0;$

         \item[(iii)]$|\hat{f}_{0}(k,\eta)|\leq C_{0} e^{-2\pi\lambda_{0}|\eta|}$ for some constant $C_{0}>0,$ where $\lambda_{0}$ is defined in (ii);

             \item[(iv)] In addition, $(\mathbf{PSC})$ holds,

             $(\mathbf{PSC}):$ for any component velocity in the $\hat{z}$ direction $v_{3}\in\mathbb{R},$
             there exists some positive constant $v_{Te}$ such that if $v_{3}=\frac{\omega}{k_{3}} \textmd{ when } k_{3}\neq0; \textmd{ or } k_{3}=0$ where $\omega,k$ are frequencies of time and space $t,x,$
             respectively,  then $|v_{3}|\gg v_{Te}.$
                  %\item[(v)] $f^{0}(\eta_{1},\eta_{2},v_{3})\leq Ce^{-\lambda_{1}|\eta_{1}|}e^{-\lambda_{1}|\eta_{2}|}e^{-\alpha_{1}|v_{3}|},$ otherwise $f^{0}(\eta_{1},\eta_{2},v_{3})<\varepsilon(v_{3}),$ for some $0<\varepsilon(v_{3})<1$ sufficient small.%$\inf_{k\in\mathbb{Z}^{3}}|\tilde{\mathcal{L}}(\omega,k)-1|>\kappa,$ for some $0<\kappa<1.$
         \end{itemize}
         Then %if $f_{0}(x,v),f^{0}(v)$ are axisymetric, then the solution $f$ of (1.1) is also axisymmetric.
for any fixed $\eta_{3},k_{3},$ and  for any $\lambda'_{0}<\lambda_{0},$ we have
         \begin{align}
         &|\hat{f}(t,k,\eta)-\hat{f}_{0}(k,\eta)|\leq  e^{-2\pi\lambda'_{0}|\eta_{3}+k_{3}t|},\quad|\hat{\rho}(t,k)-\hat{\rho}_{0}|\leq  e^{-2\pi\lambda'_{0}|\eta_{k1}|}e^{-2\pi\lambda'_{0}|\eta_{k2}|} e^{-\lambda'_{0}|k_{3}|t},\notag\\
        & |\hat{E}(t,k)|\leq e^{-2\pi\lambda'_{0}|\eta_{k1}|}e^{-2\pi\lambda'_{0}|\eta_{k2}|}e^{-2\pi\lambda'_{0}|k_{3}|t},|\hat{B}(t,k)|\leq te^{-2\pi\lambda'_{0}|\eta_{k1}|}e^{-2\pi\lambda'_{0}|\eta_{k2}|}e^{-2\pi\lambda'_{0}|k_{3}|t}.
         \end{align}
        where $\rho_{0}=\int_{\mathbb{T}^{3}}\int_{\mathbb{R}^{3}}f_{0}(x,v)dvdx,\eta_{k1}=\frac{1}{\Omega}(-k_{2}\cos\Omega t+k_{2}-k_{1}\sin\Omega t),\eta_{k2}=\frac{1}{\Omega}(-k_{2}\sin\Omega t-k_{1}+k_{1}\cos\Omega t).$
         %(这里可能不需要假设初值为轴对称在速度方向 %maybe not assume $f_{0}$ axisymmetric in $v$)
         \end{thm}

\begin{rem} In the linear case, from (1.2) and (1.24)-(1.26), it is easy to observe that  cyclotron damping is almost the  same with Landau damping
 when $B_{0}$ tends to zero. However, when $B_{0}$ is fixed, in the horizon direction there is no damping and the motion of particles is a circle;
  in the $\hat{z}$ direction, the damping still occurs. Therefore, there is no immediate tendency toward trapping. This is a crucial point for cyclotron damping.
    And in the latter case, the motion of plasma particles moves along spiral trajectories.
\end{rem}

         Before proving $Theorem$ 1.1, we give a key lemma.
         \begin{lem} Under the assumptions of Theorem 0.1, let $\Phi(t,k)=e^{2\pi\lambda'_{0}|k_{3}|t}$
         $\cdot e^{2\pi\lambda'_{0}|\eta_{k1}|}e^{2\pi\lambda'_{0}|\eta_{k2}|}\hat{\rho}(t,k),$ $A(t,k)=\int_{\mathbb{T}^{3}}\int_{\mathbb{R}^{3}}
         e^{-2\pi ix\cdot k}f_{0}(x'(0,x,v),v'(0,x,v))dv'dx,$ where $\lambda_{0}>\lambda'_{0}>0,\lambda_{0}$ is defined in Theorem 0.1, we have
         \begin{align}
\|\Phi\|_{L^{\infty}(dt)}\leq\bigg(1+\frac{C(k,W,\Omega)}{(\lambda_{0}-\lambda'_{0})^{\frac{3}{2}}}\bigg)\|e^{\lambda_{0}\nu_{k}}A\|_{L^{\infty}(dt)},
\end{align}
where $\eta_{k1}=\frac{1}{\Omega}(-k_{2}\cos\Omega t+k_{2}-k_{1}\sin\Omega t),\eta_{k2}=\frac{1}{\Omega}(-k_{2}\sin\Omega t-k_{1}+k_{1}\cos\Omega t), \nu_{k}=|\eta_{k1}|+|\eta_{k2}|+|k_{3}|t.$
         \end{lem}
          $Proof.$  First, we consider the case that $\omega\neq0,k_{3}\neq0,$
\begin{align}
&\tilde{\rho}(\omega,k)=\int_{\mathbb{R}^{+}}\int_{\mathbb{T}^{3}}\int_{\mathbb{R}^{3}}e^{2\pi it\omega}e^{-2\pi ix\cdot k}f(t,x,v)dxdvdt\notag\\
&=\int_{\mathbb{R}^{+}}\int_{\mathbb{T}^{3}}\int_{\mathbb{R}^{3}}e^{2\pi it\omega}e^{-2\pi ix\cdot k}f_{0}(x'(0,x,v),v'(0,x,v))dvdxdt-\frac{q}{m}\int_{\mathbb{R}^{+}}\int_{\mathbb{T}^{3}}\int^{t}_{0}\int_{\mathbb{R}^{3}}\notag\\
&\cdot e^{2\pi it\omega}e^{-2\pi ix\cdot k}[(E+v'(\tau,x,v)\times B)\cdot\nabla'_{v}f^{0}](\tau,x'(\tau,x,v),v'(\tau,x,v))dvdxd\tau dt.\notag\\
\end{align}
Recall $E(t,x)=W(x)\ast\rho(t,x),$ $\partial_{t}B(t,x)=\nabla_{x}\times E(t,x),$ then taking the Fourier transform in the variables $t,x,$
$$\widetilde{E}(\omega,k)=\widehat{W}(k)\tilde{\rho}(\omega,k),\quad\omega\widetilde{B}(\omega,k)=k\times\widetilde{E}(\omega,k).$$
Furthermore, we get
\begin{align}
&v\times\widetilde{B}(\omega,k)=\frac{1}{\omega}[v\times(k\times\widetilde{E}(\omega,k))]\notag\\
&=\frac{1}{\omega}\bigg(v_{2}(k_{1}\widehat{W}_{2}-k_{2}\widehat{W}_{1})-v_{3}k_{3}\widehat{W}_{1},-v_{1}(k_{1}
\widehat{W}_{2}-k_{2}\widehat{W}_{1})+v_{3}k_{3}\widehat{W}_{2},
v_{1}k_{3}\widehat{W}_{1}+v_{2}k_{3}\widehat{W}_{2}\bigg)\tilde{\rho}(\omega,k).\notag\\
\end{align}
Combining (1.3)-(1.4) and (2.3)-(2.4), and note that $dv\rightarrow dv'$ preserves the measure,  we can change between $dv$ and $dv'$ whenever we need, but in order to
simply the notations, we don't differentiate the notations $dv$ and $dv'$  in this paper,
so we have

$$\tilde{\rho}(\omega,k)=\int_{\mathbb{R}^{+}}\int_{\mathbb{T}^{3}}\int_{\mathbb{R}^{3}}e^{2\pi it\omega}e^{-2\pi ix\cdot k}f_{0}(x'(0,x,v),v'(0,x,v))dvdxdt$$
$$+\frac{q}{m}\tilde{\rho}(\omega,k)\frac{1}{\omega}\int_{\mathbb{R}^{+}}\int_{\mathbb{R}^{3}}e^{2\pi it\omega}e^{-2\pi ik_{3}v_{3}t}
e^{-2\pi i\eta_{k2}v'_{2}}e^{-2\pi i\eta_{k1}v'_{1}}(k_{3}v_{3})(\widehat{W}_{1}\partial_{v'_{1}}f^{0}-\widehat{W}_{2}\partial_{v'_{2}}f^{0})dvdt$$
$$+\frac{q}{m}\tilde{\rho}(\omega,k)\frac{1}{\omega}\int_{\mathbb{R}^{+}}\int_{\mathbb{R}^{3}}e^{2\pi it\omega}e^{-2\pi ik_{3}v_{3}t}
e^{-2\pi i\eta_{k2}v'_{2}}e^{-2\pi i\eta_{k1}v'_{1}}(k_{1}\widehat{W}_{2}-k_{2}\widehat{W}_{1})$$
$$\cdot(v_{2}\partial_{v'_{1}}f^{0}-v_{1}\partial_{v'_{2}}f^{0})dvdt-\frac{q}{m}\tilde{\rho}(\omega,k)\frac{1}{\omega}\int_{\mathbb{R}^{+}}
\int_{\mathbb{R}^{3}}e^{2\pi it\omega}e^{-2\pi ik_{3}v_{3}t}
e^{-2\pi i\eta_{k2}v'_{2}}e^{-2\pi i\eta_{k1}v'_{1}}$$
\begin{align}
&\cdot(v_{1}\widehat{W}_{1}+v_{2}\widehat{W}_{2})k_{3}\partial_{v'_{3}}f^{0}
dvdt-\frac{q}{m}\tilde{\rho}(\omega,k)\int_{\mathbb{R}^{+}}\int_{\mathbb{R}^{3}}e^{2\pi it\omega}e^{-2\pi ik_{3}v_{3}t}
e^{-2\pi i\eta_{k2}v'_{2}}e^{-2\pi i\eta_{k1}v'_{1}}\notag\\
&\cdot(\widehat{W}_{1}\partial_{v'_{1}}f^{0}+\widehat{W}_{2}\partial_{v'_{2}}f^{0})dvdt,\notag\\
\end{align}
where $\eta_{k_{1}}=\frac{1}{\Omega}(-k_{2}\cos\Omega t+k_{2}-k_{1}\sin\Omega t),$ $\eta_{k_{2}}=\frac{1}{\Omega}(-k_{2}\sin\Omega t-k_{1}+k_{1}\cos\Omega t).$

Let
%and by the definition of $\mathcal{L},$ it is easy to check
\begin{align}
&\tilde{\mathcal{L}}(\omega,k)=\frac{q}{m}\frac{1}{\omega}\int_{\mathbb{R}^{+}}\int_{\mathbb{R}^{3}}e^{2\pi it\omega}e^{-2\pi ik_{3}v_{3}t}
e^{-2\pi i\eta_{k2}v'_{2}}e^{-2\pi i\eta_{k1}v'_{1}}(k_{3}v_{3})(\widehat{W}_{1}\partial_{v'_{1}}f^{0}-\widehat{W}_{2}\partial_{v'_{2}}f^{0})dvdt\notag\\
&-\frac{q}{m}\int_{\mathbb{R}^{+}}\int_{\mathbb{R}^{3}}e^{2\pi it\omega}e^{-2\pi ik_{3}v_{3}t}
(\widehat{W}_{1}\partial_{v'_{1}}f^{0}+\widehat{W}_{2}\partial_{v'_{2}}f^{0})\cdot e^{-2\pi i\eta_{k2}v'_{2}}e^{-2\pi i\eta_{k1}v'_{1}}dvdt\notag\\
&-\frac{q}{m}\frac{1}{\omega}\int_{\mathbb{R}^{+}}\int_{\mathbb{R}^{3}}e^{2\pi it\omega}e^{-2\pi ik_{3}v_{3}t}
e^{-2\pi i\eta_{k2}v'_{2}}e^{-2\pi i\eta_{k1}v'_{1}}(v_{1}\widehat{W}_{1}+v_{2}\widehat{W}_{2})k_{3}\partial_{v'_{3}}f^{0}
dvdt\notag\\
&+\frac{q}{m}\frac{1}{\omega}\int_{\mathbb{R}^{+}}\int_{\mathbb{R}^{3}}e^{2\pi it\omega}e^{-2\pi ik_{3}v_{3}t}
e^{-2\pi i\eta_{k2}v'_{2}}e^{-2\pi i\eta_{k1}v'_{1}}(k_{1}\widehat{W}_{2}-k_{2}\widehat{W}_{1})(v_{2}\partial_{v'_{1}}f^{0}
-v_{1}\partial_{v'_{2}}f^{0})dvdt,\notag\\
\end{align}
hence
\begin{align}
\tilde{\rho}(\omega,k)=\tilde{A}(\omega,k)+\tilde{\rho}(\omega,k)\tilde{\mathcal{L}}(\omega,k).
\end{align}
%Assume $\inf_{k\in\mathbb{Z}^{3}}|\tilde{\mathcal{L}}(\omega,k)-1|>\kappa,$ for some constant $0<\kappa<1,$
Then taking the inverse Fourier transform in time $ t, $ we get
$\hat{\rho}(t,k)=\hat{A}(t,k)+\hat{\rho}(t,k)\ast\hat{\mathcal{L}}(t,k),$
and

\begin{align}
&e^{2\pi\lambda'_{0}\eta_{k1}}e^{2\pi\lambda'_{0}\eta_{k2}}e^{2\pi\lambda'_{0}|k_{3}|t}\hat{\rho}(t,k)=e^{2\pi\lambda'_{0}\eta_{k1}}e^{2\pi\lambda'_{0}\eta_{k2}}
e^{2\pi\lambda'_{0}|k_{3}|t}\hat{A}(t,k)+e^{2\pi\lambda'_{0}\eta_{k1}}e^{2\pi\lambda'_{0}\eta_{k2}}e^{2\pi\lambda'_{0}|k_{3}|t}\hat{\rho}(t,k)\ast\hat{\mathcal{L}}(t,k).
\end{align}
Let $\Phi(t,k)=e^{2\pi\lambda'_{0}\eta_{k1}}e^{2\pi\lambda'_{0}\eta_{k2}}e^{2\pi\lambda'_{0}|k_{3}|t}\hat{\rho}(t,k),$
$\mathcal{A}(t,k)=e^{2\pi\lambda'_{0}\eta_{k1}}e^{2\pi\lambda'_{0}\eta_{k2}}
e^{2\pi\lambda'_{0}|k_{3}|t}\hat{A}(t,k),$
$\mathcal{K}^{0}(t,k)$
$=e^{2\pi\lambda'_{0}\eta_{k1}}$
$e^{2\pi\lambda'_{0}\eta_{k2}}e^{2\pi\lambda'_{0}|k_{3}|t}
\hat{\mathcal{L}}(t,k),$ then from (1.8), we have
$\widetilde{\Phi}(\omega,k)=\widetilde{\mathcal{A}}(\omega,k)+\widetilde{\Phi}(\omega,k)\widetilde{\mathcal{K}}^{0}(\omega,k).$

Then
\begin{align}
&\|\Phi(t,k)\|_{L^{2}(dt)}
=\|\widetilde{\Phi}(\omega,k)\|_{L^{2}}\leq\|\widetilde{\mathcal{A}}(\omega,k)\|_{L^{2}}+\|\widetilde{\Phi}(\omega,k)\|_{L^{2}}
\|\widetilde{\mathcal{K}}^{0}(\omega,k)\|_{L^{\infty}}\notag\\
&\leq\|e^{2\pi\lambda'_{0}\nu_{k}}\hat{A}(t,k)\|_{L^{2}}+\|e^{2\pi\lambda'_{0}\nu_{k}}\hat{\rho}(t,k)\|_{L^{2}}
\|\widetilde{\mathcal{K}}^{0}(\omega,k)\|_{L^{\infty}},\notag\\
\end{align}
where $\nu_{k}=|\eta_{k1}|+|\eta_{k2}|+|k_{3}|t.$

Next we have to estimate $\|\widetilde{\mathcal{K}}^{0}(\omega,k)\|_{L^{\infty}}.$

Indeed,

$$\|\widetilde{\mathcal{K}}^{0}(\omega,k)\|_{L^{\infty}}
\leq\sup_{\omega}\frac{q}{m}\bigg[\bigg|\frac{1}{\omega}\int_{\mathbb{R}^{+}}\int_{\mathbb{R}^{3}}e^{2\pi it\omega}e^{-2\pi ik_{3}v_{3}t}e^{2\pi\lambda'_{0}\nu_{k}}
e^{-2\pi i\eta_{k2}v'_{2}}e^{-2\pi i\eta_{k1}v'_{1}}(k_{3}v_{3})\cdot(\widehat{W}_{1}\partial_{v'_{1}}f^{0}$$
$$-\widehat{W}_{2}\partial_{v'_{2}}f^{0})
dvdt\bigg|+\frac{q}{m}\bigg|\int_{\mathbb{R}^{+}}\int_{\mathbb{R}^{3}}e^{2\pi it\omega}
e^{-2\pi ik_{3}v_{3}t}e^{2\pi\lambda_{0}\nu_{k}}\cdot(\widehat{W}_{1}\partial_{v'_{1}}f^{0}+\widehat{W}_{2}\partial_{v'_{2}}f^{0})\cdot e^{-2\pi i\eta_{k2}v'_{2}}$$
$$e^{-2\pi i\eta_{k1}v'_{1}}
dvdt\bigg|-\frac{q}{m}(\omega,k)\frac{1}{\omega}\int_{\mathbb{R}^{+}}\bigg|\int_{\mathbb{R}^{3}}e^{2\pi it\omega}
e^{2\pi\lambda_{0}\nu_{k}}e^{-2\pi ik_{3}v_{3}t}e^{-2\pi i\eta_{k2}v'_{2}}e^{-2\pi i\eta_{k1}v'_{1}}\cdot(v_{1}\widehat{W}_{1}+v_{2}\widehat{W}_{2})k_{3}$$
\begin{align}
&\partial_{v'_{3}}f^{0}
dv\bigg|dt+\frac{q}{m}(\omega,k)\frac{1}{\omega}\int_{\mathbb{R}^{+}}\bigg|\int_{\mathbb{R}^{3}}e^{2\pi it\omega}e^{-2\pi ik_{3}v_{3}t}e^{-2\pi i\eta_{k2}v'_{2}}e^{2\pi \lambda_{0}\nu_{k}}e^{-2\pi i\eta_{k1}v'_{1}}
\cdot(k_{1}\widehat{W}_{2}-k_{2}\widehat{W}_{1})\notag\\
&(v_{2}\partial_{v'_{1}}f^{0}
-v_{1}\partial_{v'_{2}}f^{0})dv\bigg|dt\bigg]=I+II+III+IV.\notag\\
\end{align}
%\begin{align}
%&=\frac{q}{m}\int^{\infty}_{0}|\frac{k_{3}}{\omega}\hat{W}_{2}\int_{\mathbb{R}^{+}}e^{it\omega}(i\eta_{k1})e^{-ik_{3}v_{3}t}e^{\lambda|k_{3}|t}
%\widehat{v_{3}f^{0}}(\eta_{k1},\eta_{k2},v_{3})dv_{3}dt|d\omega\notag\\
%&+\frac{q}{m}\int^{\infty}_{0}|\frac{k_{3}}{\omega}\hat{W}_{1}\int_{\mathbb{R}^{+}}e^{it\omega}(i\eta_{k2})e^{-ik_{3}v_{3}t}e^{\lambda|k_{3}|t}
%\widehat{v_{3}f^{0}}(\eta_{k1},\eta_{k2},v_{3})dv_{3}dt|d\omega\notag\\
%&+\frac{q}{m}\int^{\infty}_{0}|\hat{W}_{1}\int_{\mathbb{R}^{+}}e^{it\omega}(i\eta_{k1})e^{-ik_{3}v_{3}t}e^{\lambda|k_{3}|t}\widehat{f^{0}}(\eta_{k1}
%,\eta_{k2},v_{3})dv_{3}dt|d\omega\notag\\
%&+\frac{q}{m}\int^{\infty}_{0}|\hat{W}_{2}\int_{\mathbb{R}^{+}}e^{it\omega}(i\eta_{k2})e^{-ik_{3}v_{3}t}e^{\lambda|k_{3}|t}\widehat{f^{0}}(\eta_{k1},\eta_{k2},
%v_{3})dv_{3}dt|d\omega\notag\\
%&+\cdots\cdots\notag\\
%&=I+II+III+IV+\cdots
%&\leq Ce^{-\alpha(\frac{\omega}{k_{3}})^{2}}+\varepsilon\leq C_{Te}e^{-\alpha(v_{Te})^{2}},
%\end{align}
%where $\eta_{k_{1}}=\frac{1}{\Omega}(-k_{2}\cos\Omega t+k_{2}-k_{1}\sin\Omega t),\eta_{k_{2}}=\frac{1}{\Omega}(-k_{2}\sin\Omega t-k_{1}+k_{1}\cos\Omega t).$

In fact, we only need to estimate one term of (1.10) because of similar processes of  other terms. Without loss of generality, we give an estimate for $ I.$ In the same way, we only estimate one term of $I,$ here we still denote $ I. $
$$ I=\sup_{\omega}\frac{q}{m}\bigg|\frac{k_{3}}{\omega}\widehat{W}_{2}\int_{\mathbb{R}^{+}}\int_{\mathbb{R}}e^{2\pi it\omega}(2\pi i\eta_{k1})
e^{2\pi\lambda_{0}|\eta_{k1}|}e^{2\pi\lambda_{0}|\eta_{k2}|}e^{-2\pi ik_{3}v_{3}t}e^{2\pi\lambda_{0}|k_{3}|t}\cdot\widehat{v_{3}f^{0}}(\eta_{k1},\eta_{k2},v_{3})
dv_{3}dt\bigg|$$
$$=\sup_{\omega}\frac{q}{m}\bigg|\frac{k_{3}}{\omega}\widehat{W}_{2}\int_{\mathbb{R}^{+}}\int_{\mathbb{R}}e^{2\pi it\omega}(2\pi i\eta_{k1})e^{2\pi\lambda_{0}|\eta_{k1}|}e^{2\pi\lambda_{0}|\eta_{k2}|}
e^{-2\pi ik_{3}v_{3}t}\cdot\sum_{n}\frac{|2\pi i\lambda_{0}|k_{3}|t|^{n}}{n!}
\widehat{v_{3}f^{0}}(\eta_{k1},\eta_{k2},v_{3})dv_{3}dt\bigg|$$
%&=\frac{q}{m}\sum_{n}\frac{\lambda^{n}}{n!}\int^{\infty}_{0}|\frac{k_{3}}{\omega}\hat{W}_{2}\int_{\mathbb{R}^{+}}e^{it\omega}(i\eta_{k1})(ik_{3}t)^{n}
%e^{-ik_{3}v_{3}t}
%\widehat{v_{3}f^{0}}(\eta_{k1},\eta_{k2},v_{3})dv_{3}dt|d\omega\notag\\
%$$=\frac{q}{m}\sum_{n}\frac{\lambda^{n}}{n!}\int^{\infty}_{0}\bigg|\frac{k_{3}}{\omega}\hat{W}_{2}\int_{\mathbb{R}^{+}}\int_{\mathbb{R}}(2\pi i\eta_{k1})(2\pi ik_{3}t)^{n}
%e^{2\pi\lambda_{0}|\eta_{k1}|}e^{2\pi\lambda_{0}|\eta_{k2}|}e^{2\pi ik_{3}t(\frac{\omega}{k_{3}}-v_{3})}\cdot\widehat{v_{3}f^{0}}(\eta_{k1},\eta_{k2},v_{3})dv_{3}dt\bigg|d\omega$$
%$$=\frac{q}{m}\sum_{n}\frac{\lambda^{n}}{n!}\int^{\infty}_{0}\bigg|\frac{k_{3}}{\omega}\hat{W}_{2}\int_{\mathbb{R}^{+}}\int_{\mathbb{R}}(2\pi i\eta_{k1})
%(-1)^{n}e^{2\pi\lambda_{0}|\eta_{k1}|}e^{2\pi\lambda_{0}|\eta_{k2}|}\nabla_{v_{3}}^{n}e^{2\pi ik_{3}t(\frac{\omega}{k_{3}}-v_{3})}\cdot\widehat{v_{3}f^{0}}(\eta_{k1},\eta_{k2},v_{3})dv_{3}dt\bigg|d\omega$$
$$=\sup_{\omega}\frac{q}{m}\sum_{n}\frac{\lambda^{n}}{n!}\bigg|\frac{k_{3}}{\omega}\widehat{W}_{2}\int_{\mathbb{R}^{+}}\int_{\mathbb{R}}(2\pi i\eta_{k1})
(-1)^{n}e^{2\pi\lambda_{0}|\eta_{k1}|}e^{2\pi\lambda_{0}|\eta_{k2}|}e^{2\pi ik_{3}t(\frac{\omega}{k_{3}}-v_{3})}
\cdot\nabla_{v_{3}}^{n}\widehat{v_{3}f^{0}}(\eta_{k1},\eta_{k2},v_{3})dv_{3} dt\bigg|$$
%$$=\frac{q}{m}\sum_{n}\frac{\lambda_{0}^{n}}{n!}\int^{\infty}_{0}\bigg|\frac{k_{3}}{\omega}\hat{W}_{2}\int_{\mathbb{R}^{+}}(2\pi i\eta_{k1})
%(-1)^{n}e^{2\pi\lambda_{0}|\eta_{k1}|}e^{2\pi\lambda_{0}|\eta_{k2}|}e^{2\pi ik_{3}t(\frac{\omega}{k_{3}})}
%(2\pi ik_{3}t)^{n}\cdot\widehat{v_{3}f^{0}}(\eta_{k1},\eta_{k2},k_{3}t) dt\bigg|d\omega$$
$$=\sup_{\omega}\frac{q}{m}\sum_{n}\frac{\lambda_{0}^{n}}{n!}\bigg|\frac{k_{3}}{\omega}\widehat{W}_{2}(2\pi i\eta_{k1})
e^{2\pi\lambda_{0}|\eta_{k1}|}e^{2\pi\lambda_{0}|\eta_{k2}|}
(-i\nabla_{\frac{\omega}{k_{3}}})^{n}\widehat{v_{3}f^{0}}(\eta_{k1},\eta_{k2},\frac{\omega}{k_{3}})\bigg |\leq\frac{q}{mv_{Te}}e^{-c_{0}v_{Te}},$$
%\begin{align}
%&\leq\frac{q}{m}\sum_{n}\frac{\lambda_{0}^{n}}{n!}\int^{\infty}_{v_{Te}}\bigg|\frac{k_{3}}{\omega}\hat{W}_{2}(2\pi i\eta_{k1})
%e^{2\pi\lambda_{0}|\eta_{k1}|}e^{2\pi\lambda_{0}|\eta_{k2}|}
%(-i\nabla_{\frac{\omega}{k_{3}}})^{n}\widehat{v_{3}f^{0}}(\eta_{k1},\eta_{k2},\frac{\omega}{k_{3}}) \bigg|d\omega\leq\frac{q}{mv_{Te}}e^{-c_{0}v_{Te}}\notag\\
%&+\frac{q}{m}\sum_{n}\frac{\lambda^{n}}{n!}\int^{v_{Te}}_{0}|\frac{k_{3}}{\omega}\hat{W}_{2}(i\eta_{k1})
%e^{\lambda|\eta_{k1}|}e^{\lambda|\eta_{k2}|}
%(-i\nabla_{\frac{\omega}{k_{3}}})^{n}\widehat{v_{3}f^{0}}(\eta_{k1},\eta_{k2},\frac{\omega}{k_{3}}) |d\omega\leq\frac{q}{mv_{Te}}e^{-c_{0}v_{Te}},
%\end{align}
where in the last inequality we use the facts that if $v_{3}=\frac{\omega}{k_{3}},$ then $v_{3}\gg v_{Te},$ and   the assumption (i) and (iv).
Then  %for any $k,\omega,$ by the assumption (iv) of $Proposition$ 1.1, we obtain
%$$|1-\tilde{\mathcal{L}}(\omega,k)|\geq\kappa,$$
 there exists  some constant $0<\kappa<1$ such that
$\|\widetilde{\mathcal{K}}^{0}(\omega,k)\|_{L^{\infty}}\leq\kappa.$

%\begin{align}
%&\mathcal{L}(t,k)=\frac{q}{m}\int^{t}_{0}k_{3}(\hat{W}_{2}(i\eta_{k_{1}})\partial_{\eta_{3}}
%\hat{f}^{0}(\eta_{k_{1}},\eta_{k_{2}},k_{3}t)+\hat{W}_{1}(i\eta_{k_{2}})\partial_{\eta_{3}}\hat{f}^{0}(\eta_{k_{1}},\eta_{k_{2}},k_{3}t))dt\notag\\
%&-\frac{q}{m}\int^{t}_{0}
%(\hat{W}_{1}(i\eta_{k_{1}})\hat{f}^{0}(\eta_{k_{1}},\eta_{k_{2}},k_{3}t)+\hat{W}_{2}(i\eta_{k_{2}})\hat{f}^{0}(\eta_{k_{1}},\eta_{k_{2}},k_{3}t))dt,\notag\\
%\end{align}

In conclusion, we have $\|e^{2\pi\lambda'_{0}\nu_{k}}\hat{\rho}(t,k)\|_{L^{2}}\leq\frac{\|e^{2\pi\lambda'_{0}\nu_{k}}A\|_{L^{2}(dt)}}{\kappa}.$

%$(\mathbf{wherther}\quad \mathbf{consider} \quad\omega=0$ ? since we only care the integration)

%If $\omega=0,$ by the assumption (iv) of $Proposition 1.1,$  it is easy to obtain that
%$$|\tilde{\rho}(\omega,k)|\leq|\tilde{A}(\omega,k)|.$$
%Then applying Plancherel's identity to find for each $ k,$
%\begin{align}
%\|\rho\|_{L^{2}(dt)}\leq\frac{\|A\|_{L^{2}(dt)}}{\kappa},
%\end{align}
%furthermore, we have $(\mathbf{right}?)$
%\begin{align}
%\|\Phi\|_{L^{2}(dt)}\leq\frac{\|e^{\lambda|k_{3}|t}A\|_{L^{2}(dt)}}{\kappa},
%\end{align}
Then we get
\begin{align}
&\|\Phi\|_{L^{\infty}(dt)}\leq\|e^{2\pi\lambda'_{0}\nu_{k}}A\|_{L^{\infty}
(d\omega)}+\frac{\|e^{2\pi\lambda'_{0}\nu_{k}}\mathcal{L}\|_{L^{2}(dt)}\|e^{2\pi\lambda'_{0}\nu_{k}}A\|_{L^{2}(dt)}}{\kappa}
\end{align}
\begin{align}
&\|e^{2\pi\lambda'_{0}\nu_{k}}\mathcal{L}\|^{2}_{L^{2}(dt)}\leq\int^{\infty}_{0}|e^{4\pi\lambda'_{0}|\eta_{k1}|}
e^{4\pi\lambda'_{0}|\eta_{k2}|}e^{4\pi\lambda'_{0}|k_{3}|t}\bigg\{\frac{q}{m}\int^{t}_{0}\bigg[k_{3}\bigg((2\pi i\eta_{k_{1}})\widehat{W}_{2}
\partial_{\eta_{3}}
\hat{f}^{0}(\eta_{k_{1}},\eta_{k_{2}},k_{3}\tau)\notag\\
&+(2\pi i\eta_{k_{2}})\widehat{W}_{1}\partial_{\eta_{3}}\hat{f}^{0}(\eta_{k_{1}},\eta_{k_{2}},k_{3}\tau)\bigg)-\frac{q}{m}
\bigg(\widehat{v'_{1}\partial_{v'_{3}}f^{0}}(\eta_{k1},\eta_{k2},k_{3}\tau)\widehat{W}_{1}
+\widehat{v'_{2}\partial_{v'_{3}}f^{0}}(\eta_{k1},\eta_{k2},k_{3}\tau)\widehat{W}_{2}\bigg)\notag\\
&\cdot k_{3}+\frac{q}{m}(k_{2}\widehat{W}_{1}+k_{1}\widehat{W}_{2})\bigg(\widehat{v'_{2}\partial_{v'_{1}}f^{0}}(\eta_{k1},\eta_{k2},k_{3}\tau)
-\widehat{v'_{1}\partial_{v'_{2}}f^{0}}(\eta_{k1},\eta_{k2},k_{3}\tau)\bigg)\bigg]d\tau\notag\\
&-\frac{q}{m}
\bigg(\widehat{W}_{1}(2\pi i\eta_{k_{1}})\hat{f}^{0}(\eta_{k_{1}},\eta_{k_{2}},k_{3}t)-\widehat{W}_{2}(2\pi i\eta_{k_{2}})
\hat{f}^{0}(\eta_{k_{1}},\eta_{k_{2}},k_{3}t)\bigg)\bigg\}^{2}dt.\notag\\
\end{align}
 Using the conditions of Theorem 0.1, by the simple  computation similar to $I,$
\begin{align}
&\|e^{2\pi\lambda'_{0}\nu_{k}}\mathcal{L}\|^{2}_{L^{2}(dt)}\leq C(W,k,\Omega)\int^{\infty}_{0}e^{-4\pi(\lambda_{0}-\lambda'_{0})|k_{3}|t}
C(W,k,\Omega)dt\leq\frac{C(W,k,\Omega)}{(\lambda_{0}-\lambda'_{0})}.
\end{align}
Now we estimate $\|e^{2\pi\lambda'_{0}|\eta_{k1}|}e^{2\pi\lambda'_{0}|\eta_{k2}|}e^{2\pi\lambda'_{0}|k_{3}|t}A\|_{L^{2}(dt)}$ as the above process,
\begin{align}
&\|e^{2\pi\lambda'_{0}|\eta_{k1}|}e^{2\pi\lambda'_{0}|\eta_{k2}|}e^{2\pi\lambda'_{0}|k_{3}|t}A\|_{L^{2}(dt)}\notag\\
%&=\bigg(\int^{\infty}_{0}|e^{2\pi\lambda_{0}|\eta_{k1}|}e^{2\pi\lambda_{0}|\eta_{k2}|}e^{2\pi\lambda_{0}|k_{3}|t}\int_{\mathbb{T}^{3}}\int_{\mathbb{R}^{3}}e^{-2\pi ix\cdot k}f_{0}(x'(0,x,v),v'(0,x,v))dvdx|^{2}dt\bigg)^{\frac{1}{2}}\notag\\
&=\bigg(\int^{\infty}_{0}|e^{2\pi\lambda'_{0}|\eta_{k1}|}e^{2\pi\lambda'_{0}|\eta_{k2}|}e^{2\pi\lambda'_{0}|k_{3}|t}
\hat{f}_{0}(k,\eta_{k1},\eta_{k2},k_{3}t)|^{2}dt\bigg)^{\frac{1}{2}}\notag\\
&\leq \frac{C(\Omega)}{(\lambda_{0}-\lambda'_{0})^{\frac{1}{2}}}\|e^{2\pi\lambda'_{0}|\eta_{k1}|}e^{2\pi\lambda'_{0}|\eta_{k2}|}
e^{2\pi\lambda_{0}|k_{3}|t}A\|_{L^{\infty}(dt)}.
%\leq\frac{C(\Omega)}{(\lambda_{1}-\lambda_{0})^{\frac{1}{2}}}\|\widetilde{e^{2\pi\lambda_{0}\nu_{k}}A}\|_{L^{1}(d\omega)}.
\end{align}
%It is easy to check that
%$\|\widetilde{e^{2\pi\lambda_{0}\nu_{k}}A}\|_{L^{1}(d\omega)}\leq\|\widetilde{e^{2\pi\lambda_{0}\nu_{k}}A}\|_{L^{1}(d\omega)},$ for $\lambda_{0}<\lambda_{1}.$

 Now we consider $k_{3}=0,$
 $ k_{1}k_{2}\neq0,$
    %$$\hat{\rho}(t,k_{1},k_{2},0)=\int_{\mathbb{T}^{3}}\int_{\mathbb{R}^{3}}e^{-2\pi i(x_{1},x_{2})\cdot (k_{1},k_{2})}f(t,x,v)dvdx=\int_{\mathbb{T}^{3}}\int_{\mathbb{R}^{3}}e^{-2\pi i(x_{1},x_{2})\cdot (k_{1},k_{2})}f_{0}(x'(0,x,v),v'(0,x,v))dvdx$$
%$$-\frac{q}{m}\int^{t}_{0}\int_{\mathbb{T}^{3}}\int_{\mathbb{R}^{3}}e^{-2\pi i(x_{1},x_{2})\cdot (k_{1},k_{2})}[(E+v'\times B)\cdot\nabla'_{v}f^{0}](\tau,x'(\tau,x,v),v'(\tau,x,v))dvdx,$$
\begin{align}
&\tilde{\rho}(\omega,k_{1},k_{2},0)=\int_{\mathbb{R}^{+}}\int_{\mathbb{T}^{3}}\int_{\mathbb{R}^{3}}e^{2\pi it\omega}e^{-2\pi i(x_{1},x_{2})\cdot (k_{1},k_{2})}f(t,x_{1},x_{2},x_{3},v)dx_{1}dx_{2}dx_{3}dvdt\notag\\
&=\int_{\mathbb{R}^{+}}\int_{\mathbb{T}^{3}}\int_{\mathbb{R}^{3}}e^{2\pi it\omega}e^{-2\pi i(x_{1},x_{2})\cdot (k_{1},k_{2})}f_{0}(x'(0,x,v),v'(0,x,v))dvdx_{1}dx_{2}dx_{3}dt-\frac{q}{m}\int_{\mathbb{R}^{+}}\int_{\mathbb{T}^{3}}\int^{t}_{0}\int_{\mathbb{R}^{3}}\notag\\
&\cdot e^{2\pi it\omega}e^{-2\pi i(x_{1},x_{2})\cdot (k_{1},k_{2})}[(E+v'(\tau,x,v)\times B)\cdot\nabla'_{v}f^{0}](\tau,x'(\tau,x,v),v'(\tau,x,v))dvdxd\tau dt.\notag\\
\end{align}

 Taking the Fourier transform in the variables $t,(x_{1},x_{2}),$
$$\widetilde{E}(\omega,k_{1},k_{2},0)=\widehat{W}(k_{1},k_{2},0)\tilde{\rho}(\omega,k_{1},k_{2},0),
\quad\omega\widetilde{B}(\omega,k)=(k_{1},k_{2},\partial_{x_{3}})\times\widetilde{E}(\omega,k_{1},k_{2},0).$$
Furthermore, we get
\begin{align}
&v\times\widetilde{B}(\omega,k_{1},k_{2},0)=\frac{1}{\omega}[v\times((k_{1},k_{2},\partial_{x_{3}})\times\widetilde{E}(\omega,k_{1},k_{2},0))]\notag\\
&=\frac{1}{\omega}\bigg(v_{2}(k_{1}\widehat{W}_{2}-k_{2}\hat{W}_{1})-v_{3}\partial_{x_{3}}\widehat{W}_{1},-v_{1}(k_{1}\widehat{W}_{2}-k_{2}\widehat{W}_{1})+v_{3}
\partial_{x_{3}}\widehat{W}_{2},
v_{1}\partial_{x_{3}}\widehat{W}_{1}+v_{2}\partial_{x_{3}}\widehat{W}_{2}\bigg)\tilde{\rho}(\omega,k_{1},k_{2},0).\notag\\
\end{align}

%$$\hat{\rho}(t,k_{1},k_{2},0)=\int_{\mathbb{T}^{3}}\int_{\mathbb{R}^{3}}e^{-2\pi i(x_{1},x_{2})\cdot (k_{1},k_{2})}f_{0}(x'(0,x,v),v'(0,x,v))dv'dx_{\perp}dx_{3}$$
%$$-\frac{q}{m}\int^{t}_{0}\int_{\mathbb{R}^{3}}e^{-2\pi i(\frac{1}{\Omega}(v'_{2}-v_{2},-v'_{1}+v_{1})\cdot (k_{1},k_{2})}[(v'\times \hat{B}(\tau,k_{1},k_{2},x_{3}))\cdot\nabla'_{v}f^{0}](v'(\tau,x,v))dv'dx_{3}d\tau.$$
%Then assume  $\hat{W}(k_{1},k_{2},0)=0,$ we get $\int_{\mathbb{R}^{3}}e^{-2\pi i(\frac{1}{\Omega}(v'_{2}-v_{2},-v'_{1}+v_{1})
%\cdot (k_{1},k_{2})}\hat{E}(\tau,k_{1},k_{2},0)\cdot\nabla'_{v}f^{0}dv'=0.$

$$\tilde{\rho}(\omega,k_{1},k_{2},0)=\int_{\mathbb{R}^{+}}\int_{\mathbb{T}^{3}}\int_{\mathbb{R}^{3}}e^{2\pi it\omega}e^{-2\pi ix_{12}\cdot k_{12}}f_{0}(x'(0,x,v),v'(0,x,v))dvdx_{12}dx_{3}dt$$
$$+\frac{q}{m}\tilde{\rho}(\omega,k_{12},0)\frac{1}{\omega}\int_{\mathbb{R}^{+}}\int_{\mathbb{R}^{3}}\int_{\mathbb{T}}e^{2\pi it\omega}
e^{-2\pi i\eta_{k2}v'_{2}}e^{-2\pi i\eta_{k1}v'_{1}}(v_{3}\partial_{x_{3}})(\widehat{W}_{1}\partial_{v'_{1}}f^{0}-\widehat{W}_{2}\partial_{v'_{2}}f^{0})dx_{3}dvdt$$
$$+\frac{q}{m}\tilde{\rho}(\omega,k_{12},0)\frac{1}{\omega}\int_{\mathbb{R}^{+}}\int_{\mathbb{R}^{3}}\int_{\mathbb{T}}e^{2\pi it\omega}
e^{-2\pi i\eta_{k2}v'_{2}}e^{-2\pi i\eta_{k1}v'_{1}}(k_{1}\widehat{W}_{2}-k_{2}\hat{W}_{1})$$
$$\cdot(v_{2}\partial_{v'_{1}}f^{0}-v_{1}\partial_{v'_{2}}f^{0})dx_{3}dvdt-\frac{q}{m}\tilde{\rho}(\omega,k_{12},0)\frac{1}{\omega}\int_{\mathbb{R}^{+}}
\int_{\mathbb{R}^{3}}\int_{\mathbb{T}}e^{2\pi it\omega}
e^{-2\pi i\eta_{k2}v'_{2}}e^{-2\pi i\eta_{k1}v'_{1}}$$
\begin{align}
&\cdot\partial_{x_{3}}(v_{1}\hat{W}_{1}+v_{2}\widehat{W}_{2})\partial_{v'_{3}}f^{0}
dx_{3}dvdt-\frac{q}{m}\tilde{\rho}(\omega,k_{12},0)\int_{\mathbb{R}^{+}}\int_{\mathbb{R}^{3}}\int_{\mathbb{T}}e^{2\pi it\omega}
e^{-2\pi i\eta_{k2}v'_{2}}e^{-2\pi i\eta_{k1}v'_{1}}\notag\\
&\cdot(\widehat{W}_{1}\partial_{v'_{1}}f^{0}+\widehat{W}_{2}\partial_{v'_{2}}f^{0})dx_{3}dvdt,\notag\\
\end{align}
where $\eta_{k_{1}}=\frac{1}{\Omega}(-k_{2}\cos\Omega t+k_{2}-k_{1}\sin\Omega t),$ $\eta_{k_{2}}=\frac{1}{\Omega}(-k_{2}\sin\Omega t-k_{1}+k_{1}\cos\Omega t).$

Let
%and by the definition of $\mathcal{L},$ it is easy to check
\begin{align}
&\tilde{\mathcal{L}}(\omega,k_{12},0)=\frac{q}{m}\frac{1}{\omega}\int_{\mathbb{R}^{+}}\int_{\mathbb{R}^{3}}\int_{\mathbb{T}}e^{2\pi it\omega}
e^{-2\pi i\eta_{k2}v'_{2}}e^{-2\pi i\eta_{k1}v'_{1}}(v_{3}\partial_{x_{3}})(\widehat{W}_{1}\partial_{v'_{1}}f^{0}-\widehat{W}_{2}\partial_{v'_{2}}f^{0})dx_{3}dvdt\notag\\
&-\frac{q}{m}\int_{\mathbb{R}^{+}}\int_{\mathbb{R}^{3}\int_{\mathbb{T}}}e^{2\pi it\omega}
(\widehat{W}_{1}\partial_{v'_{1}}f^{0}+\widehat{W}_{2}\partial_{v'_{2}}f^{0})\cdot e^{-2\pi i\eta_{k2}v'_{2}}e^{-2\pi i\eta_{k1}v'_{1}}dx_{3}dvdt\notag\\
&-\frac{q}{m}\frac{1}{\omega}\int_{\mathbb{R}^{+}}\int_{\mathbb{R}^{3}}\int_{\mathbb{T}}e^{2\pi it\omega}
e^{-2\pi i\eta_{k2}v'_{2}}e^{-2\pi i\eta_{k1}v'_{1}}\partial_{x_{3}}(v_{1}\widehat{W}_{1}+v_{2}\widehat{W}_{2})\partial_{v'_{3}}f^{0}
dx_{3}dvdt\notag\\
&+\frac{q}{m}\frac{1}{\omega}\int_{\mathbb{R}^{+}}\int_{\mathbb{R}^{3}}\int_{\mathbb{T}}e^{2\pi it\omega}
e^{-2\pi i\eta_{k2}v'_{2}}e^{-2\pi i\eta_{k1}v'_{1}}(k_{1}\widehat{W}_{2}-k_{2}\widehat{W}_{1})(v_{2}\partial_{v'_{1}}f^{0}
-v_{1}\partial_{v'_{2}}f^{0})dx_{3}dvdt,\notag\\
\end{align}
hence
\begin{align}
\tilde{\rho}(\omega,k_{12},0)=\tilde{A}(\omega,k_{12},0)+\tilde{\rho}(\omega,k_{12},0)\tilde{\mathcal{L}}(\omega,k_{12},0).
\end{align}
%Assume $\inf_{k\in\mathbb{Z}^{3}}|\tilde{\mathcal{L}}(\omega,k)-1|>\kappa,$ for some constant $0<\kappa<1,$
Then taking the inverse Fourier transform in time $ t, $ we get
$\hat{\rho}(t,k_{12},0)=\hat{A}(t,k_{12},0)+\hat{\rho}(t,k_{12},0)\ast\hat{\mathcal{L}}(t,k_{12},0),$
and

\begin{align}
&e^{2\pi\lambda'_{0}\eta_{k1}}e^{2\pi\lambda'_{0}\eta_{k2}}\hat{\rho}(t,k_{12},0)=e^{2\pi\lambda'_{0}\eta_{k1}}e^{2\pi\lambda'_{0}\eta_{k2}}
\hat{A}(t,k_{12},0)+e^{2\pi\lambda'_{0}\eta_{k1}}e^{2\pi\lambda'_{0}\eta_{k2}}\hat{\rho}(t,k_{12},0)\ast\hat{\mathcal{L}}(t,k_{12},0).
\end{align}

Then
\begin{align}
&\|\Phi(t,k_{12},0)\|_{L^{2}(dt)}
=\|\widetilde{\Phi}(\omega,k_{12},0)\|_{L^{2}}\leq\|\widetilde{\mathcal{A}}(\omega,k_{12},0)\|_{L^{2}}+\|\widetilde{\Phi}(\omega,k_{12},0)\|_{L^{2}}
\|\widetilde{\mathcal{K}}^{0}(\omega,k)\|_{L^{\infty}}\notag\\
&\leq\|e^{2\pi\lambda'_{0}\nu_{k_{12}}}\hat{A}(t,k_{12},0)\|_{L^{2}}+\|e^{2\pi\lambda'_{0}\nu_{k_{12}}}\hat{\rho}(t,k_{12},0)\|_{L^{2}}
\|\widetilde{\mathcal{K}}^{0}(\omega,k_{12},0)\|_{L^{\infty}},\notag\\
\end{align}
where $\nu_{k_{12}}=|\eta_{k1}|+|\eta_{k2}|.$

Next we have to estimate $\|\widetilde{\mathcal{K}}^{0}(\omega,k_{12},0)\|_{L^{\infty}}.$

Indeed,

$$\|\widetilde{\mathcal{K}}^{0}(\omega,k_{12},0)\|_{L^{\infty}}
\leq\sup_{\omega}\frac{q}{m}\bigg[\bigg|\frac{1}{\omega}\int_{\mathbb{R}^{+}}\int_{\mathbb{R}^{3}}\int_{\mathbb{T}}e^{2\pi it\omega}e^{2\pi\lambda'_{0}\nu_{k_{12}}}
e^{-2\pi i\eta_{k2}v'_{2}}e^{-2\pi i\eta_{k1}v'_{1}}(v_{3}\partial_{x_{3}})\cdot(\widehat{W}_{1}\partial_{v'_{1}}f^{0}$$
$$-\widehat{W}_{2}\partial_{v'_{2}}f^{0})
dx_{3}dvdt\bigg|+\frac{q}{m}\bigg|\int_{\mathbb{R}^{+}}\int_{\mathbb{R}^{3}}\int_{\mathbb{T}}e^{2\pi it\omega}
e^{2\pi\lambda_{0}\nu_{k_{12}}}\cdot(\widehat{W}_{1}\partial_{v'_{1}}f^{0}+\widehat{W}_{2}\partial_{v'_{2}}f^{0})\cdot e^{-2\pi i\eta_{k2}v'_{2}}$$
$$e^{-2\pi i\eta_{k1}v'_{1}}
dx_{3}dvdt\bigg|-\frac{q}{m}(\omega,k)\frac{1}{\omega}\int_{\mathbb{R}^{+}}\bigg|\int_{\mathbb{R}^{3}}\int_{\mathbb{T}}e^{2\pi it\omega}
e^{2\pi\lambda_{0}\nu_{k}}e^{-2\pi i\eta_{k2}v'_{2}}e^{-2\pi i\eta_{k1}v'_{1}}\cdot(v_{1}\widehat{W}_{1}+v_{2}\widehat{W}_{2})k_{3}$$
\begin{align}
&\partial_{v'_{3}}f^{0}
dv\bigg|dt+\frac{q}{m}(\omega,k)\frac{1}{\omega}\int_{\mathbb{R}^{+}}\bigg|\int_{\mathbb{R}^{3}}e^{2\pi it\omega}e^{-2\pi i\eta_{k2}v'_{2}}e^{2\pi \lambda_{0}\nu_{k}}e^{-2\pi i\eta_{k1}v'_{1}}
\cdot(k_{1}\widehat{W}_{2}-k_{2}\widehat{W}_{1})\notag\\
&(v_{2}\partial_{v'_{1}}f^{0}
-v_{1}\partial_{v'_{2}}f^{0})dv\bigg|dt\bigg]=I+II+III+IV.\notag\\
\end{align}
%\begin{align}
%&=\frac{q}{m}\int^{\infty}_{0}|\frac{k_{3}}{\omega}\hat{W}_{2}\int_{\mathbb{R}^{+}}e^{it\omega}(i\eta_{k1})e^{-ik_{3}v_{3}t}e^{\lambda|k_{3}|t}
%\widehat{v_{3}f^{0}}(\eta_{k1},\eta_{k2},v_{3})dv_{3}dt|d\omega\notag\\
%&+\frac{q}{m}\int^{\infty}_{0}|\frac{k_{3}}{\omega}\hat{W}_{1}\int_{\mathbb{R}^{+}}e^{it\omega}(i\eta_{k2})e^{-ik_{3}v_{3}t}e^{\lambda|k_{3}|t}
%\widehat{v_{3}f^{0}}(\eta_{k1},\eta_{k2},v_{3})dv_{3}dt|d\omega\notag\\
%&+\frac{q}{m}\int^{\infty}_{0}|\hat{W}_{1}\int_{\mathbb{R}^{+}}e^{it\omega}(i\eta_{k1})e^{-ik_{3}v_{3}t}e^{\lambda|k_{3}|t}\widehat{f^{0}}(\eta_{k1}
%,\eta_{k2},v_{3})dv_{3}dt|d\omega\notag\\
%&+\frac{q}{m}\int^{\infty}_{0}|\hat{W}_{2}\int_{\mathbb{R}^{+}}e^{it\omega}(i\eta_{k2})e^{-ik_{3}v_{3}t}e^{\lambda|k_{3}|t}\widehat{f^{0}}(\eta_{k1},\eta_{k2},
%v_{3})dv_{3}dt|d\omega\notag\\
%&+\cdots\cdots\notag\\
%&=I+II+III+IV+\cdots
%&\leq Ce^{-\alpha(\frac{\omega}{k_{3}})^{2}}+\varepsilon\leq C_{Te}e^{-\alpha(v_{Te})^{2}},
%\end{align}
%where $\eta_{k_{1}}=\frac{1}{\Omega}(-k_{2}\cos\Omega t+k_{2}-k_{1}\sin\Omega t),\eta_{k_{2}}=\frac{1}{\Omega}(-k_{2}\sin\Omega t-k_{1}+k_{1}\cos\Omega t).$

In fact, we only need to estimate one term of (1.10) because of similar processes of  other terms. Without loss of generality, we give an estimate for $ I.$ In the same way, we only estimate one term of $I,$ here we still denote $ I. $
$$ I=\sup_{\omega}\frac{q}{m}\bigg|\frac{1}{\omega}\int_{\mathbb{R}^{+}}\int_{\mathbb{R}}\int_{\mathbb{T}}\widehat{W}_{2}(k_{12},0)e^{2\pi it\omega}(2\pi i\eta_{k1})
e^{2\pi\lambda_{0}|\eta_{k1}|}e^{2\pi\lambda_{0}|\eta_{k2}|}\cdot\partial_{x_{3}}\widehat{v_{3}f^{0}}(\eta_{k1},\eta_{k2},v_{3})
dx_{3}dv_{3}dt\bigg|$$
%&=\frac{q}{m}\sum_{n}\frac{\lambda^{n}}{n!}\int^{\infty}_{0}|\frac{k_{3}}{\omega}\hat{W}_{2}\int_{\mathbb{R}^{+}}e^{it\omega}(i\eta_{k1})(ik_{3}t)^{n}
%e^{-ik_{3}v_{3}t}
%\widehat{v_{3}f^{0}}(\eta_{k1},\eta_{k2},v_{3})dv_{3}dt|d\omega\notag\\
%$$=\frac{q}{m}\sum_{n}\frac{\lambda^{n}}{n!}\int^{\infty}_{0}\bigg|\frac{k_{3}}{\omega}\hat{W}_{2}\int_{\mathbb{R}^{+}}\int_{\mathbb{R}}(2\pi i\eta_{k1})(2\pi ik_{3}t)^{n}
%e^{2\pi\lambda_{0}|\eta_{k1}|}e^{2\pi\lambda_{0}|\eta_{k2}|}e^{2\pi ik_{3}t(\frac{\omega}{k_{3}}-v_{3})}\cdot\widehat{v_{3}f^{0}}(\eta_{k1},\eta_{k2},v_{3})dv_{3}dt\bigg|d\omega$$
%$$=\frac{q}{m}\sum_{n}\frac{\lambda^{n}}{n!}\int^{\infty}_{0}\bigg|\frac{k_{3}}{\omega}\hat{W}_{2}\int_{\mathbb{R}^{+}}\int_{\mathbb{R}}(2\pi i\eta_{k1})
%(-1)^{n}e^{2\pi\lambda_{0}|\eta_{k1}|}e^{2\pi\lambda_{0}|\eta_{k2}|}\nabla_{v_{3}}^{n}e^{2\pi ik_{3}t(\frac{\omega}{k_{3}}-v_{3})}\cdot\widehat{v_{3}f^{0}}(\eta_{k1},\eta_{k2},v_{3})dv_{3}dt\bigg|d\omega$$
%$$=\frac{q}{m}\sum_{n}\frac{\lambda_{0}^{n}}{n!}\int^{\infty}_{0}\bigg|\frac{k_{3}}{\omega}\hat{W}_{2}\int_{\mathbb{R}^{+}}(2\pi i\eta_{k1})
%(-1)^{n}e^{2\pi\lambda_{0}|\eta_{k1}|}e^{2\pi\lambda_{0}|\eta_{k2}|}e^{2\pi ik_{3}t(\frac{\omega}{k_{3}})}
%(2\pi ik_{3}t)^{n}\cdot\widehat{v_{3}f^{0}}(\eta_{k1},\eta_{k2},k_{3}t) dt\bigg|d\omega$$
$$\leq\frac{q}{mv_{Te}}e^{-c_{0}v_{Te}},$$
%\begin{align}
%&\leq\frac{q}{m}\sum_{n}\frac{\lambda_{0}^{n}}{n!}\int^{\infty}_{v_{Te}}\bigg|\frac{k_{3}}{\omega}\hat{W}_{2}(2\pi i\eta_{k1})
%e^{2\pi\lambda_{0}|\eta_{k1}|}e^{2\pi\lambda_{0}|\eta_{k2}|}
%(-i\nabla_{\frac{\omega}{k_{3}}})^{n}\widehat{v_{3}f^{0}}(\eta_{k1},\eta_{k2},\frac{\omega}{k_{3}}) \bigg|d\omega\leq\frac{q}{mv_{Te}}e^{-c_{0}v_{Te}}\notag\\
%&+\frac{q}{m}\sum_{n}\frac{\lambda^{n}}{n!}\int^{v_{Te}}_{0}|\frac{k_{3}}{\omega}\hat{W}_{2}(i\eta_{k1})
%e^{\lambda|\eta_{k1}|}e^{\lambda|\eta_{k2}|}
%(-i\nabla_{\frac{\omega}{k_{3}}})^{n}\widehat{v_{3}f^{0}}(\eta_{k1},\eta_{k2},\frac{\omega}{k_{3}}) |d\omega\leq\frac{q}{mv_{Te}}e^{-c_{0}v_{Te}},
%\end{align}
where in the last inequality we use the facts that if $k_{3}=0,$ then $v_{3}\gg v_{Te},$ and   the assumption (i) and (iv).
Then  %for any $k,\omega,$ by the assumption (iv) of $Proposition$ 1.1, we obtain
%$$|1-\tilde{\mathcal{L}}(\omega,k)|\geq\kappa,$$
 there exists  some constant $0<\kappa<1$ such that
$\|\widetilde{\mathcal{K}}^{0}(\omega,k)\|_{L^{\infty}}\leq\kappa.$

%Since $\partial_{t}\hat{B}=(k_{1},k_{2},k_{3})\times \hat{E}(t,k_{1},k_{2},k_{3}),E(t,x)=W(x)\ast\rho(t,x),$
%then if  $ k_{3}=0,$  since $W(x)=(W_{1}(x),W_{2}(x),0),$
%$\widehat{W}(k_{1},k_{2},0)=0,$ we have
%$\hat{B}(t,k_{1},k_{2},0)\equiv0.$
%So $\int_{\mathbb{R}^{3}}e^{-2\pi i(\frac{1}{\Omega}(v'_{2}-v_{2},-v'_{1}+v_{1})\cdot (k_{1},k_{2})}\nabla'_{v}\cdot(\hat{B}\times(f^{0}v'))dv'=0.$
%Finally, we get $\hat{\rho}(t,k_{1},k_{2},0)=\int\int e^{-2\pi i(x_{1},x_{2})\cdot (k_{1},k_{2})}f_{0}dvdx'.$  And we finish the proof.

First through $Lemma $ 1.1, we can study  the asymptotic behavior of the electric field and the magnetic field  $E(t,x),B(t,x).$

\begin{cor} Under the assumptions of Theorem 1.1, and $E(t,x),B(t,x)$ satisfy the  Maxwell equations of (1.1), then for any  for any $\lambda_{0}''<\lambda_{0}'<\lambda_{0},$  we have
\begin{align}
&|\widehat{E}(t,k)|\leq e^{-2\pi\lambda'_{0}|\eta_{k1}|}e^{-2\pi\lambda'_{0}|\eta_{k2}|}e^{-2\pi\lambda'_{0}|k_{3}|t},
|\widehat{B}(t,k)|\leq te^{-2\pi\lambda''_{0}|\eta_{k1}|}e^{-2\pi\lambda''_{0}|\eta_{k2}|}e^{-2\pi\lambda'_{0}|k_{3}|t}.
 \end{align}
\end{cor}
$Proof.$ Since $\partial_{t}\widehat{B}(t,k)=k\times\widehat{E}(t,k),\widehat{E}(t,k)=\widehat{W}(k)\hat{\rho}(t,k)=(\hat{W}_{1}(k),
\widehat{W}_{2}(k),0)\hat{\rho}(t,k),$
then $\partial_{t}\widehat{B}(t,k)=(-k_{3}\widehat{W}_{2}(k),k_{3}\widehat{W}_{1}(k),k_{1}\widehat{W}_{2}(k)-k_{2}\widehat{W}_{1}(k))\hat{\rho}(t,k).$
By Lemma 1.1,
$|\partial_{t}\widehat{B}(t,k)|\leq C(k,\Omega,W)$
$e^{-2\pi\lambda'_{0}|\eta_{k1}|}$
$e^{-2\pi\lambda'_{0}|\eta_{k2}|}$
$e^{-2\pi\lambda'_{0}|k_{3}|t}.$
 Hence, from $\widehat{B}(0,k)=0,$ we have $|\widehat{B}(t,k)|\leq te^{-2\pi\lambda'_{0}|\eta_{k1}|}e^{-2\pi\lambda'_{0}|\eta_{k2}|}e^{-2\pi\lambda'_{0}|k_{3}|t}.$
 Then $|\widehat{B}(t,k)|\leq e^{-2\pi\lambda''_{0}|\eta_{k1}|}e^{-2\pi\lambda''_{0}|\eta_{k2}|}e^{-2\pi\lambda''_{0}|k_{3}|t}.$
%Furthermore,
%\begin{align}
%\|\Phi\|_{L^{\infty}(dt)}\leq(1+\frac{C(k,W,\Omega)}{(\lambda_{1}-\lambda)^{2}})\|(\widetilde{e^{\lambda|k_{3}|t}A})\|_{L^{1}(d\omega)},
%\end{align}
$$$$
         $Proof$
         $of$
         $Theorem$
         $ 0.1.$
From (1.1), we have
\begin{align}
&f(t,x,v)=f_{0}(x'(0,x,v),v'(0.x,v))-\frac{q}{m}\int^{t}_{0}[(E+v'\times B)\cdot\nabla'_{v}f^{0}](\tau.x'(\tau,x,v),v'(\tau,x,v))d\tau.
\end{align}

Taking the Fourier-Laplace transform in variables $x,v,t,$  we find
\begin{align}
&\tilde{f}(\omega,k,\eta)=\int_{\mathbb{R}^{+}}\int_{\mathbb{T}^{3}}\int_{\mathbb{R}^{3}}e^{2\pi it\cdot\omega}e^{-2\pi ix\cdot k}e^{-2\pi iv\cdot\eta}f_{0}(x'(0,x,v),v'(0,x,v))dvdxdt\notag\\
&-\frac{q}{m}\int_{\mathbb{R}^{+}}\int_{\mathbb{T}^{3}}\int_{\mathbb{R}^{3}}e^{2\pi it\cdot\omega}e^{-2\pi ix\cdot k}e^{-2\pi iv\cdot\eta}\int^{t}_{0}[B\cdot\nabla'_{v}\times(f^{0}v')](\tau,x'(\tau,x,v),v'(\tau,x,v))d\tau dvdxdt\notag\\
&-\frac{q}{m}\int_{\mathbb{R}^{+}}\int_{\mathbb{T}^{3}}\int_{\mathbb{R}^{3}}e^{2\pi it\cdot\omega}e^{-2\pi ix\cdot k}e^{-2\pi iv\cdot\eta}\int^{t}_{0}[E\cdot\nabla'_{v}f^{0}](\tau,x'(\tau,x,v),v'(\tau,x,v))d\tau dvdxdt\notag\\
&=I+II+III.
\end{align}
\begin{align}
&I=\int_{\mathbb{R}^{+}}\int_{\mathbb{T}^{3}}\int_{\mathbb{R}^{3}}e^{2\pi it\cdot\omega}e^{-2\pi iv\cdot\eta}\cdot f_{0}(x,v')\exp(-2\pi i[x+(\frac{1}{\Omega}(v'_{2}-v_{2}),
-\frac{1}{\Omega}(v'_{1}-v_{1}),v_{3}t)]\cdot k)dvdxdt\notag\\
&=\int_{\mathbb{R}^{+}}\int_{\mathbb{R}^{3}}e^{2\pi it\cdot\omega}e^{-2\pi iv_{3}(\eta_{3}+k_{3}t)}e^{-2\pi iv'_{2}(\eta_{2}+\eta_{k2})}
e^{-2\pi iv'_{1}(\eta_{1}+\eta_{k1})}\hat{f}_{0}(k,v')dvdt.\notag\\
%&=\int_{\mathbb{R}^{+}}e^{it\cdot\omega}dt\int_{\mathbb{R}}e^{-i(k_{3}t+\eta_{3})v_{3}}dv_{3}\hat{f}_{0}(k,v_{\perp},v_{3})\cdot \exp(-i\frac{v_{\perp}k_{\perp}}{\Omega}[\sin(\theta+\Omega t-\varphi)-\sin(\theta-\varphi)])dv_{3}dt\notag\\
%&=2\pi\int_{\mathbb{R}^{+}}e^{it\cdot\omega}dt\int_{\mathbb{R}}e^{-i(k_{z}t+\eta_{z})v_{z}}dv_{z}\int_{\mathbb{R}}v_{\perp}\hat{f}_{0}(k,v_{\perp},v_{z})
%\cdot(\sum^{\infty}_{n=1}J^{2}_{n}(\frac{v_{\perp}k_{\perp}}{\Omega})\cos(n\Omega t))dv_{\perp}
\end{align}
%where $J_{n}(z)$ is the Bessel function.

%In this paper we only consider $n=1.$   Therefore, from (1.6),
%\begin{align}
%I=\int_{\mathbb{R}^{+}}e^{it\cdot\omega}dt\int_{\mathbb{R}}e^{-i(k_{z}t+\eta_{z})v_{z}}dv_{z}\int_{\mathbb{R}}
%\sum^{\infty}_{n=1}v_{\perp}\hat{f}_{0}(k,v_{\perp},v_{z})J^{2}_{n}(\frac{v_{\perp}k_{\perp}}{\Omega}))\cos(n\Omega t)dv_{\perp}
%\end{align}

As the above process, we have

%$$II=-\frac{q}{m}\int_{\mathbb{R}^{+}}\int_{\mathbb{R}^{3}}e^{2\pi it\cdot\omega}e^{-2\pi i(k_{3}t+\eta_{3})v_{3}}e^{-2\pi iv'_{2}(\eta_{2}+\eta_{k2})}
%e^{-2\pi iv'_{1}(\eta_{1}+\eta_{k1})}$$
%$$ \int^{t}_{0}(\hat{B}(\tau,k)\cdot (\nabla'_{v}\times(f^{0}(v_{\perp},v_{3})v(\tau,x,v')))d\tau dvdt$$
$$II=-\frac{q}{m}\int_{\mathbb{R}^{+}}\int_{\mathbb{R}^{3}}e^{2\pi it\cdot\omega}e^{-2\pi i(k_{3}t+\eta_{3})v_{3}}e^{-2\pi iv'_{2}(\eta_{2}+\eta_{k2})}
e^{-2\pi iv'_{1}(\eta_{1}+\eta_{k1})} $$
$$\cdot\int^{t}_{0} (\hat{B}(\tau,k)\cdot (v_{3}\partial_{v'_{2}}f^{0}-v'_{2}\partial_{v_{3}}f^{0},v'_{}\partial_{v_{3}}f^{0}-v_{3}\partial_{v'_{1}}f^{0},
v'_{2}\partial_{v'_{1}}f^{0}-v'_{1}\partial_{v'_{2}}f^{0})d\tau dvdt$$
\begin{align}
&=-\frac{q}{m}\int_{\mathbb{R}^{+}}e^{2\pi i\tau\cdot\omega}\hat{B}(\tau,k)d\tau\cdot\int_{\mathbb{R}^{+}}e^{2\pi it\cdot\omega}\int_{\mathbb{R}^{3}}
e^{-2\pi i(k_{3}t+\eta_{3})v_{3}}e^{-2\pi iv'_{2}(\eta_{2}+\eta_{k2})}
e^{-2\pi iv'_{1}(\eta_{1}+\eta_{k1})}\notag\\
&\cdot(v_{3}\partial_{v'_{2}}f^{0}-v'_{2}\partial_{v_{3}}f^{0},v'_{}\partial_{v_{3}}f^{0}-v_{3}\partial_{v'_{1}}f^{0},
v'_{2}\partial_{v'_{1}}f^{0}-v'_{1}\partial_{v'_{2}}f^{0}) dvdt\notag\\
&=-\frac{q}{m}\tilde{\rho}(\omega,k)\frac{1}{\omega}\int_{\mathbb{R}^{+}}e^{2\pi it\cdot\omega}\int_{\mathbb{R}^{3}}
e^{-2\pi i(k_{z}t+\eta_{3})v_{3}}e^{-2\pi iv'_{2}(\eta_{2}+\eta_{k2})}
e^{-2\pi iv'_{1}(\eta_{1}+\eta_{k1})}(k\times\hat{W}(k))\notag\\
&\cdot (v_{3}\partial_{v'_{2}}f^{0}-v'_{2}\partial_{v_{3}}f^{0},v'_{}\partial_{v_{3}}f^{0}-v_{3}\partial_{v'_{1}}f^{0},
v'_{2}\partial_{v'_{1}}f^{0}-v'_{1}\partial_{v'_{2}}f^{0}) dvdt.\notag\\
\end{align}
%$\mathbf{Here}$
 %$\mathbf{we}$
 %$\mathbf{set}$
Since
 $\widehat{W}(k)=(\widehat{W}_{1}(k),\widehat{W}_{2}(k),0),$
then $ k\times\widehat{W}(k)=(-k_{3}\widehat{W}_{2},k_{3}\widehat{W}_{1},k_{1}\widehat{W}_{2}-k_{2}\widehat{W}_{1}),$
%$$C_{W}(t,k,v)=(k_{\perp}W_{3})(v_{\perp}\partial_{v_{z}}f^{0}-v_{z}\partial_{v_{\perp}}f^{0})\cos(\theta+\Omega t-\varphi),$$
and \begin{align}
%&II=-\frac{q}{m}\tilde{\rho}(\omega,k)\frac{1}{\omega}\int_{\mathbb{R}^{+}}e^{it\cdot\omega}\int_{\mathbb{R}}
%e^{-i(k_{z}t+\eta_{z})v_{z}}\exp(-i\frac{v_{\perp}k_{\perp}}{\Omega}
%[\sin(\theta+\Omega t-\varphi)-\sin(\theta-\varphi)])\notag\\
%&[(k_{\perp}W_{3})(v_{\perp}\partial_{v_{z}}f^{0}-v_{z}\partial_{v_{\perp}}f^{0})\cos(\theta+\Omega t-\varphi)]dv_{z}dt\notag\\
&II=\frac{q}{m}\tilde{\rho}(\omega,k)\frac{1}{\omega}\int_{\mathbb{R}^{+}}e^{2\pi it\cdot\omega}\int_{\mathbb{R}}e^{-2\pi i(k_{3}t+\eta_{3})v_{3}}
e^{-2\pi iv'_{2}(\eta_{2}+\eta_{k2})}
e^{-2\pi iv'_{1}(\eta_{1}+\eta_{k1})}
\cdot (k_{3}\widehat{W}_{2}(-v_{3}\partial_{v'_{2}}f^{0}\notag\\
&+v'_{2}\partial_{v_{3}}f^{0})+k_{3}\widehat{W}_{1}(v'_{1}\partial_{v_{3}}f^{0}-v_{3}\partial_{v'_{1}}f^{0})+
(k_{1}\widehat{W}_{2}-k_{2}\widehat{W}_{1})
(v'_{2}\partial_{v'_{1}}f^{0}-v'_{1}\partial_{v'_{2}}f^{0}))dvdt.\notag\\
%&=-\frac{q}{m}\tilde{\rho}(\omega,k)\frac{1}{\omega}\int_{\mathbb{R}^{+}}e^{it\cdot\omega}\int_{\mathbb{R}}\int_{\mathbb{R}}\sum^{\infty}_{n=1}e^{-i(k_{z}t+\eta_{z})v_{z}}
%(k_{\perp}W_{3})(v_{\perp}\partial_{v_{z}}f^{0}-v_{z}\partial_{v_{\perp}}f^{0})\frac{\Omega}{k_{\perp}}\notag\\
%& J^{2}_{n}(\frac{v_{\perp}k_{\perp}}{\Omega})\cos(n\Omega t)dv_{\perp}dv_{z}\notag\\
\end{align}

Similarly,
\begin{align}
&III=-\frac{q}{m}\tilde{\rho}(\omega,k)\int_{\mathbb{R}^{+}}e^{2\pi it\cdot\omega}\int_{\mathbb{R}^{3}}
e^{-2\pi i(k_{3}t+\eta_{3})v_{3}}e^{-2\pi iv'_{2}(\eta_{2}+\eta_{k2})}
e^{-2\pi iv'_{1}(\eta_{1}+\eta_{k1})}(\widehat{W}_{1}\partial_{v'_{1}}f^{0}+\widehat{W}_{2}\partial_{v'_{2}}f^{0})dvdt.\notag\\
%&=-\frac{q}{m}\tilde{\rho}(\omega,k)\int_{\mathbb{R}^{+}}e^{it\cdot\omega}\int_{\mathbb{R}^{2}}e^{-i(k_{z}t+\eta_{z})v_{z}}\sum^{\infty}_{n=1}J^{2}_{n}(\frac{v_{\perp}k_{\perp}}{\Omega})\cos(n\Omega %t)
%(\hat{W}_{3}\partial_{v_{z}}f^{0})v_{\perp}dv_{\perp}dv_{z}dt.
\end{align}

Therefore, from(1.17)-(1.21), we have
\begin{align}
&\omega\tilde{f}(\omega,k,v',\eta_{3})=\omega (I+II+III).
\end{align}

Furthermore,
\begin{align}
&\partial_{t}\hat{f}(t,k,\eta)=\partial_{t}\int_{\mathbb{R}^{3}}e^{-2\pi i(k_{3}t+\eta_{3})v_{3}}e^{-2\pi iv'_{2}(\eta_{2}+\eta_{k2})}
e^{-2\pi iv'_{1}(\eta_{1}+\eta_{k1})}\hat{f}_{0}(k,v')dv\notag\\
&-\frac{q}{m}\hat{\rho}(t,k)\ast_{t} \int_{\mathbb{R}^{3}}e^{-2\pi i(k_{3}t+\eta_{3})v_{3}}e^{-2\pi iv'_{2}(\eta_{2}+\eta_{k2})}
e^{-2\pi iv'_{1}(\eta_{1}+\eta_{k1})}\cdot (k_{3}\widehat{W}_{2}(-v_{3}\partial_{v'_{2}}f^{0}\notag\\
&+v'_{2}\partial_{v_{3}}f^{0})+k_{3}\widehat{W}_{1}(v'_{1}\partial_{v_{3}}f^{0}-v_{3}\partial_{v'_{1}}f^{0})+
(k_{1}\widehat{W}_{2}-k_{2}\widehat{W}_{1})
(v'_{2}\partial_{v'_{1}}f^{0}-v'_{1}\partial_{v'_{2}}f^{0}))dvdt\notag\\
&-\frac{q}{m}\hat{\rho}(t,k)\ast_{t}\partial_{t}\int_{\mathbb{R}^{3}}
e^{-2\pi i(k_{3}t+\eta_{3})v_{3}}e^{-2\pi iv'_{2}(\eta_{2}+\eta_{k2})}
e^{-2\pi iv'_{1}(\eta_{1}+\eta_{k1})}(\widehat{W}_{1}\partial_{v'_{1}}f^{0}+\widehat{W}_{2}\partial_{v'_{2}}f^{0})dvdt\notag\\
&=IV+V+VI.
\end{align}

 Now we estimate $IV,V,VI,$ respectively.
\begin{align}
&IV=\partial_{t}\int_{\mathbb{R}^{3}}e^{-2\pi i(k_{3}t+\eta_{3})v_{3}}e^{-2\pi iv'_{2}(\eta_{2}+\eta_{k2})}
e^{-2\pi iv'_{1}(\eta_{1}+\eta_{k1})}\hat{f}_{0}(k,v')dv=\partial_{t}\hat{f}_{0}(k,\eta_{1}+\eta_{k1},\eta_{2}+\eta_{k2},k_{3}t+\eta_{3})\notag\\
%&\cdot \exp(-i\frac{v_{\perp}k_{\perp}}{\Omega}[\sin(\theta+\Omega t-\varphi)-\sin(\theta-\varphi)])\notag\\
%&=\partial_{t}\int_{\mathbb{R}}e^{-i(k_{z}t+\eta_{z})v_{z}}dv_{z}\int_{\mathbb{R}}\sum^{\infty}_{n=1}J^{2}_{n}(\frac{v_{\perp}k_{\perp}}{\Omega})\cos(n\Omega t)v_{\perp}\hat{f}_{0}(k,v_{\perp},v_{z})dv_{\perp}\notag\\
%&=-\Omega \int_{\mathbb{R}}\sum^{\infty}_{n=1}nJ^{2}_{n}(\frac{v_{\perp}k_{\perp}}{\Omega})\sin(n\Omega t)v_{\perp}\hat{f}_{0}(k,v_{\perp},k_{z}t+\eta_{z})dv_{\perp}\notag\\
%&+k_{z} \int_{\mathbb{R}}\sum^{\infty}_{n=1}J^{2}_{n}(\frac{v_{\perp}k_{\perp}}{\Omega})\cos(n\Omega t)v_{\perp}\partial_{\eta_{z}}\hat{f}_{0}(k,v_{\perp})dv_{\perp}\notag\\
&\leq C(\Omega,k)e^{-2\pi\lambda_{0}|\eta_{1}+\eta_{k1}|}e^{-2\pi\lambda_{0}|\eta_{2}+\eta_{k2}|}e^{-2\pi\lambda_{0}|k_{3}t+\eta_{3}|},
\end{align}
here we use the  assumption that $|\hat{f}_{0}(k,\eta)|\leq C_{0} e^{-2\pi\lambda_{0}|\eta_{1}|}e^{-2\pi\lambda_{0}|\eta_{2}|}e^{-2\pi\lambda_{0}|\eta_{3}|}.$

$$V=-\frac{q}{m}\hat{\rho}(t,k)\ast_{t} \int_{\mathbb{R}^{3}}e^{-2\pi i(k_{3}t+\eta_{3})v_{3}}
e^{-2\pi iv'_{2}(\eta_{2}+\eta_{k2})}
e^{-2\pi iv'_{1}(\eta_{1}+\eta_{k1})}\cdot (k_{3}\widehat{W}_{2}(-v_{3}\partial_{v'_{2}}f^{0}$$
$$+v'_{2}\partial_{v_{3}}f^{0})+k_{3}\widehat{W}_{1}(v'_{1}\partial_{v_{3}}f^{0}-v_{3}\partial_{v'_{1}}f^{0})+
(k_{1}\widehat{W}_{2}-k_{2}\hat{W}_{1})
(v'_{2}\partial_{v'_{1}}f^{0}-v'_{1}\partial_{v'_{2}}f^{0}))dv$$
$$=-\frac{q}{m}\hat{\rho}(t,k)\ast_{t}[k_{3}\widehat{W}_{2}(\widehat{v'_{2}\partial_{v_{3}}f^{0}}
-\widehat{v_{3}\partial_{v'_{2}}f^{0}})(\eta_{1}+\eta_{k1},\eta_{2}+\eta_{k2},k_{3}t+\eta_{3})$$
$$+k_{3}\widehat{W}_{1}(\widehat{v'_{1}\partial_{v_{3}}f^{0}}
-\widehat{v_{3}\partial_{v'_{1}}f^{0}})(\eta_{1}+\eta_{k1},\eta_{2}+\eta_{k2},k_{3}t+\eta_{3})$$
$$+(k_{1}\widehat{W}_{2}-k_{2}\widehat{W}_{1})(\widehat{v'_{2}\partial_{v_{1}}f^{0}}
-\widehat{v_{1}\partial_{v'_{2}}f^{0}})(\eta_{1}+\eta_{k1},\eta_{2}+\eta_{k2},k_{3}t+\eta_{3})]$$
$$\leq C\int^{t}_{0}e^{-2\pi\lambda_{0}|\eta_{1}+\eta_{k1,t-\tau}|}e^{-2\pi\lambda_{0}|\eta_{2}+\eta_{k2,t-\tau}|}e^{-2\pi\lambda'_{0}|\eta_{k1,\tau}|}
e^{-2\pi\lambda'_{0}|\eta_{k2,\tau}|}e^{-2\pi\lambda'_{0}|k_{3}|\tau}\cdot|k_{3}(t-\tau)+\eta_{3}|e^{-2\pi\lambda_{0}|k_{3}(t-\tau)+\eta_{3}|}d\tau$$
$$\leq C\int^{t}_{0}e^{-2\pi\lambda_{0}|\eta_{1}+\eta_{k1,t-\tau}|}e^{-2\pi\lambda_{0}|\eta_{2}+\eta_{k2,t-\tau}|}e^{-2\pi\lambda'_{0}|\eta_{k1,\tau}|}
e^{-2\pi\lambda'_{0}|\eta_{k2,\tau}|}e^{-2\pi\lambda'_{0}|k_{3}|\tau} e^{-2\pi\frac{(\lambda_{0}+\lambda'_{0})}{2}|k_{3}(t-\tau)+\eta_{3}|}d\tau$$
$$\leq C\bigg(\int^{t}_{0}e^{-2\pi\lambda_{0}|\eta_{1}+\eta_{k1,t-\tau}|}e^{-2\pi\lambda'_{0}|\eta_{k1,\tau}|}d\tau\bigg)^{\frac{1}{3}}
\bigg(\int^{t}_{0}e^{-2\pi\lambda'_{0}|\eta_{k2,\tau}|}e^{-2\pi\lambda_{0}|\eta_{2}+\eta_{k2,t-\tau}|}d\tau\bigg)^{\frac{1}{3}}
 e^{-2\pi\lambda''_{0}|k_{3}t+\eta_{3}|}$$
\begin{align}
%&=C\bigg(\int^{t}_{0}e^{-2\pi\lambda_{0}|\eta_{1}+\frac{1}{\Omega}[-k_{2}\cos\Omega(t-\tau)+k_{2}-k_{1}\sin\Omega(t-\tau)]|}
%e^{-2\pi\lambda'_{0}|\frac{1}{\Omega}[-k_{2}\cos\Omega\tau+k_{2}-k_{1}\sin\Omega\tau]|}d\tau\bigg)^{\frac{1}{3}}\notag\\
%&\bigg(\int^{t}_{0}e^{-2\pi\lambda_{0}|\eta_{2}+\frac{1}{\Omega}[-k_{2}\sin\Omega(t-\tau)-k_{1}+k_{1}\cos\Omega(t-\tau)%]|}
%e^{-2\pi\lambda'_{0}|\frac{1}{\Omega}[-k_{2}\sin\Omega\tau-k_{1}+k_{1}\cos\Omega\tau]|}d\tau\bigg)^{\frac{1}{3}}e^{-2\pi\lambda''_{0}|k_{3}t+\eta_{3}|}\notag\\
&=C\bigg(\int^{t}_{0}e^{-2\pi\lambda_{0}|\eta_{1}+\frac{1}{\Omega}\{k_{2}-\sqrt{k^{2}_{1}+k^{2}_{2}}\sin[\Omega(t-\tau)+\varphi_{1}]\}|}
e^{-2\pi\lambda'_{0}|\frac{1}{\Omega}\{k_{2}-\sqrt{k^{2}_{1}+k^{2}_{2}}\sin[\Omega\tau+\varphi_{1}]\}|}d\tau\bigg)^{\frac{1}{3}}\notag\\
&\bigg(\int^{t}_{0}e^{-2\pi\lambda_{0}|\eta_{2}+\frac{1}{\Omega}[-\sqrt{k^{2}_{1}+k^{2}_{2}}\sin[\Omega(t-\tau)-\varphi_{2}]-k_{1}]|}
e^{-2\pi\lambda'_{0}|\frac{1}{\Omega}[-\sqrt{k^{2}_{1}+k^{2}_{2}}\sin[\Omega\tau-\varphi_{2}]-k_{1}]|}d\tau\bigg)^{\frac{1}{3}}
e^{-2\pi\lambda''_{0}|k_{3}t+\eta_{3}|},\notag\\
%&\mathbf{how}\quad\mathbf{to}\quad\mathbf{compute}\quad\mathbf{detailly}?\leq C e^{-\lambda_{2}|k_{3}t+\eta_{3}|},\notag\\
\end{align}
where $\lambda''_{0}=\lambda'_{0}-\frac{1}{2}(\lambda_{0}-\lambda'_{0}),\tan\varphi_{1}=\frac{k_{2}}{k_{1}},\tan\varphi_{2}=\frac{k_{1}}{k_{2}}.$
\begin{align}
%&VI=-\frac{q}{m}\hat{\rho}(t,k)\ast\partial_{t}\int_{\mathbb{R}^{3}}
%e^{-2\pi i(k_{3}t+\eta_{3})v_{3}}e^{-2\pi iv'_{2}(\eta_{2}+\eta_{k2})}
%e^{-2\pi iv'_{1}(\eta_{1}+\eta_{k1})}(\hat{W}_{1}\partial_{v'_{1}}f^{0}+\hat{W}_{2}\partial_{v'_{2}}f^{0})dv\notag\\
&VI=-\frac{q}{m}\hat{\rho}(t,k)\ast_{t}\partial_{t}\bigg(\hat{W}_{1}\widehat{\partial_{v'_{1}}f^{0}}+\hat{W}_{2}
\widehat{\partial_{v'_{2}}f^{0}}\bigg)(\eta_{1}+\eta_{k1},\eta_{2}+\eta_{k2},k_{3}t+\eta_{3})\notag\\
%&\leq\int^{t}_{0}C(W,k,\lambda,\lambda_{1})e^{-\lambda|k_{3}|\tau}|k_{3}(t-\tau)+\eta_{3}||k|e^{-\lambda_{1}|k_{3}(t-\tau)+\eta_{3}|}d\tau\notag\\
%&\leq C e^{-\lambda_{2}|k_{3}t+\eta_{3}|}.
&\leq C\bigg(\int^{t}_{0}e^{-2\pi\lambda_{0}|\eta_{1}+\frac{1}{\Omega}\{k_{2}-\sqrt{k^{2}_{1}+k^{2}_{2}}\sin[\Omega(t-\tau)+\varphi_{1}]\}|}
e^{-2\pi\lambda'_{0}|\frac{1}{\Omega}\{k_{2}-\sqrt{k^{2}_{1}+k^{2}_{2}}\sin[\Omega\tau+\varphi_{1}]\}|}d\tau\bigg)^{\frac{1}{3}}\notag\\
&\bigg(\int^{t}_{0}e^{-2\pi\lambda_{0}|\eta_{2}+\frac{1}{\Omega}[-\sqrt{k^{2}_{1}+k^{2}_{2}}\sin[\Omega(t-\tau)-\varphi_{2}]-k_{1}]|}
e^{-2\pi\lambda'_{0}|\frac{1}{\Omega}[-\sqrt{k^{2}_{1}+k^{2}_{2}}\sin[\Omega\tau-\varphi_{2}]-k_{1}]|}d\tau\bigg)^{\frac{1}{3}}
e^{-2\pi\lambda''_{0}|k_{3}t+\eta_{3}|}.\notag\\
\end{align}
%Now we consider $k_{z}=0,$ since $E(t,x,y,-z)=-E(t,x,y,z),$ and $\partial_{t}B(t,x,y,z)=\nabla\times E(t,x,y,z),B(0,x,y,z)=-B(0,x,y-z)$
%then$B(t,x,y,z)=B(t,x,y,-z).$
%Meanwhile, we assume $f_{0}(x,y,z,v)=-f_{0}(x,y,-z,v),$ then $\hat{\rho}(t,k_{\perp},0)=0.$
%If $k_{z}\neq0,k_{\perp}=0,$ the case is the same to Landau damping. Therefore,here we omit the details.
From (1.23)-(1.26) and Corollary 1.2, the results of Theorem 1.1 are obvious.

 $$$$
         $$$$
        % \begin{center}
%\item\section{$B_{0}\rightarrow0$}
%{\bf\large 1. \quad Introduction }
%\end{center}

%First we still write the euqation
%\begin{equation}
%\left\{\begin{array}{l}
%\partial_{t}f+v\cdot\nabla_{x}f+\frac{q}{m}(v\times B_{0})\cdot\nabla_{v}f=-\frac{q}{m}(E+v\times B_{1})\cdot\nabla_{v}f^{0},\\
%\partial_{t}B_{1}=\nabla_{x}\times E,\\
%E=W(x)\ast\rho(t,x),\rho(t,x)=\int_{\mathbb{R}^{3}}f(t,x,v)dv,\\
%f(0,x,v)=f_{0}(x,v_{\perp},v_{z}),f^{0}(v)=f^{0}(v_{\perp},v_{z}),
%\end{array}\right.
% \end{equation}

           \begin{center}
\item\section{Notation and Hybrid analytic norm}
%{\bf\large 1. \quad Introduction }
\end{center}
         %$\mathbf{Step }\mathbf{2}.$  $\mathbf{Hybrid}$ $\mathbf{analytic}$ $\mathbf{norm}$

In order to prove our result in the nonlinear case, we have to introduce the hybrid analytic norm that is one of the cornerstones of our analysis, because they will
connect well to both estimates  in $ x$  on the force field and uniform estimates in $ v. $
First, let us recall that the free transport equation in a constant magnetic field
\begin{align}
\partial_{t}f+v\cdot\nabla_{x}f+\frac{q}{m}(v\times B_{0})\cdot\nabla_{v}f=0
\end{align}
has a strong mixing property: any solution of (1.1) converges weakly in large time to a spatially homogeneous distribution equal to the space-averaging of the initial datum. Let us sketch the proof.

If $ f $ solves (2.1) in $\mathbb{T}^{3}\times\mathbb{R}^{3},$ with initial datum $f_{in}=f(\tau,\cdot),$ then
\begin{align}
f(t,x'_{t},v'_{t})=f_{in}(x'_{\tau}-v'\cdot(\frac{\sin\Omega t}{\Omega},\frac{\sin\Omega t}{\Omega},t),v'_{\tau}),
\end{align}

We give a description on the motion of charged particles in the following way:
\begin{align}
\frac{dx'}{dt}=v',\quad \frac{dv'}{dt}=\frac{q}{m}v'\times B_{0},
\end{align}
where $B_{0}=B_{0}\hat{z}.$

% that the inductive hypothesis hold for the step $n+1,$  first we study the dynamics of the motion of particles in (2.2).
We assume that $x'(t,x,v)=x=(x_{1},x_{2},x_{3}),v'(t,x,v)=v=(v_{x},v_{y},v_{3})=(v_{\perp}\cos\theta,v_{\perp}\sin\theta,v_{3})$ at $t=\tau,$
then the solution of Eq.(2.1) is obtained as follows:
\begin{align}
&v'_{x}(t)=v_{\perp}\cos(\theta+\Omega (t-\tau)),\quad v'_{y}(t)=v_{\perp}\sin(\theta+\Omega (t-\tau)),\quad v'_{3}=v_{3};\notag\\
&x'_{1}(t)=x_{1}+\frac{v_{\perp}}{\Omega}[\sin(\theta+\Omega (t-\tau))-\sin\theta],\notag\\
&x'_{2}(t)=x_{2}-\frac{v_{\perp}}{\Omega}[\cos(\theta+\Omega (t-\tau))-\cos\theta],\quad x'_{3}(t)=x_{3}+v_{3}(t-\tau),
\end{align}
where $\Omega=\frac{qB_{0}}{m},v_{\perp}=\sqrt{v^{2}_{1}+v^{2}_{2}}.$ From now on, without loss of generality, we assume $m=q=1.$
 Also through simple computation, it is easy to get
  \begin{align}
\begin{bmatrix}
 v'_{t1}\\
v'_{t2}\\
 v'_{t3}
 \end{bmatrix}=\begin{bmatrix}
 \cos\Omega(t-\tau)\quad-\sin\Omega(t-\tau)\quad 0\\
  \sin\Omega(t-\tau)\quad\cos\Omega(t-\tau)\quad\quad 0\\
  0\quad\quad\quad\quad\quad\quad 0\quad\quad\quad\quad\quad \quad\quad1
 \end{bmatrix}\begin{bmatrix}
 v'_{\tau 1}\\
 v'_{\tau 2}\\
v'_{\tau 3}
 \end{bmatrix}
 \triangleq \mathcal{R}(t-\tau)v'_{\tau},
 \end{align}

 $$\begin{bmatrix}
 x'_{t1}\\
x'_{t2}\\
 x'_{t3}
 \end{bmatrix}=\begin{bmatrix}
 x_{1}\\
x_{2}\\
 x_{3}
 \end{bmatrix}+\begin{bmatrix}
 \frac{1}{\Omega}\sin\Omega(t-\tau)\quad \frac{1}{\Omega}\cos\Omega(t-\tau)-\frac{1}{\Omega}\quad\quad 0\\
  -\frac{1}{\Omega}\cos\Omega(t-\tau)+ \frac{1}{\Omega}\quad \frac{1}{\Omega}\sin\Omega(t-\tau)\quad\quad \quad 0\\
  0\quad\quad\quad\quad\quad\quad 0\quad\quad\quad\quad\quad \quad \quad \quad \quad (t-\tau)
 \end{bmatrix}\begin{bmatrix}
  v_{1}\\
 v_{2}\\
v_{3}
 \end{bmatrix}$$
$$
 \triangleq\begin{bmatrix}
 x_{1}\\
x_{2}\\
 x_{3}
 \end{bmatrix}+\mathcal{M}(t-\tau)\begin{bmatrix}
  v_{1}\\
 v_{2}\\
v_{3}
 \end{bmatrix}$$
  \begin{align}
( x'_{t1},x'_{t2}, x'_{t3})=(x'_{\tau 1},x'_{\tau 2}, x'_{\tau 3})+
 \bigg(\frac{v'_{\tau 1}\sin\Omega(t-\tau)}{\Omega},\frac{v'_{\tau 2}\sin\Omega(t-\tau)}{\Omega},v'_{\tau 3}(t-\tau)\bigg).\end{align}

%and
%\begin{align}
%\nabla_{v'}x'=(\frac{\sin\Omega(t-\tau)}{\Omega},\frac{\sin\Omega(t-\tau)}{\Omega},(t-\tau)).
%\end{align}

We introduce an equivalence relation, that is, for any different velocity $v'_{t_{1}},v'_{t_{2}},$ there exists an orthogonal matrix $\mathcal{O}(t_{1},t_{2})$
such that $v'_{t_{1}}=\mathcal{O}(t_{1},t_{2})v'_{t_{2}},$ we say $v'_{t_{1}}\sim v'_{t_{2}}.$ And all elements satisfying the above equivalence relation are denoted by $[v'_{t}],$ for short $v'_{t}.$

           \begin{de} From the system (1.1)and the above equivalence relation, we can define the corresponding transform $$ S^{0}_{t,\tau}(x',v')\triangleq\bigg(x'_{\tau}+\bigg(\frac{v'_{\tau 1}\sin\Omega(t-\tau)}{\Omega},\frac{v'_{\tau 2}\sin\Omega(t-\tau)}{\Omega},v'_{\tau 3}(t-\tau)\bigg),v'_{\tau}\bigg),$$
           where $\mathcal{M}(t-\tau),\mathcal{R}(t-\tau)$ are defined in (2.5).
           \end{de}
           \begin{rem}
           From the above equality (2.6), we can observe clearly  the connection and difference on between Landau damping  and cyclotron damping. Indeed, it can be reduced to Landau damping as $B_{0}\rightarrow0.$

            In addition,  from the dynamics system of the above order differential system, it is known that   $S^{0}_{t,\tau}$ satisfies $$S^{0}_{t_{2},t_{3}}\circ S^{0}_{t_{1},t_{2}}=S^{0}_{t_{1},t_{3}}.$$
           % that is, $$S^{0}_{t_{2},t_{3}}\circ S^{0}_{t_{1},t_{2}}(x,v)=S^{0}_{t_{2},t_{3}}(x+\mathcal{M}(t_{2}-t_{1})v,\mathcal{R}(t_{2}-t_{1})v)$$
            %$$=(x+\mathcal{M}(t_{3}-t_{2}+t_{2}-t_{1})v,\mathcal{R}(t_{3}-t_{2}+t_{2}-t_{1})v)$$
            %$$=(x+\mathcal{M}(t_{3}-t_{1})v,\mathcal{R}(t_{3}-t_{1})v)$$
           % $$=S^{0}_{t_{1},t_{3}}(x,v).$$
           \end{rem}

           To estimate solutions and trajectories of kinetic equations, maybe we have to work on the phase space $\mathbb{T}_{x}^{3}\times\mathbb{R}_{v}^{3}.$ And we also
           use the following three parameters: $\lambda$(gliding analytic regularity), $\mu$(analytic regularity in $ x $) and $\tau$(time-shift along the free
           transport semigroup). From $Remark $ 2.2, we know that the linear Vlasov equation has the property of the free transport semigroup. This
            property is crucial to our
           analysis. In this paper, one of the cornerstones of our analysis  is to compare the solution of the nonlinear
           case at time $\tau$ with the solution of the linear case.

           Now we start to introduce the very important tools in our paper. These are time-shift pure and hybrid analytic norms. They are similar with those
           in the paper [23]
           written by Mouhot and Villani.

 First, we introduce some notations. We denote $\mathbb{T}^{3}=\mathbb{R}^{3}/\mathbb{Z}^{3}.$ For function $f(x'_{\tau},v'_{\tau}),$ we define the Fourier transform $\hat{f}(k,\eta),$
where $(k,\eta)\in\mathbb{Z}^{3}\times\mathbb{R}^{3},$ via
$$\hat{f}(k,\eta)=\int_{\mathbb{T}^{3}\times\mathbb{R}^{3}}e^{-2i\pi x'_{\tau}\cdot k}e^{-2i\pi v'_{\tau}\cdot \eta}f(x'_{\tau},v'_{\tau})dx'_{\tau}dv'_{\tau}, $$
$$\tilde{f}(\omega,k,\eta)=\int_{\mathbb{R}^{+}}e^{2i\pi t\cdot\omega}\int_{\mathbb{T}^{3}\times\mathbb{R}^{3}}
e^{-2i\pi x'_{\tau}\cdot k}e^{-2i\pi v'_{\tau}\cdot \eta}f(x'_{\tau},v'_{\tau})dx'_{\tau}dv'_{\tau}dt.$$

We also write

$$k=(k_{1},k_{2},k_{3})=(k_{\perp}\cos\varphi,k_{\perp}\sin\varphi,k_{3}),\quad
\eta=(\eta_{1},\eta_{2},\eta_{3})=(\eta_{\perp}\cos\gamma,\eta_{\perp}\sin\gamma,\eta_{3}).$$

Now we define some notations
$$\tilde{f}=\tilde{f}(\omega,k,\eta),\quad\hat{f}=\hat{f}(t,k,\eta),\quad\tilde{\rho}=\tilde{\rho}(\omega,k),\quad\hat{\rho}=\hat{\rho}(t,k).$$

            \begin{de}(Hybrid analytic norms)
            $$\|f\|_{\mathcal{C}^{\lambda,\mu}}=\sum_{m,n\in\mathbb{N}_{0}^{3}}\frac{\lambda^{n}}{n!}\frac{\mu^{m}}{m!}
            \|\nabla^{m}_{x'_{\tau}}\nabla_{v'_{\tau}}^{n}f\|_{L^{\infty}(\mathbb{T}_{x'}^{3}\times\mathbb{R}^{3}_{v'})},\quad
            \|f\|_{\mathcal{F}^{\lambda,\mu}}=\sum_{k\in\mathbb{Z}^{3}}\int_{\mathbb{R}^{3}}|\tilde{f}(k,\eta)|
           e^{2\pi\lambda|\eta|}e^{2\pi\mu|k|}d\eta,$$
            $$\|f\|_{\mathcal{Z}^{\lambda,\mu}}=\sum_{l\in\mathbb{Z}^{3}}\sum_{n\in\mathbb{N}_{0}^{3}}
            \frac{\lambda^{n}}{n!}e^{2\pi\mu|l|}\|\widehat{\nabla_{v'_{\tau}}^{n}f(l,v)}\|_{L^{\infty}(\mathbb{R}_{v'}^{3})}.$$
            \end{de}
           \begin{de} (Time-shift pure and hybrid analytic norms)
For any $\lambda,\mu\geq0,p\in[1,\infty],$ we define
            $$\|f\|_{\mathcal{C}^{\lambda,\mu}_{t,\tau}}=\|f\circ S^{0}_{t,\tau}(x',v')\|_{\mathcal{C}^{\lambda,\mu}}=\sum_{m,n\in\mathbb{N}_{0}^{3}}\frac{\lambda^{n}}{n!}\frac{\mu^{m}}{m!}
            \|\nabla^{m}_{x'_{\tau}}[\nabla_{v'_{\tau}}+\bigg(\frac{\sin\Omega(t-\tau)}{\Omega},\frac{\sin\Omega(t-\tau)}{\Omega},t-\tau\bigg)
            \cdot\nabla_{x'_{\tau}}]^{n}f\|_{L^{\infty}(\mathbb{T}_{x}^{3}\times\mathbb{R}^{3}_{v})},$$
           $$\|f\|_{\mathcal{F}^{\lambda,\mu}_{t,\tau}}=\|f\circ S^{0}_{t,\tau}(x',v')\|_{\mathcal{F}^{\lambda,\mu}}=\sum_{k\in\mathbb{Z}^{3}}\int_{\mathbb{R}^{3}}|\tilde{f}(k,\eta)|
           e^{2\pi\lambda|\eta+k\cdot(\frac{\sin\Omega(t-\tau)}{\Omega},\frac{\sin\Omega(t-\tau)}{\Omega},(t-\tau))|}e^{2\pi\mu|k|}d\eta,$$
          $$\|f\|_{\mathcal{Z}^{\lambda,\mu}_{t,\tau}}=\|f\circ S^{0}_{t,\tau}(x,v)\|_{\mathcal{Z}^{\lambda,\mu}}=\sum_{l\in\mathbb{Z}^{3}}\sum_{n\in\mathbb{N}_{0}^{3}}
            \frac{\lambda^{n}}{n!}e^{2\pi\mu|l|}\|[\nabla_{v'_{\tau}}+2i\pi\bigg(\frac{\sin\Omega(t-\tau)}{\Omega},\frac{\sin\Omega(t-\tau)}{\Omega},t-\tau\bigg)\cdot l]^{n}\hat{f}(l,v)\|
            _{L^{\infty}(\mathbb{R}_{v}^{3})},$$
               $$\|f\|_{\mathcal{Z}^{\lambda,\mu;p}_{t,\tau}}=\sum_{l\in\mathbb{Z}^{3}}\sum_{n\in\mathbb{N}_{0}^{3}}
            \frac{\lambda^{n}}{n!}e^{2\pi\mu|l|}\|[\nabla_{v'_{\tau}}+2i\pi\bigg(\frac{\sin\Omega(t-\tau)}{\Omega},\frac{\sin\Omega(t-\tau)}{\Omega},t-\tau\bigg)\cdot l]^{n}\hat{f}(l,v)\|_{L^{p}(\mathbb{R}_{v'}^{3})},$$
            $$\|f\|_{\mathcal{Y}^{\lambda,\mu}_{t,\tau}}=\|f\|_{\mathcal{F}^{\lambda,\mu;\infty}_{t,\tau}}
            =\sup_{k\in\mathbb{Z}^{3},\eta\in\mathbb{R}^{3}}e^{2\pi\mu|k|}
            e^{2\pi\lambda|\eta+k\cdot(\frac{\sin\Omega(t-\tau)}{\Omega},\frac{\sin\Omega(t-\tau)}{\Omega},(t-\tau))|}
            |\hat{f}(k,\eta)|.$$
           \end{de}

           From the above  definitions, we can state some simple and important propositions, and the related proofs can be found in [23], so we remove the proofs. From
           (2.6), we know that the damping occurs only in the $\hat{z}$ direction. Therefore, we will mainly focus on the third component of the ``phase space"
            $(x,v)$ when referring to the damping mechanics of the wave propagation.

\begin{prop} For any $\tau\in\mathbb{R},\lambda,\mu\geq 0,$
\begin{itemize}
 \item[{(i)}]   if $f$ is a function only of $x,$  that is, $t=\tau$ in the second equality of (2.5), then
           $\|f\|_{\mathcal{C}^{\lambda,\mu}_{\tau}}=\|f\|_{\mathcal{C}^{\lambda|\tau|+\mu}},\|f\|_{\mathcal{F}^{\lambda,\mu}_{\tau}}
           =\|f\|_{\mathcal{Z}^{\lambda,\mu}_{\tau}}=\|f\|_{\mathcal{F}^{\lambda|\tau|+\mu}};$
          \item[{(ii)}]    if $f$ is a function only of $v,$ then
           $\|f\|_{\mathcal{C}^{\lambda,\mu;p}_{\tau}}=\|f\|_{\mathcal{Z}^{\lambda,\mu;p}_{\tau}}=\|f\|_{\mathcal{C}^{\lambda,;p}};$
            \item[{(iii)}] for any $\lambda>0,$ then
            $\|f\circ(Id+G)\|_{\mathcal{F}^{\lambda}}\leq\|f\|_{\mathcal{F}^{\lambda+\nu}},\nu=\|G\|_{\dot{\mathcal{F}}^{\lambda}};$
            \item[{(iv)}] for any $\bar{\lambda}>\lambda,p\in[1,\infty],$ $\|\nabla f\|_{\mathcal{C}^{\lambda;p}}\leq
            \frac{1}{\lambda e\log(\frac{\bar{\lambda}}{\lambda} )}
            \| f\|_{\mathcal{C}^{\bar{\lambda};p}},$
            $\|\nabla f\|_{\mathcal{F}^{\lambda;p}}\leq\frac{1}{2\pi e(\bar{\lambda}-\lambda )}
            \| f\|_{\mathcal{F}^{\bar{\lambda};p}},$
            \item[{(v)}]  for any $\bar{\lambda}>\lambda>0,\mu>0,$ then
            $\|vf\|_{\mathcal{Z}^{\lambda,\mu;1}_{\tau}}\leq\|f\|_{\mathcal{Z}^{\bar{\lambda},\mu;1}_{\tau}};$
           \item[{(vi)}]  for any $\bar{\lambda}>\lambda,\bar{\mu}>\mu,$
          $\|\nabla_{v}f\|_{\mathcal{Z}_{\tau}^{\lambda,\mu;p}}
          \leq C(d)\bigg(\frac{1}{\lambda\log(\frac{\bar{\lambda}}{\lambda})}\|f\|_{\mathcal{Z}_{\tau}^{\bar{\lambda},\bar{\mu};p}}+\frac{\tau}{\bar{\mu}-\mu}
 \|f\|_{\dot{\mathcal{Z}}_{\tau}^{\bar{\lambda},\bar{\mu};p}}\bigg);$
 \item[{(vii)}] for any $\bar{\lambda}>\lambda,$
 $\|(\nabla_{v}+\tau\nabla_{x})f\|_{\mathcal{Z}_{\tau}^{\lambda,\mu;p}}\leq \frac{1}{C(d)\lambda\log(\frac{\bar{\lambda}}{\lambda})}
 \|f\|_{\mathcal{Z}_{\tau}^{\bar{\lambda},\mu;p}};$
  \item[{(viii)}] for any $\bar{\lambda}\geq\lambda\geq0,\bar{\mu}\geq\mu\geq0,$ then $\|f\|_{\mathcal{Z}_{\tau}^{\lambda,\mu}}
  \leq_{\mathcal{Z}_{\tau}^{\bar{\lambda},\bar{\mu}}}.$ Moreover, for any $\tau,\bar{\tau}\in\mathbb{R},$ $p\in[1,\infty],$ we have
  $\|f\|_{\mathcal{Z}_{\tau}^{\lambda,\mu;p}}\leq\|f\|_{\mathcal{Z}_{\bar{\tau}}^{\lambda,\mu+\lambda|\tau-\bar{\tau}|;p}};$
  \item[{(viiii)}] $\|f\|_{\mathcal{Y}^{\lambda,\mu}_{\tau}}\leq\|f\|_{\mathcal{Z}^{\lambda,\mu;1}_{\tau}};$
  \item[{(ix)}] for any function $f=f(x,v),$ $\|\int_{\mathbb{R}^{3}}fdv\|_{\mathcal{F}^{\lambda|\tau|+\mu}}\leq\|f\|_{\mathcal{Z}^{\lambda,\mu;1}_{\tau}}.$
 \end{itemize}
\end{prop}

$Proof.$ Here we only give the proof of (v). By the invariance under the action of free transport, it is sufficient to do the proof for $t=0.$
Applying the Fourier transform formula, we have
$$\nabla^{m}_{v}(v\hat{f})(k,v)=\int_{\mathbb{R}}\partial_{\eta}\tilde{f}(k,\eta)(2i\pi\eta)^{m}e^{2i\pi\eta\cdot v}d\eta,$$
and therefore
$$\sum_{m\in\mathbb{N}_{0}}\frac{\lambda^{m}}{m!}\int_{\mathbb{R}}\bigg|\int_{\mathbb{R}}
\partial_{\eta}\tilde{f}(k,\eta)(2i\pi\eta)^{m}e^{2i\pi\eta\cdot v}d\eta\bigg|dv
\leq\sum_{m\in\mathbb{N}_{0}}\frac{\lambda^{m}}{m!}\sup_{\eta\in\mathbb{R}}|(2\pi\eta)^{m}\partial_{\eta}\tilde{f}(k,\eta)|
\int_{\mathbb{R}}\bigg|\int_{\mathbb{R}}e^{2i\pi\eta\cdot v}d\eta\bigg|dv$$
%$$=\sum_{m\in\mathbb{N}_{0}}\frac{\lambda^{m}}{m!}\sup_{\eta\in\mathbb{R}}|(2\pi\eta)^{m}\partial_{\eta}\tilde{f}(k,\eta)|
%\leq\sum_{m\in\mathbb{N}_{0}}\frac{\lambda^{m}}{m!}\sup_{\eta\in\mathbb{R}}|\partial_{\eta}[(2\pi\eta)^{m}\tilde{f}(k,\eta)]-2\pi m(2\pi\eta)^{m-1}\tilde{f}(k,\eta)|$$
$$\leq C(\lambda,\bar{\lambda})\sum_{m\in\mathbb{N}_{0}}\frac{\bar{\lambda}^{m}}{m!}\sup_{\eta\in\mathbb{R}}|(2\pi\eta)^{m}\tilde{f}(k,\eta)|\leq
C(\lambda,\bar{\lambda})\sum_{m\in\mathbb{N}_{0}}
\frac{\bar{\lambda}^{m}}{m!}\int_{\mathbb{R}}|\nabla^{m}_{v}\hat{f}(k,v)|dv,$$
where the last second inequality uses the property (iv),
then we have $$e^{2\pi\mu|k|}\sum_{m\in\mathbb{N}_{0}}\frac{\lambda^{m}}{m!}\int_{\mathbb{R}}\bigg|\int_{\mathbb{R}}
\partial_{\eta}\tilde{f}(k,\eta)(2i\pi\eta)^{m}e^{2i\pi\eta\cdot v}d\eta\bigg|dv
\leq e^{2\pi\mu|k|}\sum_{m\in\mathbb{N}_{0}}\frac{\bar{\lambda}^{m}}{m!}\int_{\mathbb{R}}|\nabla^{m+1}_{v}\hat{f}(k,v)|dv.$$
This establishes (v).
           \begin{prop} For any $X\in\{\mathcal{C},\mathcal{F},\mathcal{Z}\}$ and any $t,\tau\in\mathbb{R},$
           $$\|f\circ S^{0}_{\tau}\|_{X^{\lambda,\mu}_{\tau}}=\|f\|_{X^{\lambda,\mu}_{t+\tau}}.$$

           \end{prop}
           \begin{lem} Let $\lambda,\mu\geq0,t\in\mathbb{R},$ and consider two functions $F,G:\mathbb{T}\times\mathbb{R}\rightarrow\mathbb{T}\times\mathbb{R}.$
           Then there is $ \varepsilon\in(0,\frac{1}{2})$ such that  if $F,G$ satisfy
\begin{align}
\|\nabla(F-Id)\|_{\mathcal{Z}^{\lambda',\mu'}_{\tau}}\leq\varepsilon,
\end{align}
where $\lambda'=\lambda +2\|F-G\|_{\mathcal{Z}^{\lambda ,\mu}_{\tau}},\quad\mu'=\mu+2(1+|\tau|)\|F-G\|_{\mathcal{Z}^{\lambda ,\mu}_{\tau}} ,$ then $F$ is invertible  and
\begin{align}
&\|F^{-1}\circ G-Id\|_{\mathcal{Z}^{\lambda ,\mu}_{\tau}}\leq2\|F-G\|_{\mathcal{Z}^{\lambda ,\mu}_{\tau}}.
\end{align}
\end{lem}

\begin{prop} For any $\lambda,\mu\geq0$ and any $p\in[1,\infty],\tau\in\mathbb{R},\sigma\in\mathbb{R},a\in\mathbb{R}\setminus\{0\}$ and $b\in\mathbb{R},$
we have
$$\|f(x+bv'+X(x,v'),av'+V(x,v'))\|_{\mathcal{Z}^{\lambda,\mu;p}_{\tau}}\leq|a|^{-\frac{3}{p}}\|f\|_{\mathcal{Z}^{\alpha,\beta;p}_{\sigma}},$$
where $\alpha=\lambda|a|+\|V\|_{\mathcal{Z}^{\lambda,\mu}_{\tau}},\quad \beta=\mu+\lambda|b+\tau-a\sigma|+\|X-\sigma V\|_{\mathcal{Z}^{\lambda,\mu}_{\tau}}.$
\end{prop}
\begin{lem} Let $G=G(x,v)$ and $R=R(x,v)$ be valued in $\mathbb{R},$ and
$\beta(x)=\int_{\mathbb{R}}(G\cdot R)(x',v')dv'.$
Then for any $\lambda,\mu,t\geq0$ and any $b>-1,$ we have
$$\|\beta\|_{\mathcal{F}^{\lambda t+\mu}}\leq3\|G\|_{\mathcal{Z}_{\tau-\frac{bt}{1+b}}^{\lambda(1+b),\mu;1}}
\|R\|_{\mathcal{Z}_{\tau-\frac{bt}{1+b}}^{\lambda(1+b),\mu}}.$$
\end{lem}
%To explain our result in this paper, we need to quantify regularity corrections to the analytic regularity in the $ x $ variable.
%\begin{de} For $\lambda,\mu,\gamma,\eta\geq0,\tau\in\mathbb{R},p\in[1,\infty],$ we define
%$$\|f\|_{\mathcal{Z}^{(\lambda,\eta),(\mu,\gamma);p}_{\tau}}=\sum_{l\in\mathbb{Z}^{d}}\sum_{n\in\mathbb{N}_{0}^{d}}
%\frac{\lambda^{n}}{n!}e^{2\pi\mu|l|}(1+|l|)^{\gamma}\|(\nabla_{v}+2i\pi\tau l)^{n}\hat{f}(l,v)(1+|v|)^{\mu}\|_{L^{p}(\mathbb{R}_{v}^{d})},$$
%$$\|f\|_{\mathcal{F}^{(\mu,\gamma);p}_{\tau}}=\sum_{l\in\mathbb{Z}^{d}}
%e^{2\pi\mu|l|}(1+|l|)^{\gamma}\|\hat{f}(l,v)\|_{L^{p}}.$$
%\end{de}

\begin{center}
\item\section{Linear cyclotron damping revisited}
%{\bf\large 1. \quad Introduction }
\end{center}

In this section, we recast the linear damping in the hybrid analytic norms. \begin{thm}
 For any $\eta,v\in\mathbb{R}^{3},k\in\mathbb{N}_{0}^{3},$  we assume that the following conditions hold in equations (0.9).
%$\quad\mathbf{ how} \quad\mathbf{ to}\quad \mathbf{make} \quad \mathbf{assumptions }$
%$\quad \mathbf{on}\quad f^{0}$ is suitable
\begin{itemize}
        \item[(i)] $ W(x)$ is an odd function on $x_{3},$  $|\widehat{W}(k)|\leq\frac{1}{1+|k|^{\gamma}},\gamma>1,$ where $ W(x)=(W_{1}(x),W_{2}(x),0);$
%\item[(i)]$f^{0}(v)$ is  Maxwellian, that is, $f^{0}(v)=f^{0}(v_{\perp},v_{z})=C_{M}exp(-\alpha(v^{2}_{\perp}+v^{2}_{z}))$ for some $\alpha>0$  ;
         \item[(ii)] $\|\nabla_{v}f^{0}\|_{\mathcal{C}^{\lambda_{0};1}}\leq C^{0};\|f_{0}\|_{\mathcal{Z}^{\lambda_{0},\mu_{0};1}_{0}}\leq\delta_{0}$
              for some constants $\lambda_{0},\alpha_{0},C^{0}>0,\delta_{0}>0;$
             \item[(iii)] In addition, $(\mathbf{PSC})$ holds,

$(\mathbf{PSC}):$ for any component velocity in the $\hat{z}$ direction $v_{3}\in\mathbb{R},$
             there exists some positive constant $v_{Te}$ such that if $v_{3}=\frac{\omega}{k_{3}} \textmd{ when } k_{3}\neq0; \textmd{ or } k_{3}=0$ where $\omega,k$ are frequencies of time and space $t,x,$
             respectively,  then $|v_{3}|\gg v_{Te}.$
                  %\item[(v)] $f^{0}(\eta_{1},\eta_{2},v_{3})\leq Ce^{-\lambda_{1}|\eta_{1}|}e^{-\lambda_{1}|\eta_{2}|}e^{-\alpha_{1}|v_{3}|},$ otherwise $f^{0}(\eta_{1},\eta_{2},v_{3})<\varepsilon(v_{3}),$ for some $0<\varepsilon(v_{3})<1$ sufficient small.%$\inf_{k\in\mathbb{Z}^{3}}|\tilde{\mathcal{L}}(\omega,k)-1|>\kappa,$ for some $0<\kappa<1.$
         \end{itemize}
         Then %if $f_{0}(x,v),f^{0}(v)$ are axisymetric, then the solution $f$ of (1.1) is also axisymmetric.
  for any $\lambda'_{0}<\lambda_{0},$ we have
         \begin{align}
         &\sup_{t\geq 0}\|\rho(t,\cdot)\|_{\mathcal{F}^{\lambda'_{0}t+\mu_{0}}}\leq C(C^{0},\lambda_{0},\lambda'_{0},\mu)\delta_{0},\notag\\
         &\sup_{t\geq0}\|f-f^{0}\|_{\mathcal{Z}^{\lambda''_{0},\mu;1}_{t}}\leq C(C^{0},\lambda_{0},\lambda'_{0},\mu)\delta_{0}.
         \end{align}
         %(这里可能不需要假设初值为轴对称在速度方向 %maybe not assume $f_{0}$ axisymmetric in $v$)
         \end{thm}

%\begin{align}
%&\mathcal{L}(t,k)=\frac{q}{m}\int^{t}_{0}k_{3}(\hat{W}_{2}(i\eta_{k_{1}})\partial_{\eta_{3}}
%\hat{f}^{0}(\eta_{k_{1}},\eta_{k_{2}},k_{3}t)+\hat{W}_{1}(i\eta_{k_{2}})\partial_{\eta_{3}}\hat{f}^{0}(\eta_{k_{1}},\eta_{k_{2}},k_{3}t))dt\notag\\
%&-\frac{q}{m}
%(\hat{W}_{1}(i\eta_{k_{1}})\hat{f}^{0}(\eta_{k_{1}},\eta_{k_{2}},k_{3}t)+\hat{W}_{2}(i\eta_{k_{2}})\hat{f}^{0}(\eta_{k_{1}},\eta_{k_{2}},k_{3}t)),\notag\\
%\end{align}
%where $\eta_{k_{1}}=\frac{1}{\Omega}(-k_{2}\cos\Omega t+k_{2}-k_{1}\sin\Omega t),\eta_{k_{2}}=\frac{1}{\Omega}(-k_{2}\sin\Omega t-k_{1}+k_{1}\cos\Omega t).$

         The revisited proof of Lemma 1.1

          $Proof.$  First $\omega\neq0,k_{3}\neq0,$
%\begin{align}
%&\tilde{\rho}(\omega,k)=\int_{\mathbb{R}^{+}}\int_{\mathbb{T}^{3}}\int_{\mathbb{R}^{3}}e^{2\pi it\omega}e^{-2\pi ix\cdot k}f(t,x,v)dxdvdt\notag\\
%&=\int_{\mathbb{R}^{+}}\int_{\mathbb{T}^{3}}\int_{\mathbb{R}^{3}}e^{2\pi it\omega}e^{-2\pi i(x_{1},x_{2})\cdot (k_{1},k_{2})}f_{0}(x'(0,x,v),v'(0,x,v))dvdxdt-\frac{q}{m}\int_{\mathbb{R}^{+}}\int_{\mathbb{T}^{3}}\int^{t}_{0}\int_{\mathbb{R}^{3}}\notag\\
%&\cdot e^{2\pi it\omega}e^{-2\pi ix\cdot k}[(E+v'(\tau,x,v)\times B)\cdot\nabla'_{v}f^{0}](\tau,x'(\tau,x,v),v'(\tau,x,v))dvdxd\tau dt.\notag\\
%\end{align}
%Recall $E(t,x)=W(x)\ast\rho(t,x),$ $\partial_{t}B(t,x)=\nabla_{x}\times E(t,x),$ then taking the Fourier transform in the variables $t,x,$
%$$\tilde{E}(\omega,k)=\widehat{W}(k)\tilde{\rho}(\omega,k),\quad\omega\tilde{B}(\omega,k)=k\times\tilde{E}(\omega,k).$$
%Furthermore, we get
%\begin{align}
%&v\times\tilde{B}(\omega,k)=\frac{1}{\omega}[v\times(k\times\tilde{E}(\omega,k))]\notag\\
%&=\frac{1}{\omega}\bigg(v_{2}(k_{1}\widehat{W}_{2}-k_{2}\widehat{W}_{1})-v_{3}k_{3}\widehat{W}_{1},-v_{1}(k_{1}
%\widehat{W}_{2}-k_{2}\widehat{W}_{1})+v_{3}k_{3}\widehat{W}_{2},
%v_{1}k_{3}\widehat{W}_{1}+v_{2}k_{3}\widehat{W}_{2}\bigg)\tilde{\rho}(\omega,k).\notag\\
%\end{align}
%Combining (1.3)-(1.4) and (2.3)-(2.4), and note that $dv\rightarrow dv'$ preserves the measure,  we can change between $dv$ and $dv'$ whenever we need, but in order to
%simply the notations, we don't differentiate the notations $dv$ and $dv'$  in this paper,
 we have

$$\tilde{\rho}(\omega,k)=\int_{\mathbb{R}^{+}}\int_{\mathbb{T}^{3}}\int_{\mathbb{R}^{3}}e^{2\pi it\omega}e^{-2\pi ix\cdot k}f_{0}(x'(0,x,v),v'(0,x,v))dvdxdt$$
$$+\frac{q}{m}\tilde{\rho}(\omega,k)\frac{1}{\omega}\int_{\mathbb{R}^{+}}\int_{\mathbb{R}^{3}}
e^{2\pi it\omega}e^{-2\pi ik\cdot(v'_{t1}\frac{\sin\Omega t}{t},v'_{t2}\frac{\sin\Omega t}{t},v_{3}t)}
(k_{3}v'_{3})(\widehat{W}_{1}\partial_{v'_{1}}f^{0}-\widehat{W}_{2}\partial_{v'_{2}}f^{0})dvdt$$
$$+\frac{q}{m}\tilde{\rho}(\omega,k)\frac{1}{\omega}\int_{\mathbb{R}^{+}}\int_{\mathbb{R}^{3}}e^{2\pi it\omega}e^{-2\pi ik\cdot(v'_{t1}\frac{\sin\Omega t}{t},v'_{t2}\frac{\sin\Omega t}{t},v_{3}t)}(k_{1}\widehat{W}_{2}-k_{2}\widehat{W}_{1})\cdot(v_{2}\partial_{v'_{1}}f^{0}-v_{1}\partial_{v'_{2}}f^{0})dvdt$$
\begin{align}
&-\frac{q}{m}\tilde{\rho}(\omega,k)\frac{1}{\omega}\int_{\mathbb{R}^{+}}
\int_{\mathbb{R}^{3}}e^{2\pi it\omega}e^{-2\pi ik\cdot(v'_{t1}\frac{\sin\Omega t}{t},v'_{t2}\frac{\sin\Omega t}{t},v_{3}t)}\cdot(v_{1}\widehat{W}_{1}+v_{2}\widehat{W}_{2})k_{3}\partial_{v'_{3}}f^{0}
dvdt\notag\\
&-\frac{q}{m}\tilde{\rho}(\omega,k)\int_{\mathbb{R}^{+}}\int_{\mathbb{R}^{3}}e^{2\pi it\omega}e^{-2\pi ik\cdot(v'_{t1}\frac{\sin\Omega t}{t},v'_{t2}\frac{\sin\Omega t}{t},v_{3}t)}(\widehat{W}_{1}\partial_{v'_{1}}f^{0}+\widehat{W}_{2}\partial_{v'_{2}}f^{0})dvdt.\notag\\
\end{align}

Let
%and by the definition of $\mathcal{L},$ it is easy to check
\begin{align}
&\tilde{\mathcal{L}}(\omega,k)=\frac{q}{m}\frac{1}{\omega}\int_{\mathbb{R}^{+}}\int_{\mathbb{R}^{3}}e^{2\pi it\omega}e^{-2\pi ik\cdot(v'_{t1}\frac{\sin\Omega t}{t},v'_{t2}\frac{\sin\Omega t}{t},v_{3}t)}(k_{3}v_{3})(\widehat{W}_{1}\partial_{v'_{1}}f^{0}-\widehat{W}_{2}\partial_{v'_{2}}f^{0})dvdt\notag\\
&-\frac{q}{m}\int_{\mathbb{R}^{+}}\int_{\mathbb{R}^{3}}e^{2\pi it\omega}e^{-2\pi ik\cdot(v'_{t1}\frac{\sin\Omega t}{t},v'_{t2}\frac{\sin\Omega t}{t},v_{3}t)}
(\widehat{W}_{1}\partial_{v'_{1}}f^{0}+\widehat{W}_{2}\partial_{v'_{2}}f^{0})\cdot dvdt\notag\\
&-\frac{q}{m}\frac{1}{\omega}\int_{\mathbb{R}^{+}}\int_{\mathbb{R}^{3}}e^{2\pi it\omega}e^{-2\pi ik\cdot(v'_{t1}\frac{\sin\Omega t}{t},v'_{t2}\frac{\sin\Omega t}{t},v_{3}t)}(v_{1}\widehat{W}_{1}+v_{2}\widehat{W}_{2})k_{3}\partial_{v'_{3}}f^{0}
dvdt\notag\\
&+\frac{q}{m}\frac{1}{\omega}\int_{\mathbb{R}^{+}}\int_{\mathbb{R}^{3}}e^{2\pi it\omega}e^{-2\pi ik\cdot(v'_{t1}\frac{\sin\Omega t}{t},v'_{t2}\frac{\sin\Omega t}{t},v_{3}t)}(k_{1}\widehat{W}_{2}-k_{2}\widehat{W}_{1})(v_{2}\partial_{v'_{1}}f^{0}
-v_{1}\partial_{v'_{2}}f^{0})dvdt,\notag\\
\end{align}
hence
\begin{align}
\tilde{\rho}(\omega,k)=\tilde{A}(\omega,k)+\tilde{\rho}(\omega,k)\tilde{\mathcal{L}}(\omega,k).
\end{align}
%Assume $\inf_{k\in\mathbb{Z}^{3}}|\tilde{\mathcal{L}}(\omega,k)-1|>\kappa,$ for some constant $0<\kappa<1,$
Then taking the inverse Fourier transform in time $ t, $ we get
$\hat{\rho}(t,k)=\hat{A}(t,k)+\hat{\rho}(t,k)\ast\hat{\mathcal{L}}(t,k),$
and

\begin{align}
&e^{2\pi\lambda'_{0} |k\cdot(\frac{\sin\Omega t}{t},\frac{\sin\Omega t}{t},t)|}e^{\mu|k|}\hat{\rho}(t,k)=e^{2\pi\lambda'_{0} |k\cdot(\frac{\sin\Omega t}{t},\frac{\sin\Omega t}{t},t)|}e^{\mu|k|}\hat{A}(t,k)+e^{2\pi\lambda'_{0} |k\cdot(\frac{\sin\Omega t}{t},\frac{\sin\Omega t}{t},t)|}e^{\mu|k|}\hat{\rho}(t,k)\ast\hat{\mathcal{L}}(t,k).
\end{align}
Let $\Phi(t,k)=e^{2\pi\lambda'_{0} |k\cdot(\frac{\sin\Omega t}{t},\frac{\sin\Omega t}{t},t)|}e^{\mu|k|}\hat{\rho}(t,k),$
$\mathcal{A}(t,k)=e^{2\pi\lambda'_{0} |k\cdot(\frac{\sin\Omega t}{t},\frac{\sin\Omega t}{t},t)|}e^{\mu|k|}\hat{A}(t,k),$
$\mathcal{K}^{0}(t,k)=e^{2\pi\lambda'_{0} |k\cdot(\frac{\sin\Omega t}{t},\frac{\sin\Omega t}{t},t)|}$
$\cdot\hat{\mathcal{L}}(t,k),$ then from (3.4), we have
$\widetilde{\Phi}(\omega,k)=\widetilde{\mathcal{A}}(\omega,k)+\widetilde{\Phi}(\omega,k)\widetilde{\mathcal{K}}^{0}(\omega,k).$

Then
\begin{align}
&\|\Phi(t,k)\|_{L^{2}(dt)}
=\|\widetilde{\Phi}(\omega,k)\|_{L^{2}}\leq\|\widetilde{\mathcal{A}}(\omega,k)\|_{L^{2}}+\|\widetilde{\Phi}(\omega,k)\|_{L^{2}}
\|\widetilde{\mathcal{K}}^{0}(\omega,k)\|_{L^{\infty}}\notag\\
&\leq\|e^{2\pi\lambda'_{0} |k\cdot(\frac{\sin\Omega t}{t},\frac{\sin\Omega t}{t},t)|}e^{\mu|k|}\hat{A}(t,k)\|_{L^{2}}+\|e^{2\pi\lambda'_{0} |k\cdot(\frac{\sin\Omega t}{t},\frac{\sin\Omega t}{t},t)|}e^{\mu|k|}\hat{\rho}(t,k)\|_{L^{2}}
\|\widetilde{\mathcal{K}}^{0}(\omega,k)\|_{L^{\infty}}.\notag\\
\end{align}

Next we have to estimate $\|\widetilde{\mathcal{K}}^{0}(\omega,k)\|_{L^{\infty}}.$

Indeed,

$$\|\widetilde{\mathcal{K}}^{0}(\omega,k)\|_{L^{\infty}}
\leq\sup_{\omega}\frac{q}{m}\bigg[\bigg|\frac{1}{\omega}\int_{\mathbb{R}^{+}}\int_{\mathbb{R}^{3}}e^{2\pi it\omega}e^{2\pi\lambda'_{0} |k\cdot(\frac{\sin\Omega t}{t},\frac{\sin\Omega t}{t},t)|}e^{-2i\pi k\cdot(\frac{v'_{1}\sin\Omega t}{t},\frac{v'_{2}\sin\Omega t}{t},v'_{3}t)}\cdot(\widehat{W}_{1}\partial_{v'_{1}}f^{0}$$
$$-\widehat{W}_{2}\partial_{v'_{2}}f^{0})
dvdt\bigg|+\frac{q}{m}\bigg|\int_{\mathbb{R}^{+}}\int_{\mathbb{R}^{3}}e^{2\pi it\omega}e^{2\pi\lambda'_{0} |k\cdot(\frac{\sin\Omega t}{t},\frac{\sin\Omega t}{t},t)|}e^{-2i\pi k\cdot(\frac{v'_{1}\sin\Omega t}{t},\frac{v'_{2}\sin\Omega t}{t},v'_{3}t)}
(\widehat{W}_{1}\partial_{v'_{1}}f^{0}+\widehat{W}_{2}\partial_{v'_{2}}f^{0})dvdt\bigg|$$
$$-\frac{q}{m}(\omega,k)\frac{1}{\omega}\int_{\mathbb{R}^{+}}\bigg|\int_{\mathbb{R}^{3}}e^{2\pi it\omega}e^{2\pi\lambda'_{0} |k\cdot(\frac{\sin\Omega t}{t},\frac{\sin\Omega t}{t},t)|}e^{-2i\pi k\cdot(\frac{v'_{1}\sin\Omega t}{t},\frac{v'_{2}\sin\Omega t}{t},v'_{3}t)}
(v_{1}\widehat{W}_{1}+v_{2}\widehat{W}_{2})k_{3}\partial_{v'_{3}}f^{0}
dv\bigg|dt$$
\begin{align}
&+\frac{q}{m}(\omega,k)\frac{1}{\omega}\int_{\mathbb{R}^{+}}\bigg|\int_{\mathbb{R}^{3}}e^{2\pi it\omega}e^{2\pi\lambda'_{0} |k\cdot(\frac{\sin\Omega t}{t},\frac{\sin\Omega t}{t},t)|}e^{-2i\pi k\cdot(\frac{v'_{1}\sin\Omega t}{t},\frac{v'_{2}\sin\Omega t}{t},v'_{3}t)}
(k_{1}\widehat{W}_{2}-k_{2}\widehat{W}_{1})(v_{2}\partial_{v'_{1}}f^{0}
-v_{1}\partial_{v'_{2}}f^{0})dv\bigg|dt\bigg]\notag\\
&=I+II+III+IV.\notag\\
\end{align}
%\begin{align}
%&=\frac{q}{m}\int^{\infty}_{0}|\frac{k_{3}}{\omega}\hat{W}_{2}\int_{\mathbb{R}^{+}}e^{it\omega}(i\eta_{k1})e^{-ik_{3}v_{3}t}e^{\lambda|k_{3}|t}
%\widehat{v_{3}f^{0}}(\eta_{k1},\eta_{k2},v_{3})dv_{3}dt|d\omega\notag\\
%&+\frac{q}{m}\int^{\infty}_{0}|\frac{k_{3}}{\omega}\hat{W}_{1}\int_{\mathbb{R}^{+}}e^{it\omega}(i\eta_{k2})e^{-ik_{3}v_{3}t}e^{\lambda|k_{3}|t}
%\widehat{v_{3}f^{0}}(\eta_{k1},\eta_{k2},v_{3})dv_{3}dt|d\omega\notag\\
%&+\frac{q}{m}\int^{\infty}_{0}|\hat{W}_{1}\int_{\mathbb{R}^{+}}e^{it\omega}(i\eta_{k1})e^{-ik_{3}v_{3}t}e^{\lambda|k_{3}|t}\widehat{f^{0}}(\eta_{k1}
%,\eta_{k2},v_{3})dv_{3}dt|d\omega\notag\\
%&+\frac{q}{m}\int^{\infty}_{0}|\hat{W}_{2}\int_{\mathbb{R}^{+}}e^{it\omega}(i\eta_{k2})e^{-ik_{3}v_{3}t}e^{\lambda|k_{3}|t}\widehat{f^{0}}(\eta_{k1},\eta_{k2},
%v_{3})dv_{3}dt|d\omega\notag\\
%&+\cdots\cdots\notag\\
%&=I+II+III+IV+\cdots
%&\leq Ce^{-\alpha(\frac{\omega}{k_{3}})^{2}}+\varepsilon\leq C_{Te}e^{-\alpha(v_{Te})^{2}},
%\end{align}
%where $\eta_{k_{1}}=\frac{1}{\Omega}(-k_{2}\cos\Omega t+k_{2}-k_{1}\sin\Omega t),\eta_{k_{2}}=\frac{1}{\Omega}(-k_{2}\sin\Omega t-k_{1}+k_{1}\cos\Omega t).$

In fact, we only need to estimate one term of (3.6) because of similar processes of  other terms. Without loss of generality, we give an estimate for $ I.$ In the same way, we only estimate one term of $I,$ here we still denote $ I. $
$$ I=\sup_{\omega}\frac{q}{m}\bigg|\frac{k_{3}}{\omega}\widehat{W}_{2}\int_{\mathbb{R}^{+}}\int_{\mathbb{R}}e^{2\pi it\omega}\bigg(2\pi i\frac{k_{1}\sin\Omega t}{\Omega}\bigg)
e^{2\pi\lambda_{0} |k\cdot(\frac{\sin\Omega t}{t},\frac{\sin\Omega t}{t},t)|}e^{-2\pi ik_{3}v_{3}t}\cdot\widehat{v_{3}f^{0}}(\frac{k_{1}\sin\Omega t}{\Omega},\frac{k_{2}\sin\Omega t}{\Omega},v_{3})
dv_{3}dt\bigg|$$
$$=\sup_{\omega}\frac{q}{m}\bigg|\frac{k_{3}}{\omega}\widehat{W}_{2}\int_{\mathbb{R}^{+}}\int_{\mathbb{R}}e^{2\pi it\omega}\bigg(2\pi i\frac{k_{1}\sin\Omega t}{\Omega}\bigg)
e^{2\pi\lambda_{0} |k\cdot(\frac{\sin\Omega t}{t},\frac{\sin\Omega t}{t},t)|}e^{-2\pi ik_{3}v_{3}t}e^{-2\pi\lambda_{0}|k_{3}|t}\cdot\sum_{n}\frac{|2\pi i\lambda_{0}|k_{3}|t|^{n}}{n!}$$
$$\cdot\widehat{v_{3}f^{0}}(\frac{k_{1}\sin\Omega t}{\Omega},\frac{k_{2}\sin\Omega t}{\Omega},v_{3})dv_{3}dt\bigg|$$
$$=\sup_{\omega}\frac{q}{m}\sum_{n}\frac{\lambda_{0}^{n}}{n!}\bigg|\frac{k_{3}}{\omega}\widehat{W}_{2}\int_{\mathbb{R}^{+}}\int_{\mathbb{R}}\bigg(2\pi i\frac{k_{1}\sin\Omega t}{\Omega}\bigg)
e^{2\pi\lambda_{0} |k\cdot(\frac{\sin\Omega t}{t},\frac{\sin\Omega t}{t},t)|}e^{-2\pi\lambda_{0}|k_{3}|t}e^{2\pi ik_{3}t(\frac{\omega}{k_{3}}-v_{3})}$$
$$\cdot\nabla_{v_{3}}^{n}\widehat{v_{3}f^{0}}(\frac{k_{1}\sin\Omega t}{\Omega},\frac{k_{2}\sin\Omega t}{\Omega},v_{3})dv_{3}dt\bigg|$$
%&=\frac{q}{m}\sum_{n}\frac{\lambda^{n}}{n!}\int^{\infty}_{0}|\frac{k_{3}}{\omega}\hat{W}_{2}\int_{\mathbb{R}^{+}}e^{it\omega}(i\eta_{k1})(ik_{3}t)^{n}
%e^{-ik_{3}v_{3}t}
%\widehat{v_{3}f^{0}}(\eta_{k1},\eta_{k2},v_{3})dv_{3}dt|d\omega\notag\\
%$$=\frac{q}{m}\sum_{n}\frac{\lambda^{n}}{n!}\int^{\infty}_{0}\bigg|\frac{k_{3}}{\omega}\hat{W}_{2}\int_{\mathbb{R}^{+}}\int_{\mathbb{R}}(2\pi i\eta_{k1})(2\pi ik_{3}t)^{n}
%e^{2\pi\lambda_{0}|\eta_{k1}|}e^{2\pi\lambda_{0}|\eta_{k2}|}e^{2\pi ik_{3}t(\frac{\omega}{k_{3}}-v_{3})}\cdot\widehat{v_{3}f^{0}}(\eta_{k1},\eta_{k2},v_{3})dv_{3}dt\bigg|d\omega$$
%$$=\frac{q}{m}\sum_{n}\frac{\lambda^{n}}{n!}\int^{\infty}_{0}\bigg|\frac{k_{3}}{\omega}\hat{W}_{2}\int_{\mathbb{R}^{+}}\int_{\mathbb{R}}(2\pi i\eta_{k1})
%(-1)^{n}e^{2\pi\lambda_{0}|\eta_{k1}|}e^{2\pi\lambda_{0}|\eta_{k2}|}\nabla_{v_{3}}^{n}e^{2\pi ik_{3}t(\frac{\omega}{k_{3}}-v_{3})}\cdot\widehat{v_{3}f^{0}}(\eta_{k1},\eta_{k2},v_{3})dv_{3}dt\bigg|d\omega$$
%$$=\frac{q}{m}\sum_{n}\frac{\lambda_{0}^{n}}{n!}\int^{\infty}_{0}\bigg|\frac{k_{3}}{\omega}\hat{W}_{2}\int_{\mathbb{R}^{+}}(2\pi i\eta_{k1})
%(-1)^{n}e^{2\pi\lambda_{0}|\eta_{k1}|}e^{2\pi\lambda_{0}|\eta_{k2}|}e^{2\pi ik_{3}t(\frac{\omega}{k_{3}})}
%(2\pi ik_{3}t)^{n}\cdot\widehat{v_{3}f^{0}}(\eta_{k1},\eta_{k2},k_{3}t) dt\bigg|d\omega$$
$$=\sup_{\omega}\frac{q}{m}\sum_{n}\frac{\lambda_{0}^{n}}{n!}\bigg|\frac{k_{3}}{\omega}\widehat{W}_{2}\bigg(2\pi i\frac{k_{1}\sin\Omega t}{\Omega}\bigg)
e^{2\pi\lambda_{0} |k\cdot(\frac{\sin\Omega t}{t},\frac{\sin\Omega t}{t},t)|}e^{-2\pi\lambda_{0}|k_{3}|t}
(-i\nabla_{\frac{\omega}{k_{3}}})^{n}\widehat{v_{3}f^{0}}(\frac{k_{1}\sin\Omega t}{\Omega},\frac{k_{2}\sin\Omega t}{\Omega},\frac{\omega}{k_{3}})\bigg |$$
$$\leq\frac{q}{mv_{Te}}e^{-c_{0}v_{Te}},$$
%\begin{align}
%&\leq\frac{q}{m}\sum_{n}\frac{\lambda_{0}^{n}}{n!}\int^{\infty}_{v_{Te}}\bigg|\frac{k_{3}}{\omega}\hat{W}_{2}(2\pi i\eta_{k1})
%e^{2\pi\lambda_{0}|\eta_{k1}|}e^{2\pi\lambda_{0}|\eta_{k2}|}
%(-i\nabla_{\frac{\omega}{k_{3}}})^{n}\widehat{v_{3}f^{0}}(\eta_{k1},\eta_{k2},\frac{\omega}{k_{3}}) \bigg|d\omega\leq\frac{q}{mv_{Te}}e^{-c_{0}v_{Te}}\notag\\
%&+\frac{q}{m}\sum_{n}\frac{\lambda^{n}}{n!}\int^{v_{Te}}_{0}|\frac{k_{3}}{\omega}\hat{W}_{2}(i\eta_{k1})
%e^{\lambda|\eta_{k1}|}e^{\lambda|\eta_{k2}|}
%(-i\nabla_{\frac{\omega}{k_{3}}})^{n}\widehat{v_{3}f^{0}}(\eta_{k1},\eta_{k2},\frac{\omega}{k_{3}}) |d\omega\leq\frac{q}{mv_{Te}}e^{-c_{0}v_{Te}},
%\end{align}
where in the last inequality we use the facts that if $v_{3}=\frac{\omega}{k_{3}},$ then $v_{3}\gg v_{Te},$ and   the assumption (i) and (iv).
Then  %for any $k,\omega,$ by the assumption (iv) of $Proposition$ 1.1, we obtain
%$$|1-\tilde{\mathcal{L}}(\omega,k)|\geq\kappa,$$
 there exists  some constant $0<\kappa<1$ such that
$\|\widetilde{\mathcal{K}}^{0}(\omega,k)\|_{L^{\infty}}\leq\kappa.$

%\begin{align}
%&\mathcal{L}(t,k)=\frac{q}{m}\int^{t}_{0}k_{3}(\hat{W}_{2}(i\eta_{k_{1}})\partial_{\eta_{3}}
%\hat{f}^{0}(\eta_{k_{1}},\eta_{k_{2}},k_{3}t)+\hat{W}_{1}(i\eta_{k_{2}})\partial_{\eta_{3}}\hat{f}^{0}(\eta_{k_{1}},\eta_{k_{2}},k_{3}t))dt\notag\\
%&-\frac{q}{m}\int^{t}_{0}
%(\hat{W}_{1}(i\eta_{k_{1}})\hat{f}^{0}(\eta_{k_{1}},\eta_{k_{2}},k_{3}t)+\hat{W}_{2}(i\eta_{k_{2}})\hat{f}^{0}(\eta_{k_{1}},\eta_{k_{2}},k_{3}t))dt,\notag\\
%\end{align}

In conclusion, we have \begin{align}
\sup_{t\geq0}\|\rho(t,\cdot)\|_{\mathcal{F}^{\lambda'_{0}t+\mu'}}\leq C\sup_{t\geq0}\bigg\|f_{0}\circ S^{0}_{-t}dv\bigg\|_{\mathcal{F}^{\lambda_{0}t+\mu}}\leq \|f_{0}\circ S^{0}_{-t}\|_{\mathcal{Z}_{t}^{\lambda_{0},\mu;1}}=\|f_{0}\|_{\mathcal{Z}_{0}^{\lambda_{0},\mu;1}}\leq\delta_{0}.\end{align}

%$(\mathbf{wherther}\quad \mathbf{consider} \quad\omega=0$ ? since we only care the integration)

    %$$\hat{\rho}(t,k_{1},k_{2},0)=\int_{\mathbb{T}^{3}}\int_{\mathbb{R}^{3}}e^{-2\pi i(x_{1},x_{2})\cdot (k_{1},k_{2})}f(t,x,v)dvdx=\int_{\mathbb{T}^{3}}\int_{\mathbb{R}^{3}}e^{-2\pi i(x_{1},x_{2})\cdot (k_{1},k_{2})}f_{0}(x'(0,x,v),v'(0,x,v))dvdx$$
%$$-\frac{q}{m}\int^{t}_{0}\int_{\mathbb{T}^{3}}\int_{\mathbb{R}^{3}}e^{-2\pi i(x_{1},x_{2})\cdot (k_{1},k_{2})}[(E+v'\times B)\cdot\nabla'_{v}f^{0}](\tau,x'(\tau,x,v),v'(\tau,x,v))dvdx,$$

%Since $\partial_{t}\hat{B}=(k_{1},k_{2},k_{3})\times \hat{E}(t,k_{1},k_{2},k_{3}),E(t,x)=W(x)\ast\rho(t,x),$
%then if  $ k_{3}=0,$  since $W(x)=(W_{1}(x),W_{2}(x),0),$
%$\widehat{W}(k_{1},k_{2},0)=0,$ we have
%$\hat{B}(t,k_{1},k_{2},0)\equiv0.$
%So $\int_{\mathbb{R}^{3}}e^{-2\pi i(\frac{1}{\Omega}(v'_{2}-v_{2},-v'_{1}+v_{1})\cdot (k_{1},k_{2})}\nabla'_{v}\cdot(\hat{B}\times(f^{0}v'))dv'=0.$
%Finally, we get $\hat{\rho}(t,k_{1},k_{2},0)=\int\int e^{-2\pi i(x_{1},x_{2})\cdot (k_{1},k_{2})}f_{0}dvdx'.$  And we finish the proof.

%Furthermore,
%\begin{align}
%\|\Phi\|_{L^{\infty}(dt)}\leq(1+\frac{C(k,W,\Omega)}{(\lambda_{1}-\lambda)^{2}})\|(\widetilde{e^{\lambda|k_{3}|t}A})\|_{L^{1}(d\omega)},
%\end{align}
$$$$
         $Proof$
         $of$
         $Theorem$
         $ 3.1.$
From (2.1), we have
\begin{align}
&f(t,x,v)=f_{0}(x'(0,x,v),v'(0.x,v))-\frac{q}{m}\int^{t}_{0}[(E+v'\times B)\cdot\nabla'_{v}f^{0}](\tau.x'(\tau,x,v),v'(\tau,x,v))d\tau.
\end{align}
Then for any $\lambda''_{0}<\lambda'_{0},$ we have, for all $t\geq0,$
\begin{align}
&\|f\|_{\mathcal{Z}^{\lambda''_{0},\mu';1}_{t}}\leq \|f_{0}\circ S_{-t}^{0}\|_{\mathcal{Z}^{\lambda''_{0},\mu';1}_{t}}+\int^{t}_{0}\|(W(x)\ast\rho_{\tau})\circ S_{-(t-\tau)}^{0}\|_{\mathcal{Z}^{\lambda''_{0},\mu';\infty}_{t}}\|\nabla_{v}f^{0}\|_{\mathcal{Z}^{\lambda''_{0},\mu';1}_{t}}d\tau\notag\\
&+\int^{t}_{0}\|B\circ S_{-(t-\tau)}^{0}\|_{\mathcal{Z}^{\lambda''_{0},\mu';\infty}_{t}}\|\nabla_{v}vf^{0}\|_{\mathcal{Z}^{\lambda''_{0},\mu';1}_{t}}d\tau\notag\\
&=\|f_{0}\|_{\mathcal{Z}^{\lambda''_{0},\mu';1}}+
\int^{t}_{0}\|(W(x)\ast\rho_{\tau})\|_{\mathcal{F}^{\lambda''_{0}\tau+\mu'}}\|\nabla_{v}f^{0}\|_{\mathcal{C}^{\lambda''_{0};1}}d\tau
+\int^{t}_{0}\|B\|_{\mathcal{F}^{\lambda''_{0}\tau+\mu'}}\|\nabla_{v}vf^{0}\|_{\mathcal{C}^{\lambda''_{0};1}}d\tau\notag\\
&\leq\delta_{0}+C\delta_{0}/(\lambda_{0}-\lambda''_{0}).
\end{align}

 $$$$
         $$$$
        % \begin{center}
%\item\section{$B_{0}\rightarrow0$}
%{\bf\large 1. \quad Introduction }
%\end{center}

%First we still write the euqation
%\begin{equation}
%\left\{\begin{array}{l}
%\partial_{t}f+v\cdot\nabla_{x}f+\frac{q}{m}(v\times B_{0})\cdot\nabla_{v}f=-\frac{q}{m}(E+v\times B_{1})\cdot\nabla_{v}f^{0},\\
%\partial_{t}B_{1}=\nabla_{x}\times E,\\
%E=W(x)\ast\rho(t,x),\rho(t,x)=\int_{\mathbb{R}^{3}}f(t,x,v)dv,\\
%f(0,x,v)=f_{0}(x,v_{\perp},v_{z}),f^{0}(v)=f^{0}(v_{\perp},v_{z}),
%\end{array}\right.
% \end{equation}

         \begin{center}
\item\section{Nonlinear Cyclotron damping}
%{\bf\large 1. \quad Introduction }
\end{center}

Next we give the proof of the main Theorem 0.3, stating the primary steps as propositions which are proved in subsections.

 \begin{center}
\item\subsection{The improved Newton iteration}
%{\bf\large 1. \quad Introduction }
\end{center}
The first idea which may come to mind is a classical Newton iteration as done by Mouhot and Villani [23]:
 Let
$$f^{0}=f^{0}(v) \quad\textrm{be\quad given},$$
and $$f^{n}=f^{0}+h^{1}+\ldots+h^{n},$$
where
 \begin{align}
  \left\{\begin{array}{l}
      \partial_{t}h^{1}+v\cdot\nabla_{x}h^{1}+v\times B_{0}\cdot\nabla_{v}h^{1}+(E[h^{1}]+v\times B[h^{1}])\cdot\nabla_{v}f^{0}=0,\\
     h^{1}(0,x,v) =f_{0}-f^{0}, \\
   \end{array}\right.
           \end{align}
 and now we consider the Vlasov equation in step $n+1,n\geq1,$
 \begin{align}
  \left\{\begin{array}{l}
      \partial_{t}h^{n+1}+v\cdot\nabla_{x}h^{n+1}+v\times B_{0}\cdot\nabla_{v}h^{n+1}+E[f^{n}]\cdot\nabla_{v}h^{n+1} +v\times B[f^{n}]\cdot\nabla_{v}h^{n+1}\\
      =-E[h^{n+1}]\cdot\nabla_{v}f^{n}-v\times B[h^{n+1}]\cdot\nabla_{v}f^{n}-E[h^{n}]\cdot\nabla_{v}h^{n}-v\times B[h^{n}]\cdot\nabla_{v}h^{n},\\
      h^{n+1}(0,x,v) =0, \\
   \end{array}\right.
           \end{align}
the corresponding dynamical system is described by the equations:
for any $(x,v)\in\mathbb{T}^{3}\times\mathbb{R}^{3},$
let $(X^{n}_{t,\tau},V^{n}_{t,\tau})$ as the solution of the following ordinary differential equations
  $$
   \left\{\begin{array}{l}
   \frac{d}{dt}X^{n}_{t,\tau}(x,v)=V^{n}_{t,\tau}(x,v),\\
 X^{n}_{\tau,\tau}(x,v)=x,
  \end{array}\right.
           $$

         \begin{align}
   \left\{\begin{array}{l}
   \frac{d}{dt}V^{n}_{t,\tau}(x,v)=V^{n}_{t,\tau}(x,v)\times (B_{0}+B[f^{n}](t,X^{n}_{t,\tau}(x,v)))
   +E[f^{n}](t,X^{n}_{t,\tau}(x,v)),\\
 V^{n}_{\tau,\tau}(x,v)=v.
  \end{array}\right.
           \end{align}
            At the same time, we consider the corresponding  linear dynamics system as follows,
            \begin{align}
   \left\{\begin{array}{l}
  \frac{d}{dt}X^{0}_{t,\tau}(x,v)=V^{0}_{t,\tau}(x,v),\quad  \frac{d}{dt}V^{0}_{t,\tau}(x,v)=V^{0}_{t,\tau}(x,v)\times B_{0},\\
 X^{0}_{\tau,\tau}(x,v)=x,\quad V^{0}_{\tau,\tau}(x,v)=v.
  \end{array}\right.
           \end{align}
             It is easy to check that
             $$\Omega^{n}_{t,\tau}-Id\triangleq(\delta X^{n}_{t,\tau},\delta V^{n}_{t,\tau})\circ (X^{0}_{\tau,t},V^{0}_{\tau,t})
             = ( X^{n}_{t,\tau}\circ(X^{0}_{\tau,t},V^{0}_{\tau,t})-Id, V^{n}_{t,\tau}\circ(X^{0}_{\tau,t},V^{0}_{\tau,t})-Id).$$
            Therefore,  in order to estimate $( X^{n}_{t,\tau}\circ(X^{0}_{\tau,t},V^{0}_{\tau,t})-Id, V^{n}_{t,\tau}\circ(X^{0}_{\tau,t},V^{0}_{\tau,t})-Id),$
             we only need to study $(\delta X^{n}_{t,\tau},\delta V^{n}_{t,\tau})\circ (X^{0}_{\tau,t},V^{0}_{\tau,t}).$

             From Eqs.(4.3) and (4.4),
             $$
   \left\{\begin{array}{l}
   \frac{d}{dt}\delta X^{n}_{\tau,t}(x,v)=\delta V^{n}_{\tau,t}(x,v),\\
 \delta X^{n}_{\tau,\tau}(x,v)=0,
  \end{array}\right.
          $$
            \begin{align}
   \left\{\begin{array}{l}
   \frac{d}{dt}\delta V^{n}_{\tau,t}(x,v)=\delta V^{n}_{\tau,t}(x,v)\times B_{0}+E[f^{n}](t,X^{n}_{\tau,t}(x,v))+(\delta V^{n}_{\tau,t}(x,v)+ V^{0}_{\tau,t}(x,v))\times B[f^{n}](t,X^{n}_{\tau,t}(x,v)),\\
\delta V^{n}_{\tau,\tau}(x,v)=0,
  \end{array}\right.
           \end{align}
           and $|V^{0}_{\tau,t}(x,v)|=|v|.$ Since $v\in\mathbb{R}^{3}$ and $B[f^{n}](t,X^{n}_{\tau,t}(x,v))$ independent of $v,$
           there is almost no hope to get a ``good " estimates of $\Omega^{n}_{t,\tau}-Id.$ Furthermore, when $k_{3}=0,$ because of $\widehat{E[f^{n}]}(s,k_{1},k_{2},0)=0,\widehat{B[f^{n}]}(s,k_{1},k_{2},0)\neq0,$  the deflection estimates are in the absence of a decaying perturbed magnetic field.

           To circumvent these difficulties,  we recall the basic physical Law on Lenz's Law:

          $ The$ $direction$ $of$ $current$ $induced$ $in$ $a$ $conductor$ $by$ $a$ $changing$ $magnetic$ $field$ $due$ $to$
           $induction$ $is$ $such$ $that$ $it$ $creates$ $a$ $magnetic$ $field$ $that$ $opposes$ $the$
           $change$ $that$ $produced$ $it.$

           According to  the statement of Lenz's Law and Maxwell equations, based on the approximation equations (4.2),
           it is easy to know that we only need to consider the following dynamical system

           \begin{align}
   \left\{\begin{array}{l}
   \frac{d}{dt}X^{n}_{t,\tau}(x,v)=V^{n}_{t,\tau}(x,v), \frac{d}{dt}V^{n}_{t,\tau}(x,v)=V^{n}_{t,\tau}(x,v)\times B_{0}
   +E[f^{n}](t,X^{n}_{t,\tau}(x,v)),\\
 X^{n}_{\tau,\tau}(x,v)=x, V^{n}_{\tau,\tau}(x,v)=v,
  \end{array}\right.
           \end{align}

then we get
           \begin{align}
   \left\{\begin{array}{l}
   \frac{d}{dt}\delta X^{n}_{t,\tau}(x,v)=\delta V^{n}_{t,\tau}(x,v), \frac{d}{dt}\delta V^{n}_{t,\tau}(x,v)=\delta V^{n}_{t,\tau}(x,v)\times B_{0}
  +E[f^{n}](t,X^{n}_{t,\tau}(x,v)),\\
 \delta X^{n}_{\tau,\tau}(x,v)=0,\delta V^{n}_{\tau,\tau}(x,v)=0,
  \end{array}\right.
           \end{align}
            and we write  the approximation equations (4.2)  into the following form,
           \begin{align}
  \left\{\begin{array}{l}
      \partial_{t}h^{n+1}+v\cdot\nabla_{x}h^{n+1}+v\times B_{0}\cdot\nabla_{v}h^{n+1}+E[f^{n}]\cdot\nabla_{v}h^{n+1} \\
      =-E[h^{n+1}]\cdot\nabla_{v}f^{n}-v\times B[h^{n+1}]\cdot\nabla_{v}f^{n}-E[h^{n}]\cdot\nabla_{v}h^{n}-v\times B[h^{n}]\cdot\nabla_{v}h^{n}\\
      -v\times B[f^{n}]\cdot\nabla_{v}h^{n+1},\\
      h^{n+1}(0,x,v) =0.\\
   \end{array}\right.
           \end{align}

 \begin{center}
\item\subsection{Main challenges}
%{\bf\large 1. \quad Introduction }
\end{center}
           Integrating (4.8) in time and $ h^{n+1}(0,x,v) =0,$ we get
            \begin{align}
            h^{n+1}(t,X^{n}_{t,0}(x,v),V^{n}_{t,0}(x,v))=\int^{t}_{0}\Sigma^{n+1}(s,X^{n}_{s,0}(x,v),V^{n}_{s,0}(x,v))ds,
            \end{align}
           where$$\Sigma^{n+1}(t,x,v)=-E[h^{n+1}]\cdot\nabla_{v}f^{n}-v\times B[h^{n+1}]\cdot\nabla_{v}f^{n}-E[h^{n}]\cdot\nabla_{v}h^{n}
           -v\times B[h^{n}]\cdot\nabla_{v}h^{n}
      -v\times B[f^{n}]\cdot\nabla_{v}h^{n+1}.$$

           By the definition of $(X^{n}_{t,\tau}(x,v),V^{n}_{t,\tau}(x,v)),$ we have
           $$h^{n+1}(t,x,v)=\int^{t}_{0}\Sigma^{n+1}(s,X^{n}_{s,t}(x,v),V^{n}_{s,t}(x,v))ds
           =\int^{t}_{0}\Sigma^{n+1}(s,\delta X^{n}_{s,t}(x,v)+ X^{0}_{s,t}(x,v),\delta V^{n}_{s,t}(x,v)+ V^{0}_{s,t}(x,v))ds.$$

           Since the unknown $h^{n+1}$ appears on both sides of (4.9), we hope to get a self-consistent estimate. For this, we have little choice but to
           integrate in $v$ and get an integral equation on $\rho[h^{n+1}]=\int_{\mathbb{R}^{3}}h^{n+1}dv,$ namely

 \begin{align}
           &\rho[h^{n+1}](t,x)=\int^{t}_{0}\int_{\mathbb{R}^{3}}(\Sigma^{n+1}\circ\Omega^{n}_{s,t}(x,v))(s, X^{0}_{s,t}(x,v),V^{0}_{s,t}(x,v))dvds\notag\\
           &=\int^{t}_{0}\int_{\mathbb{R}^{3}}-\bigg[(\mathcal{E}^{n+1}_{s,t}\cdot G_{s,t}^{n})-(F^{n+1}_{s,t}\cdot G_{s,t}^{n,v})-(\mathcal{E}^{n}_{s,t}\cdot H^{n}_{s,t})-(F^{n}_{s,t}\cdot H^{n,v}_{s,t})\notag\\
           &-(B[f^{n}]\circ\Omega^{n}_{s,t}(x,v))\cdot H_{s,t}^{n+1,v}\bigg](s, X^{0}_{s,t}(x,v),V^{0}_{s,t}(x,v))dvds\notag\\
          &=I^{n+1,n}+II^{n+1,n}+III^{n,n}+IV^{n,n}+V^{n,n+1},
           \end{align}
  where
  $$
  \left\{\begin{array}{l}
   \mathcal{E}^{n+1}_{s,t}=E[h^{n+1}]\circ\Omega^{n}_{s,t}(x,v),\quad \mathcal{E}^{n}_{s,t}=E[h^{n}]\circ\Omega^{n}_{s,t}(x,v),\\
   G_{s,t}^{n}=(\nabla'_{v}f^{n})\circ\Omega^{n}_{s,t}(x,v),\quad G_{s,t}^{n,v}=(\nabla'_{v}f^{n}\times V^{0}_{s,t}(x,v))\circ\Omega^{n}_{s,t}(x,v),\\
   F^{n+1}_{s,t}=B[h^{n+1}]\circ\Omega^{n}_{s,t}(x,v),\quad F^{n}_{s,t}=B[h^{n}]\circ\Omega^{n}_{s,t}(x,v),\\
    H^{n}_{s,t}=(\nabla'_{v}h^{n})\circ\Omega^{n}_{s,t}(x,v),\quad H_{s,t}^{n,v}=(\nabla'_{v}h^{n}\times V^{0}_{s,t}(x,v))\circ\Omega^{n}_{s,t}(x,v).\\
   \end{array}\right.
   $$

  It is obvious that Eq.(4.10) is not a closed equation, while $\rho[h^{n+1}](t,x)$ satisfies a closed equation written
   in [23] with the electric field case. To obtain a self-consistent estimate, we go back the Vlasov equations (4.9), composing with
    $((X^{n}_{0,\tau}(x,v),V^{n}_{0,\tau}(x,v))),$ where $0\leq\tau\leq t,$ this gives
  \begin{align}
  h^{n+1}(t,X^{n}_{\tau,t}(x,v),V^{n}_{\tau,t}(x,v))=\int^{t}_{0}\Sigma^{n+1}(s,X^{n}_{\tau,s}(x,v),V^{n}_{\tau,s}(x,v))ds.
  \end{align}
   In order to achieve the goal,  we have to combine Eq.(4.10) and Eq.(4.11) to form a iteration, then we obtain a self-consistent estimate on $\rho[h^{n+1}].$
   At the technical level, it is more difficult than that in  Mouhot and Villani'paper [23],
   since the new term $v\times B[f^{n}]\cdot\nabla_{v}h^{n+1}$ in Eqs.(4.9)-(4.10) brings different kinds of resonances (in terms of different norms).

   \begin{center}
   \item\subsection{Inductive hypothesis}
%{\bf\large 1. \quad Introduction }
\end{center}

For n=1, from (1.1), it is easy to see  that (4.1) is a linear Vlasov equation. From section 3, we know that  the conclusions of Theorem 0.3 hold.

 Now for any $i\leq n,i\in\mathbb{N}_{0},$ we assume that the following estimates hold,
 \begin{align}
 &\sup_{t\geq0}\|\rho[h^{i}](t,\cdot)\|_{\mathcal{F}^{(\lambda_{i}-B_{0})t+\mu_{i}}}\leq\delta_{i},\notag\\
 &\sup_{0\leq\tau\leq t}\|h^{i}_{\tau}\circ\Omega^{i-1}_{t,\tau}\|
_{\mathcal{Z}_{\tau-\frac{bt}{1+b}}^{(\lambda_{i}-B_{0})(1+b),\mu_{i};1}}\leq\delta_{i},
 \sup_{0\leq\tau\leq t}\|(h^{i}_{\tau}v)\circ\Omega^{i-1}_{t,\tau}\|_{\mathcal{Z}^{(\lambda_{i}-B_{0})(1+b),
\mu_{i};1}_{\tau-\frac{bt}{1+b}}}\leq\delta_{i},\notag\\
\end{align}
then we have the following inequalities, denote $(\mathbf{E}^{n}):$
$$\sup_{t\geq 0}\|E[h^{i}](t,\cdot)\|_{\mathcal{F}^{(\lambda_{i}-B_{0})t+\mu_{i}}}<\delta_{i},\quad
           \sup_{t\geq0}\|B[h^{i}](t,\cdot)\|_{\mathcal{F}^{(\lambda_{i}-B_{0})t+\mu_{i}}}<\delta_{i},$$
$$\sup_{0\leq\tau\leq t}\| \nabla_{x}((h^{i}_{\tau}v)\circ\Omega^{i-1}_{t,\tau})\|_{\mathcal{Z}^{(\lambda_{i}-B_{0})(1+b),
\mu_{i};1}_{\tau-\frac{bt}{1+b}}}\leq\delta_{i},\sup_{0\leq\tau\leq t}\| (\nabla_{x}(h^{i}_{\tau}v))\circ\Omega^{i-1}_{t,\tau}\|_{\mathcal{Z}^{(\lambda_{i}-B_{0})(1+b),
\mu_{i};1}_{\tau-\frac{bt}{1+b}}}\leq\delta_{i},$$
$$\sup_{0\leq\tau\leq t}\|\nabla_{x}(h^{i}_{\tau}\circ\Omega^{i-1}_{t,\tau})\|
_{\mathcal{Z}_{\tau-\frac{bt}{1+b}}^{(\lambda_{i}-B_{0})(1+b),\mu_{i};1}}\leq\delta_{i},\sup_{0\leq\tau\leq t}\|(\nabla_{x}h^{i}_{\tau})\circ\Omega^{i-1}_{t,\tau}\|
_{\mathcal{Z}_{\tau-\frac{bt}{1+b}}^{(\lambda_{i}-B_{0})(1+b),\mu_{i};1}}\leq\delta_{i},$$
$$\|(\nabla_{v}+\tau\nabla_{x})(h^{i}_{\tau}\circ\Omega^{i-1}_{t,\tau})\|
_{\mathcal{Z}_{\tau-\frac{bt}{1+b}}^{(\lambda_{i}-B_{0})(1+b),\mu_{i};1}}\leq\delta_{i},\|((\nabla_{v}+\tau\nabla_{x})h^{i}_{\tau})\circ\Omega^{i-1}_{t,\tau}\|
_{\mathcal{Z}_{\tau-\frac{bt}{1+b}}^{(\lambda_{i}-B_{0})(1+b),\mu_{i};1}}\leq\delta_{i},$$
$$\|(\nabla_{v}+\tau\nabla_{x})(h^{i}_{\tau}v\circ\Omega^{i-1}_{t,\tau})\|
_{\mathcal{Z}_{\tau-\frac{bt}{1+b}}^{(\lambda_{i}-B_{0})(1+b),\mu_{i};1}}\leq\delta_{i},\|((\nabla_{v}+\tau\nabla_{x})h^{i}_{\tau}v)\circ\Omega^{i-1}_{t,\tau}\|
_{\mathcal{Z}_{\tau-\frac{bt}{1+b}}^{(\lambda_{i}-B_{0})(1+b),\mu_{i};1}}\leq\delta_{i},$$
$$\sup_{0\leq\tau\leq t}\|(\nabla\nabla_{v} h^{i}_{\tau})\circ\Omega^{i-1}_{t,\tau}\|
_{\mathcal{Z}_{\tau-\frac{bt}{1+b}}^{(\lambda_{i}-B_{0})(1+b),\mu_{i};1}}\leq\delta_{i},\sup_{0\leq\tau\leq t}\|(\nabla\nabla_{v}
(h^{i}_{\tau}v))\circ\Omega^{i-1}_{t,\tau}\|
_{\mathcal{Z}_{\tau-\frac{bt}{1+b}}^{(\lambda_{i}-B_{0})(1+b),\mu_{i};1}}\leq\delta_{i},$$
$$\sup_{0\leq\tau\leq t}(1+\tau)^{2}\|(\nabla_{v}\times h^{i})\circ\Omega^{i-1}_{t,\tau}-\nabla_{v}\times(h^{i}\circ\Omega^{i-1}_{t,\tau})\|
_{\mathcal{Z}_{\tau-\frac{bt}{1+b}}^{(\lambda_{i}-B_{0})(1+b),\mu_{i};1}}\leq\delta_{i},$$
$$\sup_{0\leq\tau\leq t}(1+\tau)^{2}\|(\nabla_{v}\times(h^{i}v))\circ\Omega^{i-1}_{t,\tau}-\nabla_{v}\times((h^{i}v)\circ\Omega^{i-1}_{t,\tau})\|
_{\mathcal{Z}_{\tau-\frac{bt}{1+b}}^{(\lambda_{i}-B_{0})(1+b),\mu_{i};1}}\leq\delta_{i}.$$

It is easy to check  that the first two inequalities of $(\mathbf{E}^{n})$ hold under our assumptions and (4.12), so we need to show that the other equalities of
$(\mathbf{E}^{n})$ also hold, the related proofs are found in section 5.
  \begin{center}
   \item\subsection{Local time iteration}
%{\bf\large 1. \quad Introduction }
\end{center}

Before working out the core of the proof of Theorem 0.1, we shall show a short time estimate, which will play a role as an initial data layer for the Newton scheme.
The main tool in this section is given by the following lemma, which is through  the direct computation from the definition of the corresponding norms. Therefore, we
omit the proof.
\begin{lem} Let $f$ be an analytic function, $\lambda(t)=\lambda-Kt$ and $\mu(t)=\mu-Kt,$ $ K>0,$ let $T>0$ be so small that $\lambda(t)>\lambda'(t)>0,
\mu(t)>\mu'(t)>0$ for $0\leq t\leq T.$ Then
for any $\tau\in[0,T]$ and any $p\geq1,$
$$\frac{d^{+}}{dt}|_{t=\tau}\|f\|_{\mathcal{Z}_{\tau}^{\lambda(t),\mu(t);1}}\leq-\frac{K}{2(1+\tau)}\|\nabla f\|_{\mathcal{Z}_{\tau}^{\lambda(\tau),\mu(\tau);1}}
-\frac{K}{2(1+\tau)}\|v\nabla f\|_{\mathcal{Z}_{\tau}^{\lambda'(\tau),\mu'(\tau);1}},$$
where $\frac{d^{+}}{dt}$ stands for the upper right derivative.
\end{lem}

 For $n\geq1,$ now let us solve
 $$ \partial_{t}h^{n+1}+v\cdot\nabla_{x}h^{n+1}+v\times B_{0}\cdot\nabla_{v}h^{n+1} =\widetilde{\Sigma}^{n+1},$$
 where $$\widetilde{\Sigma}^{n+1}=-E[f^{n}]\cdot\nabla_{v}h^{n+1}-E[h^{n+1}]\cdot\nabla_{v}f^{n}
 -v\times B[h^{n+1}]\cdot\nabla_{v}f^{n}-E[h^{n}]\cdot\nabla_{v}h^{n}-v\times B[h^{n}]\cdot\nabla_{v}h^{n} -v\times B[f^{n}]\cdot\nabla_{v}h^{n+1}.$$

 Hence $$\|h^{n+1}\|_{\mathcal{Z}_{t}^{\lambda_{n+1}(t),\mu_{n+1}(t);1}}\leq\int^{t}_{0}\|\widetilde{\Sigma}_{\tau}^{n+1}\circ S_{-(t-\tau)}^{0}\|
 _{\mathcal{Z}_{t}^{\lambda_{n+1}(t),\mu_{n+1}(t);1}}d\tau\leq\int^{t}_{0}\|\widetilde{\Sigma}_{\tau}^{n+1}\|
 _{\mathcal{Z}_{\tau}^{\lambda_{n+1}(t),\mu_{n+1}(t);1}}d\tau,$$
 then by Lemma 4.1,
 $$\frac{d^{+}}{dt}\|h^{n+1}\|_{\mathcal{Z}_{t}^{\lambda_{n+1}(t),\mu_{n+1}(t);1}}\leq-\frac{K}{2}\|\nabla_{x}h^{n+1}\|_{\mathcal{Z}_{t}
 ^{\lambda_{n+1},\mu_{n+1};1}}
 -\frac{K}{2}\|\nabla_{v}h^{n+1}\|_{\mathcal{Z}_{t}^{\lambda_{n+1},\mu_{n+1};1}}$$
 $$-\frac{K}{2}\|v\nabla_{x}h^{n+1}\|_{\mathcal{Z}_{t}^{\lambda'_{n+1},\mu'_{n+1};1}}
 -\frac{K}{2}\|v\nabla_{v}h^{n+1}\|_{\mathcal{Z}_{t}^{\lambda'_{n+1},\mu'_{n+1};1}}$$
 $$+\|E[f^{n}]\|_{\mathcal{F}^{\lambda_{n+1}t+\mu_{n+1}}}
 \|\nabla_{v}h^{n+1}\|_{\mathcal{Z}_{t}^{\lambda_{n+1},\mu_{n+1};1}}+\|B[f^{n}]\|_{\mathcal{F}^{\bar{\lambda}_{n+1}t+\bar{\mu}_{n+1}}}
 \|\nabla_{v}h^{n+1}\times v\|_{\mathcal{Z}_{t}^{\lambda'_{n+1},\mu'_{n+1};1}}$$
 $$+\|E[h^{n+1}]\|_{\mathcal{F}^{\lambda_{n+1}t+\mu_{n+1}}}
 \|\nabla_{v}f^{n}\|_{\mathcal{Z}_{t}^{\lambda_{n+1},\mu_{n+1};1}}+\|B[h^{n+1}]\|_{\mathcal{F}^{\lambda_{n+1}t+\mu_{n+1}}}
 \|\nabla_{v}f^{n}\times v\|_{\mathcal{Z}_{t}^{\lambda_{n+1},\mu_{n+1};1}}$$
 $$+\|E[h^{n}]\|_{\mathcal{F}^{\lambda_{n+1}t+\mu_{n+1}}}
 \|\nabla_{v}h^{n}\|_{\mathcal{Z}_{t}^{\lambda_{n+1},\mu_{n+1};1}}+\|B[h^{n}]\|_{\mathcal{F}^{\lambda_{n+1}t+\mu_{n+1}}}
 \|\nabla_{v}h^{n}\times v\|_{\mathcal{Z}_{t}^{\lambda_{n+1},\mu_{n+1};1}},$$
 where the third term of the right-hand inequality uses (6.2) of Theorem 6.1 in  the following section 6.
 First, we easily get $\|E[h^{n}]\|_{\mathcal{F}^{\lambda_{n+1}t+\mu_{n+1}}}\leq C\|\nabla h^{n}\|_{\mathcal{Z}_{t}^{\lambda_{n+1},\mu_{n+1};1}},$
 $\|B[h^{n}]\|_{\mathcal{F}^{\lambda_{n+1}t+\mu_{n+1}}}\leq C\|\nabla h^{n}\|_{\mathcal{Z}_{t}^{\lambda_{n+1},\mu_{n+1};1}}.$ Moreover,
 $$ \|\nabla_{v}f^{n}\|_{\mathcal{Z}_{t}^{\lambda_{n+1},\mu_{n+1};1}}\leq\sum^{n}_{i=1}\|\nabla_{v}h^{i}\|_{\mathcal{Z}_{t}^{\lambda_{n+1},\mu_{n+1};1}}
 \leq C\sum^{n}_{i=1}\frac{\|h^{i}\|_{\mathcal{Z}_{t}^{\lambda_{i+1},\mu_{i+1};1}}}{\min\{\lambda_{i}-\lambda_{n+1},\mu_{i}-\mu_{n+1}\}},$$
 $$ \|\nabla_{v}f^{n}\times v\|_{\mathcal{Z}_{t}^{\lambda_{n+1},\mu_{n+1};1}}\leq\sum^{n}_{i=1}\|\nabla_{v}h^{i}\times v\|_{\mathcal{Z}_{t}^{\lambda_{n+1},\mu_{n+1};1}}
 \leq C\sum^{n}_{i=1}\frac{\|h^{i}\|_{\mathcal{Z}_{t}^{\lambda_{i},\mu_{i};1}}}{\min\{\lambda_{i}-\lambda_{n+1},\mu_{i}-\mu_{n+1}\}}.$$

 We gather the above estimates,
 $$\frac{d^{+}}{dt}\|h^{n+1}\|_{\mathcal{Z}_{t}^{\lambda_{n+1}(t),\mu_{n+1}(t);1}}\leq \bigg( C\sum^{n}_{i=1}\frac{\delta_{i}}
 {\min\{\lambda_{i}-\lambda_{n+1},\mu_{i}-\mu_{n+1}\}}-\frac{K}{2}\bigg)\|\nabla h^{n+1}\|_{\mathcal{Z}_{t}^{\lambda_{n+1},\mu_{n+1};1}}$$
 $$+\bigg( C\sum^{n}_{i=1}\frac{\delta_{i}}
 {\min\{\lambda_{i}-\bar{\lambda}_{n+1},\mu_{i}-\bar{\mu}_{n+1}\}}-\frac{K}{2}\bigg)\|\nabla h^{n+1}\|_{\mathcal{Z}_{t}^{\lambda'_{n+1},\mu'_{n+1};1}}
 +\frac{\delta^{2}_{n}}{\min\{\lambda_{n}-\bar{\lambda}_{n+1},\mu_{n}-\bar{\mu}_{n+1}\}}.$$

We may choose $$\delta_{n+1}=\frac{\delta^{2}_{n}}{\min\{\lambda_{n}-\bar{\lambda}_{n+1},\mu_{n}-\bar{\mu}_{n+1}\}},$$
if \begin{align}
C\max\bigg\{\sum^{n}_{i=1}\frac{\delta_{i}}
 {\min\{\lambda_{i}-\lambda_{n+1},\mu_{i}-\mu_{n+1}\}},\sum^{n}_{i=1}\frac{\delta_{i}}
 {\min\{\lambda_{i}-\bar{\lambda}_{n+1},\mu_{i}-\bar{\mu}_{n+1}\}}\bigg\}\leq\frac{K}{2}\end{align}
  holds.

 We choose $\lambda_{i}-\lambda_{i+1}=\mu_{i}-\mu_{i+1}=\frac{\Lambda}{i^{2}},$ where $\Lambda>0$ is arbitrarily small. Then for $i\leq n,\lambda_{i}-\lambda_{n+1}\geq
 \frac{\Lambda}{i^{2}},$ and $\delta_{n+1}\leq\delta^{2}_{n}n^{2}/\Lambda.$ Next we need to check that $\sum^{\infty}_{n=1}\delta_{n}n^{2}<\infty.$
 In fact, we choose $ K $ large enough and T small enough such that $\lambda_{0}-KT\geq \lambda_{\ast},\mu_{0}-KT\geq \mu_{\ast},$ and (4.13) holds, where
  $\lambda_{0}>\lambda_{\ast},\mu_{0}>\mu_{\ast}$ are fixed.

  If $\delta_{1}=\delta,$ then $\delta_{n}=n^{2}\frac{\delta^{2^{n}}}{\Lambda^{n}}(2^{2})^{2^{n-2}}(4^{2})^{2^{n-2}}\ldots((n-1)^{2})^{2}n^{2}.$
  To prove the sequence convergence for $\delta$ small enough, by induction that $\delta_{n}\leq z^{a^{n}},$ where $z$ small enough and $a\in(1,2).$ We
  claim that the conclusion holds for $n+1.$ Indeed, $\delta_{n+1}\leq\frac{z^{2a^{n}}}{\Lambda}n^{2}\leq z^{a^{n+1}}\frac{z^{(2-a)a^{n}}n^{2}}{\Lambda}.$ If $ z$ is
  so small that $z^{(2-a)a^{n}}\leq\frac{\Lambda}{n^{2}}$ for all $n\in\mathbb{N},$ then $\delta_{n+1}\leq z^{a^{n+1}},$ this concludes the local-time argument.
  \begin{rem} It is worthy to note that there are  resonances to occur in local time  which are  caused by
  the action of the magnetic field, in detail, the resonances are from the term $v\times B[f^{n}]\cdot\nabla_{v}h^{n+1}$. This phenomenon is very different from Landau damping in [23] in local time.
  \end{rem}

   \begin{center}
\item\subsection{Global time iteration }
%{\bf\large 1. \quad Introduction }
\end{center}

Based on the estimates of local-time iteration, without loss of generality, sometimes  we only consider the case  $\tau\geq\frac{bt}{1+b},$ where $b$
is small enough.

First, we give deflection estimates to compare the free evolution with the true evolution from the particles trajectories.
  \begin{prop}
           Assume  for any $i\in\mathbb{N},$ $0<i\leq n, $
           $$\sup_{t\geq 0}\|E[h^{i}](t,\cdot)\|_{\mathcal{F}^{(\lambda_{i}-B_{0})t+\mu_{i}}}<\delta_{i},\quad
           \sup_{t\geq0}\|B[h^{i}](t,\cdot)\|_{\mathcal{F}^{(\lambda_{i}-B_{0})t+\mu_{i}}}<\delta_{i}.$$

           And there exist  constants $\lambda_{\star}>B_{0},\mu_{\star}>0$ such that $\lambda_{0}-B_{0}>\lambda'_{0}-B_{0}>\lambda_{1}-B_{0}
           >\lambda'_{1}-B_{0}>\ldots>\lambda_{i}-B_{0}>\lambda'_{i}-B_{0}>\ldots>\lambda_{\star}-B_{0},$
            $\mu_{0}>\mu_{1}>\mu'_{1}>\ldots>\mu_{i}>\mu'_{i}>\ldots>\mu_{\star}.$
           % satisfy the following conditions:

            %$ C(\lambda'_{k})\varepsilon e^{-2\pi(\lambda_{k}-\lambda'_{k})s}\min\bigg\{(s-\tau),\frac{1}{2\pi(\lambda_{k}-\lambda'_{k})}\bigg\}\leq\frac{1}{2}(\mu_{k+1}-\mu'_{k+1}),$

         % and

          %$|B_{0}|C(\lambda'_{k}) e^{-2\pi(\lambda_{k}-\lambda'_{k})s}\min\bigg\{(s-\tau),\frac{1}{2\pi(\lambda_{k}-\lambda'_{k})}\bigg\}\leq C(\lambda'_{k+1})e^{-2\pi(\lambda_{k+1}-\lambda'_{k+1})s},$   for all
          %$\tau\leq s\leq t.$

           % Meanwhile, $\mathbf{here}$  $\mathbf{assume}$  $\mathbf{that}$ $\mathbf{there}$ $\mathbf{exist}$ $\mathbf{a}$  $\mathbf{zero}$  $\mathbf{measure}$
            % $set$  $A$  and sufficiently large constant $V_{M}>0$ such that $$\sup_{\mathbb{R}^{3}\setminus A}|v|\leq V_{M}.$$
          Then we have %for any $k\in\mathbb{Z}^{3},i \in\mathbb{N}$
          $$\|\delta X^{n+1}_{t,\tau}\circ (X^{0}_{\tau,t},V^{0}_{\tau,t})\|_{\mathcal{Z}_{\tau-\frac{bt}{1+b}}^{\lambda'_{n}-B_{0},\mu'_{n}}}
          \leq C\sum^{n}_{i=1}\delta_{i}  e^{-\pi(\lambda_{i}-\lambda'_{i})\tau}\min\bigg\{\frac{(t-\tau)^{2}}{2},
           \frac{1}{2\pi(\lambda_{i}-\lambda'_{i})^{2}}\bigg\},$$

            $$\|\delta V^{n+1}_{t,\tau}\circ (X^{0}_{\tau,t},V^{0}_{\tau,t})\|_{\mathcal{Z}_{\tau-\frac{bt}{1+b}}^{\lambda'_{n}-B_{0},\mu'_{n}}}
            \leq C\sum^{n}_{i=1}\delta_{i} e^{-\pi(\lambda_{i}-\lambda'_{i})\tau}\min\bigg\{\frac{(t-\tau)}{2},
           \frac{1}{2\pi(\lambda_{i}-\lambda'_{i})}\bigg\},$$
           for $0<\tau<t,$ $b=b(t,\tau)$ sufficiently small.
           %$$\|\delta X^{k}_{t,\tau}\|_{\mathcal{Z}_{t+\sigma}^{\lambda'_{k},\mu'_{k}}}\leq C(\lambda'_{k})\varepsilon e^{-2\pi(\lambda_{k}-\lambda'_{k})t}\min\bigg\{\frac{(t-\tau)^{2}}{2},\frac{1}{2\pi(\lambda_{k}-\lambda'_{k})^{2}}\bigg\},$$
    %$$\|\delta V^{k}_{t,\tau}\|_{\mathcal{Z}_{t+\sigma}^{\lambda'_{k},\mu'_{k}}}\leq C(\lambda'_{k})\varepsilon e^{-2\pi(\lambda_{k}-\lambda'_{k})t}\min\bigg\{(t-\tau),\frac{1}{2\pi(\lambda_{k}-\lambda'_{k})}\bigg\}.$$
           \end{prop}
           \begin{rem}
From Proposition 4.3, it is easy to know that the gliding analytic regularity $\lambda> B_{0}$ in the cyclotron damping, comparing with $\lambda>0$ in Landau damping,
here $B_{0}$ is called cyclotron frequency. In other words, resonances in cyclotron damping occur at the cyclotron frequency, not zero frequency (Landau resonances
occur at zero frquency).
\end{rem}

           \begin{prop} Under the assumptions of Proposition 4.3, then
           $$\bigg\|\nabla\Omega^{n+1}X_{t,\tau}-(Id,0)\bigg\|_{\mathcal{Z}_{s+\frac{bt}{1-b}}^{(\lambda'_{n+1}-B_{0})(1-b),\mu'_{n+1}}}<\mathcal{C}_{1}^{n},
           \bigg\|\nabla\Omega^{n+1}V_{t,\tau}-(0,Id)\bigg\|_{\mathcal{Z}_{s+\frac{bt}{1-b}}^{(\lambda'_{n+1}-B_{0})(1-b),\mu'_{n+1}}}<\mathcal{C}_{1}^{n}+
           \mathcal{C}_{2}^{n},$$
           where $\mathcal{C}_{1}^{n}= C\sum^{n}_{i=1}\frac{ e^{-\pi(\lambda_{i}-\lambda'_{i})\tau}\delta_{i}}
           {2\pi(\lambda_{i}-\lambda'_{i})^{2}}\min\bigg\{\frac{(t-\tau)^{2}}{2},
           1\bigg\},\mathcal{C}_{2}^{n}= C\sum^{n}_{i=1}\frac{ e^{-\pi(\lambda_{i}-\lambda'_{j})\tau}\delta_{i}}{2\pi(\lambda_{i}-\lambda'_{i})}
           \min\bigg\{t-\tau,1\bigg\}.$
           \end{prop}
           \begin{prop} Under the assumptions of Proposition 4.3, then
           $$\bigg\|\Omega^{i}X_{t,\tau}-\Omega^{n}X_{t,\tau}\bigg\|_{\mathcal{Z}_{s-\frac{bt}{1+b}}^{(\lambda'_{n}-B_{0})(1-b),\mu'_{n}}}
           < \mathcal{C}_{1}^{i,n},\bigg\|\Omega^{i}V_{t,\tau}-\Omega^{n}V_{t,\tau}\bigg\|_{\mathcal{Z}_{s-\frac{bt}{1+b}}^{(\lambda'_{n}-B_{0})(1-b),\mu'_{n}}}
           <\mathcal{C}_{1}^{i,n}+ \mathcal{C}_{2}^{i,n},$$
           where $\mathcal{C}_{1}^{i,n}=C\sum^{n}_{j=i+1} \frac{ e^{-\pi(\lambda_{j}-\lambda'_{j})\tau}\delta_{j}}
           {2\pi(\lambda_{j}-\lambda'_{j})^{2}}\min\bigg\{\frac{(t-\tau)^{2}}{2},
           1\bigg\},\mathcal{C}_{2}^{i,n}=C\sum^{n}_{j=i+1}
           \frac{ e^{-\pi(\lambda_{j}-\lambda'_{j})\tau}\delta_{j}}{2\pi(\lambda_{j}-\lambda'_{j})}
           \min\bigg\{t-\tau,1\bigg\}.$
           \end{prop}
           \begin{rem} Note that $\mathcal{C}_{1}^{i,n},\mathcal{C}_{2}^{i,n}$  decay fast as $\tau\rightarrow\infty,i\rightarrow\infty,$ and uniformly in
           $n\geq i,$ since the sequence $\{\delta_{n}\}^{\infty}_{n=1}$ has fast convergence. Hence, if $r\in\mathbb{N}$ is given, we shall have
           \begin{align}
           \mathcal{C}_{1}^{i,n}\leq\omega_{i,n}^{r,1},\quad \textmd{and}\quad \mathcal{C}_{2}^{i,n}\leq\omega_{i,n}^{r,2},\quad \textmd{all}\quad r\geq1,
           \end{align}
           with $\omega_{i,n}^{r,1}=C^{r}_{\omega}\sum^{n}_{j=i+1} \frac{ \delta_{j}}
           {2\pi(\lambda_{j}-\lambda'_{j})^{2+r}}\frac{\min\{\frac{(t-\tau)^{2}}{2},
           1\}}{(1+\tau)^{r}}$ and $\omega_{i,n}^{r,2}=C^{r}_{\omega}\sum^{n}_{j=i+1} \frac{ \delta_{j}}
           {2\pi(\lambda_{j}-\lambda'_{j})^{1+r}}\frac{\min\{\frac{(t-\tau)^{2}}{2},
           1\}}{(1+\tau)^{r}},$ for some absolute constant $C^{r}_{\omega}$ depending only on $ r.$
           \end{rem}
           \begin{prop} Under the assumptions of Proposition 4.3, then
           $$\bigg\|(\Omega^{i}_{t,\tau})^{-1}\circ\Omega_{t,\tau}^{n}-Id\bigg\|_{\mathcal{Z}_{s-\frac{bt}{1+b}}^{(\lambda'_{n}-B_{0})(1-b),\mu'_{n}}}
            <\mathcal{C}_{1}^{i,n}+ \mathcal{C}_{2}^{i,n}.$$
          \end{prop}

            To give a self-consistent estimate, we have to control each term of Eq.(4.10): I,II,III,IV,V. And the most difficult terms are $I,II,V,$ respectively,
            because there is some resonance phenomena occurring in these terms that makes the  propagated wave away from equilibrium.

            Let us first consider the first two terms I, II, because they have  the same  proofs.
 \begin{align}
 &I^{n+1,n}(t,x)+II^{n+1,n}(t,x)=\int^{t}_{0}\int_{\mathbb{R}^{3}}-(\mathcal{E}^{n+1}_{s,t}\cdot G_{s,t}^{n})(s, X^{0}_{s,t}(x,v),V^{0}_{s,t}(x,v))\notag\\
 &-(F^{n+1}_{s,t}\cdot G_{s,t}^{n,v})(s, X^{0}_{s,t}(x,v),V^{0}_{s,t}(x,v))dvds.
 \end{align}

 To handle these terms, we start by introducing
\begin{align}
&\bar{G}^{n}_{s,t}=\nabla_{v}f^{0}+\sum^{n}_{i=1}\nabla_{v}(h^{i}\circ\Omega^{i-1}_{s,t}),
\quad\bar{G}^{n,v}_{s,t}=\nabla_{v}\times(f^{0}v)+\sum^{n}_{i=1}\nabla_{v}\times((h^{i}v)\circ\Omega^{i-1}_{s,t}),
\end{align}
and the error terms $\mathcal{R}_{0},\tilde{\mathcal{R}}_{0},\mathcal{R}_{1},\tilde{\mathcal{R}}_{1}$ are defined by
\begin{align}
\mathcal{R}_{0}=\int^{t}_{0}\int_{\mathbb{R}^{3}}((B[h^{n+1}]\circ\Omega^{n}_{s,t}(x,v)-B[h^{n+1}])\cdot G_{s,t}^{n,v})(s, X^{0}_{s,t}(x,v),V^{0}_{s,t}(x,v))dvds,
\end{align}
\begin{align}
\tilde{\mathcal{R}}_{0}=\int^{t}_{0}\int_{\mathbb{R}^{3}}(B[h^{n+1}]\cdot (G_{s,t}^{n,v}-\bar{G}_{s,t}^{n,v}))(s, X^{0}_{s,t}(x,v),V^{0}_{s,t}(x,v))dvds,
\end{align}
\begin{align}
\mathcal{R}_{1}=\int^{t}_{0}\int_{\mathbb{R}^{3}}((E[h^{n+1}]\circ\Omega^{n}_{s,t}(x,v)-E[h^{n+1}])\cdot G_{s,t}^{n,v})(s, X^{0}_{s,t}(x,v),V^{0}_{s,t}(x,v))dvds,
\end{align}
\begin{align}
\tilde{\mathcal{R}}_{1}=\int^{t}_{0}\int_{\mathbb{R}^{3}}(E[h^{n+1}]\cdot (G_{s,t}^{n,v}-\bar{G}_{s,t}^{n,v}))(s, X^{0}_{s,t}(x,v),V^{0}_{s,t}(x,v))dvds,
\end{align}
then we can decompose $$I^{n+1,n}=\bar{I}^{n+1,n}+\mathcal{R}_{1}+\tilde{\mathcal{R}}_{1},\quad II^{n+1,n}=\bar{II}^{n+1,n}+\mathcal{R}_{0}+\tilde{\mathcal{R}}_{0}.$$
Because dealing with the first term  $I$ is the same to the second term $II,$
to simply the proof, here we only prove the second term $II.$
Now first we consider $\bar{II}^{n+1,n},$ which we decompose as
$$\bar{II}^{n+1,n}=\bar{II}_{0}^{n+1,n}+\sum^{n}_{i=1}\bar{II}_{i}^{n+1,n},$$
where
$$\bar{II}_{0}^{n+1,n}(t,x)=\int^{t}_{0}\int_{\mathbb{R}^{3}}B[h^{n+1}](\tau,x'(\tau,x,v),
v'(\tau,x,v))\cdot(\nabla'_{v}\times(f^{0}v))(v')dvd\tau,$$
$$\bar{II}_{i}^{n+1,n}(t,x)=\int^{t}_{0}\int_{\mathbb{R}^{3}}B[h^{n+1}](\tau,x'(\tau,x,v),
v'(\tau,x,v))\cdot(\nabla'_{v}\times(h_{\tau}^{i}v)\circ\Omega^{i-1}_{t,\tau})(\tau,x'(\tau,x,v),v'(\tau,x,v))dvd\tau.$$
Since both $B[h^{n+1}](\tau,x'(\tau,x,v),
v'(\tau,x,v))$ and $(\nabla'_{v}\times(h_{\tau}^{i}v)\circ\Omega^{i-1}_{t,\tau})(\tau,x'(\tau,x,v),v'(\tau,x,v))$ have the variable $x,$
then applying Fourier transform in $ x, $ we get
\begin{align}
|[B[h^{n+1}]\cdot(\nabla'_{v}\times(h_{\tau}^{i}v)\circ\Omega^{i-1}_{t,\tau})(x'(\tau,x,v))]^{\wedge}(k)|
=|\sum_{l}\widehat{B[h^{n+1}]}(k-l)(\nabla'_{v}\times(h_{\tau}^{i}v)\circ\Omega^{i-1}_{t,\tau})^{\wedge}(l)|.
\end{align}
It is easy to see that  Eq.(4.21) has two waves of distinct frequencies $k-l,l,$ which may interact. When interacting at certain particular times, the influence
of the waves becomes very strong: this is known in plasma physics as  plasma echo (we explain it in  detail in next section),  and can be thought of as a kind of
resonance. It is the key point in our paper. Based on the iteration scheme different from that in Mouhot and Villani'papaer [23], we have to deal with the new term
$v\times B[f^{n}]\cdot\nabla_{v}h^{n+1}$ which also generates resonances.
\begin{prop} ($\textbf{M}$ain term $I$)
Assume %$\lambda^{\star}>\lambda_{1}>\lambda_{2}>\ldots>\lambda_{n}>\ldots>0,\mu^{\star}>\mu_{1}>\mu_{2}>\ldots>\mu_{n}>\ldots>0,\mu'_{n+1}=\mu_{n+1}+\eta(\frac{t-\tau}{1+t}),
%\nu'_{n+1}=\lambda_{n+1}(1+b)|\tau-\frac{bt}{1+b}|+\mu'_{n+1},$ where
 $b(t,\tau,\Omega)\geq0,\eta>0$ small. And there exist  constants $\lambda_{\star}>B_{0},\mu_{\star}>0$ such that $\lambda_{0}-B_{0}>\lambda'_{0}-B_{0}
 >\lambda_{1}-B_{0}>\lambda'_{1}-B_{0}>\ldots>\lambda_{i}-B_{0}>\lambda'_{i}-B_{0}>\ldots>\lambda_{\star}-B_{0},$
 $\mu_{0}>\mu_{1}>\mu''_{1}>\mu'_{1}>\ldots>\mu_{i}>\mu''_{i}>\mu'_{i}>\ldots>\mu_{\star}.$

We have
$$\| \bar{II}_{i}^{n+1,n}(t,\cdot)\|_{\mathcal{F}^{(\lambda'_{n}-B_{0}) t+\mu'_{n}}}$$
 $$\leq C\int^{t}_{0} K^{n+1}_{1}(t,\tau)\| \nabla'_{v}\times
((h^{i}_{\tau}v)\circ\Omega^{i-1}_{t,\tau}))-\langle\nabla'_{v}\times
((h^{i}_{\tau}v)\circ\Omega^{i-1}_{t,\tau}))\rangle\|_{\mathcal{Z}^{(\lambda'_{i}-B_{0})(1+b),
\mu'_{i};1}_{\tau-\frac{bt}{1+b}}}$$
$$
\| B[h^{n+1}]\|_{_{\mathcal{F}^{\nu}}}d\tau+\int^{t}_{0} K^{n+1}_{0}(t,\tau)\| \langle\nabla'_{v}\times
((h^{i}_{\tau}v)\circ\Omega^{i-1}_{t,\tau}))\rangle\|_{\mathcal{C}^{(\lambda_{i}'-B_{0})(1+b);1}}\cdot\| B[h^{n+1}]\|_{_{\mathcal{F}^{\nu}}}d\tau,$$
 where
$$\nu=\max\bigg\{(\lambda_{n}'-B_{0})\tau+\mu''_{n}-\frac{1}{2}(\lambda_{n}'-B_{0})b(t-\tau),0\bigg\},$$
$$K^{n}_{0}(t,\tau)=e^{-\pi(\lambda'_{i}-\lambda'_{n})(t-\tau)},$$
$$K^{n+1}_{1}(t,\tau)=\sup_{k_{3},l_{3}\in\mathbb{Z}}e^{-2\pi(\mu'_{i}-\mu'_{n})|l_{3}|}e^{-\pi(\lambda'_{i}-\lambda'_{n})|k_{3}(t-\tau)+l_{3}\tau|}
  e^{-2\pi(\frac{\lambda'_{n}}{2}(\tau-\tau')+\mu''_{n}-\mu'_{n})|k_{3}-l_{3}|}.$$
\end{prop}
\begin{lem}
We have $\frac{1}{\tau}\| B[h^{n+1}]\|_{\mathcal{F}^{\nu}}\leq C \sup_{0\leq s\leq \tau}\| \rho[h^{n+1}]\|_{\mathcal{F}^{(\lambda'_{n}-B_{0})s+\mu'_{n}}},\quad
\tau\|h^{n}\|_{\mathcal{Z}^{(\lambda'_{n}-B_{0})(1+b),\mu'_{i};1}_{\tau-\frac{bt}{1+b}}}\leq C
\|h^{n}\|$
$_{\mathcal{Z}^{(\bar{\lambda}'_{n}-B_{0})(1+b),\bar{\mu}'_{i};1}_{\tau-\frac{bt}{1+b}}}, $ where $\lambda_{i}'<\bar{\lambda}_{i}'<\lambda_{i},
\mu_{i}'<\bar{\mu}_{i}'<\mu_{i},\tau>0.$
\end{lem}
$Proof.$ Since $\partial_{t}B=\nabla\times E,\quad E=W(x)\ast\rho(f),$
$B=\int^{t}_{0}\nabla\times E(s,x)ds=\int^{t}_{0}\nabla\times( W(x)\ast\rho(f))(s,x)ds,$
then $\|B[h^{n+1}]\|_{\mathcal{F}^{\nu}}\leq\int^{\tau}_{0}\|\nabla\times E[h^{n+1}]\|_{\mathcal{F}^{\nu}}ds$
$\leq\int^{\tau}_{0}\|\nabla E[h^{n+1}]\|_{\mathcal{F}^{\nu}}ds\leq \tau \sup_{0\leq s\leq \tau}\| \rho[h^{n+1}]\|_{\mathcal{F}^{(\lambda'_{n}-B_{0})s+\mu'_{n}}}, $
%$$\leq\int^{\tau}_{0}e^{-2\pi(\nu''_{n+1}\tau-\nu'_{n+1}\tau)}\| \rho[h^{n+1}]\|_{\mathcal{\dot{F}}^{\nu''_{n+1}}}d\tau,$$
%where we set $\nu''_{n+1}=\lambda_{n+1}(1+b)|\tau-\frac{bt}{1+b}|+\mu''_{n+1},\mu''_{n+1}=\mu_{n+1}+2\eta(\frac{t-\tau}{1+t}).$
%$\mathbf{maybe}$ $\mathbf{we}$ $\mathbf{need}$ $\mathbf{take}$ $\mathcal{F}$ $\mathbf{norm}$ on $\mathbf{x_{\perp}}$ in order to remove $\nabla_{x_{\perp}}.$
where we use $\nu< (\lambda'_{n}-B_{0})\tau+\mu'_{n}.$
The last inequality can be obtained from the definition of the norm $\mathcal{Z}^{\lambda,\mu;1}_{t}$ and (vi)-(vii) of Proposition 2.5.
%Indeed, leaving apart the small-time case, we assume $\tau\geq\frac{bt}{1+b},$ then
%$$\nu''_{n+1}=\lambda_{n+1}(1+b)|\tau-\frac{bt}{1+b}|+\mu_{n+1}+2\eta(\frac{t-\tau}{1+t})$$
%$$=\lambda_{n+1}\tau+\mu_{n+1}+(\frac{2\eta}{1+t}-\lambda_{n+1}b)(t-\tau),$$
%which is bounded by $\lambda_{n+1}\tau+\mu_{n+1}$ if we choose $\eta\leq\frac{1}{2}\lambda^{\star}(1+t).$

%So we obtain $$\|B[h^{n+1}]\|_{\mathcal{F}^{\nu'_{n+1}}}\leq\int^{t}_{0}e^{-2\pi(\nu''_{n+1}\tau-\nu'_{n+1}\tau)}\| \rho[h^{n+1}]\|_{\mathcal{\dot{F}}^{\lambda_{n+1}\tau+\mu_{n+1}}}d\tau.$$
\begin{cor}
From the above statement, we have
$$
\|\bar{II}_{i}^{n+1,n}(t,\cdot)\|_{\mathcal{F}^{(\lambda'_{n}-B_{0}) t+\mu'_{n}}}
\leq\int^{t}_{0}K^{n+1}_{0}(t,\tau)\delta_{i}\| \rho[h^{n+1}]\|_{\mathcal{F}^{(\lambda'_{n}-B_{0}) \tau+\mu'_{n}}}d\tau$$
$$+\int^{t}_{0}K^{n+1}_{1}(t,\tau)(1+\tau)\delta_{i}\|\rho[h^{n+1}]\|_{\mathcal{F}^{(\lambda'_{n}-B_{0})\tau+\mu'_{n}}}d\tau,
$$
where $K^{n}_{0}(t,\tau)=e^{-\pi(\lambda'_{i}-\lambda'_{n})(t-\tau)},$
and
$$K^{n+1}_{1}(t,\tau)=e^{-2\pi(\mu'_{i}-\mu'_{n})|l_{3}|}e^{-\pi(\lambda'_{i}-\lambda'_{n})|k_{3}(t-\tau)+l_{3}\tau|}
  e^{-2\pi(\frac{(\lambda'_{n}-B_{0})}{2}(\tau-\tau')+\mu''_{n}-\mu'_{n})|k_{3}-l_{3}|}.$$
\end{cor}
\begin{prop}($\textbf{E}$rror term $I$)
 $$\|\mathcal{R}_{0}(t,\cdot)\|_{\mathcal{F}^{(\lambda'_{n}-B_{0}) t+\mu'_{n}}}
 \leq C\bigg(C'_{0}+\sum^{n}_{i=1}\delta_{i}\bigg)\bigg(\sum^{n}_{i=1}\frac{\delta_{i}}{(\lambda_{i}-\lambda'_{i})^{5}}\bigg)
\int^{t}_{0}\|\rho[h^{n+1}]\|_{\mathcal{F}^{(\lambda'_{n}-B_{0})\tau+\mu'_{n}}}\frac{d\tau}{(1+\tau)^{2}}.$$
 \end{prop}
 \begin{prop}($\textbf{E}$rror term $II$)
$$
\|\tilde{\mathcal{R}}_{0}(t,\cdot)\|_{\mathcal{F}^{(\lambda'_{n}-B_{0}) \tau+\mu'_{n}}}
\leq\bigg(C^{4}_{\omega}\bigg(C'_{0}+\sum^{n}_{i=1}\delta_{i}\bigg)\bigg(\sum^{n}_{j=1}\frac{\delta_{j}}{2\pi(\lambda_{j}-\lambda'_{j})^{6}}\bigg)
+\sum^{n}_{i=1}\delta_{i}\bigg)
\int^{t}_{0}\|\rho\|_{\mathcal{F}^{(\lambda'_{n}-B_{0}) \tau+\mu'_{n}}}\frac{1}{(1+\tau)^{2}}d\tau$$
$$=\int^{t}_{0}\tilde{K}_{1}^{n+1}\|\rho\|_{\mathcal{F}^{(\lambda'_{n}-B_{0}) \tau+\mu'_{n}}}\frac{1}{(1+\tau)^{2}}d\tau.$$
\end{prop}
\begin{prop}($\textbf{M}$ain term $V$)
$$
\|V\|_{\mathcal{F}^{(\lambda'_{n}-B_{0})t+\mu'_{n}}}\leq\int^{t}_{0}e^{-\pi(\lambda_{n}-\lambda'_{n})|k(t-s)+ls|}\bigg(\sum^{n}_{i=1}\delta_{i}\bigg)
\| h^{n+1}\circ\Omega^{n}_{s,t}(x,v)\|_{\mathcal{Z}_{s+\frac{bt}{1-b}}^{(\lambda'_{n}-B_{0})(1-\frac{1}{2}b),\mu'_{n};1}}ds.
$$
\end{prop}
\begin{prop}($\textbf{M}$ain term $III$)
$$
\sup_{0\leq s\leq t}\|
h^{n+1}\circ\Omega^{n}_{t,\tau}\|_{\mathcal{Z}^{(\lambda'_{n}-B_{0})(1-\frac{1}{2}b),
\mu'_{n};1}_{s+\frac{bt}{1-b}}}\leq\delta^{2}_{n}+\bigg(\sum^{n}_{i=1}\delta_{i}\bigg)
\sup_{0\leq s\leq t}\|\rho^{n+1}\|_{\mathcal{F}^{(\lambda'_{n}-B_{0})s+\mu'_{n}}}.
$$
\end{prop}

   \begin{center}
\item\subsection{  The proof of main theorem}
%{\bf\large 1. \quad Introduction }
\end{center}

%$$K^{n}_{1,2}(t,\tau)=\sup_{k_{1,2},l_{1,2}} e^{-2\pi(\mu_{i} -\mu_{n+1})|l_{1,2}|}$$
 %$$\cdot e^{-2\pi(\lambda_{i}-\lambda_{n+1})|k_{1,2}\frac{\sin\Omega(t-\tau)}
%{\Omega}+l_{1,2}\frac{\sin\Omega\tau}{\Omega}|}e^{-2\pi(\frac{\lambda_{n+1}}{2}(|\frac{\sin\Omega\tau}{\Omega}|
%-|\frac{\sin\Omega\tau'}{\Omega}|)+\mu'_{n+1}-\mu_{n+1})|k_{1,2}-l_{1,2}|},$$
  %$$K^{n}_{3}(t,\tau)=\sup_{k_{3},l_{3}}e^{-2\pi(\mu_{i}-\mu_{n+1})|l_{3}|}e^{-\pi(\lambda_{i}-\lambda_{n+1})|k_{3}(t-\tau)+l_{3}\tau|}
 % e^{-2\pi(\frac{\lambda_{n+1}}{2}(\tau-\tau')+\mu'_{n+1}-\mu_{n+1})|k_{3}-l_{3}|}.$$
$$$$
$\mathbf{Step }$ $\mathbf{2}.$
Note if $\varepsilon$ in (0.13) is small enough, up to slightly lowering $\lambda_{1},$ we may choose all parameters in such a way that
$\lambda_{k}-B_{0},\lambda'_{k}-B_{0}\rightarrow\lambda_{\infty}-B_{0}>\underline{\lambda}-B_{0}\quad \textmd{and}\quad\mu_{k},
\mu'_{k}\rightarrow\mu_{\infty}>\underline{\mu},\quad \textmd{as}
\quad k\rightarrow\infty;$
then we pick up $D>0$ such that
$\mu_{\infty}-\lambda_{\infty}(1+D)D\geq\mu'_{\infty}>\underline{\mu},$
and we let $b(t)=\frac{D}{1+t}.$
From the iteration, we have,  for all $k\geq2,$
\begin{align}
\sup_{0\leq\tau\leq t}\|h^{k}_{\tau}\circ\Omega^{k-1}_{t,\tau}\|_{\mathcal{Z}^{\lambda_{\infty}(1+b),\mu_{\infty};1}_{\tau-\frac{bt}{1+b}}}\leq\delta_{k},
\end{align}
where $\sum^{\infty}_{k=2}\delta_{k}\leq C\delta.$
Choosing $\tau=t$ in (3.22) yields
$\sup_{0\leq\tau\leq t}\|h^{k}_{\tau}\|_{\mathcal{Z}^{(\lambda_{\infty}-B_{0})(1+D),\mu_{\infty};1}_{t-\frac{Dt}{1+D+t}}}\leq\delta_{k}.$
This implies that
$\sup_{ t\geq0}\|h^{k}_{t}\|_{\mathcal{Z}^{(\lambda_{\infty}-B_{0})(1+D),\mu_{\infty}-\lambda_{\infty}(1+D)D;1}_{t}}\leq\delta_{k}.$
In particular, we have a uniform estimate on $h^{k}_{t}$ in $\mathcal{Z}^{(\lambda_{\infty}-B_{0}),\mu'_{\infty};1}_{t}.$ Summing up over $ k $ yields for
$f=f^{0}+\sum^{\infty}_{k=1}h^{k}$ the estimate
\begin{align}
\sup_{t\geq0}\|f(t,\cdot)-f^{0}\|_{\mathcal{Z}^{\lambda_{\infty}-B_{0},\mu'_{\infty};1}_{t}}\leq C\delta.
\end{align}
From (viii) of Proposition 2.5, we can deduce from (4.21) that
\begin{align}
\sup_{t\geq0}\|f(t,\cdot)-f^{0}\|_{\mathcal{Y}^{\underline{\lambda}-B_{0},\underline{\mu}}_{t}}\leq C\delta.
\end{align}
Moreover, $\rho=\int_{\mathbb{R}^{3}}fdv$ satisfies similarly
$\sup_{t\geq0}\|\rho(t,\cdot)\|_{\mathcal{F}^{(\lambda_{\infty}-B_{0})t+\mu_{\infty}}}\leq C\delta.$
It follows that $|\hat{\rho}(t,k)|\leq C\delta $
$e^{-2\pi(\lambda_{\infty}-B_{0})|k_{3}|t}e^{-2\pi\mu_{\infty}|k|}$ for any $k\neq 0.$ On the one hand, by Sobolev embedding,
we deduce that for any $r\in\mathbb{N},$
$$\|\rho(t,\cdot)-\langle\rho\rangle\|_{C^{r}(\mathbb{T}^{3})}\leq C_{r}\delta e^{-2\pi(\lambda'-B_{0})t};$$
on the other hand, multiplying $\hat{\rho}$ by the Fourier transform of $W,$ and $\partial_{t}B=\nabla\times E,$ we see that the electric and magnetic fields
$E,B$ satisfy
\begin{align}
\sup_{t\geq0}\|E(t,\cdot)\|_{\mathcal{F}^{(\lambda'-B_{0})t+\mu'}}\leq\delta,\quad \sup_{t\geq0}\|E(t,\cdot)\|_{\mathcal{F}^{(\lambda'-B_{0})t+\mu'}}\leq\delta,
%|\hat{E}(t,k)|\leq C\delta e^{-2\pi(\lambda'-B_{0})|k_{3}|t}e^{-2\pi\mu'|k|},\quad|\hat{B}(t,k)|\leq C\delta e^{-2\pi(\lambda'-B_{0})|k_{3}|t}e^{-2\pi\mu'|k|},
\end{align}
for some $\lambda_{0}>\lambda'>\underline{\lambda},$ $\mu_{0}>\mu'>\underline{\mu}.$

Now, from (4.22), we have, for any fixed $(k_{3},\eta_{3})\in\mathbb{Z}\times\mathbb{R}$ and any $t\geq0,$
\begin{align}
|\hat{f}(t,k,\eta_{1}+k_{1}\frac{\sin\Omega t}{\Omega},\eta_{2}+k_{2}\frac{\sin\Omega t}{\Omega},\eta_{3}+k_{3}t)-\hat{f}^{0}(\eta)|\leq C\delta e^{-2\pi\mu'|k_{3}|}e^{-2\pi(\lambda'-B_{0})|\eta_{3}|},
\end{align}
this finishes the proof of Theorem 0.1.

         $$$$
         \begin{center}
\item\section{ Dynamical behavior of the particles' trajectory}
%{\bf\large 1. \quad Introduction }
\end{center}

          To prove Proposition 4.3,  by the classical Picard iteration, we only
          need to  consider the following equivalent equations
           \begin{align}
   \left\{\begin{array}{l}
   \frac{d}{dt}\delta X^{n+1}_{t,\tau}(x,v)=\delta V^{n+1}_{t,\tau}(x,v),\\
  \frac{d}{dt}\delta V^{n+1}_{t,\tau}(x,v)=\delta V^{n+1}_{t,\tau}(x,v)\times B_{0} +E[f^{n}](t,\delta X^{n}_{t,\tau}(x,v)+X^{0}_{t,\tau}(x,v)),\\
 \delta X^{n+1}_{\tau,\tau}(x,v)=0,\delta V^{n+1}_{\tau,\tau}(x,v)=0.
  \end{array}\right.
           \end{align}
            It is easy to check that $$\Omega^{n+1}_{t,\tau}-Id\triangleq(\delta X^{n+1}_{t,\tau},\delta V^{n+1}_{t,\tau})\circ (X^{0}_{\tau,t},V^{0}_{\tau,t})
            = ( X^{n+1}_{t,\tau}\circ(X^{0}_{\tau,t},V^{0}_{\tau,t})-Id, V^{n+1}_{t,\tau}\circ(X^{0}_{\tau,t},V^{0}_{\tau,t})-Id).$$
            Therefore, in order to estimate $( X^{n+1}_{t,\tau}\circ(X^{0}_{\tau,t},V^{0}_{\tau,t})-Id, V^{n+1}_{t,\tau}\circ(X^{0}_{\tau,t},V^{0}_{\tau,t})-Id),$
             we only need to study $(\delta X^{n+1}_{t,\tau},\delta V^{n+1}_{t,\tau})\circ (X^{0}_{\tau,t},V^{0}_{\tau,t}).$

             Now we give a detailed proof of Proposition 4.3.

           $Proof.$ If $n=0,$ first, it is trivial that $\delta V^{0}_{t,\tau}(x,v)=0,$ then Eqs.(5.1) reduces to the following equations
            \begin{align}
   \left\{\begin{array}{l}
   \frac{d}{dt}\delta X^{1}_{t,\tau}(x,v)=\delta V^{1}_{t,\tau}(x,v),\frac{d}{dt}\delta V^{1}_{t,\tau}(x,v)=\delta V^{1}_{t,\tau}(x,v)\times B_{0}+E[f^{0}](t,X^{0}_{t,\tau}(x,v)),\\
\delta X^{1}_{\tau,\tau}(x,v)=0,\delta V^{1}_{\tau,\tau}(x,v)=0.
  \end{array}\right.
           \end{align}
           Then we have
           $$\delta X^{1}_{t,\tau}\circ (X^{0}_{\tau,t},V^{0}_{\tau,t})(x,v)=\int^{t}_{\tau}\delta V^{1}_{s,\tau}\circ (X^{0}_{\tau,s},V^{0}_{\tau,s})(x,v)ds,$$
          $$ \delta V^{1}_{t,\tau}\circ (X^{0}_{\tau,t},V^{0}_{\tau,t})(x,v)=\int^{t}_{\tau}e^{B_{0}(t-s)}E[f^{0}](s,X^{0}_{s,t}(x,v))ds.$$
%so $$[\delta X^{1}_{t,\tau}\circ (X^{0}_{\tau,t},V^{0}_{\tau,t})]^{\wedge}(k,v)=\int^{t}_{\tau}(t-s)\widehat{E[f^{0}]}(s,k)ds,$$
%$$\|\delta X^{1}_{t,\tau}\circ (X^{0}_{\tau,t},V^{0}_{\tau,t})\|_{\mathcal{Z}^{(\lambda_{0}'-B_{0})(1+b),\mu_{0}'}_{\tau-\frac{bt}{1+b}}}
%=\int^{t}_{\tau}(t-s)\|E[f^{0}](s,x+\mathcal{M}(t-s)v)\|_{\mathcal{Z}^{(\lambda_{0}'-B_{0})(1+b),\mu_{0}'}_{\tau-\frac{bt}{1+b}}}ds$$
%$$=\int^{t}_{\tau}(t-s)\|E[f^{0}](s,\cdot)\|_{\mathcal{Z}^{(\lambda_{0}'-B_{0})(1+b),\mu_{0}'}_{s-\frac{bt}{1+b}}}ds.$$
By the definition of $E[f^{0}],$ we know that $\|E[f^{0}](s,\cdot)\|_{\mathcal{Z}^{(\lambda_{0}'-B_{0})(1+b),\mu_{0}'}_{s-\frac{bt}{1+b}}}=0,$ it is trivial that
$\|\delta V^{1}_{t,\tau}\circ (X^{0}_{\tau,t},V^{0}_{\tau,t})\|_{\mathcal{Z}^{(\lambda_{0}'-B_{0})(1+b),\mu_{0}'}_{\tau-\frac{bt}{1+b}}}
   \leq C\delta_{0} e^{-2\pi(\lambda_{0}-\lambda'_{0})\tau}\min\bigg\{\frac{t-\tau}{2},\frac{1}{2\pi(\lambda_{0}-\lambda'_{0})}
        \bigg\}.$

   Similarly,
        $\|\delta X^{1}_{t,\tau}\circ (X^{0}_{\tau,t},V^{0}_{\tau,t})$
$\|_{\mathcal{Z}^{(\lambda_{0}'-B_{0})(1+b),\mu_{0}'}_{\tau-\frac{bt}{1+b}}}
\leq C\delta_{0} e^{-2\pi(\lambda_{0}-\lambda'_{0})\tau}
       \min\bigg\{\frac{(t-\tau)^{2}}{2},\frac{1}{2\pi(\lambda_{0}-\lambda'_{0})^{2}}\bigg\}.$

    Suppose for  $n>1,$ both
    $$\|\delta X^{n}_{t,\tau}\circ (X^{0}_{\tau,t},V^{0}_{\tau,t})\|_{\mathcal{Z}^{(\lambda_{n-1}'-B_{0})(1+b),\mu_{n-1}'}_{\tau-\frac{bt}{1+b}}}
    \leq C\sum^{n-1}_{i=1}\delta_{i} e^{-2\pi(\lambda_{i}-\lambda'_{i})\tau}
        \min\bigg\{\frac{(t-\tau)^{2}}{2},\frac{1}{2\pi(\lambda_{i}-\lambda'_{i})^{2}}\bigg\},$$
        and
        $$\|\delta V^{n}_{t,\tau}\circ (X^{0}_{\tau,t},V^{0}_{\tau,t})\|_{\mathcal{Z}^{(\lambda_{n-1}'-B_{0})(1+b),\mu_{n-1}'}_{\tau-\frac{bt}{1+b}}}
        \leq C\sum^{n-1}_{i=1}\delta_{i} e^{-2\pi(\lambda_{i}-\lambda'_{i})\tau}
        \min\bigg\{\frac{(t-\tau)}{2},\frac{1}{2\pi(\lambda_{i}-\lambda'_{i})}\bigg\}.$$
    %$$\|\delta X^{k}_{t,\tau}\|_{\mathcal{Z}_{t+\sigma}^{\lambda'_{k},\mu'_{k}}}\leq C(\lambda'_{k},V_{M})\varepsilon e^{-2\pi(\lambda_{k}-\lambda'_{k})t}\min\bigg\{\frac{(t-\tau)^{2}}{2},\frac{1}{2\pi(\lambda_{k}-\lambda'_{k})^{2}}\bigg\},$$
    %$$\|\delta V^{k}_{t,\tau}\|_{\mathcal{Z}_{t+\sigma}^{\lambda'_{k},\mu'_{k}}}\leq C(\lambda'_{k},V_{M})\varepsilon e^{-2\pi(\lambda_{k}-\lambda'_{k})t}\min\bigg\{(t-\tau),\frac{1}{2\pi(\lambda_{k}-\lambda'_{k})}\bigg\}.$$
   % and
   % $$\|\delta V^{k}_{\tau,t}(\cdot,v)\|_{\mathcal{F}^{\lambda_{k}t+\mu_{k}}}\leq\varepsilon\int^{t}_{\tau}e^{(\lambda_{k}-\lambda_{k-1})s}
    %(1 + |V^{0}_{\tau,s}(x,v)|)ds$$ hold.
    Then for $n+1,$ since $ (\delta X^{n+1}_{t,\tau},\delta V^{n+1}_{t,\tau})$ satisfy
     \begin{align}
   \left\{\begin{array}{l}
   \frac{d}{dt}\delta X^{n+1}_{t,\tau}(x,v)=\delta V^{n+1}_{t,\tau}(x,v),\\
  \frac{d}{dt}\delta V^{n+1}_{t,\tau}(x,v)=\delta V^{n+1}_{t,\tau}(x,v)\times B_{0} +E[f^{n}](t,\delta X^{n}_{t,\tau}(x,v) +X^{0}_{t,\tau}(x,v)),\\
 \delta X^{n+1}_{\tau,\tau}(x,v)=0,\delta V^{n+1}_{\tau,\tau}(x,v)=0.
  \end{array}\right.
           \end{align}
          Then we have
           $\delta V^{n+1}_{t,\tau}=\int^{t}_{\tau}e^{\mathbf{B}_{0}(t-s)}E[f^{n}](s,\delta X^{n}_{s,\tau}(x,v) +X^{0}_{s,\tau}(x,v))ds,$
           and
            $\delta V^{n+1}_{t,\tau}\circ (X^{0}_{\tau,t},V^{0}_{\tau,t})=\int^{t}_{\tau}e^{\mathbf{B}_{0}(t-s)}
            [E[f^{n}]\circ(\delta X^{n}_{s,\tau}\circ (X^{0}_{\tau,s},V^{0}_{\tau,s}))](s,X^{0}_{s,t}(x,v) )ds.$

           Hence
           $$\|\delta V^{n+1}_{t,\tau}\circ (X^{0}_{\tau,t},V^{0}_{\tau,t})\|_{\mathcal{Z}^{(\lambda_{n}'-B_{0})(1+b),\mu_{n}'}_{\tau-\frac{bt}{1+b}}}$$
           $$\leq\int^{t}_{\tau}e^{B_{0}(t-s)}\|[E[f^{n}]\circ(\delta X^{n}_{s,\tau}\circ (X^{0}_{\tau,s},V^{0}_{\tau,s}))](s,X^{0}_{s,t}(x,v) )\|
           _{\mathcal{Z}^{(\lambda_{n}'-B_{0})(1+b),\mu_{n}'}_{\tau-\frac{bt}{1+b}}}ds$$
           $$=\int^{t}_{\tau}\|E[f^{n}]\circ(\delta X^{n}_{s,\tau}\circ (X^{0}_{\tau,s},V^{0}_{\tau,s}))\|
           _{\mathcal{Z}^{(\lambda_{n}'-B_{0})(1+b),\mu_{n}'}_{s-\frac{bt}{1+b}}}ds=\int^{t}_{\tau}\|E[f^{n}](s,\cdot)\|_{\mathcal{F}^{\nu'_{n}}}ds
           \leq\sum^{n}_{i=1}\int^{t}_{\tau}\|E[h^{i}](s,\cdot)\|_{\mathcal{F}^{\nu'_{n}}}ds,$$
where $\nu'_{n}=(\lambda'_{n}-B_{0})|s-b(t-s)|+\mu'_{n}+\|\delta X^{n}_{s,\tau}\circ (X^{0}_{\tau,s},V^{0}_{\tau,s})\|
_{\mathcal{Z}^{(\lambda_{n}'-B_{0})(1+b),\mu_{n}'}_{\tau-\frac{bt}{1+b}}}.$

Note that $\|\delta X^{n}_{s,\tau}\circ (X^{0}_{\tau,s},V^{0}_{\tau,s})\|
_{\mathcal{Z}^{(\lambda_{n}'-B_{0})(1+b),\mu_{n}'}_{\tau-\frac{bt}{1+b}}}\leq \mathcal{C}^{n}_{1},$ if $s\geq\frac{bt}{1+b},$ then
 $\nu'_{n}\leq(\lambda'_{n}-B_{0})s+\mu'_{n}+\mathcal{C}^{n}_{1}\leq (\lambda_{i}-B_{0})s+\mu_{i}-(\lambda_{i}-\lambda'_{n})s$ as soon as
 $\mathcal{C}^{n}_{1}\leq\frac{(\lambda_{i}-B_{0})b(t-s)}{2}(I);$
if $s\leq\frac{bt}{1+b},$ then $\nu'_{n}\leq(\lambda'_{n}-B_{0})bt+\mu'_{n}-(\lambda'_{n}-B_{0})(1+b)s+\mathcal{C}^{n}_{1}
\leq(\lambda'_{n}-B_{0})D+\mu'_{n}-(\lambda_{i}-\lambda'_{n})s+\mathcal{C}^{n}_{1}\leq \mu_{0}-(\lambda_{i}-\lambda'_{n})s$ as soon as
$\mathcal{C}^{n}_{1}\leq\frac{\mu_{0}-\mu'_{n}}{2}(II).$ In order to the feasibility of the conditions $(I)$ and $(II),$ we only need to check that the following
assumption $(\mathbf{I})$ holds
$$2C^{1}_{\omega}\bigg(\sum^{n}_{i=1}\frac{\delta_{i}}{(2\pi(\lambda_{i}-\lambda'_{n}))^{3}}\bigg)\leq\min\bigg\{\frac{(\lambda_{i}-B_{0})b(t-s)}{6},
\frac{\mu_{0}-\mu'_{n}}{2}\bigg\},$$
since $\mathcal{C}^{n}_{1}\leq\omega^{1,2}_{0,n}=2C^{1}_{\omega}\bigg(\sum^{n}_{i=1}\frac{\delta_{i}}{(2\pi(\lambda_{i}-\lambda'_{n}))^{3}}\bigg)
\frac{\min\{\frac{1}{2}(t-s)^{2},1\}}{1+s}.$

          %$$ \leq\sum_{l}e^{2\pi|l_{12}|(\lambda'_{i}|\frac{\sin\Omega t}{\Omega}|+\mu'_{i})} e^{2\pi|l_{3}|(\lambda'_{i}t+\mu'_{i})}|\hat{E[f^{i}]}(t,l)|
            %$$
           %$$\cdot e^{\varepsilon_{i-1} e^{-\pi|l_{3}-k_{3}|(\lambda_{i-1}-\lambda'_{i})t}(1+|V^{0}_{t,\tau}|)}.$$
         % ( $\mathbf{how}\quad \mathbf{to}\quad\mathbf{get}? $ $If $ $we$ $prove $ $e^{2\pi|k_{12}-l_{12}|(\lambda'_{i}|\frac{\sin\Omega t}{\Omega}|+\mu'_{i})} e^{2\pi|k_{3}-l_{3}|(\lambda'_{i}t+\mu'_{i})}|(\Omega^{i}_{\tau,t})^{-1}-Id|^{\wedge}(t,l-k)\leq 2e^{2\pi|k_{12}-l_{12}|(\lambda'_{i}|\frac{\sin\Omega t}{\Omega}|+\mu'_{i})} e^{2\pi|k_{3}-l_{3}|(\lambda'_{i}t+\mu'_{i})}|\Omega^{i}_{\tau,t}-Id|^{\wedge}(t,l-k),$)

         % $$e^{2\pi|k_{12}|(\lambda_{i}'|\frac{\sin\Omega t}{\Omega}|+\mu_{i}')}e^{2\pi|k_{3}|(\lambda_{i}'t+\mu_{i}')}|[(\Omega^{n+1}_{t,\tau})^{-1}-Id]^{\wedge}(k,v)|\leq 2\kappa_{i},\kappa_{i}=(\kappa_{xi},\kappa_{vi}),$$
         % and $$e^{2\pi|k_{12}-l_{12}|(\lambda'_{i}|\frac{\sin\Omega t}{\Omega}|+\mu'_{i})}\cdot  e^{2\pi|k_{3}-l_{3}|(\lambda'_{i}t+\mu'_{i})}|[e^{2i\pi k\cdot(\Omega^{i}_{\tau,t}-Id)}]^{\wedge}(l-k)|$$
         % $$\leq e^{2\pi|k_{12}-l_{12}|(\lambda'_{i}|\frac{\sin\Omega t}{\Omega}|+\mu'_{i})}\cdot  e^{2\pi|k_{3}-l_{3}|(\lambda'_{i}t+\mu'_{i})}|[e^{2i\pi k\cdot(\Omega^{i}_{\tau,t}-Id)^{-1}(l-k)}]||(\Omega^{i}_{\tau,t})^{-1}-Id)|$$
          We can obtain the following conclusion,
            $$\|\delta V^{n+1}_{t,\tau}\circ (X^{0}_{\tau,t},V^{0}_{\tau,t})\|_{\mathcal{Z}^{(\lambda_{n}'-B_{0})(1+b),\mu_{n}'}
            _{\tau-\frac{bt}{1+b}}}\leq C\sum^{n}_{i=1}\delta_{i}\int^{t}_{\tau}e^{-2\pi(\lambda_{i}-\lambda'_{i})s}ds
            \leq C\sum^{n}_{i=1}\delta_{i} e^{-2\pi(\lambda_{i}-\lambda'_{i})\tau}
             \min\bigg\{\frac{(t-\tau)}{2},\frac{1}{2\pi(\lambda_{i}-\lambda'_{i})}\bigg\},$$

then we have  $$\|\delta X^{n+1}_{t,\tau}\circ (X^{0}_{\tau,t},V^{0}_{\tau,t})\|_{\mathcal{Z}^{(\lambda_{n}'-B_{0})(1+b),\mu_{n}'}_{\tau-\frac{bt}{1+b}}}
\leq C\sum^{n}_{i=1}\delta_{i}\int^{t}_{\tau}e^{-2\pi(\lambda_{i}-\lambda'_{i})s}ds
\leq C\sum^{n}_{i=1}\delta_{i} e^{-2\pi(\lambda_{i}-\lambda'_{i})\tau}\min\bigg\{\frac{(t-\tau)^{2}}{2},\frac{1}{2\pi(\lambda_{i}
             -\lambda'_{i})^{2}}\bigg\}.$$
             We finish the proof of  Proposition 4.3.

            In the following we  estimate $\nabla\Omega_{t,\tau}^{n}-Id.$
In fact, we  write
$(\Omega_{t,\tau}^{n}-Id)(x,v)=(\delta X^{n}_{t,\tau},\delta V^{n}_{t,\tau})\circ (X^{0}_{\tau,t},V^{0}_{\tau,t}),$
we get by differentiation
$\nabla_{x}\Omega_{t,\tau}^{n+1}-(I,0)=\nabla_{x}(\delta X^{n}_{t,\tau}\circ (X^{0}_{\tau,t},V^{0}_{\tau,t}),$
 $\delta V^{n}_{t,\tau}\circ (X^{0}_{\tau,t},V^{0}_{\tau,t})),$
$\nabla_{v}\Omega_{t,\tau}^{n}-(0,I)=(\nabla_{v}+\mathcal{M}(t-\tau)\nabla_{x})(\delta X^{n}_{t,\tau}\circ (X^{0}_{\tau,t},V^{0}_{\tau,t}),
\delta V^{n}_{t,\tau}\circ (X^{0}_{\tau,t},V^{0}_{\tau,t})).$
 \begin{align}
   \left\{\begin{array}{l}
   \frac{d}{dt}\nabla_{x}\delta X^{i}_{t,\tau}(x,v)=\nabla_{x}\delta V^{i}_{t,\tau}(x,v),\\
  \frac{d}{dt}\nabla_{x}\delta V^{i}_{t,\tau}(x,v)=\nabla_{x}\delta V^{i}_{t,\tau}(x,v)\times B_{0}
  +\nabla_{x}E[f^{i}](t,\delta X^{i}_{t,\tau}(x,v) +X^{0}_{t,\tau}(x,v)),\\
 \delta X^{i}_{\tau,\tau}(x,v)=0,\quad\delta V^{i}_{\tau,\tau}(x,v)=0.
  \end{array}\right.
           \end{align}

           Using the same process in the proof of Proposition 4.3, we can obtain Proposition 4.5.

           To establish a control of $\Omega^{i}_{t,\tau}-\Omega^{n}_{t,\tau}$ in norm $\mathcal{Z}^{(\lambda'_{n}-B_{0})(1+b),\mu'_{n}}_{\tau-\frac{bt}{1+b}},$
           we start again from the differential equation satisfied by $\delta V^{i}_{t,\tau}$ and $\delta V^{n}_{t,\tau}:$
          \begin{align}
   \left\{\begin{array}{l}
   \frac{d}{dt}(\delta X^{i}_{t,\tau}-\delta X^{n}_{t,\tau})(x,v)=\delta V^{i}_{t,\tau}(x,v)-\delta V^{n}_{t,\tau}(x,v),\\
  \frac{d}{dt}(\delta V^{i}_{t,\tau}-\delta V^{n}_{t,\tau})(x,v)=(\delta V^{i}_{t,\tau}(x,v)-\delta V^{n}_{t,\tau}(x,v))\times B_{0}
  +E[f^{i-1}](t,\delta X^{i-1}_{t,\tau}(x,v) +X^{0}_{t,\tau}(x,v))\\
  -E[f^{n-1}](t,\delta X^{n-1}_{t,\tau}(x,v) +X^{0}_{t,\tau}(x,v)),\\
 (\delta X^{i}_{t,\tau}-\delta X^{n}_{t,\tau})(x,v)=0,\quad (\delta V^{i}_{t,\tau}-\delta V^{n}_{t,\tau})(x,v)=0.
  \end{array}\right.
           \end{align}

           So from (5.5), $\delta V^{i}_{t,\tau}-\delta V^{n}_{t,\tau}$ satisfies the equation:
           $$\frac{d}{dt}(\delta V^{i}_{t,\tau}-\delta V^{n}_{t,\tau})(x,v)=(\delta V^{i}_{t,\tau}(x,v)-\delta V^{n}_{t,\tau}(x,v))\times B_{0}
  +E[f^{i-1}](t,\delta X^{i-1}_{t,\tau}(x,v) +X^{0}_{t,\tau}(x,v))$$
 $$-E[f^{n-1}](t,\delta X^{i-1}_{t,\tau}(x,v) +X^{0}_{t,\tau}(x,v))+E[f^{n-1}](t,\delta X^{i-1}_{t,\tau}(x,v) +X^{0}_{t,\tau}(x,v))
  -E[f^{n-1}](t,\delta X^{n-1}_{t,\tau}(x,v) +X^{0}_{t,\tau}(x,v)).$$

  Under the assumption $(\mathbf{I}),$  we can use the similar proof of Proposition 4.3 to finish Proposition 4.6.

  Let $\varepsilon$ be the small constant appearing in Lemma 2.7.

   If
  $$3\mathcal{C}_{1}^{i}+\mathcal{C}_{2}^{i}\leq\varepsilon,\quad \textmd{for}\quad \textmd{all}\quad i\geq1,\quad (\mathbf{II})$$
  then $\|\nabla\Omega_{t,\tau}^{i}\|_{\mathcal{Z}^{(\lambda'_{i}-B_{0})(1+b),\mu'_{i}}_{\tau-\frac{bt}{1+b}}}\leq\varepsilon;$ if in addition
  $$2(1+\tau)(1+B)(3\mathcal{C}_{1}^{i,n}+\mathcal{C}_{2}^{i,n})(\tau,t)\leq\max\{\lambda'_{i}-\lambda'_{n},\mu'_{i}-\mu'_{n}\},\quad (\mathbf{III})$$
  for all $i\in\{1,\ldots,n-1\}$ and all $t\geq\tau,$
  then  \begin{align}
   \left\{\begin{array}{l}
   \lambda'_{n}(1+b)+2\|\Omega^{n}-\Omega^{i}\|_{\mathcal{Z}^{(\lambda'_{n}-B_{0})(1+b),\mu'_{n}}_{\tau-\frac{bt}{1+b}}}\leq \lambda'_{i}(1+b),\\
  \mu'_{n}+2(1+|\tau-\frac{bt}{1+b}|)\|\Omega^{n}-\Omega^{i}\|_{\mathcal{Z}^{(\lambda'_{n}-B_{0})(1+b),\mu'_{n}}_{\tau-\frac{bt}{1+b}}}\leq\mu'_{i}.
 \end{array}\right.
           \end{align}

Then Lemma 2.7 and (5.6) yield Proposition 4.8.

As  a corollary of Proposition 4.8 and Proposition 2.8, under the assumption $(\mathbf{IV}):
$
$$4(1+\tau)(\mathcal{C}_{1}^{i,n}+\mathcal{C}_{2}^{i,n})\leq\min\{\lambda_{i}-\lambda'_{n},\mu_{i}-\mu'_{n}\},$$
for all $i\in\{1,\ldots,n\}$ and all $\tau\in[0,t],$
we have
\begin{cor} under the assumption (4.12), we have
$$\|h^{i}_{\tau}\circ\Omega^{n}_{t,\tau}\|
_{\mathcal{Z}_{\tau-\frac{bt}{1+b}}^{(\lambda'_{n}-B_{0})(1+b),\mu'_{n};1}}\leq\delta_{i},\quad
 \sup_{0\leq\tau\leq t}\|(h^{i}_{\tau}v)\circ\Omega^{n}_{t,\tau}\|_{\mathcal{Z}^{(\lambda'_{n}-B_{0})(1+b),
\mu'_{n};1}_{\tau-\frac{bt}{1+b}}}\leq\delta_{i},$$
$$\sup_{0\leq\tau\leq t}\| (\nabla_{x}(h^{i}_{\tau}v))\circ\Omega^{n}_{t,\tau}\|_{\mathcal{Z}^{(\lambda'_{n}-B_{0})(1+b),
\mu'_{n};1}_{\tau-\frac{bt}{1+b}}}\leq\delta_{i},\quad\sup_{0\leq\tau\leq t}\| (\nabla_{x}h^{i}_{\tau})\circ\Omega^{n}_{t,\tau}\|_{\mathcal{Z}^{(\lambda'_{n}-B_{0})(1+b),
\mu'_{n};1}_{\tau-\frac{bt}{1+b}}}\leq\delta_{i},$$
$$\|(\nabla'_{v}+\tau\nabla_{x})(h^{i}_{\tau})\circ\Omega^{n}_{t,\tau}\|
_{\mathcal{Z}_{\tau-\frac{bt}{1+b}}^{(\lambda'_{n}-B_{0})(1+b),\mu'_{n};1}}\leq\delta_{i},\quad\|((\nabla'_{v}+\tau\nabla_{x})h^{i}_{\tau})\circ\Omega^{n}_{t,\tau}\|
_{\mathcal{Z}_{\tau-\frac{bt}{1+b}}^{(\lambda'_{n}-B_{0})(1+b),\mu'_{n};1}}\leq\delta_{i}.$$
\end{cor}

$$$$
\begin{center}
\item\section{The estimates of main terms}
%{\bf\large 1. \quad Introduction }
\end{center}

In order to estimate  these terms $I,II,V,$ we have to make good understanding of plasma echoes. Therefore, we will try to explain  plasma echoes and give
the relevant  mathematical results.
\begin{center}
\item\subsection{Plasma echoes}
%{\bf\large 1. \quad Introduction }
\end{center}

This section is one of  the key sections in our paper. And from Theorem 6.1 in this section, we can see, it is reasonable that
 the influence of the magnetic field is regarded as an error term.  First, we plan to  briefly  explain plasma echoes though  a simple example from [17], then we control plasma echoes
 (or obtain plasma echoes) in time-shift pure and hybrid analytic norms.

The unusual non-linear phenomena that results from  the undamped oscillations of the distribution function $f$ satisfying
the nonlinear Vlasov equation is called plasma echo. Let a perturbation be specified at the initial instant, such that
 the distribution function $\delta f$ is the perturbation of that of Maxwellian plasma $f^{0}(v)\sim \exp(-\alpha v^{2}),\alpha>0$  is a constant and varies periodically in the $x-$direction.
 Without loss of generality, we assume
 $
 \delta f=A_{1}f^{0}(v)\cos k_{1}x \quad \textmd{at}\quad t=0;$
in this section,  $A_{i} $ denotes the amplitude and $k_{i}$ denotes the wave number for $i=1,2.$
The perturbation of the density, i.e. the integral $\int\delta f dv,$ varies in the same manner in the $x-$direction at $t=0.$ Subsequently, the perturbation
of the distribution function varies at time $t$ according to
$\delta f=A_{1}f^{0}(v)\cos k_{1}(x-vt),$
which corresponds to a free movement of each particle in the $x-$direction with its own speed $v.$ But the density perturbation is damped
(in a time $\sim\frac{1}{k_{1}v_{T}}$), because
 $\int\delta f dv$ is made small by the speed-oscillatory factor $\cos k_{1}(x-vt).$ The asymptotic form of the damping at times $t\gg\frac{1}{k_{1}v_{T}}$
 is given by
 \begin{align}
 \delta\rho=\int\delta f dv\varpropto\exp(-\frac{1}{2}k^{2}_{1}v^{2}_{T}t^{2}),
 \end{align}
where the proof of (6.1) can be found in [17].

 Now let the distribution function be again modulates at a time $t=\tau\gg\frac{1}{k_{1}v_{T}},$ with amplitude $A_{2}$ and a new wave number $k_{2}> k_{1}.$
 The resulting density perturbation is damped in time $t\sim\frac{1}{k_{2}v_{T}},$ but reappears at a time $\tau'=\frac{k_{2}\tau}{k_{2}-k_{1}},$ since the second modulation creates
  in the distribution function for $t=\tau$ a second-order term of the form
  $$\delta f^{(2)}=A_{1}A_{2}f^{0}(v)\cos k_{1}(x-v\tau)\cos k_{2}x,$$
  whose further development at $t>\tau$ changes  into
  $$\delta f^{(2)}=A_{1}A_{2}f^{0}(v)\cos k_{1}(x-v\tau)\cos k_{2}[x-v(t-\tau)]$$
  $$=A_{1}A_{2}f^{0}(v)\{\cos [(k_{1}-k_{2})(x-vt)+k_{2}v\tau]+\cos[( k_{1}+k_{2})(x-vt)+k_{2}v\tau]\}.$$
  We see that at $t=\tau'$ the oscillatory dependence of the first term on $ v $ disappears, so that this term makes a finite contribution to the perturbation of
  the density with wave number $k_{2}-k_{1}.$ The resulting echo is then damped in a time $\sim\frac{1}{v_{T}(k_{2}-k_{1})},$ and the final stage of this
  damping follows a law similar to (6.1).

From the above physical point of view, under the assumption of the stability condition, we are discovering that, even in magnetic field case, echoes occurring
 at distinct frequencies are asymptotically well separated. In the following, through  complicate computation, we give a detailed description by
  using mathematical tool. The same to Section 1, since  resonances only  occur in the $\hat{z}$ direction,
  in order to simplify the statement of the proof of the following theorem, we assume $(x,v)=(x_{3},v_{3})\in\mathbb{T}\times\mathbb{R}.$

 \begin{thm} Let $\lambda,\bar{\lambda},\mu,\bar{\mu},\mu',\hat{\mu}$ be such that
$2\lambda\geq\bar{\lambda}>\lambda>0,$ $ \bar{\mu}\geq\mu'>\mu>\hat{\mu}>0,$ and let $b=b(t,s)>0,$
  $ R=R(t,x),G=G(t,x,v)$ and assume $\widehat{G}(t,k_{1},k_{2},0,v)=0,$ we have if
$$\sigma(t,x,v)=\int^{t}_{0}R(s,x+\mathcal{M}(t-s)v)G(s,x+\mathcal{M}(t-s)v,v)ds,$$
$$\sigma_{1}(t,x)=\int^{t}_{0}\int_{\mathbb{R}^{3}}R(s,x+\mathcal{M}(t-s)v)G(s,x+\mathcal{M}(t-s)v,v)dvds.$$
 Then
\begin{align}
 &\|\sigma(t,\cdot)\|_{\mathcal{Z}^{\lambda,\mu;1}_{t}}\leq\int^{t}_{0}\sup_{k\neq l,k,l\in\mathbb{Z}_{\ast}}e^{-2\pi(\bar{\mu}-\mu)|k-l|}
 e^{-2\pi(\bar{\lambda}-\lambda)|k-l|s}\|R\|_{\mathcal{F}^{\bar{\lambda}s+\bar{\mu}}}\|G\|_{\mathcal{Z}^{\lambda(1-b),\hat{\mu};1}_{s}}ds,\\
 & \|\sigma(t,\cdot)\|_{\mathcal{Z}^{\lambda,\mu;1}_{t}}\leq\int^{t}_{0}\sup_{k,l\in\mathbb{Z}_{\ast}}e^{-\pi(\bar{\mu}-\mu)|l|}
e^{-\pi(\bar{\lambda}-\lambda)|k(t-s)+ls|}e^{-2\pi[\mu'-\mu+\lambda b(t-s)]|k-l|}
\cdot \|R\|_{\mathcal{F}^{\lambda s+\mu'-\lambda b(t-s)}}\|G\|_{\mathcal{Z}^{\bar{\lambda}(1+b),\bar{\mu};1}_{s-\frac{bt}{1+b}}}ds,
\end{align}
\begin{align}
& \|\sigma_{1}(t,\cdot)\|_{\mathcal{F}^{\lambda t+\mu}}\leq\int^{t}_{0}\sup_{k,l\in\mathbb{Z}_{\ast}}e^{-\pi(\bar{\mu}-\mu)|l|}
e^{-\pi(\bar{\lambda}-\lambda)|k(t-s)+ls|}e^{-2\pi[\mu'-\mu+\lambda b(t-s)]|k-l|}
\cdot \|R\|_{\mathcal{F}^{\lambda s+\mu'-\lambda b(t-s)}}\|G\|_{\mathcal{Z}^{\bar{\lambda}(1+b),\bar{\mu};1}_{s-\frac{bt}{1+b}}}ds,\\
&\|\sigma_{1}(t,\cdot)\|_{\mathcal{F}^{\lambda t+\mu}}\leq\int^{t}_{0}\sup_{k\neq l,k,l\in\mathbb{Z}_{\ast}} e^{-2\pi(\bar{\lambda}-\lambda)|k-l|s}
\|R\|_{\mathcal{F}^{\bar{\lambda}s+\mu+\lambda b(t-s)}}\|G\|_{\mathcal{Z}^{\lambda(1-b),\mu;1}_{s+\frac{bt}{1-b}}}ds.
\end{align}
\end{thm}
$Proof.$
$$ \|\sigma(t,x,v)\|_{\mathcal{Z}^{\lambda,\mu;1}_{t}}\leq\int^{t}_{0}\|(RG)\circ S^{0}_{s-t}(s,\cdot)\|_{\mathcal{Z}^{\lambda,\mu;1}_{t}}ds=\int^{t}_{0}\|(RG)(s,\cdot)\|_{\mathcal{Z}^{\lambda,\mu;1}_{s}}ds.$$

Let $s'=s+b(t-s), b=\frac{Ds}{t(1+t)},$ where  some constant $D>0$ small enough.  Note that
$$\|(RG)(s,\cdot)\|_{\mathcal{Z}^{\lambda,\mu;1}_{s}}=\sum_{k\in\mathbb{Z}}\sum_{n\in\mathbb{N}_{0}}\frac{\lambda^{n}}{n!}e^{2\pi\mu|k|}
\bigg\|\bigg[\nabla_{v}+2i\pi ks\bigg]^{n}\widehat{(RG)}(s,k,v)\bigg\|_{L_{dv}^{1}}$$
$$=\sum_{k\in\mathbb{Z}}\sum_{n\in\mathbb{N}_{0}}\frac{\lambda^{n}}{n!}e^{2\pi\mu|k|}
\bigg\|e^{2i\pi v\cdot k(t-s)}\bigg[\nabla_{v}+2i\pi ks\bigg]^{n}\widehat{(RG)}(s,k,v)\bigg\|_{L_{dv}^{1}}$$
$$=\sum_{k\in\mathbb{Z}}\sum_{n\in\mathbb{N}_{0}}\frac{\lambda^{n}}{n!}e^{2\pi\mu|k|}
\bigg\|\bigg[2i\pi k(t-s)+2i\pi ks\bigg]^{n}e^{2i\pi v\cdot k(t-s)}\widehat{(RG)}(s,k,v)\bigg\|_{L_{dv}^{1}}$$
$$=\sum_{k\in\mathbb{Z}}\sum_{n\in\mathbb{N}_{0}}\frac{\lambda^{n}}{n!}e^{2\pi\mu|k|}
\bigg\|\bigg[2i\pi k(t-s')+2i\pi k(s'-s)+2i\pi (k-l(1-b))s+2i\pi ls(1-b)\bigg]^{n}$$
$$\cdot e^{2i\pi v\cdot k(t-s)}\sum_{l}\hat{R}(s,k-l)\hat{G}(s,l,v)\bigg\|_{L_{dv}^{1}}$$
$$\leq\sum_{k\in\mathbb{Z}}\sum_{n\in\mathbb{N}_{0}}\frac{\lambda^{n}}{n!}e^{2\pi\mu|k|}\sum_{l}|\hat{R}(s,k-l)|
\sum^{n}_{\gamma=0}C^{\gamma}_{n}\bigg|\bigg[2i\pi k(s'-s)+2i\pi (k-l(1-b))s\bigg]^{\gamma}\bigg|$$
$$\bigg\|\bigg[2i\pi k(t-s')+2i\pi ls(1-b)\bigg]^{n-\gamma}\cdot e^{2i\pi v\cdot k(t-s)}\hat{G}(s,l,v)\bigg\|_{L_{dv}^{1}}$$
$$\leq\sum_{k\in\mathbb{Z}}\sum_{n\in\mathbb{N}_{0}}e^{2\pi\mu|k|}\sum_{l}|\hat{R}(s,k-l)|
\sum^{n}_{\gamma=0}\frac{\lambda^{\gamma}}{\gamma!}\bigg|\bigg[2i\pi k(s'-s)+2i\pi (k-l(1-b))s\bigg]^{\gamma}\bigg|$$
$$\frac{\lambda^{n-\gamma}}{(n-\gamma)!}\bigg\|\bigg[2i\pi k(t-s')+2i\pi ls(1-b)\bigg]^{n-\gamma}\cdot e^{2i\pi v\cdot k(t-s)}\hat{G}(s,l,v)\bigg\|_{L_{dv}^{1}}$$
$$\leq\sum_{k\in\mathbb{Z}}\sum_{n\in\mathbb{N}_{0}}e^{2\pi\mu|k|}\sum_{l}|\hat{R}(s,k-l)|
\sum^{n}_{\gamma=0}\frac{\lambda^{\gamma}}{\gamma!}\bigg|\bigg[2i\pi k(s'-s)+2i\pi (k-l(1-b))s\bigg]^{\gamma}\bigg|$$
$$\frac{\lambda^{n-\gamma}(1-b)^{n-\gamma}}{(n-\gamma)!}\bigg\|\bigg[2i\pi k(t-s)+2i\pi ls\bigg]^{n-\gamma}\cdot e^{2i\pi v\cdot k(t-s)}\hat{G}(s,l,v)\bigg\|_{L_{dv}^{1}}$$
$$\leq\sum_{k\in\mathbb{Z}}\sum_{n\in\mathbb{N}_{0}}e^{2\pi\mu|k|}\sum_{l}|\hat{R}(s,k-l)|
\sum^{n}_{\gamma=0}\frac{\lambda^{\gamma}}{\gamma!}\bigg|\bigg[2i\pi k(s'-s)+2i\pi (k-l(1-b))s\bigg]^{\gamma}\bigg|$$
$$\cdot\frac{\lambda^{n-\gamma}(1-b)^{n-\gamma}}{(n-\gamma)!}\|[\nabla_{v}+2i\pi ls]^{n-\gamma}\hat{G}(s,l,v)\|_{L_{dv}^{1}}.$$
Now  we will divide $k,l$ into the following cases:

Case 1. $\min\{|k|,|l|\}>k-l>0.$ We still decompose this case into two steps.

Step 1. If $\min\{|k|,|l|\}>k-l>0,$ $k<0,$
$$\sum^{n}_{\gamma=0}\frac{\lambda^{\gamma}}{\gamma!}\bigg|\bigg[2i\pi k(s'-s)+2i\pi (k-l(1-b))s\bigg]^{\gamma}\bigg|$$
$$=\sum^{n}_{\gamma=0}\frac{\lambda^{\gamma}}{\gamma!}\bigg|\bigg[2\pi (k-l)s+2i\pi kb(t-s)+2\pi lbs\bigg]^{\gamma}\bigg|
=\sum^{n}_{\gamma=0}\frac{\lambda^{\gamma}}{\gamma!}\bigg|\bigg[2\pi (k-l)s(1-b)+2\pi kbt\bigg]^{\gamma}\bigg|.$$
If $s\geq\frac{-kDs}{(1+t)(k-l)(1-b)},$ we have
$\sum^{n}_{\gamma=0}\frac{\lambda^{\gamma}}{\gamma!}\bigg|\bigg[2i\pi k(s'-s)+2i\pi (k-l(1-b))s\bigg]^{\gamma}\bigg|
=\sum^{n}_{\gamma=0}\frac{\lambda^{\gamma}}{\gamma!}\bigg[2\pi (k-l)s(1-b)+2\pi kbt\bigg]^{\gamma}.$

Then
 $$\sum_{n\in\mathbb{N}_{0}}e^{2\pi\mu|k|}\sum_{l}|\hat{R}(s,k-l)|
\sum^{n}_{\gamma=0}\frac{\lambda^{\gamma}}{\gamma!}\bigg|\bigg[2i\pi k(s'-s)+2i\pi (k-l(1-b))s\bigg]^{\gamma}\bigg|
\cdot\frac{\lambda^{n-\gamma}(1-b)^{n-\gamma}}{(n-\gamma)!}\|[\nabla_{v}+2i\pi ls]^{n-\gamma}\hat{G}(s,l,v)\|_{L_{dv}^{1}}$$
$$\leq\sum_{l}\sum_{n\in\mathbb{N}_{0}}e^{2\pi\mu|k|}e^{2\pi (k-l)s(1-b)+2\pi kbt}|\hat{R}(s,k-l)|
\cdot\frac{\lambda^{n}(1-b)^{n}}{n!}\|[\nabla_{v}+2i\pi ls]^{n}\hat{G}(s,l,v)\|_{L_{dv}^{1}}$$
$$\leq\sum_{l}\sum_{n\in\mathbb{N}_{0}}e^{2\pi(\mu+\lambda bt)|k-l|}e^{2\pi \lambda|k-l|s(1-b)}|\hat{R}(s,k-l)|\cdot e^{2\pi(\mu-\lambda bt)|l|}
\frac{\lambda^{n}(1-b)^{n}}{n!}\|[\nabla_{v}+2i\pi ls]^{n}\hat{G}(s,l,v)\|_{L_{dv}^{1}}.$$

%$$\|(RG)(s,\cdot)\|_{\mathcal{Z}^{\lambda,\mu;1}_{s}}\leq \sup_{k,l\in\mathbb{Z}^{3},k\neq l} e^{2\pi\lambda b|k(t-s)+ls|}e^{-2\pi(\bar{\lambda}-\lambda)|k-l|s}\|R(s,\cdot)\|_{\mathcal{F}^{\bar{\lambda}s+\mu}}
%\|G(s,\cdot)\|_{\mathcal{Z}^{\lambda(1-b),\mu;1}_{s}},$$

If $s\leq\frac{-kDs}{(1+t)(k-l)(1-b)}\leq\frac{-kD}{(k-l)(1-b)}\leq t,$ for some constant $0<\epsilon_{0}<\frac{b}{1-b},$
$$\sum^{n}_{\gamma=0}\frac{\lambda^{\gamma}}{\gamma!}\bigg|\bigg[2i\pi k(s'-s)+2i\pi (k-l(1-b))s\bigg]^{\gamma}\bigg|
=\sum^{n}_{\gamma=0}\frac{\lambda^{\gamma}}{\gamma!}\bigg[-2\pi (k-l)s(1-b)-2\pi kbt\bigg]^{\gamma}$$
$$\leq\sum^{n}_{\gamma=0}\frac{\lambda^{\gamma}}{\gamma!}\bigg[-(2+\epsilon_{0})\pi kbt+\pi\epsilon_{0} kbt\bigg]^{\gamma}
=\sum^{n}_{\gamma=0}\frac{\lambda^{\gamma}}{\gamma!}\bigg[-(2+\epsilon_{0})\pi k\frac{Ds}{t(1+t)}t+\pi\epsilon_{0} kbt\bigg]^{\gamma}$$
$$=\sum^{n}_{\gamma=0}\frac{\lambda^{\gamma}}{\gamma!}\bigg[(2+\epsilon_{0})\pi |k-l|\frac{|k|Ds}{|k-l|(1+t)}-\pi\epsilon_{0} |l|\frac{|k|Ds}{|l|(1+t)}\bigg]^{\gamma}$$
$$\leq\sum^{n}_{\gamma=0}\frac{\lambda^{\gamma}}{\gamma!}\bigg[(2+\epsilon_{0})\pi |k-l|\frac{|k|Bs}{|k-l|}\frac{|k-l|(1-b)}{|k|D}-\pi \epsilon_{0}|l|\frac{|k|Ds}{|l|(1+t)}\bigg]^{\gamma}$$
$$\leq\sum^{n}_{\gamma=0}\frac{\lambda^{\gamma}}{\gamma!}\bigg[(2+\epsilon_{0})\pi |k-l|(1-b)s-\pi\epsilon_{0} |l|\frac{Bs}{(1+t)}\bigg]^{\gamma},$$
so we get
$$\sum_{n\in\mathbb{N}_{0}}e^{2\pi\mu|k|}\sum_{l}|\hat{R}(s,k-l)|
\sum^{n}_{\gamma=0}\frac{\lambda^{\gamma}}{\gamma!}\bigg|\bigg[2i\pi k(s'-s)+2i\pi (k-l(1-b))s\bigg]^{\gamma}\bigg|
\cdot\frac{\lambda^{n-\gamma}(1-b)^{n-\gamma}}{(n-\gamma)!}\|[\nabla_{v}+2i\pi ls]^{n-\gamma}\hat{G}(s,l,v)\|_{L_{dv}^{1}}$$
$$\leq\sum_{l}\sum_{n\in\mathbb{N}_{0}}e^{2\pi\mu|k-l|}e^{2\pi \lambda|k-l|s(1-\frac{b}{2})}|\hat{R}(s,k-l)|
\cdot e^{2\pi(\mu-\frac{\lambda\epsilon_{0} Ds}{2(1+t)})|l|}\frac{\lambda^{n}(1-b)^{n}}{n!}\|[\nabla_{v}+2i\pi ls]^{n}\hat{G}(s,l,v)\|_{L_{dv}^{1}}.$$

Step 2. $\min\{|k|,|l|\}>k-l>0,k>0.$
Indeed, we only consider $k-l>0,l>0,$  since $\min\{|k|,|l|\}>k-l>-\min\{|k|,|l|\}.$

 It is easy to check that
$$\sum^{n}_{\gamma=0}\frac{\lambda^{\gamma}}{\gamma!}\bigg|\bigg[2i\pi k(s'-s)+2i\pi (k-l(1-b))s\bigg]^{\gamma}\bigg|
=\sum^{n}_{\gamma=0}\frac{\lambda^{\gamma}}{\gamma!}[2\pi (k-l)s'-2\pi lbt]^{\gamma}.$$
This can be reduced to case 1, here we omit the details.
 We can obtain
$$\sum_{n\in\mathbb{N}_{0}}e^{2\pi\mu|k|}\sum_{l}|\hat{R}(s,k-l)|
\sum^{n}_{\gamma=0}\frac{\lambda^{\gamma}}{\gamma!}\bigg|\bigg[2i\pi k(s'-s)+2i\pi (k-l(1-b))s\bigg]^{\gamma}\bigg|
\cdot\frac{\lambda^{n-\gamma}(1-b)^{n-\gamma}}{(n-\gamma)!}\|[\nabla_{v}+2i\pi ls]^{n-\gamma}\hat{G}(s,l,v)\|_{L_{dv}^{1}}$$
$$\leq\sum_{l}\sum_{n\in\mathbb{N}_{0}}e^{2\pi\mu|k-l|}e^{2\pi \lambda|k-l|s'}|\hat{R}(s,k-l)|\cdot e^{2\pi\hat{\mu}|l|}\frac{\lambda^{n}(1-b)^{n}}{n!}\|[\nabla_{v}+2i\pi ls]^{n}\hat{G}(s,l,v)\|_{L_{dv}^{1}}.$$

Case 2. If  $-\min\{|k|,|l|\}<k-l<0,$ the method of this case is the same to  Case 1.
$$\sum_{n\in\mathbb{N}_{0}}e^{2\pi\mu|k|}\sum_{l}|\hat{R}(s,k-l)|
\sum^{n}_{\gamma=0}\frac{\lambda^{\gamma}}{\gamma!}\bigg|\bigg[2i\pi k(s'-s)+2i\pi (k-l(1-b))s\bigg]^{\gamma}\bigg|
\cdot\frac{\lambda^{n-\gamma}(1-b)^{n-\gamma}}{(n-\gamma)!}\|[\nabla_{v}+2i\pi ls]^{n-\gamma}\hat{G}(s,l,v)\|_{L_{dv}^{1}}$$
$$\leq\sum_{l}\sum_{n\in\mathbb{N}_{0}}e^{2\pi\mu|k-l|}e^{2\pi \lambda|k-l|s'}|\hat{R}(s,k-l)|
\cdot e^{2\pi\hat{\mu|}l|}\frac{\lambda^{n}(1-b)^{n}}{n!}\|[\nabla_{v}+2i\pi ls]^{n}\hat{G}(s,l,v)\|_{L_{dv}^{1}}.$$

Case 3.  When $k-l>\min\{|k|,|l|\},$ or  $k-l<-\min\{|k|,|l|\},$  we only need to consider one of two cases. Without loss of generality, we  assume
 $k-l>\min\{|k|,|l|\}.$
 $$\sum^{n}_{\gamma=0}\frac{\lambda^{\gamma}}{\gamma!}\bigg|\bigg[2i\pi k(s'-s)+2i\pi (k-l(1-b))s\bigg]^{\gamma}\bigg|=\sum^{n}_{\gamma=0}\frac{\lambda^{\gamma}}{\gamma!}\bigg|2\pi (k-l)s(1-b)+2\pi kbt\bigg|^{\gamma}$$
$$\leq\sum^{n}_{\gamma=0}\frac{\lambda^{\gamma}}{\gamma!}\bigg[\bigg|2\pi (k-l)s(1-b)\bigg|+\bigg|2\pi kbt\bigg|\bigg]^{\gamma}
\leq\sum^{n}_{\gamma=0}\frac{\lambda^{\gamma}}{\gamma!}\bigg[\bigg|2\pi (k-l)s(1-b)\bigg|+\bigg|2\pi (k-l)bt\bigg|\bigg]^{\gamma}$$
$$\leq\sum^{n}_{\gamma=0}\frac{\lambda^{\gamma}}{\gamma!}\bigg[2\pi(s(1-b)+bt)| k-l|\bigg]^{\gamma}
=\sum^{n}_{\gamma=0}\frac{\lambda^{\gamma}}{\gamma!}\bigg[2\pi\bigg(1-b+\frac{D}{1+t}\bigg)| k-l|s\bigg]^{\gamma}.$$
%$$\leq\sum^{n}_{\gamma=0}\frac{\lambda^{\gamma}}{\gamma!}[-2\pi (k-l+l)\frac{Bs}{(1+t)}]^{\gamma}.$$

$$\sum_{n\in\mathbb{N}_{0}}e^{2\pi\mu|k|}\sum_{l}|\hat{R}(s,k-l)|
\sum^{n}_{\gamma=0}\frac{\lambda^{\gamma}}{\gamma!}\bigg|\bigg[2i\pi k(s'-s)+2i\pi (k-l(1-b))s\bigg]^{\gamma}\bigg|
\cdot\frac{\lambda^{n-\gamma}(1-b)^{n-\gamma}}{(n-\gamma)!}\|[\nabla_{v}+2i\pi ls]^{n-\gamma}\hat{G}(s,l,v)\|_{L_{dv}^{1}}$$
$$\leq\sum_{l}\sum_{n\in\mathbb{N}_{0}}e^{2\pi\mu|k-l|}e^{2\pi \lambda(1-b+\frac{2D}{1+t})|k-l|s}|\hat{R}(s,k-l)|
\cdot e^{2\pi\hat{\mu}|l|}\frac{\lambda^{n}(1-b)^{n}}{n!}\|[\nabla_{v}+2i\pi ls]^{n}\hat{G}(s,l,v)\|_{L_{dv}^{1}}.$$

Case 4. If $t<\frac{|k|Ds}{(1+t)|k-l|(1-b)}\leq\frac{|k|D}{|k-l|(1-b)}$ with $\min\{|k|,|l|\}>k-l>-\min\{|k|,|l|\},k\neq l.$ Without loss of generality,
 we assume that $\min\{|k|,|l|\}>k-l>0,$
%Note that
%\begin{align}
%|k(t-s)+ls|\sim
%\left\{
%\begin{array}
%c|k-l|t,\quad if \quad |k|\gg|l|\quad or\quad |l|\gg|k| \\
%|l|t,\quad \quad if\quad\quad \quad |k|\sim|l| ,\\
%\end{array}
%\right.
%\end{align}
$$\sum^{n}_{\gamma=0}\frac{\lambda^{\gamma}}{\gamma!}\bigg|\bigg[2i\pi k(s'-s)+2i\pi (k-l(1-b))s\bigg]^{\gamma}\bigg|
=\sum^{n}_{\gamma=0}\frac{\lambda^{\gamma}}{\gamma!}|[2\pi (k-l)s'-2\pi lbt]^{\gamma}|.$$
If $k-l>0,$ $l>0,$

$$\sum^{n}_{\gamma=0}\frac{\lambda^{\gamma}}{\gamma!}\bigg|\bigg[2\pi (k-l)s'-2\pi lbt\bigg]^{\gamma}\bigg|
\leq\sum^{n}_{\gamma=0}\frac{\lambda^{\gamma}}{\gamma!}\bigg|\bigg[2\pi (k-l)\cdot\frac{|k|Ds}{(1+t)|k-l|(1-b)}-2\pi lbt\bigg]^{\gamma}\bigg|$$
$$\leq\sum^{n}_{\gamma=0}\frac{\lambda^{\gamma}}{\gamma!}\bigg|\bigg[2\pi \cdot\frac{|k|Ds}{(1+t)(1-b)}-2\pi l\frac{Ds}{1+t}\bigg]^{\gamma}\bigg|
\leq\sum^{n}_{\gamma=0}\frac{\lambda^{\gamma}}{\gamma!}\bigg[2\pi \frac{D|k-l|}{(1+t)(1-b)}s-2\pi\epsilon_{0} l\frac{Ds}{1+t}\bigg]^{\gamma}$$
$$\leq\sum^{n}_{\gamma=0}\frac{\lambda^{\gamma}}{\gamma!}\bigg[2\pi |k-l|s-2\pi\epsilon_{0} |l|\frac{Ds}{1+t}\bigg]^{\gamma}.$$
If $k-l>0,l<0,$ this is equal to $k-l>0,k<0,$ since $\min\{|k|,|l|\}>k-l>-\min\{|k|,|l|\}.$
$$\sum^{n}_{\gamma=0}\frac{\lambda^{\gamma}}{\gamma!}\bigg|\bigg[2i\pi k(s'-s)+2i\pi (k-l(1-b))s\bigg]^{\gamma}\bigg|$$
$$=\sum^{n}_{\gamma=0}\frac{\lambda^{\gamma}}{\gamma!}\bigg[2\pi (k-l)s(1-b)+2\pi kbt\bigg]^{\gamma}
\leq\sum^{n}_{\gamma=0}\frac{\lambda^{\gamma}}{\gamma!}\bigg[2\pi |k-l|s-2\pi\epsilon_{0} |l|\frac{Ds}{1+t}\bigg]^{\gamma}.$$
In summary,
$$\sum_{k,l\in\mathbb{Z}_{\ast}}\sum_{n\in\mathbb{N}_{0}}e^{2\pi\mu|k|}|\hat{R}(s,k-l)|
\sum^{n}_{\gamma=0}\frac{\lambda^{\gamma}}{\gamma!}\bigg|\bigg[2i\pi k(s'-s)+2i\pi (k-l(1-b))s\bigg]^{\gamma}\bigg|
\cdot\frac{\lambda^{n-\gamma}(1-b)^{n-\gamma}}{(n-\gamma)!}\|[\nabla_{v}+2i\pi ls]^{n-\gamma}\hat{G}(s,l,v)\|_{L_{dv}^{1}}$$
$$\leq \sup_{k\neq l,k,l\in\mathbb{Z}_{\ast}}e^{-2\pi(\bar{\mu}-\mu)|k-l|} e^{-2\pi(\bar{\lambda}-\lambda)|k-l|s}\|R\|_{\mathcal{F}^{\bar{\lambda}s+\bar{\mu}}}\|G\|_{\mathcal{Z}^{\lambda(1-b),\hat{\mu};1}_{s}},$$
then
$$ \|\sigma(t,\cdot)\|_{\mathcal{Z}^{\lambda,\mu;1}_{t}}\leq\int^{t}_{0}\sup_{k\neq l,k,l\in\mathbb{Z}_{\ast}^{3}}e^{-2\pi(\bar{\mu}-\mu)|k-l|} e^{-2\pi(\bar{\lambda}-\lambda)|k-l|s}\|R\|_{\mathcal{F}^{\bar{\lambda}s+\bar{\mu}}}\|G\|_{\mathcal{Z}^{\lambda(1-b),\hat{\mu};1}_{s}}ds.$$
Now we estimate the second inequality of Theorem 6.1.
$$ \|\sigma(t,x,v)\|_{\mathcal{Z}^{\lambda,\mu;1}_{t}}\leq\int^{t}_{0}\|(RG)\circ S^{0}_{s-t}(s,\cdot)\|_{\mathcal{Z}^{\lambda,\mu;1}_{t}}ds=\int^{t}_{0}\|(RG)(s,\cdot)\|_{\mathcal{Z}^{\lambda,\mu;1}_{s}}ds.$$
Note that
$$\|(RG)(s,\cdot)\|_{\mathcal{Z}^{\lambda,\mu;1}_{s}}=\sum_{k\in\mathbb{Z}^{3}}\sum_{n\in\mathbb{N}_{0}^{3}}\frac{\lambda^{n}}{n!}e^{2\pi\mu|k|}
\|[\nabla_{v}+2i\pi ks]^{n}\widehat{(RG)}(s,k,v)\|_{L_{dv}^{1}}$$
$$=\sum_{k\in\mathbb{Z}}\sum_{n\in\mathbb{N}_{0}}\frac{\lambda^{n}}{n!}e^{2\pi\mu|k|}
\|e^{2i\pi v\cdot k(t-s)}[\nabla_{v}+2i\pi ks]^{n}\widehat{(RG)}(s,k,v)\|_{L_{dv}^{1}}$$
$$=\sum_{k\in\mathbb{Z}}\sum_{n\in\mathbb{N}_{0}}\frac{\lambda^{n}}{n!}e^{2\pi\mu|k|}
\|[2i\pi k(t-s)+2i\pi ks]^{n}e^{2i\pi v\cdot k(t-s)}\widehat{(RG)}(s,k,v)\|_{L_{dv}^{1}}$$
$$\leq\sum_{k,l\in\mathbb{Z}}e^{2\pi\lambda|k|t}e^{2\pi\mu|k|}|\hat{R}(s,k-l)|
\|e^{2i\pi v\cdot k(t-s)}\hat{G}(s,l,v)\|_{L_{dv}^{1}}.$$
The proof of the inequality (6.3) is as follows, here we sketch the main steps. We let $s'=s-b(t-s)$ and write
$$e^{2\pi(\lambda t+\mu)|k|}\leq e^{-2\pi(\bar{\mu}-\mu)|l|}e^{-2\pi\lambda(s-s')|k-l|}e^{-2\pi(\mu'-\mu)|k-l|}
e^{-2\pi(\bar{\lambda}-\lambda)|k(t-s')+ls'|}$$
$$\times e^{2\pi\bar{\mu}|l|}e^{2\pi(\lambda s+\mu')|k-l|}e^{2\pi\bar{\lambda}|k(t-s')+ls'|}$$
$$\leq e^{-2\pi(\bar{\mu}-\mu)|l|}e^{-\pi(\bar{\lambda}-\lambda)|k(t-s)+ls|}e^{-2\pi[\mu'-\mu+\lambda(s-s')/2]|k-l|}$$
$$\times e^{2\pi\bar{\mu}|l|}e^{2\pi[\lambda s+\mu'-\lambda(s-s')/2]|k-l|}\sum_{n\in\mathbb{N}_{0}}
\frac{|2i\pi\bar{\lambda}(k(t-s')+ls')|^{n}}{n!},$$
then we can deduce that
$$ \|\sigma(t,\cdot)\|_{\mathcal{Z}^{\lambda,\mu;1}_{t}}\leq\int^{t}_{0}\sup_{k,l\in\mathbb{Z}_{\ast}}e^{-\pi(\bar{\mu}-\mu)|l|}
e^{-\pi(\bar{\lambda}-\lambda)|k(t-s)+ls|}e^{-2\pi[\mu'-\mu+\lambda b(t-s)]|k-l|}\cdot \|R\|_{\mathcal{F}^{\lambda s+\mu'-\lambda b(t-s)}}\|G\|_{\mathcal{Z}^{\bar{\lambda}(1+b),\bar{\mu};1}_{s-\frac{bt}{1+b}}}ds.$$
In the following we estimate the norm $\mathcal{F}^{\lambda t+\mu}$ of the function $\sigma_{1}(t,x),$
$$\sigma_{1}(t,x)=\int^{t}_{0}\int_{\mathbb{R}}R(s,x+\mathcal{M}(t-s)v)G(s,x+\mathcal{M}(t-s)v,v)dvds.$$

We write
$$\hat{\sigma}_{1}(t,k)=\int^{t}_{0}\sum_{l\in\mathbb{Z}}\int_{\mathbb{R}}\hat{R}(s,k-l)\hat{G}(s,l,v)e^{2\pi i v\cdot k(t-s)}dvds,$$
$$|\hat{\sigma}_{1}(t,k)|\leq\int^{t}_{0}\bigg(\sum_{l\in\mathbb{Z}}\bigg|\int_{\mathbb{R}}\hat{G}(s,l,v)e^{2\pi i v\cdot k(t-s)}dv\bigg||\hat{R}(s,k-l)|\bigg)ds.$$
Next, $$e^{2\pi(\lambda t+\mu)|k|}\leq e^{2\pi\lambda( s+b(t-s))|k-l|}e^{2\pi\lambda(1-b)|k(t-s)+l(s+\frac{bt}{1-b})|}e^{2\pi\mu|k-l|}e^{2\pi\mu|l|}$$
$$\leq e^{2\pi\bar{\lambda}|k-l|s}e^{2\pi\lambda(1-b)|k(t-s)+l(s+\frac{bt}{1-b})|}e^{2\pi(\mu+\lambda b(t-s))|k-l|}e^{2\pi\mu|l|}e^{-2\pi(\bar{\lambda}-\lambda)|k-l|s}.$$
Hence
$$\|\sigma_{1}(t,\cdot)\|_{\mathcal{F}^{\lambda t+\mu}}\leq\int^{t}_{0}\sup_{k\neq l,k,l\in\mathbb{Z}_{\ast}^{3}} e^{-2\pi(\bar{\lambda}-\lambda)|k-l|s}
\|R\|_{\mathcal{F}^{\bar{\lambda}s+\bar{\mu}}}\|G\|_{\mathcal{Z}^{\lambda(1-b),\mu;1}_{s+\frac{bt}{1-b}}}ds,$$
where $\bar{\mu}=\mu+\lambda b(t-s).$

\begin{center}
\item\subsection{ Estimates of main terms}
%{\bf\large 1. \quad Introduction }
\end{center}

% $ \mathbf{ Note}$ that,
  % \begin{align}
  % &\int_{\mathbb{T}^{3}}\int_{\mathbb{R}^{3}}-(\mathcal{E}^{n+1}_{t,s}\cdot G_{t,s}^{n})(s, X^{0}_{t,s}(x,v),V^{0}_{t,s}(x,v))dxdv\neq0\notag\\
   %&\int_{\mathbb{T}^{3}}\int_{\mathbb{R}^{3}}-(F^{n+1}_{t,s}\cdot G_{t,s}^{n,v})(s, X^{0}_{t,s}(x,v),V^{0}_{t,s}(x,v)))dxdv\neq0\notag\\
   %&\int_{\mathbb{T}^{3}}\int_{\mathbb{R}^{3}}-(\mathcal{E}^{n}_{t,s}\cdot H^{n}_{t,s})(s, X^{0}_{t,s}(x,v),V^{0}_{t,s}(x,v)))dxdv\neq0\notag\\
 % &\int_{\mathbb{T}^{3}}\int_{\mathbb{R}^{3}}-(F^{n}_{t,s}\cdot H_{t,s}^{n,v})(s, X^{0}_{t,s}(x,v),V^{0}_{t,s}(x,v)))dxdv\neq0\notag\\
  % \end{align}

   %$\mathbf{Problem:}$ how to add condition such that $\hat{\rho}^{n+1}(t,k_{\perp},0)=0$?

%Since $\hat{W}(k_{1},k_{2},0)=0,$ it is trivial to check that $\hat{\rho}^{n+1}(t,k_{\perp},0)=0.$

In the following we estimate $\bar{II}_{i}^{n+1,n}(t,x).$ Note that their zero modes vanish. For any $n\geq i\geq1,$
$$\widehat{\bar{II}_{i}^{n+1,n}}(t,k)
=\int^{t}_{0}\int_{\mathbb{T}^{3}}\int_{\mathbb{R}^{3}}e^{-2\pi ik\cdot x} \bigg(B[h^{n+1}]
\cdot (\nabla'_{v}\times((h^{i}_{\tau}v)\circ\Omega^{i-1}_{t,\tau})\bigg)(\tau,x'(\tau,x,v),v'(\tau,x,v))dv'dxd\tau,$$
\begin{align}
&|\widehat{\bar{II}_{i}^{n+1,n}}(t,k)|\leq\int^{t}_{0}\bigg(\sum_{l\in\mathbb{Z}^{3}_{\ast}}\bigg|\int_{\mathbb{R}^{3}}e^{-2\pi ik_{3}\cdot(v_{3}(t-\tau))}e^{-2\pi i\eta_{k1}\cdot v'_{1}}e^{-2\pi i\eta_{k2}\cdot v'_{2}}\notag\\
&\cdot (\widehat{\nabla'_{v}\times((h^{i}_{\tau}v)\circ\Omega^{i-1}_{t,\tau}))}(\tau,l,v')dv'\bigg||\widehat{B[h^{n+1}]}(\tau,k-l)|\bigg)d\tau\notag\\
&=\int^{t}_{0}\bigg(\sum_{l\in\mathbb{Z}^{3}_{\ast}}\bigg|\int_{\mathbb{R}}e^{-2\pi ik_{3}\cdot(v_{3}(t-\tau))}e^{-2\pi i\eta_{k1}\cdot v'_{1}}e^{-2\pi i\eta_{k2}\cdot v'_{2}}\notag\\
&\cdot (\widehat{\nabla'_{v}\times((h^{i}_{\tau}v)\circ\Omega^{i-1}_{t,\tau}))}(\tau,l,\eta_{k1},\eta_{k2},v_{3})dv_{3}\bigg||\widehat{B[h^{n+1}]}(\tau,k-l)|
\bigg)d\tau.\notag\\
\end{align}

From (6.3) of Theorem 6.1 and (6.6), we can get Proposition 4.9.

To finish Propositions 4.11-4.12,
 we shall again use the Vlasov equation. We rewrite it as
$$h^{n+1}(t,X^{n}_{\tau,t}(x,v),V^{n}_{\tau,t}(x,v))=\int^{t}_{0}\Sigma^{n+1}(s,X^{n}_{\tau,s}(x,v),V^{n}_{\tau,s}(x,v))ds.$$
Then we get
$$\|h^{n+1}(t,X^{n}_{\tau,t}(x,v),V^{n}_{\tau,t}(x,v))\|_{\mathcal{Z}^{(\lambda'_{n}-B_{0})(1-b),\mu'_{n};1}_{t+\frac{bt}{1-b}}}$$
$$\leq\int^{t}_{0}\|\Sigma^{n+1}(s,X^{n}_{\tau,s}(x,v),V^{n}_{\tau,s}(x,v))\|_{\mathcal{Z}^{(\lambda'_{n}-B_{0})(1-b),\mu'_{n};1}
_{t+\frac{bt}{1-b}}}(\eta_{k_{1}},\eta_{k_{2}})ds$$
$$=\int^{t}_{0}\|\Sigma^{n+1}(s,\Omega^{n}_{\tau,s}(x,v))\|_{\mathcal{Z}^{(\lambda'_{n}-B_{0})(1-b),\mu'_{n};1}_{s+\frac{bt}{1-b}}}ds$$
$$\leq\int^{t}_{0}\sup_{k,l\in\mathbb{Z}_{\ast}^{3}}e^{-\pi(\mu''_{n}-\mu'_{n})|l|}
e^{-\pi(\lambda_{n}-\lambda'_{n})|k(t-s)+l(s+\frac{bt}{1-b})|}e^{-2\pi[\mu''_{n}-\mu'_{n}+(\lambda'_{n}-B_{0}) b(t-s)]|k-l|}$$
$$\cdot\|\mathcal{E}^{n+1}_{s,t}\|_{\mathcal{F}
^{\varsigma'_{n}}}\|G^{n}_{s,t}\|_{\mathcal{Z}^{(\lambda_{n}-B_{0})(1+b),\mu''_{n};1}_{s+\frac{bt}{1-b}-\frac{bt}{1+b}}}ds+\int^{t}_{0}
\sup_{k,l\in\mathbb{Z}_{\ast}^{3}}e^{-\pi(\mu''_{n}-\mu'_{n})|l|}
e^{-\pi(\lambda''_{n}-\lambda'_{n})|k(t-s)+l(s+\frac{bt}{1-b})|}$$
$$\cdot e^{-2\pi[\mu''_{n}-\mu'_{n}+(\lambda'_{n}-B_{0}) b(t-s)]|k-l|}\|F^{n+1}_{s,t}\|_{\mathcal{F}
^{\varsigma'_{n}}}\|G^{n,v}_{s,t}\|_{\mathcal{Z}^{(\lambda''_{n}-B_{0})(1+b),\mu''_{n};1}_{s+\frac{bt}{1-b}-\frac{bt}{1+b}}}ds$$
$$+\int^{t}_{0}
\bigg[\|H^{n}_{s,t}\|_{\mathcal{Z}^{(\lambda'_{n}-B_{0})(1-b),\mu'_{n};1}_{s+\frac{bt}{1-b}}}
\|\mathcal{E}^{n}_{s,t}\|_{\mathcal{Z}^{(\lambda'_{n}-B_{0})(1-b),\mu'_{n}}_{s+\frac{bt}{1-b}}}
+\|H^{n,v}_{s,t}\|_{\mathcal{Z}^{(\lambda'_{n}-B_{0})(1-b),\mu'_{n};1}_{s+\frac{bt}{1-b}}}$$
$$\cdot\|F^{n}_{s,t}\|_{\mathcal{Z}^{(\lambda'_{n}-B_{0})(1-b),\mu'_{n}}
_{s+\frac{bt}{1-b}}}\bigg]ds+\int^{t}_{0} \sup_{k\neq l,k,l\in\mathbb{Z}_{\ast}^{3}}e^{-2\pi(\mu_{n}-\mu'_{n})|k-l|}
e^{-2\pi(\lambda_{n}-\lambda'_{n})|k-l|s}$$
$$\cdot\|B[f^{n}]\circ\Omega_{t,s}^{n}\|_{\mathcal{F}^{\varsigma''_{n}}}
\|H^{n+1,v}_{s,t}\|_{\mathcal{Z}^{(\lambda'_{n}-B_{0})(1-b)^{2},\hat{\mu}'_{n};1}_{s+\frac{bt}{1-b}}}ds,$$
where $\varsigma'_{n}=(\lambda'_{n}-B_{0})(1-b)s+\mu''_{n}-(\lambda'_{n}-B_{0})(1-b)b(t-s),$ $\varsigma''_{n}=\lambda''_{n}s+\mu'+(\lambda'_{n}-B_{0})b(t-s).$

Then  $$\sup_{0<s\leq t}\|h^{n+1}(t,X^{n}_{\tau,t}(x,v),V^{n}_{\tau,t}(x,v))\|
_{\mathcal{Z}^{(\lambda'_{n}-B_{0})(1-b),\mu'_{n};1}_{s+\frac{bt}{1-b}}}$$
$$\leq\int^{t}_{0}\sup_{k,l\in\mathbb{Z}_{\ast}^{3}}e^{-\pi(\mu_{n}-\mu'_{n})|l|}
e^{-\pi(\lambda_{n}-\lambda'_{n})|k(t-s)+ls|}e^{-2\pi[\mu''_{n}-\mu'_{n}+(\lambda'-B_{0}) b(t-s)]|k-l|}$$
$$\cdot e^{-2\pi(\lambda_{n}-\lambda'_{n})|k-l|s}\bigg(C^{'0}+\sum^{n}_{i=1}\delta_{i}\bigg)\|\rho[h^{n+1}]\|_{\mathcal{F}
^{(\lambda'_{n}-B_{0})s+\mu'_{n}}}ds+\frac{\delta^{2}_{n}}{(\lambda_{n}-\lambda'_{n})^{2}}$$
$$\leq C\bigg(C^{'0}+\sum^{n}_{i=1}\delta_{i}\bigg)\sup_{0<s\leq t}\|\rho[h^{n+1}]\|_{\mathcal{F}
^{(\lambda'_{n}-B_{0})s+\mu'_{n}}}+\frac{\delta^{2}_{n}}{(\lambda_{n}-\lambda'_{n})^{2}}.$$

Therefore, we obtain
 \begin{align}
\sup_{0\leq s\leq t}\|
h^{n+1}\circ\Omega^{n}_{t,\tau}\|_{\mathcal{Z}^{(\lambda'_{n}-B_{0})(1-\frac{1}{2}b),
\mu'_{n};1}_{s+\frac{bt}{1-b}}}\leq\delta^{2}_{n}+\bigg(\sum^{n}_{i=1}\delta_{i}\bigg)
\sup_{0\leq s\leq t}\|\rho^{n+1}\|_{\mathcal{F}^{(\lambda'_{n}-B_{0})s+\mu'_{n}}},
\end{align}
this is the conclusion of Proposition 4.11.

 Finally,  we estimate the last term
$$V(t,x)=-\int^{t}_{0}\int_{\mathbb{R}^{3}}\bigg[\bigg(B[f^{n}]\circ\Omega^{n}_{s,t}(x,v)\bigg)\cdot H_{s,t}^{n+1,v}\bigg](s, X^{0}_{s,t}(x,v),V^{0}_{s,t}(x,v))dvds$$
$$=-\int^{t}_{0}\int_{\mathbb{R}^{3}}\bigg[\bigg(B[f^{n}]\circ\Omega^{n}_{s,t}(x,v)\bigg)\cdot \bigg((\nabla_{v}h^{n+1}\times v)\circ\Omega^{n}_{s,t}(x,v)\bigg)\bigg]
(s, X^{0}_{s,t}(x,v),V^{0}_{s,t}(x,v))dvds$$
$$=-\int^{t}_{0}\int_{\mathbb{R}^{3}}\bigg(B[f^{n}]\circ\Omega^{n}_{s,t}(x,v)\bigg)\cdot \bigg[\bigg((\nabla_{v}\times(h^{n+1}v))\circ\Omega^{n}_{s,t}(x,v)\bigg)$$
$$-\bigg(\nabla_{v}\times\bigg[(h^{n+1}v)\circ\Omega^{n}_{s,t}(x,v)\bigg]\bigg](s, X^{0}_{s,t}(x,v),V^{0}_{s,t}(x,v))dvds$$
%$$+\int^{t}_{0}\int_{\mathbb{R}^{3}}\bigg[(B[f^{n}]\circ\Omega^{n}_{s,t}(x,v))\cdot \bigg]\bigg)\bigg](s, X^{0}_{s,t}(x,v),V^{0}_{s,t}(x,v))dvds$$
$$-\int^{t}_{0}\int_{\mathbb{R}^{3}}\bigg[(B[f^{n}]\circ\Omega^{n}_{s,t}(x,v))\cdot
 \bigg(\nabla_{v}\times\bigg[(h^{n+1}v)\circ\Omega^{n}_{s,t}(x,v)\bigg]\bigg)\bigg](s, X^{0}_{s,t}(x,v),V^{0}_{s,t}(x,v))dvds.$$

 We claim that for $\varepsilon>0$ sufficiently small, $$\|(\nabla\Omega_{t,s}^{n})^{-1}-Id\|
 _{\mathcal{Z}_{s+\frac{bt}{1-b}}^{(\lambda_{n}-B_{0})(1-b),\mu'_{n}}}<\varepsilon.$$

 If the claim holds, then
 $$\|V\|_{\mathcal{F}^{(\lambda'_{n}-B_{0})t+\mu'_{n}}}$$
$$\leq\int^{t}_{0}\sup_{k\neq l,k,l\in\mathbb{Z}_{\ast}^{3}}e^{-2\pi(\lambda_{n}-\lambda''_{n})|k-l|s}
\|B[f^{n}]\|_{\mathcal{F}^{\nu'_{n}}}\cdot \bigg\|\bigg[\bigg((\nabla_{v}\times(h^{n+1}v))\circ\Omega^{n}_{s,t}(x,v)\bigg)$$
$$-\bigg((\nabla_{v}\times\bigg[(h^{n+1}v)\circ\Omega^{n}_{s,t}(x,v)\bigg]\bigg)\bigg]\bigg\|
_{\mathcal{Z}_{s+\frac{bt}{1-b}}^{(\lambda'_{n}-B_{0})(1-b),\mu'_{n};1}}
(\eta_{k_{1}},\eta_{k_{2}})ds+\int^{t}_{0}\sup_{k\neq l,k,l\in\mathbb{Z}_{\ast}^{3}}e^{-2\pi(\lambda_{n}-\lambda''_{n})|k-l|s}$$
%$$+\int^{t}_{0}\int_{\mathbb{R}^{3}}\bigg[(B[f^{n}]\circ\Omega^{n}_{s,t}(x,v))\cdot \bigg]\bigg)\bigg](s, X^{0}_{s,t}(x,v),V^{0}_{s,t}(x,v))dvds$$
$$\cdot\|B[f^{n}]\|_{\mathcal{F}^{\nu'_{n}}}\cdot \bigg\|
\nabla_{v}\times\bigg[(h^{n+1}v)\circ\Omega^{n}_{s,t}(x,v)\bigg]\bigg\|_{\mathcal{Z}_{s+\frac{bt}{1-b}}^{(\lambda'_{n}-B_{0})(1-b),\mu'_{n};1}}
(\eta_{k_{1}},\eta_{k_{2}})ds$$
$$\leq\int^{t}_{\tau}e^{-\pi(\lambda_{n}-\lambda'_{n})|k(t-s)+ls|}
\|B[f^{n}]\|
_{\mathcal{F}^{\nu'_{n}}}\| h^{n+1}\circ\Omega^{n}_{s,t}(x,v)\|_{\mathcal{Z}_{s+\frac{bt}{1-b}}^{(\lambda'_{n}-B_{0})(1-\frac{1}{2}b),\mu'_{n};1}}
(\eta_{k_{1}},\eta_{k_{2}})ds,$$
%$$\bigg(\nabla_{v}h^{n+1}\times v)\circ\Omega^{n}_{s,t}(x,v)\bigg)\bigg](s, X^{0}_{s,t}(x,v),V^{0}_{s,t}(x,v))$$
%$$+\int^{t}_{\tau}e^{-\pi(\bar{\lambda}_{n}-\lambda^{\ast}_{n})|k(t-s)+ls|}
%\|B[f^{n}]\circ\Omega^{n}_{s,t}(x,v)\|
%_{\mathcal{F}^{\nu^{\ast}_{n}}}\|(\nabla\Omega_{t,s}^{n})^{-1}-Id\|_{\mathcal{Z}_{s+\frac{bt}{1+b}}^{\bar{\lambda}_{n}(1-b),\mu^{\ast}_{n}}}$$
%$$\cdot\|((\nabla_{v}+(t-s)\nabla_{x})\cdot( h^{n+1}\circ\Omega^{n}_{s,t}(x,v)))\times (v\circ\Omega^{n}_{s,t}(x,v))\|_{\mathcal{Z}_{s+\frac{bt}{1+b}}^{\bar{\lambda}_{n}(1-b),\mu^{\ast}_{n};1}}ds$$
where $$\nu'_{n}=(\lambda_{n}''-B_{0})s+\mu'_{n}+(\lambda'_{n}-B_{0})b(t-s)+
\|\Omega_{t,s}^{n}-Id\|_{\mathcal{Z}_{s}^{\lambda''_{n}-B_{0},\mu'_{n}+(\lambda'_{n}-B_{0})b(t-s)}}$$
$$\leq(\lambda''_{n}-B_{0})s+\mu'_{n}+(\lambda'_{n}-B_{0})b(t-s)+
\|\Omega_{t,s}^{n}-Id\|_{\mathcal{Z}_{s-\frac{bt}{1+b}}^{(\lambda''_{n}-B_{0}),\mu'_{n}+2(\lambda'_{n}-B_{0})b(t-s)}}$$
$$\leq(\lambda''_{n}-B_{0})s+\mu'_{n}+(\lambda'_{n}-B_{0})b(t-s)+
\|\Omega_{t,s}^{n}-Id\|_{\mathcal{Z}_{s-\frac{bt}{1+b}}^{\lambda''_{n}-B_{0},\mu''_{n}}},$$
by $Proposition $ 4.3, we have
$\nu'_{n}\leq(\lambda_{n}-B_{0})s+\mu_{n}-(\lambda_{n}-\lambda''_{n})s$ as soon as
$$
4C^{1}_{\omega}\sum^{n}_{i=1}\frac{\delta_{i}}{(2\pi(\lambda_{i}-\lambda'_{i}))^{3}
       }\leq\frac{\lambda'_{\infty}D}{3}\quad (\mathbf{V}).$$

So $$\|V\|_{\mathcal{F}^{(\lambda'_{n}-B_{0})t+\mu'_{n}}}\leq\int^{t}_{\tau}e^{-\pi(\lambda_{n}-\lambda'_{n})|k(t-s)+ls|}\bigg(\sum^{n}_{i=1}\delta_{i}\bigg)
\| h^{n+1}\circ\Omega^{n}_{s,t}(x,v)\|_{\mathcal{Z}_{s+\frac{bt}{1-b}}^{(\lambda'_{n}-B_{0})(1-\frac{1}{2}b),\mu'_{n};1}}ds.$$
%$$\|V\|_{\mathcal{F}^{\lambda_{n+1} t+\mu_{n+1}}}$$
 %$$\leq C\int^{t}_{0} K^{n+1}_{1}(t,\tau)\|\nabla'_{v}\times
We finish the proof of Proposition 4.12.

$$$$

\begin{center}
\item\section{ Estimates of error terms}
%{\bf\large 1. \quad Introduction }
\end{center}

In the following we estimate one of the error terms $\mathcal{R}_{0}.$

 Recall
$$\mathcal{R}_{0}(t,x)=\int^{t}_{0}\int_{\mathbb{R}^{3}}\bigg(\bigg(B[h^{n+1}]\circ\Omega^{n}_{t,s}(x,v)-B[h^{n+1}]\bigg)\cdot G_{t,s}^{n,v}\bigg)(s, X^{0}_{t,s}(x,v),V^{0}_{t,s}(x,v)))dvds.$$
First,
$$\|\mathcal{R}_{0}(t,\cdot)\|_{\mathcal{F}^{(\lambda'_{n}-B_{0}) t+\mu'_{n}}}$$
$$\leq\int^{t}_{0}\|B[h^{n+1}]\circ\Omega^{n}_{t,s}(x,v)-B[h^{n+1}]\|_{\mathcal{Z}_{\tau-\frac{bt}{1+b}}^{(\lambda'_{n}-B_{0})(1+b),\mu'_{n}}}
\|G_{t,s}^{n,v}\|_{\mathcal{Z}_{\tau-\frac{bt}{1+b}}^{(\lambda'_{n}-B_{0})(1+b),\mu'_{n};1}}ds.$$
Next,
\begin{align}
&\|(B[h^{n+1}]\circ\Omega^{n}_{t,s}(x,v)-B[h^{n+1}]\|_{\mathcal{Z}_{\tau-\frac{bt}{1+b}}^{(\lambda'_{n}-B_{0})(1+b),\mu'_{n}}}\notag\\
&\leq\|\Omega^{n}_{t,\tau}-Id\|_{\mathcal{Z}_{\tau-\frac{bt}{1+b}}^{(\lambda'_{n}-B_{0})(1+b),\mu'_{n}}}
\int^{1}_{0}\|\nabla B[h^{n+1}]((1-\theta)Id+\theta\Omega^{n}_{t,s})\|_{\mathcal{Z}_{\tau-\frac{bt}{1+b}}^{(\lambda'_{n}-B_{0})(1+b),\mu'_{n}}}d\theta\notag\\
&\leq\|\nabla B[h^{n+1}]\|_{\mathcal{F}^{\nu'_{n}}}
\|\Omega^{n}_{t,\tau}-Id\|_{\mathcal{Z}_{\tau-\frac{bt}{1+b}}^{(\lambda'_{n}-B_{0})(1+b),\mu'_{n}}},\notag\\
\end{align}
where
$\nu'_{n}=(\lambda'_{n}-B_{0})(1+b)\bigg|\tau-\frac{bt}{1+b}\bigg|+\mu'_{n}+\|\Omega^{n}X_{t,\tau}-x'\|_{\mathcal{Z}_{\tau-\frac{bt}{1+b}}
^{(\lambda'_{n}-B_{0})(1+b),\mu'_{n}}}.$

Here we only focus on the case $\tau\geq\frac{bt}{1+b},$ then we need to show
$\|\nabla B[h^{n+1}]\|_{\mathcal{F}^{\nu'_{n}}}\leq\|\rho[h^{n+1}]\|_{\mathcal{F}^{(\lambda'_{n}-B_{0})\tau+\mu'_{n}}}.$
For that, we have to prove
$\nu'_{n}<(\lambda'_{n}-B_{0})\tau+\mu'_{n}-\iota,$
for some constant $\iota>0$ sufficiently small.

Indeed,$$\nu'_{n}\leq(\lambda'_{n}-B_{0})\tau+\mu'_{n}-(\lambda'_{n}-B_{0})b(t-\tau)+C\sum^{n}_{i=1}
\delta_{i}e^{-\pi|k_{3}|(\lambda_{i}-(\lambda'_{i}))t}\cdot\min\bigg\{\frac{(t-\tau)^{2}}{2},\frac{1}{2\pi(\lambda_{i}-\lambda'_{i})^{2}}\bigg\}$$
$$\leq(\lambda'_{n}-B_{0})\tau+\mu'_{n}-(\lambda'_{n}-B_{0})\frac{B(t-\tau)}{1+t}+C\bigg(\sum^{n}_{i=1}\frac{\delta_{i}}{(\lambda_{i}
-\lambda'_{i})^{3}}\bigg)
\frac{\min\{t-\tau,1\}}{1+\tau}.$$
Note that
$\frac{\min\{t-\tau,1\}}{1+\tau}\leq3\frac{t-\tau}{1+t}.$
In the following we also need to show that
$$
C\sum^{n}_{i=1}\frac{\delta_{i}}{(\lambda_{i}-\lambda'_{i})^{3}}\leq\frac{\lambda^{\ast}B}{3}-\iota, \quad (\mathbf{VI})
$$

Now we assume that (7.1) holds, then
$$\|(B[h^{n+1}]\circ\Omega^{n}_{t,s}(x,v)-B[h^{n+1}]\|_{\mathcal{Z}_{\tau-\frac{bt}{1+b}}^{(\lambda'_{n}-B_{0})(1+b),\mu'_{n}}}
\leq C\bigg(\sum^{n}_{i=1}\frac{\delta_{i}}{(\lambda_{i}-\lambda'_{i})^{5}}\bigg)
\frac{1}{(1+\tau)^{3}}\|\rho[h^{n+1}]\|_{\mathcal{F}^{(\lambda'_{n}-B_{0})\tau+\mu'_{n}}}.$$

 Since $G_{t,s}^{n,v}=(\nabla'_{v}f^{n}\times V^{0}_{t,s}(x,v))\circ\Omega^{n}_{t,s}(x,v),$
 $$\|G_{t,s}^{n,v}\|_{\mathcal{Z}_{\tau-\frac{bt}{1+b}}^{(\lambda'_{n}-B_{0})(1+b),\mu'_{n};1}}
 \leq\|(\nabla'_{v}f^{0}\times V^{0}_{t,s}(x,v))\circ\Omega^{n}_{t,s}(x,v)\|
 _{\mathcal{Z}_{\tau-\frac{bt}{1+b}}^{(\lambda'_{n}-B_{0})(1+b),\mu'_{n};1}}$$
$$+\sum^{n}_{i=1}\|(\nabla'_{v}h_{\tau}^{i}\times V^{0}_{t,s}(x,v))\circ\Omega^{n}_{t,s}(x,v)\|_{\mathcal{Z}_{\tau-\frac{bt}{1+b}}^{(\lambda'_{n}-B_{0})(1+b),\mu'_{n};1}}$$
$$\leq C'_{0}+\bigg(\sum^{n}_{i=1}\delta_{i}\bigg)(1+\tau).$$

  We can conclude
$$\|\mathcal{R}_{0}(t,\cdot)\|_{\mathcal{F}^{(\lambda'_{n}-B_{0}) t+\mu'_{n}}}\leq C\bigg(C'_{0}+\sum^{n}_{i=1}\delta_{i}\bigg)\bigg(\sum^{n}_{i=1}\frac{\delta_{i}}{(\lambda_{i}-\lambda'_{i})^{5}}\bigg)
\int^{t}_{0}\|\rho[h^{n+1}]\|_{\mathcal{F}^{(\lambda'_{n}-B_{0})\tau+\mu'_{n}}}\frac{d\tau}{(1+\tau)^{2}}.$$

In order to finish the control of $\tilde{\mathcal{R}}_{0},$
 we still need the estimate of $\|G_{t,s}^{n,v}-\bar{G}_{t,s}^{n,v}\|_{\mathcal{Z}_{\tau-\frac{bt}{1+b}}^{(\lambda'_{n}-B_{0})
 (1+b),\mu'_{n};1}}.$

  In fact,
$$\|G_{t,s}^{n,v}-\bar{G}_{t,s}^{n,v}\|_{\mathcal{Z}_{\tau-\frac{bt}{1+b}}^{(\lambda'_{n}-B_{0})(1+b),\mu'_{n};1}}$$
$$\leq
\|(\nabla_{v}\times(f^{0}v))\circ\Omega^{n}_{t,\tau}-\nabla_{v}\times(f^{0}v)\|_{\mathcal{Z}_{\tau-\frac{bt}{1+b}}^{(\lambda'_{n}-B_{0})(1+b),\mu'_{n};1}}
$$
$$+\sum^{n}_{i=1}\|(\nabla_{v}\times(h^{i}v))\circ\Omega^{n}_{t,\tau}-(\nabla_{v}\times(h^{i}v))\circ\Omega^{i-1}_{t,\tau}\|
_{\mathcal{Z}_{\tau-\frac{bt}{1+b}}^{(\lambda'_{n}-B_{0})(1+b),\mu'_{n};1}}$$
$$+\sum^{n}_{i=1}\|(\nabla_{v}\times(h^{i}v))\circ\Omega^{i-1}_{t,\tau}-\nabla_{v}\times((h^{i}v)\circ\Omega^{i-1}_{t,\tau})\|
_{\mathcal{Z}_{\tau-\frac{bt}{1+b}}^{(\lambda'_{n}-B_{0})(1+b),\mu'_{n};1}}.$$

  Now on the one hand, we treat the second term
 $$\sum^{n}_{i=1}\|(\nabla_{v}\times(h^{i}v))\circ\Omega^{n}_{t,\tau}-(\nabla_{v}\times(h^{i}v))\circ\Omega^{i-1}_{t,\tau}\|_{\mathcal{Z}_{\tau-\frac{bt}{1+b}}
 ^{(\lambda'_{n}-B_{0})(1+b),\mu'_{n};1}}$$
 $$\leq\int^{1}_{0}\|\nabla\nabla_{v}h^{i}_{\tau}((1-\theta)\Omega^{n}_{t,\tau}+\theta\Omega^{i-1}_{t,\tau})\|_{\mathcal{Z}
 _{\tau-\frac{bt}{1+b}}^{(\lambda'_{n}-B_{0})(1+b),\mu'_{n};1}}
 \cdot\|\Omega^{n}_{t,\tau}-\Omega^{i-1}_{t,\tau}\|
 _{\mathcal{Z}_{\tau-\frac{bt}{1+b}}^{(\lambda'_{n}-B_{0})(1+b),\mu'_{n}}}d\theta$$
 $$\leq2\|\nabla\nabla_{v}h^{i}_{\tau}\circ\Omega^{i-1}_{t,\tau}\|_{\mathcal{Z}
 _{\tau-\frac{bt}{1+b}}^{(\lambda'_{i}-B_{0})(1+b),\mu'_{i};1}}\|\Omega^{n}_{t,\tau}-\Omega^{i-1}_{t,\tau}\|
 _{\mathcal{Z}_{\tau-\frac{bt}{1+b}}^{(\lambda'_{n}-B_{0})(1+b),\mu'_{n}}}$$
 $$\leq4C\delta_{i}\bigg(\sum^{n}_{j=i}\frac{\delta_{j}}{(\lambda_{j}-\lambda'_{j})^{6}}\bigg)\frac{1}{(1+\tau)^{2}},$$
 where from Proposition 4.6, we know
 $$\|\Omega^{n}X_{t,\tau}-\Omega^{i-1}X_{t,\tau}\|
 _{\mathcal{Z}_{\tau-\frac{bt}{1+b}}^{(\lambda'_{n}-B_{0})(1+b),\mu'_{n}}}\leq2\mathcal{R}^{i-1,n}_{2}(t,\tau),$$
 $$\|\Omega^{n}V_{t,\tau}-\Omega^{i-1}V_{t,\tau}\|
 _{\mathcal{Z}_{\tau-\frac{bt}{1+b}}^{(\lambda'_{n}-B_{0})(1+b),\mu'_{n}}}
 \leq\mathcal{R}^{i-1,n}_{1}(t,\tau)+\mathcal{R}^{i-1,n}_{2}(t,\tau)$$
 with $$\mathcal{R}^{i-1,n}_{1}(t,\tau)=\bigg(\sum^{n}_{j=i}\frac{\delta_{j}e^{-2\pi(\lambda_{j}-\lambda'_{j})\tau}}{2\pi(\lambda_{j}-\lambda'_{j})}\bigg)
 \min\{t-\tau,1\},\mathcal{R}^{i-1,n}_{2}(t,\tau)=\bigg(\sum^{n}_{j=i}\frac{\delta_{j}e^{-2\pi(\lambda_{j}-\lambda'_{j})\tau}}
 {(2\pi(\lambda_{j}-\lambda'_{j}))^{2}}\bigg)
 \min\bigg\{\frac{(t-\tau)^{2}}{2},1\bigg\}.$$
On the other hand, by the  induction hypothesis, since  $\mathcal{Z}^{\lambda,\mu}_{\tau}$ norms are increasing as a function of $\lambda$ and $\mu,$ if fixed $\tau,$
$$\sum^{n}_{i=1}\|(\nabla_{v}h^{i}_{\tau})\circ\Omega^{i-1}_{t,\tau}-\nabla_{v}(h^{i}_{\tau}\circ\Omega^{i-1}_{t,\tau})\|
_{\mathcal{Z}_{\tau-\frac{bt}{1+b}}^{(\lambda'_{n}-B_{0})(1+b),\mu'_{n}}}
\leq\bigg(\sum^{n}_{i=1}\delta_{i}\bigg)\frac{1}{(1+\tau)^{2}}.$$
%The proof of $\|\Omega^{n}_{t,\tau}-\Omega^{i-1}_{t,\tau}\|
% _{\mathcal{Z}^{\lambda_{i}(1+b),
%\mu_{i}}_{\tau-\frac{bt}{1+b}}}$ is similar to the proof of $Proposition$ 3.2, so here we omit the details  of the proof.
%$$\|\Omega^{n}_{t,\tau}-\Omega^{i-1}_{t,\tau}\|
% _{\mathcal{Z}^{\lambda_{i}(1+b),
%\mu_{i}}_{\tau-\frac{bt}{1+b}}}$$
%$$\leq\|\Omega^{i-1}_{t,\tau}\| _{\mathcal{Z}^{\lambda_{i}(1+b),
%\mu_{i}}_{\tau-\frac{bt}{1+b}}}\|(\Omega^{i-1}_{t,\tau})^{-1}\Omega^{n}_{t,\tau}-Id\|
% _{\mathcal{Z}^{\lambda_{i}(1+b),
%\mu_{i}}_{\tau-\frac{bt}{1+b}}}$$
%$$\leq\sum^{n}_{j=i}\frac{\delta_{j}}{(2\pi(\lambda_{j}-\lambda'_{j})^{6})}\frac{1}{(1+\tau)^{5}}$$
%$\mathbf{Assumption:}$ $$\|\nabla\nabla_{v}h^{i}_{\tau}(\Omega^{i-1}_{t,\tau})\|_{\mathcal{Z}^{\lambda_{i}(1+b),
%\mu_{i};1}_{\tau-\frac{bt}{1+b}}}\leq (1+\tau)^{2}\delta_{i}.$$
%So we have
%$$\|\nabla\nabla_{v}h^{i}_{\tau}((1-\theta)\Omega^{n}_{t,\tau}+\theta\Omega^{i-1}_{t,\tau})\|_{\mathcal{Z}^{\lambda_{i}(1+b),
%\mu_{i};1}_{\tau-\frac{bt}{1+b}}}\leq2\|\nabla\nabla_{v}h^{i}_{\tau}(\Omega^{i-1}_{t,\tau})\|_{\mathcal{Z}^{\lambda_{i}(1+b),
%\mu_{i};1}_{\tau-\frac{bt}{1+b}}}\leq (1+\tau)^{2}\delta_{i}.$$

So we have
\begin{align}
&\|\tilde{\mathcal{R}}_{0}(t,\cdot)\|_{\mathcal{F}^{(\lambda'_{n}-B_{0}) \tau+\mu'_{n}}}
\leq\bigg(C^{4}_{\omega}\bigg(C'_{0}+\sum^{n}_{i=1}\delta_{i}\bigg)\bigg(\sum^{n}_{j=1}\frac{\delta_{j}}{2\pi(\lambda_{j}-\lambda'_{j})^{6}}\bigg)
+\sum^{n}_{i=1}\delta_{i}\bigg)
\int^{t}_{0}\|\rho\|_{\mathcal{F}^{(\lambda'_{n}-B_{0}) \tau+\mu'_{n}}}\frac{1}{(1+\tau)^{2}}d\tau\notag\\
&=\int^{t}_{0}\tilde{K}_{1}^{n+1}\|\rho\|_{\mathcal{F}^{(\lambda'_{n}-B_{0}) \tau+\mu'_{n}}}\frac{1}{(1+\tau)^{2}}d\tau.\notag\\
\end{align}
Up to now, we finish the estimates of error terms.
%In the following we continue to estimate $\tilde{\mathcal{R}}_{0}$ to finish the estimate of $II^{n+1,n}.$ Now we recall
%$$\tilde{\mathcal{R}}_{0}=\int^{t}_{0}\int_{\mathbb{R}^{3}}(B[h^{n+1}]\cdot (G_{t,s}^{n,v}-\bar{G}_{t,s}^{n,v}))(s, X^{0}_{t,s}(x,v),V^{0}_{t,s}(x,v))dvds,$$
$$$$

\begin{center}
\item\section{ Iteration}
%{\bf\large 1. \quad Introduction }
\end{center}

Now let us first deal with the source term
\begin{align}
 &\hat{III}^{n,n}(t,k)+\hat{IV}^{n,n}(t,k)=-\int^{t}_{0}\int_{\mathbb{T}^{3}}\int_{\mathbb{R}^{3}}e^{-2\pi ik\cdot x}(\mathcal{E}^{n}_{t,s}\cdot H^{n}_{t,s})(s, X^{0}_{t,s}(x,v),V^{0}_{t,s}(x,v))\notag\\
           &dvdxds-\int^{t}_{0}\int_{\mathbb{T}^{3}}\int_{\mathbb{R}^{3}}e^{-2\pi ik\cdot x}(F^{n}_{t,s}\cdot H^{n,v}_{t,s})(s, X^{0}_{t,s}(x,v),V^{0}_{t,s}(x,v))dvdxds,\notag\\
          % &=-\int^{t}_{0}\int_{\mathbb{T}^{3}}\int_{\mathbb{R}^{3}}\exp(-2\pi ik\cdot(x-v\circ\mathcal{I}_{0}))(\mathcal{E}^{n}_{t,s}\cdot H^{n}_{t,s})(s, x,V^{0}_{t,s}(x,v))dvdxds\notag\\
            % &-\int^{t}_{0}\int_{\mathbb{T}^{3}}\int_{\mathbb{R}^{3}}e^{-2\pi ik\cdot (x-v\circ\mathcal{I}_{0})}(F^{n}_{t,s}\cdot H^{n,v}_{t,s})(s, x,V^{0}_{t,s}(x,v))dvdxds\notag\\
 \end{align}
 %where $X^{0}_{t,s}(x,v)=x+v\circ\mathcal{I}_{0},V^{0}_{t,s}(x,v)=v\circ\mathcal{I}_{1}$
 %$$(v\circ\mathcal{I}_{0})(t,\tau)=(-\frac{v_{\perp}}{\Omega}[\sin(\theta+\Omega(t-\tau))-\sin\theta],\frac{v_{\perp}}{\Omega}[\cos(\theta+\Omega(t-\tau))-\cos\theta],-v_{z}(t-\tau))$$
%$$(v\circ\mathcal{I}_{1})(t,\tau)=(v_{\perp}\cos(\theta+\Omega(t-\tau)),v_{\perp}\sin(\theta+\Omega(t-\tau)),v_{z}),$$
then
%\begin{align}
 %&\hat{III}^{n,n}(t,k)+\hat{IV}^{n,n}(t,k)\notag\\
% &=-\int^{t}_{0}\sum_{l}\int_{\mathbb{R}^{3}}\exp(-2\pi ik\cdot(-v\circ\mathcal{I}_{0}))\mathcal{E}^{n}_{t,s}(s,k-l) H^{n}_{t,s}(s, l,V^{0}_{t,s}(x,v))dvds\notag\\
             %&-\int^{t}_{0}\sum_{l}\int_{\mathbb{R}^{3}}e^{-2\pi ik\cdot (-v\circ\mathcal{I}_{0})}(F^{n}_{t,s}(s,k-l) H^{n,v}_{t,s}(s, l,V^{0}_{t,s}(x,v))dvds\notag\\
             %&=-\int^{t}_{0}\sum_{l}\int_{\mathbb{R}^{3}}\exp(-2\pi ik\cdot(-v\circ\mathcal{I}_{0}))\mathcal{E}^{n}_{t,s}(s,k-l) H^{n}_{t,s}(s, l,V^{0}_{t,s}(x,v))dvds\notag\\
 %\end{align}

 $$\|III(t,\cdot)\|_{\mathcal{F}^{(\lambda'_{n}-B_{0}) t+\mu'_{n}}}+\|IV(t,\cdot)\|_{\mathcal{F}^{(\lambda'_{n}-B_{0}) t+\mu'_{n}}}$$
 $$\leq\int^{t}_{0}\|\mathcal{E}^{n}_{\tau,s}\|_{\mathcal{Z}^{(\lambda'_{n}-B_{0})(1+b),
\mu'_{n}}_{\tau-\frac{bt}{1+b}}}\|H^{n}_{\tau,s}\|_{\mathcal{Z}^{(\lambda'_{n}-B_{0})(1+b),
\mu'_{n};1}_{\tau-\frac{bt}{1+b}}}+\|F^{n}_{\tau,s}\|_{\mathcal{Z}^{(\lambda'_{n}-B_{0})(1+b),
\mu'_{n}}_{\tau-\frac{bt}{1+b}}}$$
$$\cdot\|H^{n,v}_{\tau,s}\|_{\mathcal{Z}^{(\lambda'_{n}-B_{0})(1+b),
\mu'_{n};1}_{\tau-\frac{bt}{1+b}}}d\tau$$
\begin{align}
&\leq\int^{t}_{0}\|\rho^{n}_{\tau,s}\|_{\mathcal{F}^{\nu_{n+1}}}\|H^{n}_{\tau,s}\|_{\mathcal{Z}^{(\lambda'_{n}-B_{0})(1+b),
\mu'_{n};1}_{\tau-\frac{bt}{1+b}}}d\tau+
\int^{t}_{0}\|\rho^{n}_{\tau,s}\|_{\mathcal{F}^{\nu'_{n}}}
\cdot\|H^{n,v}_{\tau,s}\|_{\mathcal{Z}^{(\lambda'_{n}-B_{0})(1+b),
\mu'_{n};1}_{\tau-\frac{bt}{1+b}}}d\tau\notag\\
&\leq\int^{t}_{0}\|\rho^{n}_{\tau,s}\|_{\mathcal{F}^{\nu'_{n}}}(1+\tau)\delta_{n}d\tau\leq\int^{t}_{0}\|\rho^{n}_{\tau,s}\|
_{\mathcal{F}^{(\lambda'_{n}-B_{0}) \tau+\mu'_{n}}}e^{-2\pi\tau(\lambda_{n} -\lambda'_{n} )}(1+\tau)\delta_{n}d\tau
\leq\frac{C\delta^{2}_{n}}{(\lambda_{n} -\lambda'_{n})^{2}}.\notag\\
\end{align}

From Propositions 4.9-4.12, combining (4.10), we conclude
$$\|\rho[h^{n+1}](t,\cdot)\|_{\mathcal{F}^{(\lambda'_{n}-B_{0}) t+\mu'_{n}}}$$
$$\leq\frac{C\delta^{2}_{n}}{(\lambda_{n} -\lambda'_{n})^{2}}+\int^{t}_{0}|K^{n}_{1}(t,\tau)|(1+\tau)\sum^{n}_{i=1}\delta_{i}\| \rho[h^{n+1}]\|
_{\mathcal{F}^{(\lambda'_{n}-B_{0})\tau+\mu'_{n}}}d\tau+2\int^{t}_{0}|K^{n}_{0}(t,\tau)|$$
$$\cdot\sum^{n}_{i=1}\delta_{i}\| \rho[h^{n+1}]\|_{\mathcal{F}^{(\lambda'_{n}-B_{0}) \tau+\mu'_{n}}}d\tau+\int^{t}_{0}
\frac{(\tilde{K}_{0}^{n+1}+\tilde{K}_{1}^{n+1})}{(1+\tau)^{2}}\|\rho\|_{\mathcal{F}^{(\lambda'_{n}-B_{0}) \tau+\mu'_{n}}}d\tau,$$
where $K^{n}_{0}(t,\tau),K^{n}_{1}(t,\tau)$ are defined in Proposition 4.9, and
 $$\tilde{K}_{0}^{n+1}\triangleq 2C\bigg(C_{0}+\sum^{n}_{i=1}\delta_{i}\bigg)\bigg(\sum^{n}_{i=1}\frac{\delta_{i}}{(2\pi(\lambda_{i}-\lambda'_{i}))^{5}}\bigg),$$
$$\tilde{K}_{1}^{n+1}\triangleq\bigg(C^{4}_{\omega}\bigg(C'_{0}+\sum^{n}_{i=1}\delta_{i}\bigg)\bigg(\sum^{n}_{j=1}\frac{\delta_{j}}
{2\pi(\lambda_{j}-\lambda'_{j})^{6}}\bigg)
+\sum^{n}_{i=1}\delta_{i}\bigg).
$$
\begin{prop}
From the above inequality, we obtain the following integral inequality:
\begin{align}
&\bigg\|\rho[h^{n+1}](t,x)-\int^{t}_{0}\int_{\mathbb{R}^{3}}\bigg(E[h^{n+1}]+v\times B[h^{n+1}]\bigg)(\tau,x+\mathcal{M}(t-\tau)v)\cdot\nabla_{v}f^{0}dvd\tau\bigg\|
_{\mathcal{F}^{(\lambda'_{n}-B_{0})t+\mu'_{n}}}\notag\\
&\leq\frac{C\delta^{2}_{n}}{(\lambda_{n}-\lambda'_{n})^{2}}+\int^{t}_{0}\bigg(K_{1}^{'n}(t,\tau)+K_{0}^{'n}(t,\tau)+\frac{c^{n}_{0}}{(1+\tau)^{2}}\bigg)
\|\rho[h^{n+1}](\tau,\cdot)\|
_{\mathcal{F}^{(\lambda'_{n}-B_{0})\tau+\mu'_{n}}}d\tau,
\end{align}
where $K_{1}^{'n}(t,\tau)=|K^{n}_{1}(t,\tau)|(1+\tau)\sum^{n}_{i=1}\delta_{i},\quad K_{0}^{'n}(t,\tau)=|K^{n}_{0}(t,\tau)|\sum^{n}_{i=1}\delta_{i},$
$c^{n}_{0}=\tilde{K}_{0}^{n+1}+\tilde{K}_{1}^{n+1}.$
\end{prop}
%$  \mathbf{Step}$ $\mathbf{6}.$ $(\mathbf{Iteration})$
%$  \mathbf{Step}$ $\mathbf{7}.$ $(\mathbf{Growth}$ $\mathbf{control})$

\begin{center}
\item\subsection{ Exponential moments of the kernel}
%{\bf\large 1. \quad Introduction }
\end{center}

Before doing the iteration, based on the analysis and computation in Sections 4-7, we explain the connection and difference in  between electromagnetic field case and
 electric field case. First, we state the connection, from (6.3)-(6.4) of Theorem 6.1 in section 6.1, the kernel is the same under  two different norms $\mathcal{Z}_{t}^{\lambda,\mu;1}$ and
  $\mathcal{F}^{\lambda t+\mu}$,
 in other words, there are no new echoes to generate in different norms $\mathcal{Z}_{t}^{\lambda,\mu;1}$ and
  $\mathcal{F}^{\lambda t+\mu},$ this is a key point for us, it implies that the first-order term of the magnetic field has no influence
  on the resonances of the plasma particles, that
 is, the first-order term of the magnetic field can be regarded as an error term and doesn't play an important role in dynamical behavior of the particles' trajectory.
The difference is to add the new term $V=-\int^{t}_{0}\int_{\mathbb{R}^{3}}[(B[f^{n}]\circ\Omega^{n}_{s,t}(x,v))\cdot
((\nabla_{v}h^{n+1}\times v)\circ\Omega^{n}_{s,t}(x,v))](s, X^{0}_{s,t}(x,v),V^{0}_{s,t}(x,v))dvds$ in the density equation on $\rho[h^{n+1}].$  In order to
deal with this term,
we have to go back to the Vlasov equation on $h^{n+1}(s, X^{n}_{s,t}(x,v),V^{n}_{s,t}(x,v)),$ and then there are new echoes to appear during estimating
$h^{n+1}(s, X^{n}_{s,t}(x,v),$ $V^{n}_{s,t}(x,v)),$ but the echoes' form is still the same with that generated by  the Vlasov-Poisson equation, the details can be found
in section 5.1. From the inequalities (6.3)-(6.4), we know
that the reason  is that  echoes don't change as the norm changes. In summary, in order to iterate on the density $\rho[h^{n+1}],$ we only need to estimate the
same kernel with Landau damping in electric field case. The following theorems are the same with the results in [23] and the proofs can also be found in Section 7 in [23], so we
sketch the proof.

\begin{prop}(Exponential moments of the kernel) Let $\gamma\in(1,\infty)$ be given. For any $\alpha\in(0,1),$ let $K^{(\alpha),\gamma}$ be defined
$$K^{(\alpha),\gamma}(t,\tau)=(1+\tau)\sup_{k,l\in\mathbb{Z}_{\ast}}\frac{e^{-\alpha|l|}e^{-\alpha(t-\tau)\frac{|k-l|}{t}}e^{-\alpha|k(t-\tau)+l\tau|}}{1+|k-l|^{\gamma}}.$$
Then for any $\gamma<\infty,$ there is $\bar{\alpha}=\bar{\alpha}(\gamma)>0$ such that if $\alpha\leq\bar{\alpha}$ and $\varepsilon\in(0,1),$ then for any
$t>0,$
$$e^{-\varepsilon t}\int^{t}_{0}K^{(\alpha),\gamma}(t,\tau)e^{\varepsilon \tau}d\tau
\leq C\bigg(\frac{1}{\alpha\varepsilon^{\gamma}t^{\gamma-1}}+\frac{1}{\alpha\varepsilon^{\gamma}t^{\gamma}}\log\frac{1}{\alpha}
+\frac{1}{\alpha^{2}\varepsilon^{1+\gamma}t^{1+\gamma}}+\bigg(\frac{1}{\alpha^{3}}+\frac{1}{\alpha^{2}\varepsilon}\log\frac{1}{\alpha}\bigg)
e^{-\frac{\varepsilon t}{4}}+\frac{e^{-\frac{\alpha t}{2}}}{\alpha^{3}}\bigg),$$
where $C=C(\gamma).$

In particular, if $\varepsilon\leq\alpha,$ then
$e^{-\varepsilon t}\int^{t}_{0}K^{(\alpha),\gamma}(t,\tau)e^{\varepsilon \tau}d\tau\leq\frac{C(\gamma)}{\alpha^{3}\varepsilon^{1+\gamma}t^{\gamma-1}}.$
\end{prop}
$Proof.$ Without loss of generality, we shall set $d=1$ and first consider $\tau\leq\frac{1}{2}t.$ We can write
$$K^{(\alpha)}(t,\tau)\leq(1+\tau)\sup_{l\in\mathbb{Z},k\in\mathbb{Z}}e^{-\alpha|l|}e^{-\alpha|k-l|/2}e^{-\alpha|k(t-\tau)+l\tau|}.$$
By symmetry, we may also assume that $k>0.$

Explicit computations yield
\begin{align}
\int^{\frac{t}{2}}_{0}e^{-\alpha|k(t-\tau)+l\tau|}(1+\tau)d\tau\leq
\left\{\begin {array}{l}
\frac{1}{\alpha(l-k)}+\frac{1}{\alpha^{2}(l-k)^{2}},\quad \textmd{if}\quad l>k,\\
e^{-\alpha kt}(\frac{t}{2}+\frac{t^{2}}{8}), \quad \textmd{if}\quad l=k,\\
\frac{e^{-\alpha(k+l)t/2}}{\alpha|k-l|}(1+\frac{t}{2}), \quad \textmd{if}\quad -k\leq l<k,\\
\frac{2}{\alpha|k-l|}+\frac{2kt}{\alpha|k-l|^{2}}+\frac{1}{\alpha^{2}|k-l|^{2}}, \quad \textmd{if}\quad l<-k.\\
\end{array}\right.
\end{align}

So from (8.4), we have
$$e^{-\varepsilon t}\int^{\frac{t}{2}}_{0}e^{-\alpha|k(t-\tau)+l\tau|}(1+\tau)e^{\varepsilon\tau}d\tau$$
$$\leq e^{-\frac{\varepsilon t}{4}}\bigg(\frac{3}{\alpha|k-l|}+\frac{1}{\alpha^{2}|k-l|^{2}}+\frac{8z}{\alpha\varepsilon|k-l|}\bigg)
1_{k\neq l}+e^{-\frac{t\alpha}{2}}\bigg(\frac{z}{\alpha}+\frac{8z^{2}}{\alpha^{2}}\bigg)1_{l=k},$$
where $z=\sup_{x}xe^{-x}=e^{-1}.$

Using the bounds (for $\alpha\sim 0^{+}$)
$$\sum_{l\in\mathbb{Z}}e^{-\alpha l}=O\bigg(\frac{1}{\alpha}\bigg),\quad \sum_{l\in\mathbb{Z}}\frac{e^{-\alpha l}}{l}=O\bigg(\log\frac{1}{\alpha}\bigg),\quad
\sum_{l\in\mathbb{Z}}\frac{e^{-\alpha l}}{l^{2}}=O(1),$$
we end up, for $\alpha\leq\frac{1}{4},$ with a bound like
$$e^{-\varepsilon t}\int^{\frac{t}{2}}_{0}K^{(\alpha)}(t,\tau)e^{\varepsilon\tau}d\tau
\leq C\bigg[e^{-\frac{\varepsilon t}{4}}\bigg(\frac{1}{\alpha^{3}}+\frac{1}{\alpha^{2}\varepsilon}\bigg)
+\frac{e^{-\alpha t/2}}{\alpha^{3}}\bigg].$$

Next we turn to the more delicate contribution of $\tau\geq\frac{1}{2}t.$ We write
\begin{align}
K^{(\alpha)}(t,\tau)\leq(1+\tau)\sup_{l\in\mathbb{Z}_{\ast}}e^{-\alpha|l|}\sup_{k\in\mathbb{Z}}\frac{e^{-\alpha|k(t-\tau)+l\tau|}}{1+|k-l|^{\gamma}}.
\end{align}

Without loss of generality, we restrict the supremum $l>0.$ The function $x\rightarrow(1+|x-l|^{\gamma})^{-1}e^{-\alpha|x(t-\tau)+l\tau|}$
is decreasing for $x\geq l,$ increasing for $x\leq-l\tau/(t-\tau);$ and on the interval $[-l\tau/(t-\tau),l],$ its logarithmic derivative goes
from $\bigg(-\alpha+\frac{\gamma/lt}{1+((t-\tau)/lt)^{\gamma}}\bigg)(t-\tau)$ to $-\alpha(t-\tau).$ It is easy to check that a given integer $k$
occurs in the supremum only for some times $\tau$ satisfying $k-1<-l\tau/(t-\tau)<k+1.$ We can assume $k\geq0,$ then $k-1<\frac{l\tau}{t-\tau}<k+1$ holds, and it
is equivalent to $\frac{k-1}{k+l-1}t<\tau<\frac{k+1}{k+l+1}t.$ More importantly, $\tau>\frac{1}{2}t$ implies that $k\geq l.$ Thus, for $t\geq\frac{\gamma}{\alpha},$
we have \begin{align}
e^{-\varepsilon t}\int^{t}_{\frac{t}{2}}K^{(\alpha)}(t,\tau)e^{\varepsilon\tau}d\tau
\leq e^{-\varepsilon t}\sum^{\infty}_{l=1}e^{-\alpha l}\sum^{\infty}_{k=l}\int^{\frac{(k+1)t}{k+l+1}}_{\frac{(k-1)t}{k+l-1}}(1+\tau)
\frac{e^{\alpha|k(t-\tau)-l\tau|}e^{\varepsilon\tau}}{1+|k+l|^{\gamma}}d\tau.
\end{align}
 For $t\leq\frac{\gamma}{\alpha},$ it is easy to check that  $e^{-\varepsilon t}\int^{t}_{\frac{t}{2}}K^{(\alpha)}(t,\tau)
e^{\varepsilon\tau}d\tau\leq\frac{\gamma}{2\alpha}$ holds. Next we shall focus on $(8.6).$
According to $\tau$ smaller or larger than $kt/(k+l),$ we separate the integral in the right-hand side of (8.6) into two parts, and by simple computation, we get
the explicit bounds
$$ e^{-\varepsilon t}\int^{kt/(k+l)}_{(k-l)t/(k+l-1)}(1+\tau)e^{-\alpha|k(t-\tau)-l\tau|}e^{\varepsilon\tau}d\tau\leq e^{-\frac{\varepsilon lt}{k+l}}
\bigg(\frac{1}{\alpha(k+l)}+\frac{kt}{\alpha(k+l)^{2}}\bigg),$$
$$e^{-\varepsilon t}\int^{\frac{(k+1)t}{k+l+1}}_{\frac{kt}{k+l}}(1+\tau)e^{-\alpha|k(t-\tau)-l\tau|}e^{\varepsilon\tau}d\tau
\leq e^{-\frac{\varepsilon lt}{k+l+1}}
\bigg(\frac{1}{\alpha(k+l)}+\frac{kt}{\alpha(k+l)^{2}}+\frac{1}{\alpha^{2}(k+l)^{2}}\bigg).$$
Hence, (8.6) is bounded above by
\begin{align}
C\sum^{\infty}_{l=1}e^{-\alpha l}\sum^{\infty}_{k=l}\bigg(\frac{1}{\alpha^{2}(k+l)^{2+\gamma}}+\frac{1}{\alpha(k+l)^{1+\gamma}}
+\frac{kt}{\alpha(k+l)^{2+\gamma}}\bigg)e^{-\frac{\varepsilon lt}{k+l}}.
\end{align}

We consider the first term $I(t)$ of (8.7) and change variables $(x,y)\mapsto (x,u),$ where $u(x,y)=\frac{\varepsilon xt}{x+y},$ then we can find that
$$I(t)=\frac{1}{\alpha^{2}\varepsilon^{1+\gamma}t^{1+\gamma}}\int^{\infty}_{1}\frac{e^{-\alpha x}}{x^{1+\gamma}}dx\int^{\varepsilon t/2}_{0}e^{-u}u^{\gamma}du
=O\bigg(\frac{1}{\alpha^{2}\varepsilon^{1+\gamma}t^{1+\gamma}}\bigg).$$

The same computation for the second integral of (8.7) yields
$$\frac{1}{\alpha\varepsilon^{\gamma}t^{\gamma}}\int^{\infty}_{1}\frac{e^{-\alpha x}}{x^{\gamma}}dx\int^{\varepsilon t/2}_{0}e^{-u}u^{\gamma-1}du
=O\bigg(\frac{1}{\alpha\varepsilon^{\gamma}t^{\gamma}}\bigg).$$

Finally, we estimate the last term of (8.7) that is the worst. It yields a contribution $\frac{t}{\alpha}\sum^{\infty}_{l=1}e^{-\alpha l}\sum^{\infty}_{k=l}
\frac{e^{-\varepsilon ltk}/(k+l)}{(k+l)^{2+\gamma}}.$ We compare this with the integral $\frac{t}{\alpha}\int^{\infty}_{1}e^{-\alpha x}\int^{\infty}_{x}
\frac{e^{-\varepsilon ltx/(x+y)}}{(x+y)^{2+\gamma}}dydx,$
and the same change of variables as before equates this with
$$\frac{1}{\alpha\varepsilon^{\gamma}t^{\gamma-1}}\int^{\infty}_{1}\frac{e^{-\alpha x}}{x^{\gamma}}dx\int^{\frac{\varepsilon t}{2}}_{0}e^{-u}u^{\gamma-1}du
-\frac{1}{\alpha\varepsilon^{1+\gamma}t^{\gamma}}\int^{\infty}_{1}\frac{e^{-\alpha x}}{x^{\gamma}}dx\int^{\frac{\varepsilon t}{2}}_{0}e^{-u}u^{\gamma}du
=O\bigg(\frac{1}{\alpha\varepsilon^{\gamma}t^{\gamma-1}}\bigg).$$
The proof of Proposition 8.2 follows by collecting all these bounds and keeping only the worst one.
To finish the growth control,  we have to prove the following result.
\begin{prop} With the same notations as in Proposition 8.2, for any $\gamma>1,$ we have
\begin{align}
\sup_{\tau\geq0}e^{\varepsilon\tau}\int^{\infty}_{\tau}e^{-\varepsilon t}K^{(\alpha),\gamma}(t,\tau)dt\leq C(\gamma)\bigg(\frac{1}{\alpha^{2}\varepsilon}+
\frac{1}{\alpha\varepsilon^{\gamma}}\bigg).
\end{align}
\end{prop}

$Proof.$ We first still reduce to $d=1,$ and split the integral as
$$e^{\varepsilon\tau}\int^{\infty}_{\tau}e^{-\varepsilon t}K^{(\alpha),\gamma}(t,\tau)dt
=e^{\varepsilon\tau}\int^{\infty}_{2\tau}e^{-\varepsilon t}K^{(\alpha),\gamma}(t,\tau)dt
+e^{\varepsilon\tau}\int^{2\tau}_{\tau}e^{-\varepsilon t}K^{(\alpha),\gamma}(t,\tau)dt=I_{1}+I_{2}.$$
For the first term $I_{1},$ we have $K^{(\alpha),\gamma}(t,\tau)\leq(1+\tau)\sum^{\infty}_{k=2}\sum_{l\in\mathbb{Z}_{\ast}}
e^{-\alpha|l|-\alpha|k-l|/2}\leq\frac{C(1+\tau)}{\alpha^{2}},$ and thus
$e^{\varepsilon\tau}$
$\int^{\infty}_{\tau}e^{-\varepsilon t}$
$K^{(\alpha),\gamma}(t,\tau)dt\leq\frac{C}{\varepsilon\alpha^{2}}.$

We treat the second term $I_{2}$ as in the proof of Proposition 8.2:
$$e^{\varepsilon\tau}\int^{\infty}_{\tau}e^{-\varepsilon t}K^{(\alpha),\gamma}(t,\tau)dt\leq e^{\varepsilon\tau}(1+\tau)\sum^{\infty}_{l=1}e^{-\alpha l}
\sum^{\infty}_{k=l}\int^{\frac{(k+l-1)\tau}{k-1}}_{\frac{(k+l+1)\tau}{k+1}}\frac{e^{-\alpha|k(t-\tau)-l\tau|}}{1+(k+l)^{\gamma}}e^{-\varepsilon t}dt
\leq\frac{C}{\alpha\varepsilon^{\gamma}},$$
where the last inequality is obtained by the same method in Proposition 8.2 with the change of variable $u=\frac{\varepsilon x\tau}{y}.$
\begin{center}
\item\subsection{ Growth control}
%{\bf\large 1. \quad Introduction }
\end{center}

From now on, we will state the main result of this section that is the same with section 7.4 in [23].
We define $\|\Phi(t)\|_{\lambda}=\sum_{k\in\mathbb{Z}_{\ast}^{3}}|\Phi(k,t)|e^{2\pi\lambda|k|}.$

\begin{thm} Assume that $f^{0}(v),W=W(x)$ satisfy the conditions of Theorem 0.1, and the $(\mathbf{PSC})$ condition holds.
Let $A\geq0,\mu\geq0$ and $\lambda\in(0,\lambda^{\ast}]$
with $0<\lambda^{\ast}<\lambda_{0}.$ Let $(\Phi(k,t))_{k\in\mathbb{Z}_{\ast}^{3},t\geq0}$ be a continuous functions of $t\geq0,$ valued in
$\mathbb{C}^{\mathbb{Z}^{3}_{\ast}},$ such that for all $t\geq0,$
\begin{align}\|\Phi(t)-\int^{t}_{0}K^{0}(t-\tau)\Phi(\tau)d\tau\|_{\lambda t+\mu}\leq A+\int^{t}_{0}(K_{0}(t,\tau)+K_{1}(t,\tau)+\frac{c_{0}}{(1+\tau)^{m}})
\|\Phi(\tau)\|_{\lambda\tau+\mu} d\tau,
\end{align}
where $c_{0}\geq0,m>1,$ and $K_{0}(t,\tau),K_{1}(t,\tau)$ are non-negative kernels. Let $\varphi(t)=\|\Phi(t)\|_{\lambda t+\mu}.$
Then we have the following:

(i) Assume that $\gamma>1$ and $K_{1}=cK^{(\alpha),\gamma}$ for some $c>0,\alpha\in(0,\bar{\alpha}(\gamma)),$ where $K^{(\alpha),\gamma},\bar{\alpha}(\gamma)$
are the same with that
 defined by Proposition 8.2. Then there are positive constants $C,\chi,$ depending only on $\gamma,\lambda^{\ast},\lambda_{0},\kappa,c_{0}, C_{W}$ and $ m,$ uniform
 as $\gamma\rightarrow1,$ such that if
  $\sup_{t\geq0}\int^{t}_{0}K_{0}(t,\tau)d\tau\leq\chi$ and
 $\sup_{t\geq0}$
 $(\int^{t}_{0}K_{0}(t,\tau)^{2}d\tau)^{\frac{1}{2}}+\sup_{t\geq0}\int^{\infty}_{t}K_{0}(t,\tau)dt\leq1,$ then for any $\varepsilon\in(0,\alpha),$
  for all $t\geq0,$
 \begin{align}\varphi(t)\leq CA\frac{1+c^{2}_{0}}{\sqrt{\varepsilon}}e^{Cc_{0}}(1+\frac{c}{\alpha\varepsilon})e^{CT}e^{Cc(1+T^{2})}e^{\varepsilon t},
 \end{align}
 where $T_{\varepsilon}=C\max\bigg\{\bigg(\frac{c^{2}}{\alpha^{5}}\varepsilon^{2+\gamma}\bigg)^{\frac{1}{\gamma-1}},
\bigg(\frac{c}{\alpha^{2}}\varepsilon^{\frac{1}{2}+\gamma}\bigg)^{\frac{1}{\gamma-1}},\bigg(\frac{c^{2}_{0}}{\varepsilon}\bigg)^{\frac{1}{2m-1}}\bigg\}.$

(ii) Assume that $K_{1}=\sum^{N}_{j=1}c_{j}K^{(\alpha_{j},1)}$ for some $\alpha_{j}\in(0,\bar{\alpha}(\gamma)),$ where $\bar{\alpha}(\gamma)$ also appears in proposition 7.2;
then there is a numeric constant $\Gamma>0$ such that whenever $1\geq\varepsilon\geq\Gamma\sum^{N}_{j=1}\frac{c_{j}}{\alpha^{3}_{j}},$ with the same notation
 as in (I), for all $t\geq0,$  one has,

  \begin{align}\varphi(t)\leq CA\frac{1+c^{2}_{0}}{\sqrt{\varepsilon}}e^{Cc_{0}}(1+\frac{c}{\alpha\varepsilon})e^{CT}e^{Cc(1+T^{2})}e^{\varepsilon t},
 \end{align}
 where $c=\sum^{N}_{j=1}c_{j}$ and $T=\max\bigg\{\frac{1}{\varepsilon^{2}}\sum^{N}_{j=1}\frac{c_{j}}{\alpha^{3}_{j}},
 \bigg(\frac{c^{2}_{0}}{\varepsilon}\bigg)^{\frac{1}{2m-1}}\bigg\}.$

\end{thm}
$Proof.$ Here we only prove (i), the  proof of (ii) is similar. We decompose the proof into three steps.

$Step$ 1. Crude pointwise bounds. From (8.9), we have
$$\varphi(t)=\sum_{k\in\mathbb{Z}^{3}_{\ast}}|\Phi(k,t|e^{2\pi(\lambda t+\mu)|k|}
\leq A+\sum_{k\in\mathbb{Z}^{3}_{\ast}}\int^{t}_{0}|K^{0}(k,t-\tau)|e^{2\pi(\lambda t+\mu)|k|}|\Phi(t,\tau)|d\tau$$
$$+\int^{t}_{0}(K_{0}(t,\tau)+K_{1}(t,\tau)+\frac{c_{0}}{(1+\tau)^{m}})\varphi(\tau)d\tau$$
$$\leq A+\int^{t}_{0}(K_{0}(t,\tau)+K_{1}(t,\tau)+\frac{c_{0}}{(1+\tau)^{m}}+\sup_{k\in\mathbb{Z}^{3}_{\ast}}|K^{0}(k,t-\tau)|e^{2\pi\lambda(t-\tau)|k|})
\varphi(\tau)d\tau.$$
We note that for any $k\in\mathbb{Z}^{3}_{\ast}$ and $t\geq0,$
$$|K^{0}(k,t-\tau)|e^{2\pi\lambda|k|(t-\tau)}\leq 4\pi^{2}|\widehat{W}(k)|C_{0}e^{-2\pi(\lambda_{0}-\lambda)|k|t}|k|^{2}t
\leq\frac{CC_{0}C_{W}}{\lambda_{0}-\lambda},$$
where (here and below) $C$ stands for a numeric constant which may change from line to line. Assuming that $\int^{t}_{0}K_{0}(t,\tau)d\tau\leq\frac{1}{2},$
we deduce that
$$\varphi(t)\leq A+\frac{1}{2}\sup_{0\leq\tau\leq t}\varphi(\tau)+C\int^{t}_{0}\bigg(\frac{C_{0}C_{W}}{\lambda_{0}-\lambda}+c(1+\tau)
+\frac{c_{0}}{(1+\tau)^{m}}\bigg)\varphi(\tau)d\tau,$$
and, by Gr$\ddot{\textmd{o}}$nwall's lemma,
\begin{align}
\varphi(t)\leq 2Ae^{C(C_{0}C_{W}t/(\lambda_{0}-\lambda)+c(t+t^{2})+c_{0}C_{m})},
\end{align}
where $C_{m}=\int^{\infty}_{0}(1+\tau)^{-m}d\tau.$

$Step$ 2. $L^{2}$ bound. For all $k\in\mathbb{Z}^{3}_{\ast}$ and $t\geq0,$ we define
$\Psi_{k}(t)=e^{-\varepsilon t}\Phi(k,t)e^{2\pi(\lambda t+\mu)|k|}, \mathcal{K}^{0}_{k}(t)=e^{-\varepsilon t}K^{0}(k,t)$
$e^{2\pi(\lambda t+\mu)|k|},$
$R_{k}(t)$
$=e^{-\varepsilon t}\bigg(\Phi(k,t)-\int^{t}_{0}K^{0}(k,t-\tau)\Phi(k,\tau)d\tau\bigg)e^{2\pi(\lambda t+\mu)|k|}=(\Psi_{k}-\Psi_{k}\ast
\mathcal{K}^{0}_{k})(t),$ and we extend all these functions by $0$ for negative values of $t.$ Taking Fourier transform in the time-variable yields
$\hat{R}_{k}=(1-\widehat{\mathcal{K}}^{0}_{k})\widehat{\Psi_{k}}.$ Since the $(\mathbf{PSC})$ condition implies that $|1-\widehat{\mathcal{K}}^{0}_{k}|\geq\kappa,$
 we can deduce that $\|\hat{\Psi}_{k}\|_{L^{2}}\leq\kappa^{-1}\|\hat{R}_{k}\|_{L^{2}},$ i.e., $\|\Psi_{k}\|_{L^{2}}\leq\kappa^{-1}\|R_{k}\|_{L^{2}}.$ So we have
 \begin{align}\|\Psi_{k}-R_{k}\|_{L^{2}(dt)}\leq\kappa^{-1}\|\mathcal{K}^{0}_{k}\|_{L^{1}(dt)}\|R_{k}\|_{L^{2}(dt)} \quad \textmd{for}\quad \textmd{all}
 \quad k\in\mathbb{Z}_{\ast}^{3}.
 \end{align}

 Then \begin{align}
& \|\varphi(t)e^{-\varepsilon t}\|_{L^{2}(dt)}=\|\sum_{k\in\mathbb{Z}^{3}}|\Psi_{k}|\|_{L^{2}(dt)}\leq\|\sum_{k\in\mathbb{Z}^{3}}|R_{k}|\|_{L^{2}(dt)}
 +\sum_{k\in\mathbb{Z}^{3}}\|R_{k}-\Psi_{k}\|_{L^{2}(dt)}\notag\\
 &\leq\|\sum_{k\in\mathbb{Z}^{3}}|R_{k}|\|_{L^{2}(dt)}(1+\frac{1}{\kappa}).
 \end{align}

 Next, we note that
 $$\|\mathcal{K}^{0}_{k}\|_{L^{1}(dt)}\leq 4\pi^{2}|\widehat{W}(k)|\int^{\infty}_{0}C_{0}e^{-2\pi(\lambda_{0}-\lambda)|k|t}|k|^{2}t dt
 \leq 4\pi|\widehat{W}(k)|\frac{C_{0}}{(\lambda_{0}-\lambda)^{2}},$$
 so $\sum_{ k\in\mathbb{Z}_{\ast}^{3}}\|\mathcal{K}^{0}_{k}\|_{L^{1}(dt)}\leq 4\pi(\sum_{ k\in\mathbb{Z}_{\ast}^{3}}|\widehat{W}(k)|)\frac{C_{0}}
 {(\lambda_{0}-\lambda)^{2}}.$
 Furthermore, we get
 \begin{align}
 &\|\varphi(t)e^{-\varepsilon t}\|_{L^{2}(dt)}\leq\bigg(1+\frac{CC_{0}C_{W}}{\kappa(\lambda_{0}-\lambda)^{2}}\bigg)
 \|\sum_{k\in\mathbb{Z}_{\ast}^{3}}\|_{L^{2}(dt)}\notag\\
 &\leq \bigg(1+\frac{CC_{0}C_{W}}{\kappa(\lambda_{0}-\lambda)^{2}}\bigg)\bigg(\int^{\infty}_{0}e^{-2\varepsilon t}
 \bigg(A+\int^{t}_{0}\bigg(K_{1}+K_{0}+\frac{c_{0}}{(1+\tau)^{m}}\bigg)\varphi(\tau)d\tau\bigg)^{2}\bigg)^{\frac{1}{2}}.
 \end{align}

 By Minkowski's inequality, we separate (8.15) into various contributions which we estimate separately. First,
 $\bigg(\int^{\infty}_{0}e^{-2\varepsilon t}A^{2}dt\bigg)^{\frac{1}{2}}=\frac{A}{\sqrt{2\varepsilon}}.$ Next, for any $T\geq1,$ by Step 1 and
 $\int^{t}_{0}K_{1}(t,\tau)d\tau\leq\frac{Cc(1+t)}{\alpha},$ we have
 \begin{align}
 &\bigg(\int^{T}_{0}e^{-2\varepsilon t}\bigg(\int^{t}_{0}K_{1}(t,\tau)\varphi(\tau)\bigg)^{2}\bigg)^{\frac{1}{2}}
 \leq (\sup_{0\leq t\leq T}\varphi(t)) \bigg(\int^{T}_{0}e^{-2\varepsilon t}\bigg(\int^{t}_{0}K_{1}(t,\tau)\bigg)^{2}\bigg)^{\frac{1}{2}}\notag\\
 &\leq CAe^{C(C_{0}C_{W}T/(\lambda_{0}-\lambda)+c(T+T^{2}))}\frac{c}{\alpha}\bigg(\int^{\infty}_{0}e^{-2\varepsilon t}(1+t)^{2}dt\bigg)^{\frac{1}{2}}\notag\\
 &\leq CA\frac{c}{a\varepsilon^{\frac{3}{2}}}e^{C(C_{0}C_{W}T/(\lambda_{0}-\lambda)+c(T+T^{2}))}.
 \end{align}

 Invoking Jensen's inequality and Fubini's theorem, we also have
 \begin{align}
 &\int^{\infty}_{T}e^{-2\varepsilon t}\bigg(\int^{t}_{0}K_{1}(t,\tau)\varphi(\tau)d\tau\bigg)^{2}dt\bigg)^{\frac{1}{2}}
 =\int^{\infty}_{T}\bigg(\int^{t}_{0}K_{1}(t,\tau)e^{-2\varepsilon (t-\tau)}e^{-2\varepsilon\tau}\varphi(\tau)d\tau\bigg)^{2}dt\bigg)^{\frac{1}{2}}\notag\\
 &\leq\int^{\infty}_{T}\bigg(\int^{t}_{0}K_{1}(t,\tau)e^{-\varepsilon (t-\tau)}d\tau\bigg)\bigg(
\int^{\infty}_{T}\bigg(\int^{t}_{0}K_{1}(t,\tau)e^{-\varepsilon (t-\tau)}e^{-2\varepsilon\tau}\varphi(\tau)^{2}d\tau\bigg)dt\bigg)^{\frac{1}{2}}\notag\\\notag\\
&\leq\bigg(\sup_{t\geq T}\int^{t}_{0}K_{1}(t,\tau)e^{-\varepsilon (t-\tau)}d\tau\bigg)^{\frac{1}{2}}\bigg(
\int^{\infty}_{T}\bigg(\int^{t}_{0}K_{1}(t,\tau)e^{-\varepsilon (t-\tau)}e^{-2\varepsilon\tau}\varphi(\tau)^{2}d\tau\bigg)dt\bigg)^{\frac{1}{2}}\notag\\\notag\\
&=\bigg(\sup_{t\geq T}\int^{t}_{0}K_{1}(t,\tau)e^{-\varepsilon (t-\tau)}d\tau\bigg)^{\frac{1}{2}}\bigg(
\int^{\infty}_{0}\int^{\infty}_{\max\{\tau,T\}}K_{1}(t,\tau)e^{-\varepsilon (t-\tau)}e^{-2\varepsilon\tau}\varphi(\tau)^{2}dtd\tau \bigg)^{\frac{1}{2}}\notag\\
&\leq\bigg(\sup_{t\geq T}\int^{t}_{0}K_{1}(t,\tau)e^{-\varepsilon (t-\tau)}d\tau\bigg)^{\frac{1}{2}}\bigg(
\sup_{\tau\geq0}\int^{\infty}_{\tau}K_{1}(t,\tau)e^{-\varepsilon (t-\tau)}e^{-2\varepsilon\tau}dt\bigg)^{\frac{1}{2}}
\bigg(\int^{\infty}_{0}\varphi(\tau)^{2}e^{-2\varepsilon\tau}d\tau \bigg)^{\frac{1}{2}}.\notag\\
\end{align}
Similarly,
\begin{align}
 &\int^{\infty}_{T}e^{-2\varepsilon t}\bigg(\int^{t}_{0}K_{0}(t,\tau)\varphi(\tau)d\tau\bigg)^{2}dt\bigg)^{\frac{1}{2}}
 \leq\bigg(\sup_{t\geq T}\int^{t}_{0}K_{0}(t,\tau)d\tau\bigg)^{\frac{1}{2}}\bigg(
\sup_{\tau\geq0}\int^{\infty}_{\tau}K_{0}(t,\tau)dt\bigg)^{\frac{1}{2}}
\bigg(\int^{\infty}_{0}\varphi(\tau)^{2}d\tau \bigg)^{\frac{1}{2}}.
\end{align}

The last term is also split, this time according to $\tau\leq T$ or $\tau> T:$
\begin{align}&\bigg(\int^{\infty}_{0}e^{-2\varepsilon t}\bigg(\int^{T}_{0}\frac{c_{0}\varphi(\tau)}{(1+\tau)^{m}}d\tau\bigg)^{2}dt\bigg)^{\frac{1}{2}}
\leq c_{0}(\sup_{0\leq\tau\leq T}\varphi(\tau))
\bigg(\int^{\infty}_{0}e^{-2\varepsilon t}\bigg(\int^{T}_{0}\frac{d\tau}{(1+\tau)^{m}}\bigg)^{2}dt\bigg)^{\frac{1}{2}}\notag\\
 &\leq c_{0}\frac{CA}{\sqrt{\varepsilon}}e^{C(C_{0}C_{W}T/(\lambda_{0}-\lambda)+c(T+T^{2}))}C_{m},
\end{align}
and \begin{align}
&\bigg(\int^{\infty}_{0}e^{-2\varepsilon t}\bigg(\int^{t}_{T}\frac{c_{0}\varphi(\tau)}{(1+\tau)^{m}}d\tau\bigg)^{2}dt\bigg)^{\frac{1}{2}}
\leq c_{0}\bigg(\int^{\infty}_{0}e^{-2\varepsilon t}\varphi(t)^{2}\bigg)^{\frac{1}{2}}
\bigg(\int^{\infty}_{0}\int^{t}_{T}\frac{e^{-2\varepsilon (t-\tau)}}{(1+\tau)^{2m}}d\tau dt\bigg)^{\frac{1}{2}}\notag\\
&=c_{0}\bigg(\int^{\infty}_{0}e^{-2\varepsilon t}\varphi(t)^{2}\bigg)^{\frac{1}{2}}
\bigg(\int^{\infty}_{T}\frac{d\tau}{(1+\tau)^{2m}}\bigg)^{\frac{1}{2}}
\bigg(\int^{\infty}_{0}e^{-2\varepsilon s} ds\bigg)^{\frac{1}{2}}
=\frac{C_{2m}^{\frac{1}{2}}c_{0}}{T^{m-\frac{1}{2}}\sqrt{\varepsilon}}\bigg(\int^{\infty}_{0}e^{-2\varepsilon t}\varphi(t)^{2}\bigg)^{\frac{1}{2}}.
\end{align}
Gathering estimates (8.16)-(8.20), we deduce from (8.15) that
\begin{align}&\|\varphi(t)e^{-\varepsilon t}\|_{L^{2}(dt)}\leq\bigg(1+\frac{CC_{0}C_{W}}{\kappa(\lambda_{0}-\lambda)^{2}}\bigg)\frac{CA}{\sqrt{\varepsilon}}
\bigg(1+\frac{c}{a\varepsilon}+c_{0}C_{m}\bigg)e^{C(C_{0}C_{W}T/(\lambda_{0}-\lambda)+c(T+T^{2}))}\notag\\
&\quad \quad \quad \quad\quad\quad\quad+a\|\varphi(t)e^{-\varepsilon t}\|_{L^{2}(dt)},
\end{align}
where $$a=\bigg(1+\frac{CC_{0}C_{W}}{\kappa(\lambda_{0}-\lambda)^{2}}\bigg)\bigg[
\bigg(\sup_{t\geq T}\int^{t}_{0}e^{-\varepsilon t}K_{1}(t,\tau)e^{\varepsilon\tau}d\tau\bigg)^{\frac{1}{2}}
\bigg(\sup_{\tau\geq 0}\int^{\infty}_{\tau}e^{\varepsilon \tau}K_{1}(t,\tau)e^{-\varepsilon t}dt\bigg)^{\frac{1}{2}}$$
$$+\bigg(\sup_{t\geq T}\int^{t}_{0}K_{0}(t,\tau)d\tau\bigg)^{\frac{1}{2}}\bigg(
\sup_{\tau\geq0}\int^{\infty}_{\tau}K_{0}(t,\tau)dt\bigg)^{\frac{1}{2}}+\frac{C_{2m}^{\frac{1}{2}}c_{0}}{T^{m-\frac{1}{2}}\sqrt{\varepsilon}}\bigg].$$

Using Propositions 8.2 and 8.3, together with the assumptions of Theorem 8.4, we see that $a\leq\frac{1}{2}$ for $\chi$ sufficiently small. Then we have
$$\|\varphi(t)e^{-\varepsilon t}\|_{L^{2}(dt)}\leq\bigg(1+\frac{CC_{0}C_{W}}{\kappa(\lambda_{0}-\lambda)^{2}}\bigg)\frac{CA}{\sqrt{\varepsilon}}
\bigg(1+\frac{c}{a\varepsilon}+c_{0}C_{m}\bigg)e^{C(C_{0}C_{W}T/(\lambda_{0}-\lambda)+c(T+T^{2}))}.$$

$Step$ 3. For $t\geq T,$ using (8.9), we get
\begin{align}& e^{-\varepsilon t}\varphi(t)\leq Ae^{-\varepsilon t}+\bigg[\bigg(\int^{t}_{0}\bigg(\sup_{k\in\mathbb{Z}^{3}_{\ast}}|K^{0}(k,t-\tau)|
e^{2\pi\lambda(t-\tau)|k|}\bigg)^{2}d\tau\bigg)^{\frac{1}{2}}\notag\\
&+\bigg(\int^{t}_{0}K_{0}(t,\tau)^{2}d\tau\bigg)^{\frac{1}{2}}+\bigg(\int^{\infty}_{0}\frac{c^{2}_{0}}{(1+\tau)^{2m}}d\tau\bigg)^{\frac{1}{2}}
+\bigg(\int^{t}_{0}e^{-2\varepsilon t}K_{1}(t,\tau)^{2}e^{2\varepsilon \tau}d\tau\bigg)^{\frac{1}{2}}\bigg]
\bigg(\int^{\infty}_{0}\varphi(\tau)e^{-\varepsilon\tau}d\tau\bigg)^{\frac{1}{2}}.
\end{align}

We note that, for any $k\in\mathbb{Z}^{3}_{\ast},$ $(|K^{0}(k,t)|e^{2\pi\lambda|k|t})^{2}\leq C\pi^{4}|\widehat{W}(k)|^{2}|\hat{f}^{0}(kt)|^{2}|k|^{4}t^{2}
\leq \frac{CC_{0}}{(\lambda_{0}-\lambda)^{2}}C^{2}_{W}e^{-2\pi(\lambda_{0}-\lambda)t},$ so we get $\int^{t}_{0}\bigg(\sup_{k\in\mathbb{Z}_{\ast}^{3}}
|K^{0}(k,t-\tau)|e^{2\pi\lambda(t-\tau)|k|}\bigg)^{2}d\tau\leq\frac{CC^{2}_{0}C^{2}_{W}}{(\lambda_{0}-\lambda)^{3}}.$

From Proposition 8.2, (8.22), the conditions of Theorem 8.4 and Step 2, the conclusion is finished.
\begin{cor} Assume that $f^{0}=f^{0}(v),$ under the assumptions of Theorem 0.1,
 we pick up $\lambda^{\ast}_{n}-B_{0}<\lambda'_{n}-B_{0}$ such that $2\pi(\lambda'_{n}-\lambda^{\ast}_{n})\leq\alpha_{n},$ choose $\varepsilon
=2\pi(\lambda'_{n}-\lambda^{\ast}_{n});$
recalling that $\hat{\rho}(t,0)=0,$  our conditions imply an upper bound on $c_{n}$ and $c^{n}_{0},$ we have the uniform control,
$$\|\rho[h^{n+1}](t,x)\|
_{\mathcal{F}^{(\lambda^{\ast}_{n}-B_{0})t+\mu'_{n}}}\leq\frac{C\delta^{2}_{n}(1+c^{n}_{0})^{2}}{\sqrt{\varepsilon}(\lambda_{n}-\lambda'_{n})^{2}}
\bigg(1+\frac{1}{\alpha_{n}(\lambda'_{n}-\lambda^{\ast}_{n})^{\frac{3}{2}}}\bigg)e^{CT^{2}_{n}},$$
where $T_{n}=C\bigg(\frac{1}{\alpha^{5}_{n}(\lambda'_{n}-\lambda^{\ast}_{n})}\bigg)^{\frac{1}{\gamma-1}}.$
\end{cor}

$Proof.$ From Propositions 8.1-8.3, we know that
$$\int^{t}_{0}K^{n}_{0}(t,\tau)d\tau\leq C_{W}\sum^{n}_{i=1}\frac{\delta_{i}}{\pi(\lambda_{i}-\lambda'_{n})},\quad
\int^{\infty}_{\tau}K^{n}_{0}(t,\tau)d\tau\leq C_{W}\sum^{n}_{i=1}\frac{\delta_{i}}{\pi(\lambda_{i}-\lambda'_{n})},$$
$$\bigg(\int^{t}_{0}K^{n}_{0}(t,\tau)^{2}d\tau\bigg)^{\frac{1}{2}}\leq C_{W}\sum^{n}_{i=1}\frac{\delta_{i}}{\sqrt{2\pi(\lambda_{i}-\lambda'_{n})}}.$$
Here $\alpha_{n}=\pi\min\{(\mu_{n}-\mu'_{n}),(\lambda_{n}-\lambda'_{n})\},$ and assume that $\alpha_{n}$ is smaller than $\bar{\alpha}(\gamma)$ in Theorem 8.4, and
that
$$\bigg(C^{4}_{\omega}\bigg(C'_{0}+\sum^{n}_{i=1}\delta_{i}+1\bigg)\bigg(\sum^{n}_{j=1}\frac{\delta_{j}}
{2\pi(\lambda_{j}-\lambda'_{j})^{6}}\bigg)\leq\frac{1}{8},\quad\quad (\mathbf{VII})$$
$$C_{W}\sum^{n}_{i=1}\frac{\delta_{i}}{\sqrt{2\pi(\lambda_{i}-\lambda'_{n})}}\leq\frac{1}{4},\quad
\sum^{n}_{i=1}\frac{\delta_{i}}{\pi(\lambda_{i}-\lambda'_{n})}\leq\max\{\chi,\frac{1}{8}\}.\quad (\mathbf{VIII})$$
Applying Theorem 8.4, we can deduce that for any $\varepsilon\in(0,\alpha_{n})$ and $t\geq0,$
$$\|\rho[h^{n+1}](t,x)\|
_{\mathcal{F}^{(\lambda^{\ast}_{n}-B_{0})t+\mu'_{n}}}\leq e^{-2\pi(\lambda'_{n}-\lambda^{\ast}_{n})t}\|\rho[h^{n+1}](t,x)\|
_{\mathcal{F}^{(\lambda'_{n}-B_{0})t+\mu'_{n}}}$$
$$\leq\frac{C\delta^{2}_{n}(1+c^{n}_{0})^{2}}{\sqrt{\varepsilon}(\lambda_{n}-\lambda'_{n})^{2}}
\bigg(1+\frac{1}{\alpha_{n}(\lambda'_{n}-\lambda^{\ast}_{n})^{\frac{3}{2}}}\bigg)e^{CT^{2}_{n}},$$
where $T_{n}=C\bigg(\frac{1}{\alpha^{5}_{n}(\lambda'_{n}-\lambda^{\ast}_{n})}\bigg)^{\frac{1}{\gamma-1}}.$
%\begin{align}
%&\|\rho[h^{n+1}](t,x)\|
%_{\mathcal{F}^{(\lambda'_{n}-B_{0})t+\mu'_{n}}}\leq \frac{C\delta^{2}_{n}(1+c^{n}_{0})^{2}}{\sqrt{\varepsilon}(\lambda_{n}-\lambda'_{n})^{2}}
%e^{Cc^{n}_{0}}\bigg(1+\frac{c_{n}}{\alpha_{n}\varepsilon}\bigg)e^{CT_{\varepsilon,n}}e^{Cc_{n}(1+T^{2}_{\varepsilon,n})}e^{\varepsilon t},\notag\\
%\end{align}
%where $c_{n}=2C\sum^{n}_{i=1}\delta_{i},$
%and $T_{\varepsilon,n}=C_{\gamma}\max\bigg\{\bigg(\frac{c^{2}_{n}}{\alpha^{5}_{n}}\varepsilon^{2+\gamma}\bigg)^{\frac{1}{\gamma-1}},
%\bigg(\frac{c_{n}}{\alpha^{2}_{n}}\varepsilon^{\frac{1}{2}+\gamma}\bigg)^{\frac{1}{\gamma-1}},\frac{(c^{n}_{0})^{\frac{2}{3}}}{\varepsilon^{\frac{1}{3}}}\bigg\}.$

\begin{center}
\item\subsection{ Estimates related to $h^{n+1}(t,X^{n}_{\tau,t}(x,v),V^{n}_{\tau,t}(x,v)))$}
%{\bf\large 1. \quad Introduction }
\end{center}

Next we show the control on $h^{i}.$
\begin{lem} For  any $ n \geq i\geq1,$
$$\| (\nabla'_{v}\times((h^{i}_{\tau}v)\circ\Omega^{i-1}_{t,\tau})-\langle\nabla'_{v}\times((h^{i}_{\tau}v)\circ\Omega^{i-1}_{t,\tau})
\rangle)\|_{\mathcal{Z}^{(\lambda'_{i}-B_{0})(1+b),
\mu'_{i};1}_{\tau-\frac{bt}{1+b}}}\leq (1+\tau)\delta_{i}.$$
 \end{lem}
$Proof.$  First, we consider $i=1.$

 In fact,
$$\|\nabla'_{v}\times(h^{1}_{\tau}v)-\langle\nabla'_{v}\times((h^{1}_{\tau}v)
\rangle\|_{\mathcal{Z}^{(\lambda'_{1}-B_{0})(1+b),
\mu'_{1};1}_{\tau-\frac{bt}{1+b}}}$$
$$\leq\|\nabla'_{v}\times(h^{1}
_{\tau}v)\|_{\mathcal{Z}^{(\lambda'_{1}-B_{0})(1+b),
\mu'_{1};1}_{\tau-\frac{bt}{1+b}}}+\bigg\|\int_{\mathbb{T}^{3}}\nabla'_{v}\times((h^{1}
_{\tau}v)dx\bigg\|_{\mathcal{C}^{(\lambda'_{1}-B_{0})(1+b);1}}$$
$$\leq\|\nabla'_{v}(h^{1}
_{\tau}v)\|_{\mathcal{Z}^{(\lambda'_{1}-B_{0})(1+b),
\mu'_{1};1}_{\tau-\frac{bt}{1+b}}}\leq\|h^{1}
_{\tau}\|_{\mathcal{Z}^{(\lambda_{1}-B_{0})(1+b),
\mu_{1};1}_{\tau-\frac{bt}{1+b}}}\leq (1+\tau)\delta_{1},$$
 where we use the property $(v)$ of Proposition 2.5.

$$\bigg\|\int_{\mathbb{T}^{3}}\nabla'_{v}\times(h^{1}
_{\tau}v)dx\bigg\|_{\mathcal{C}^{(\lambda'_{1}-B_{0})(1+b);1}}=\|\langle\nabla'_{v}\times(h^{1}
_{\tau}v)\rangle\|_{\mathcal{C}^{(\lambda'_{1}-B_{0})(1+b);1}}$$
$$=\|\langle(\nabla'_{v}+\tau\nabla_{x})\times(h^{1}
_{\tau}v)\rangle\|_{\mathcal{C}^{(\lambda'_{1}-B_{0})(1+b);1}}
\leq\|(\nabla'_{v}+\tau\nabla_{x})\times(h^{1}
_{\tau}v)\|_{\mathcal{Z}^{(\lambda'_{1}-B_{0})(1+b),\mu'_{1};1}_{\tau-\frac{bt}{1+b}}}\leq  \delta_{1},$$
where we use $(vi)$ of Proposition 2.5.

Suppose that $i=k,$ the conclusion holds, that is,
$$\bigg\|\nabla'_{v}\times((h^{k}_{\tau}v)\circ\Omega^{k-1}_{t,\tau})-\langle\nabla'_{v}\times((h^{k}_{\tau}v)\circ\Omega^{k-1}_{t,\tau})
\rangle\bigg\|_{\mathcal{Z}^{(\lambda'_{k}-B_{0})(1+b),
\mu'_{k};1}_{\tau-\frac{bt}{1+b}}}$$
$$\leq\bigg\|h^{k}_{\tau}\circ\Omega^{k-1}_{t,\tau}-\langle h^{k}_{\tau}\circ\Omega^{k-1}_{t,\tau}
\rangle\bigg\|_{\mathcal{Z}^{(\lambda_{k}-B_{0})(1+b),
\mu_{k};1}_{\tau-\frac{bt}{1+b}}}\leq (1+\tau)\delta_{k}.$$

We need to show that the conclusion still holds for $i=k+1.$ We can get the estimate for
$h^{k+1}(t,X^{k}_{t,\tau}(x,v),V^{k}_{t,\tau}(x,v))$ from (4.11).
%For this, we have to reconsider $$h^{k+1}(t,X^{k}_{0,t}(x,v),V^{k}_{0,t}(x,v))=\int^{t}_{0}\Sigma^{k+1}(s,X^{k}_{0,s}(x,v),V^{k}_{0,s}(x,v))ds,$$
%then for $\tau\geq t,$ composing with $(X^{k}_{\tau,0}(x,v),V^{k}_{\tau,0}(x,v)),$ this gives
%$$h^{k+1}(t,X^{k}_{\tau,t}(x,v),V^{k}_{\tau,t}(x,v))=\int^{t}_{0}\Sigma^{k+1}(s,X^{k}_{\tau,s}(x,v),V^{k}_{\tau,s}(x,v))ds.$$

           Note that \begin{align}
           \left\{
           \begin{array}{l}
           (\nabla h)\circ\Omega=(\nabla\Omega)^{-1}\nabla(h\circ\Omega),\\
           (\nabla^{2} h)\circ\Omega=(\nabla\Omega)^{-2}\nabla^{2}(h\circ\Omega)-(\nabla\Omega)^{-1}\nabla^{2}\Omega(\nabla\Omega)^{-1}(\nabla h\circ\Omega).
           \end{array}
           \right.
           \end{align}
         Therefore, from (8.23),   we get
           \begin{align}
          &\|(\nabla h_{\tau}^{n+1})\circ\Omega_{t,\tau}^{n}\|_{\mathcal{Z}_{\tau-\frac{bt}{1+b}}
           ^{(\lambda'_{n+1}-B_{0})(1+b),\mu'_{n+1};1}}
           \leq C(d)\|\nabla (h_{\tau}^{n+1}\circ\Omega_{t,\tau}^{n})\|_{\mathcal{Z}_{\tau-\frac{bt}{1+b}}
           ^{(\lambda'_{n+1}-B_{0})(1+b),\mu'_{n+1};1}}\notag\\
           &\leq\frac{C(d)(1+\tau)}{\min\{\lambda^{\mathfrak{b}}_{n+1}-\lambda'_{n+1},\mu_{n+1}-\mu'_{n+1}\}}\|h_{\tau}^{n+1}\circ\Omega_{t,\tau}^{n}\|_{\mathcal{Z}_{\tau-\frac{bt}{1+b}}
           ^{(\lambda^{\mathfrak{b}}_{n+1}-B_{0})(1+b),\mu_{n+1};1}},
           \end{align}
           and

           \begin{align}
          &\|(\nabla^{2} h_{\tau}^{n+1})\circ\Omega_{t,\tau}^{n}\|_{\mathcal{Z}_{\tau-\frac{bt}{1+b}}
           ^{(\lambda'_{n+1}-B_{0})(1+b),\mu'_{n+1};1}}\notag\\
           &\leq C(d)\bigg[\|\nabla^{2} (h_{\tau}^{n+1}\circ\Omega_{t,\tau}^{n})\|_{\mathcal{Z}_{\tau-\frac{bt}{1+b}}
           ^{(\lambda'_{n+1}-B_{0})(1+b),\mu'_{n+1};1}}\notag\\
           & +\|\nabla^{2}\Omega^{n}_{t,\tau}\|_{\mathcal{Z}_{\tau-\frac{bt}{1+b}}
           ^{(\lambda'_{n+1}-B_{0})(1+b),\mu'_{n+1}}}\|(\nabla h_{\tau}^{n+1})\circ\Omega_{t,\tau}^{n}\|_{\mathcal{Z}_{\tau-\frac{bt}{1+b}}
           ^{(\lambda'_{n+1}-B_{0})(1+b),\mu'_{n+1};1}}\bigg]\notag\\
           &\leq\frac{C(d)(1+\tau)}{\min\{\lambda^{\mathfrak{b}}_{n+1}-\lambda'_{n+1},\mu_{n+1}-\mu'_{n+1}\}}\|h_{\tau}^{n+1}\circ\Omega_{t,\tau}^{n}\|_{\mathcal{Z}_{\tau-\frac{bt}{1+b}}
           ^{(\lambda^{\mathfrak{b}}_{n+1}-B_{0})(1+b),\mu_{n+1};1}}\notag\\
            &\leq\frac{C(d)(1+\tau)}{\min\{\lambda^{\mathfrak{b}}_{n+1}-\lambda'_{n+1},\mu_{n+1}-\mu'_{n+1}\}}\|h_{\tau}^{n+1}\circ\Omega_{t,\tau}^{n}\|
            _{\mathcal{Z}_{\tau+\frac{bt}{1-b}}
           ^{(\lambda'_{n}-B_{0})(1-b),\mu'_{n};1}}.
           \end{align}

           We first write $$\nabla(h_{\tau}^{n+1}\circ\Omega_{t,\tau}^{n})-(\nabla h_{\tau}^{n+1})\circ\Omega_{t,\tau}^{n}=\nabla(\Omega_{t,\tau}^{n}-Id)
           \cdot[(\nabla h_{\tau}^{n+1})\circ\Omega_{t,\tau}^{n}],$$
           and we get $$\|\nabla(h_{\tau}^{n+1}\circ\Omega_{t,\tau}^{n})-(\nabla h_{\tau}^{n+1})\circ\Omega_{t,\tau}^{n}\|_{\mathcal{Z}^{(\lambda_{n}'^{\dag}-B_{0})(1+b)
           ,\mu_{n}'^{\dag};1}_{\tau-\frac{bt}{1+b}}}$$
           $$\leq\|\nabla(\Omega_{t,\tau}^{n}-Id)\|_{\mathcal{Z}^{(\lambda_{n}'^{\dag}-B_{0})(1+b)
           ,\mu_{n}'^{\dag}}_{\tau-\frac{bt}{1+b}}}\|(\nabla h_{\tau}^{n+1})\circ\Omega_{t,\tau}^{n}\|_{\mathcal{Z}^{(\lambda_{n}'^{\dag}-B_{0})(1+b)
           ,\mu_{n}'^{\dag};1}_{\tau-\frac{bt}{1+b}}}$$
           $$\leq C\bigg(\frac{1+\tau}{\min\{\lambda_{n}'-\lambda_{n}'^{\dag},\mu_{n}'-\mu_{n}'^{\dag}\}}\bigg)^{2}
           \|\Omega_{t,\tau}^{n}-Id\|_{\mathcal{Z}^{(\lambda_{n}'-B_{0})(1+b)
           ,\mu_{n}'}_{\tau-\frac{bt}{1+b}}}\| h_{\tau}^{n+1}\circ\Omega_{t,\tau}^{n}\|_{\mathcal{Z}^{(\lambda_{n}'-B_{0})(1+b)
           ,\mu_{n}';1}_{\tau-\frac{bt}{1+b}}}$$
           $$\leq\frac{CC^{4}_{\omega}}{\min\{\lambda_{n}'-\lambda_{n}'^{\dag},\mu_{n}'-\mu_{n}'^{\dag}\}^{2}}\bigg(
           \sum^{n}_{k=1}\frac{\delta_{k}}{(2\pi(\lambda_{k}-\lambda'_{k}))^{6}}\bigg)(1+\tau)^{-2}\| h_{\tau}^{n+1}\circ\Omega_{t,\tau}^{n}\|_{\mathcal{Z}^{(\lambda_{n}'-B_{0})(1+b)
           ,\mu_{n}';1}_{\tau-\frac{bt}{1+b}}},$$
           the above inequality implies $\nabla(h_{\tau}^{n+1}\circ\Omega_{t,\tau}^{n})\simeq(\nabla h_{\tau}^{n+1})\circ\Omega_{t,\tau}^{n}$ as $\tau\rightarrow\infty.$

           Since $$\|\nabla(h_{\tau}^{n+1}\circ\Omega_{t,\tau}^{n})\|_{\mathcal{Z}^{(\lambda_{n}'^{\dag}-B_{0})(1+b)
           ,\mu_{n}'^{\dag};1}_{\tau-\frac{bt}{1+b}}}\leq C\bigg(\frac{1+\tau}{\min\{\lambda_{n}'-\lambda_{n}'^{\dag},\mu_{n}'-\mu_{n}'^{\dag}\}}\bigg)^{2}
           \| h_{\tau}^{n+1}\circ\Omega_{t,\tau}^{n}\|_{\mathcal{Z}^{(\lambda_{n}'-B_{0})(1+b)
           ,\mu_{n}';1}_{\tau-\frac{bt}{1+b}}},$$
           and
 $$\|\nabla_{x}(h_{\tau}^{n+1}\circ\Omega_{t,\tau}^{n})\|_{\mathcal{Z}^{(\lambda_{n}'^{\dag}-B_{0})(1+b)
           ,\mu_{n}'^{\dag};1}_{\tau-\frac{bt}{1+b}}}+\|(\nabla_{x}+\tau\nabla_{v})(h_{\tau}^{n+1}\circ\Omega_{t,\tau}^{n})\|_{\mathcal{Z}^{(\lambda_{n}'^{\dag}-B_{0})(1+b)
           ,\mu_{n}'^{\dag};1}_{\tau-\frac{bt}{1+b}}}$$
          $$ \leq \frac{C}{\min\{\lambda_{n}'-\lambda_{n}'^{\dag},\mu_{n}'-\mu_{n}'^{\dag}\}}
           \| h_{\tau}^{n+1}\circ\Omega_{t,\tau}^{n}\|_{\mathcal{Z}^{(\lambda_{n}'-B_{0})(1+b)
           ,\mu_{n}';1}_{\tau-\frac{bt}{1+b}}},$$
           we have $$\|(\nabla_{x}h_{\tau}^{n+1})\circ\Omega_{t,\tau}^{n}\|_{\mathcal{Z}^{(\lambda_{n}'^{\dag}-B_{0})(1+b)
           ,\mu_{n}'^{\dag};1}_{\tau-\frac{bt}{1+b}}}+\|((\nabla_{x}+\tau\nabla_{v})(h_{\tau}^{n+1})\circ\Omega_{t,\tau}^{n}\|_{\mathcal{Z}^{(\lambda_{n}'^{\dag}-B_{0})(1+b)
           ,\mu_{n}'^{\dag};1}_{\tau-\frac{bt}{1+b}}}$$
           $$\leq C\bigg(\frac{C^{4}_{\omega}}{\min\{\lambda_{n}'-\lambda_{n}'^{\dag},\mu_{n}'-\mu_{n}'^{\dag}\}^{2}}\bigg(
           \sum^{n}_{k=1}\frac{\delta_{k}}{(2\pi(\lambda_{k}-\lambda'_{k}))^{6}}\bigg)
           +\frac{1}{\min\{\lambda_{n}'-\lambda_{n}'^{\dag},\mu_{n}'-\mu_{n}'^{\dag}\}}\bigg) \| h_{\tau}^{n+1}\circ\Omega_{t,\tau}^{n}\|_{\mathcal{Z}^{(\lambda_{n}'-B_{0})(1+b)
           ,\mu_{n}';1}_{\tau-\frac{bt}{1+b}}}.$$

\begin{center}
\item\subsection{The Choice of $\delta_{n+1}$}
%{\bf\large 1. \quad Introduction }
\end{center}

Now we see that the  corresponding $n+1$th step conclusion of the inductive hypothesis  has all been established with
$$\delta_{n+1}=\frac{C_{F}(1+C_{F})(1+C^{4}_{\omega})e^{CT^{2}_{n}}}{\min\{\lambda^{\ast}_{n}-\lambda_{n+1},\mu^{\ast}_{n}-\mu_{n+1}^{9}\}}
\max\bigg\{1,\sum^{n}_{i=1}\delta_{k}\bigg\}\bigg(1+\sum^{n}_{i=1}\frac{\delta_{i}}{(2\pi(\lambda_{i}-\lambda^{\ast}_{i}))^{6}}\bigg)\delta^{2}_{n}.$$

For any $n\geq1,$ we set
$\lambda_{n}-\lambda^{\ast}_{n}=\lambda^{\ast}_{n}-\lambda_{n+1}=\mu_{n}-\mu^{\ast}_{n}=\mu^{\ast}_{n}-\mu_{n+1}=\frac{\Lambda}{n^{2}}$
for some $\Lambda>0.$ By choosing $\Lambda$ small enough, we can make sure that the conditions
$2\pi(\lambda_{k}-\lambda^{\ast}_{k})<1\quad \textmd{and}\quad 2\pi(\mu_{k}-\mu^{\ast}_{k})<1$
are satisfied for all $k,$ as well as the other smallness assumptions made throughout this section.  We also have
$\lambda_{k}-\lambda^{\ast}_{k}\geq\frac{\Lambda}{k^{2}}.$  $(\mathbf{I})-(\mathbf{VIII})$  will be satisfied if  we  choose constants $\Lambda,\omega>0$
such that $\sum^{\infty}_{i=1}i^{12}\delta_{i}\leq\Lambda^{6}\omega.$

Then we have that $T_{n}\leq C_{\gamma}(n^{2}/\Lambda)^{\frac{7+\gamma}{\gamma-1}},$ so  the induction relation on $\delta_{n}$ gives
$\delta_{1}\leq C\delta \quad \textmd{\textmd{and}} \quad\delta_{n+1}=C(\frac{n^{2}}{\Lambda})^{9}e^{C(n^{2}/\Lambda)^{(14+2\gamma)/(\gamma-1)}}\delta^{2}_{n}.$

To make this relation hold, we also assume that $\delta_{n}$ is bounded below by the error coming from the short-time iteration; but this
follows easily by construction, since the constraints imposed on $\delta_{n}$ are much worse than those on the short-time iteration.

Having fixed $\Lambda,$ we will check that for $\delta$ small enough, the above relation for $\delta_{n+1}$ holds and there is  a fast convergence of $\{\delta_{i}\}^{\infty}_{i=1}.$
 The details are similar to that of the local-time case,and it can be also found in [23], here we omit it.
 $$$$

 Acknowledgements: I would like to  thank  professor Pin Yu in YMSC, Tsinghua University who gave me lots of advice. At the same time,
  I am deeply grateful to Professor Jinxin Xue  in YMSC, Tsinghua University who spent lots of time to discuss with me. At last, I'm also very grateful to the anonymous editor who give me lots of useful advice.

\end{document}